\documentclass[11pt]{article}
\usepackage[a4paper,margin=25mm]{geometry}
\usepackage{amsmath,amssymb,amsthm,mathtools,mathrsfs,microtype,enumitem,booktabs,xcolor,needspace}
\usepackage{tikz}
\usetikzlibrary{arrows.meta,calc,decorations.markings}
\definecolor{navy}{RGB}{25,53,82}
\usepackage[colorlinks=true,linkcolor=navy,urlcolor=navy,citecolor=navy]{hyperref}
\newtheorem{theorem}{Theorem}[section]
\newtheorem{lemma}[theorem]{Lemma}
\newtheorem{proposition}[theorem]{Proposition}
\newtheorem{corollary}[theorem]{Corollary}
\newtheorem{definition}[theorem]{Definition}
\theoremstyle{remark}\newtheorem{remark}[theorem]{Remark}
\newcommand{\T}{\mathbb T}\newcommand{\R}{\mathbb R}\newcommand{\C}{\mathbb C}
\newcommand{\Z}{\mathbb Z}

\newcommand{\norm}[1]{\left\|#1\right\|}
\newcommand{\dd}{\,d}
\DeclareMathOperator{\sn}{sn}\DeclareMathOperator{\cn}{cn}\DeclareMathOperator{\dn}{dn}

\hypersetup{pdftitle={Interior caustics and rigidity near noncircular ellipses},pdfauthor={Guan Huang}}
\title{Interior caustics and rigidity near noncircular ellipses}
\author{Guan Huang}
\date{}
\begin{document}
\maketitle
\begin{abstract}
We prove local rigidity near every fixed noncircular ellipse from integrable
rational caustics in an arbitrary nonempty open rotation interval
$I\subset(0,1/2)$. For each fixed $r>11/2$ and bound on the $H^r$ norm of
the elliptic graph displacement, every sufficiently $C^1$-small positively
curved perturbation admitting such caustics at all rational rotations in $I$
is an ellipse. A fixed $C^6$ bound may replace the Sobolev bound, with $C^1$
smallness retained. The caustic hypothesis can be restricted to a suitably
chosen sequence of rational rotations in $I$ with any prescribed rational
$\rho_*\in I$ as its unique accumulation point, together with finitely many
supplementary rotations in $I$; this family depends only on the reference
ellipse, interval, and anchor. The main linear result classifies all real
$L^2$ displacements satisfying the constant rational orbit-sum equations
throughout $I$ as infinitesimal ellipse deformations, and establishes their
analyticity and rational continuation in the complexified boundary parameter.
The proof combines analytic continuation and elimination of monodromy with
quantitative caustic estimates, finite detection, and Fourier interpolation.
\end{abstract}

\section{Introduction and main results}
Let $\Omega\subset\R^2$ be a bounded strictly convex domain with a $C^2$
boundary of positive curvature. The billiard map sends a boundary point and
an inward unit direction to the next reflection, according to the equal-angle
law. A convex caustic is a simple closed convex curve inside $\Omega$ whose
supporting lines remain supporting lines after reflection. We choose the
orientation of tangency for which the associated rotational invariant curve
has rotation number $\rho\in(0,1/2)$. For an ellipse, the confocal elliptic
caustics realize all these rotations; the limits $\rho\to0$ and
$\rho\to1/2$ correspond to the boundary and the focal segment, respectively
\cite{CF}.

The Birkhoff--Poritsky problem has several formulations, depending on the
region of phase space and the notion of integrability; see the accounts
\cite{FKS,BMSurvey2026}. The form considered here asks whether a sufficiently
smooth strictly convex table must be an ellipse, with the circle included,
if an open annulus in its phase cylinder is foliated by rotational invariant
curves corresponding to convex caustics. Lazutkin \cite{Lazutkin} proved that
sufficiently smooth positively curved tables already possess a positive-measure
family of caustics near the boundary; this alone does not give such a foliation.
We ask whether, near a fixed noncircular ellipse, integrable rational caustics
with rotation numbers in an arbitrary prescribed interior interval already
force elliptic rigidity. In particular, the prescribed rotations need not
approach either limiting value.

\begin{definition}We call a caustic
\emph{integrable rational} if every point on its invariant curve is
periodic. 
\end{definition}

For reduced rotation $p/d$, every point on this invariant curve has primitive
period $d$; the \hyperref[geo:par:common-period]{common-period argument}
is given in Section~\ref{geo:sec:calibration}.

Local rigidity from near-boundary rational caustics was established by
Avila--De Simoi--Kaloshin \cite{ADK} for ellipses of small eccentricity and
by Kaloshin--Sorrentino \cite{KS} near every fixed ellipse, using the
rotations $1/q$, $q\ge3$. Huang--Kaloshin--Sorrentino
\cite{HKS} studied smaller near-boundary rotation ranges
at small positive eccentricity: their general-cutoff result is conditional
on explicit nondegeneracy tests. Koval \cite{Koval}
proved rigidity from rational integrability in $0<\rho\le1/q_0$, for each
fixed $q_0>2$, outside a cutoff-dependent locally finite exceptional set of
eccentricities. His theorem uses a fixed $C^k$ bound, $k\ge39$, and $C^1$
smallness; its caustic hypothesis can be reduced to a $1/q$ tail and
finitely many supplementary rotations. These prescribed families accumulate
at zero, rather than at a positive interior anchor.

\subsection*{Main rigidity result}
For a noncircular ellipse $E$ with semiaxes $a>b>0$, use principal-axis
coordinates and write a nearby boundary as
$\Gamma=\{\mathcal X_E(\mu(\varphi),\varphi):\varphi\in\T\}$, where
$\T=\R/(2\pi\Z)$. The transverse elliptic chart $\mathcal X_E$ is defined
in \eqref{geo:eq:chart} and satisfies
$\mathcal X_E(0,\varphi)=(a\cos\varphi,b\sin\varphi)$, so the zero graph
is $E$. All graph norms refer to the elliptic boundary parameter $\varphi$,
not the Euclidean polar angle. Our Sobolev convention is
\begin{equation}\label{main:eq:sobolev}
\|u\|_{H^s}^2=\sum_{n\in\Z}(1+|n|)^{2s}|\widehat u(n)|^2,
\qquad
\widehat u(n)=\frac1{2\pi}\int_0^{2\pi}u(\varphi)e^{-in\varphi}\,d\varphi,
\end{equation}
with distributional Fourier coefficients when necessary.

\begin{theorem}[Rigidity from any interior band]\label{main:band}
Fix a noncircular ellipse $E_0$, a nonempty open interval
$I\subset(0,1/2)$, a real $r>11/2$, and $M>0$. There exists
$\varepsilon=\varepsilon(E_0,I,r,M)>0$ such that every positively curved
$H^r$ elliptic graph $\Gamma$ over $E_0$ satisfying
\[
 \|\mu\|_{H^r}\le M,\qquad \|\mu\|_{C^1}<\varepsilon
\]
is an ellipse if its billiard has an integrable rational caustic at
every rotation in $\mathbb Q\cap I$.
\end{theorem}

Only the $C^1$ displacement is required to be small; the high-norm bound $M$
is fixed but need not be small. The resulting ellipse may differ from $E_0$
in its center, semiaxes, and orientation. By
Corollary~\ref{main:cor:holder-regularity}, a fixed $C^{5,\alpha}$ bound,
$1/2<\alpha\le1$, can replace the Sobolev bound; in particular, a fixed
$C^6$ bound suffices. The threshold is local to the chosen noncircular
ellipse: no uniformity at eccentricity zero or one, and no optimality of
the regularity exponent, is asserted.

The full-band hypothesis can be weakened substantially. For any prescribed
rational anchor $\rho_*\in I$, one suitably chosen sequence accumulating at
$\rho_*$, together with finitely many additional rotations in $I$, suffices.
Specifically, Theorem~\ref{band:thm:sparse-rigidity} requires integrable
rational caustics only at
\[
 \mathcal S_{E_0,I,\rho_*}
 =\left\{\frac{\lceil q\rho_*\rceil}{q}:q\ge N\right\}\cup F,
 \qquad F\subset\mathbb Q\cap I\ \text{finite}.
\]
The cutoff $N$ and the set $F$ depend on the reference ellipse, the band,
and the anchor, but not on $r$ or $M$. Every rotation in this family lies
in $I$, and its only accumulation point is $\rho_*$. The displayed
fractions need not be reduced; the hypothesis concerns their rotation
values. The finite set supplements an infinite sequence; rigidity from
finitely many caustics alone is not asserted.

\begin{remark} Theorem \ref{main:band} does not include the circular case. A relative  local rigidity result based on interior rotations was obtained
by Huang~\cite{Huang}. For $P\ge3$ and $0<\delta<1/6$, it is established there  rigidity near ellipses of sufficiently small
eccentricity, including the circle, under a fixed $C^{16}$ bound and
sufficiently small $C^1$ displacement. The caustic hypothesis requires
integrable rational caustics at every rational rotation in
$[1/P,1/P+\delta)$ and at finitely many prescribed rotations of the form
$1/q$, which need not lie in that interval.
\end{remark}
Other rigidity results use different hypotheses. A foliation of the entire
phase cylinder by rotational invariant curves forces a disk by Bialy's
theorem \cite{BialyHopf}, while convex caustics with rotations tending to
$1/2$ force an ellipse by Innami's theorem \cite{innami} (see also
Arnold--Bialy \cite[Theorem~2]{AB}). Thus the full-band conclusion of
Theorem~\ref{main:band} also follows from Innami's theorem when $\sup I=1/2$.
This observation does not cover the sparse hypothesis at an arbitrary
interior anchor. Bialy--Mironov \cite{BM} obtain elliptic rigidity under
central symmetry and a foliation between the boundary and a $4$-periodic
invariant curve. Kaloshin--Koudjinan--Zhang \cite{KKZ} treat a
high-regularity perturbative setting near a centrally symmetric domain,
assuming rational integrability for $0<p/q\le1/4$ and an additional
integrable $1/3$-caustic. Our assumptions impose no central symmetry on the
perturbed table.

Glutsyuk \cite{GlutsyukComplex} obtains conic rigidity from a complex
algebraic caustic, without assuming closeness to an ellipse. Its complex
tangent lines must be preserved by complexified reflection, an additional
hypothesis beyond the existence of a real convex caustic. Here neither the
boundary nor its caustics are assumed algebraic: the rational continuation
proved below concerns the linearized displacement, not a caustic supplied
as input.

\subsection*{Linearized rigidity}
For an integrable rational caustic, the perimeter of the associated periodic
polygons is independent of the starting point. Linearizing this phase
independence at an ellipse gives a constant orbit-sum condition for the
boundary displacement in the graph chart above. Infinitesimal changes of
the center, the two semiaxes, and the orientation necessarily satisfy these
equations. The central linear question is whether there are any other
solutions.

Pinto-de-Carvalho--Ram\'irez-Ros \cite{PRR} express the first-order
persistence condition as constancy of an elliptic Melnikov sum and use
complex singularities to obstruct odd-period persistence for nonconstant
first variations that are entire in the complexified boundary angle.
Huang--Kaloshin \cite{HK} obtained a finite-dimensional
infinitesimal reduction for smooth deformations preserving near-boundary
integrability.

Recent persistence results address complementary questions. In an analytic
one-parameter family of strictly convex analytic billiard domains,
Fierobe--Sorrentino \cite[Theorem~2.16]{FS} show that a periodic invariant
curve at a fixed rational rotation exists either for every parameter or only
for a discrete set. Koudjinan--Ram\'irez-Ros \cite{KRR} develop higher-order
persistence conditions for resonant caustics in perturbations of a circle.
Here we instead vary the rotation over a band for a fixed reference ellipse,
and require exact caustics, not finite-order persistence, in the nonlinear
problem. Our linear classification starts with an $L^2$ displacement:
analyticity and rational continuation follow from its trace equations.

Let $x$ be the action coordinate in which the reference caustic map is
translation by $2\pi\rho$, and let $\Phi_{E,\rho}$ convert $x$ to the
boundary parameter $\varphi$; its formula is given in \eqref{geo:eq:Phi}.
The constant in each trace below may depend on $\rho$, but not on the
starting phase $x$. For an analytic representative, $V(z)$ denotes its
holomorphic germ in $z=e^{i\varphi}$. This complexified boundary parameter
is not a Cartesian coordinate of the billiard table.

\begin{theorem}[Interior-band classification]\label{main:analytic}
Fix a noncircular ellipse $E$ with semiaxes $a>b>0$ and a nonempty open
interval $I\subset(0,1/2)$. Let $v\in L^2(\T;\R)$. Suppose that, for
every reduced rational $\rho=p/d\in I$,
\[
 \sum_{j=0}^{d-1}v\bigl(\Phi_{E,\rho}(x+2\pi jp/d)\bigr)
\]
is constant for almost every $x$. Here $\Phi_{E,\rho}$ is the caustic
action-angle parametrization in \eqref{geo:eq:Phi}. Then $v$ agrees
almost everywhere with an element of
\[
 T_E=D_E^{-1}\operatorname{span}_{\R}
 \{1,\cos\varphi,\sin\varphi,\cos2\varphi,\sin2\varphi\},
 \qquad D_E=b^2\cos^2\varphi+a^2\sin^2\varphi.
\]
In particular, it has a real analytic representative. The physical germ
$V(z)$ of that representative, $z=e^{i\varphi}$, is rational, regular
at zero and infinity, with at most simple poles at the four zeros of
\[
D_E(z)=\frac{a^2+b^2}{2}-\frac{a^2-b^2}{4}(z^2+z^{-2}).
\]
\end{theorem}

The space $T_E$ is the tangent space to the five-parameter family of ellipses
in the graph chart. Conversely, each of its elements has constant rational
traces by Corollary~\ref{geo:cor:tangentkill}, so the classification is exact.
Corollary~\ref{band:cor:distributional-classification} extends it to periodic
distributions, interpreting each summand as a pullback.

The classification enters the nonlinear argument through a quantitative
inverse estimate. A sequence of equations selected near one rational anchor
leaves a finite-dimensional analytic kernel. After fixing the five ellipse
parameters, full-band classification of the reference operator allows
finitely many additional resonant rows from $I$ to separate the remaining
kernel. Only the selected sequence and these finitely many rotations are
then required of the perturbed table. The inverse estimate combines with
quadratic caustic errors to yield rigidity.

\subsection*{Proof strategy}
Our proof follows the perturbative approach developed by
Kaloshin--Sorrentino \cite[Section~4.2 and Sections~5--8]{KS} and
Koval \cite[Sections~1.3--1.6 and 5--7]{Koval}. Preservation of rational
caustics yields linear constraints on the boundary displacement, with
quadratic remainders. Nondegeneracy modulo infinitesimal ellipse
deformations is then combined with quantitative estimates to obtain
nonlinear rigidity. Complex extensions and singularity analysis of the
dynamical modes play a central role in \cite{KS}. Koval combines
finite-dimensional nondegeneracy with holomorphic dependence on
eccentricity, obtaining rigidity outside a locally finite exceptional set.

Here the selected rotations approach a prescribed positive rational
anchor rather than zero. The resulting high-frequency operator has a
nontrivial principal part that depends on the Fourier index modulo the
anchor denominator. Twist makes this principal part invertible, but the
associated parametrix gives only a finite-dimensional analytic kernel.
To identify the full-band solution space with $T_E$, we continue the
unknown displacement and eliminate its monodromy. The argument proceeds
in four stages.

\smallskip\noindent
\emph{Step 1: normalization and nonlinear estimates.}
In Sections~\ref{geo:sec:convention}--\ref{geo:sec:calibration}, we choose
a nearby reference ellipse so that the residual displacement has zero
moments in the five ellipse directions. On a fixed compact rotation
interval $J\Subset I$ containing the chosen sparse caustic family,
calibration gives an orbit-sum error $O(d^3\|\mu\|_{C^1}^2)$ and a
resonant Fourier error $O(d^2\|\mu\|_{C^1}^2)$, with a smallness threshold
independent of the primitive denominator $d$.
Appendix~\ref{app:sec:norms} controls the changes of graph coordinates.

\smallskip\noindent
\emph{Step 2: automatic regularity.}
Section~\ref{op:sec:operator} constructs the operator associated with the
selected sequence. Twist separates the nodes of its principal Vandermonde
matrix, and a localized parametrix shows that the kernel is finite
dimensional and that every distributional kernel element is analytic.
This establishes regularity, but not yet the full-band classification.

\smallskip\noindent
\smallskip\noindent
\emph{Step 3: global classification.}
Sections~\ref{arith:section}--\ref{cont:sec:logarithmic-closure} identify
an analytic solution of the full-band trace equations by continuing its
physical germ $V(z)$ and controlling its possible singularities.
Section~\ref{arith:section} supplies the geometric and arithmetic tools:
a proper degree-two Jacobi map relates the action and physical variables,
the two preimages of each exterior focal value differ by one caustic
step, and a finite-grid lemma and Picard--Fuchs identities constrain how
marked preimages can move with the caustic parameter. The geometric
degree does not bound the number of branches of the unknown $V$.

In Section~\ref{vol:sec:borel-volterra}, a convergent expansion of the
selected row symbols gives an exact Borel--Volterra equation. Retaining
the full integration history when restarting this equation gives
continuation along every path in the finite physical exterior outside a
fixed closed countable set with only finitely many accumulation points
there. Section~\ref{cont:sec:multiple-bases} chooses two rational bases
with disjoint accumulation bounds and glues their continuations while
preserving the germ already reached. The remaining possible obstructions
have no finite accumulation, although the continuation may still be
multivalued.

Sections~\ref{cont:sec:task6-finite-bad-set}--\ref{cont:sec:task6-nonfocal}
reduce these obstructions to the exterior focal pair. A point is retained
only if some accessible germ fails to extend holomorphically across it.
For one odd-denominator trace, the nonconstant companion maps are finite
at both Jacobi poles; their punctured neighborhoods then make every
sufficiently large finite exception removable for all accessible germs.
Together with the exception-free collar, this makes the retained set
finite. A retained nonfocal value would require a companion collision at
every rotation in the band with odd prime denominator. The finite-grid
lemma forces one normalized gap between geometric preimages to be
affine in the rotation, whereas the Picard--Fuchs identities rule out
such persistent nonfocal collisions. Once only the focal values remain,
an odd trace also extends every accessible branch holomorphically
through infinity.

Section~\ref{cont:sec:logarithmic-closure} removes the remaining monodromy.
At a focal point, a companion branch $g$ is locked to a prescribed source
branch $f$ before the caustic parameter $m$ varies. In a fixed local
physical coordinate $y$, their sum $f(y)+g(T_m(y))$ extends holomorphically
across $y=0$ whenever the rotation is rational with odd denominator;
here $T_m$ is the physical involution induced by the focal pair.
Density in the parameter first gives a single-valued paired sum;
fixed-contour Laurent coefficients then remove its puncture jointly.
Its parameter derivative is $g'(T_m(y))\partial_mT_m(y)$, where
$\partial_mT_m(y)$ vanishes to exactly second order at zero. This forces
each branch derivative to be single-valued meromorphic locally, with
pole order at most two, but still allows logarithms in the branches
themselves.

We then descend $dV$ globally, cancel each focal residue separately,
and show that all periods vanish. The resulting single-valued
continuation has at most simple focal poles. Reflection using reality
on the unit circle yields a rational function, regular at zero and
infinity, with possible simple poles only at the four zeros of $D_E$;
its boundary restriction belongs to $T_E$. Neither finite monodromy
nor initial meromorphicity is assumed.
\smallskip\noindent

\emph{Step 4: nonlinear closure.}
Finally, Section~\ref{close:sec:closure} uses the classification to select
finitely many additional resonant rows from $I$ that detect the transverse
kernel. The resulting inverse estimate, the quadratic errors, and Fourier
interpolation force the normalized residual to vanish for $r>11/2$.
This proves the sparse-family theorem, and hence the full-band theorem.

\section{A fixed real action-angle convention}\label{geo:sec:convention}
We recall the classical Jacobi parametrization of elliptic billiards;
see Chang--Friedberg \cite{CF} and, for the explicit formulas,
Damasceno--Dias Carneiro--Ram\'irez-Ros
\cite[Lemma~2.1 and equations~(2.3)--(2.6)]{DCRR}.
Fix a noncircular ellipse $E$ with semiaxes $a>b>0$ and use its
positively oriented parametrization
$r(\varphi)=(a\cos\varphi,b\sin\varphi)$, $\varphi\in\T$. Put
\[
 D_E(\varphi)=|r'(\varphi)|^2
             =b^2\cos^2\varphi+a^2\sin^2\varphi.
\]
We fix the phase to match this parametrization; $\varphi$ is the
elliptic boundary parameter, not the Euclidean polar angle.

We write elliptic functions in terms of the \emph{parameter} $m=k^2$,
where $k$ is the modulus:
\[
 F(\alpha\mid m)=\int_0^\alpha\frac{\dd t}{\sqrt{1-m\sin^2t}},
 \qquad K(m)=F(\pi/2\mid m),\qquad 0\le m<1.
\]
The amplitude $\operatorname{am}(u\mid m)$ is the inverse of
$F(\cdot\mid m)$ as an increasing map from $\R$ onto $\R$.
For real arguments,
$\operatorname{sn}(u\mid m)=\sin\operatorname{am}(u\mid m)$,
$\operatorname{cn}(u\mid m)=\cos\operatorname{am}(u\mid m)$, and
$\operatorname{dn}(u\mid m)=\sqrt{1-m\operatorname{sn}^2(u\mid m)}$.
For $0<\lambda<b$, the confocal elliptic caustic and its elliptic
parameter are
\[
C_\lambda:\quad \frac{x^2}{a^2-\lambda^2}+\frac{y^2}{b^2-\lambda^2}=1,
\qquad m_\lambda=\frac{a^2-b^2}{a^2-\lambda^2}.
\]
The parametrization in \cite{DCRR} is
$(a\operatorname{sn}u,b\operatorname{cn}u)$, with modulus
$k=\sqrt{m_\lambda}$. Reflecting in the horizontal axis preserves the
billiard law and gives the positively oriented parametrization
\[
q_\lambda(u)=\bigl(a\operatorname{sn}(u\mid m_\lambda),
                  -b\operatorname{cn}(u\mid m_\lambda)\bigr).
\]
In this variable, the forward billiard map associated with $C_\lambda$
acts by $u\mapsto u+\delta_\lambda$, where
\[
\delta_\lambda=2F(\arcsin(\lambda/b)\mid m_\lambda)
\in(0,2K(m_\lambda)).
\]
To express it in our boundary parameter and normalize the action angle,
set
\[
\varphi=\operatorname{am}(u\mid m_\lambda)-\frac\pi2,
\qquad u=K(m_\lambda)+\frac{2K(m_\lambda)}\pi\theta.
\]
Then $q_\lambda(u)=(a\cos\varphi,b\sin\varphi)$, and the resulting
coordinate and its inverse are
\begin{equation}\label{geo:eq:G}
G_\lambda(\varphi)=\frac{\pi}{2K(m_\lambda)}
 \big[F(\varphi+\pi/2\mid m_\lambda)-K(m_\lambda)\big],
\end{equation}
\begin{equation}\label{geo:eq:Phi}
\Phi_\lambda(\theta)=
\operatorname{am}\!\left(K(m_\lambda)+\frac{2K(m_\lambda)}\pi\theta\mid m_\lambda\right)
-\frac\pi2.
\end{equation}
Both are increasing degree-one lifts fixing zero. Differentiation gives
\[
G_\lambda'(\varphi)=
\frac{\pi\sqrt{a^2-\lambda^2}}{2K(m_\lambda)\sqrt{D_E(\varphi)-\lambda^2}}>0.
\]
Since one real Jacobi period is $4K(m_\lambda)$, the shift
$\delta_\lambda$ becomes translation through
\begin{equation}\label{geo:eq:rotation}
2\pi\rho(\lambda),\qquad
\rho(\lambda)=\frac{F(\arcsin(\lambda/b)\mid m_\lambda)}{2K(m_\lambda)}.
\end{equation}

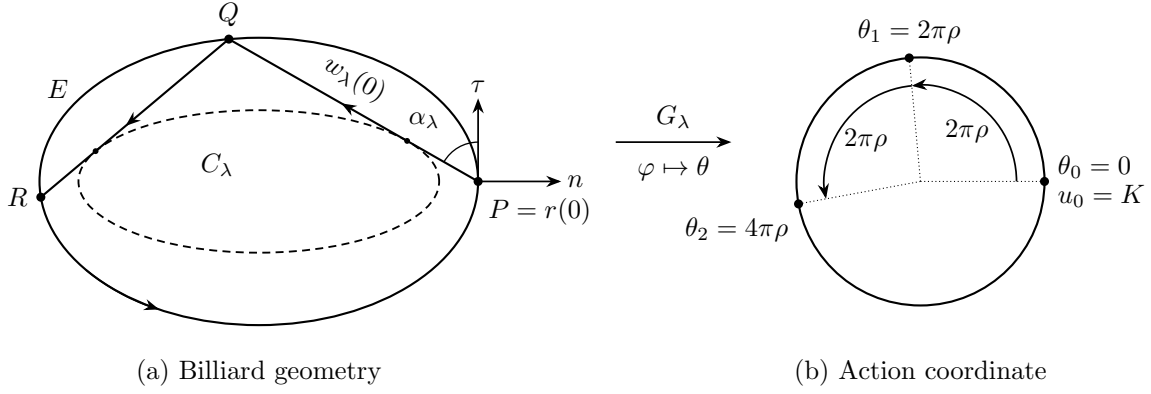
\begin{figure}[htbp]
\centering
\begin{tikzpicture}[
  x=1cm,y=1cm,font=\small,
  boundary/.style={line width=.8pt},
  caustic/.style={densely dashed,line width=.7pt},
  orbit/.style={line width=.8pt,postaction={decorate},
    decoration={markings,mark=at position .55 with {\arrow{Stealth[length=2mm]}}}},
  vector/.style={-{Stealth[length=2mm]},line width=.65pt},
  every node/.style={inner sep=2pt}]
\coordinate (P) at (2.900000000000,0.000000000000);
\coordinate (Q) at (-0.396517195714,1.882155851147);
\coordinate (R) at (-2.882533893686,-0.208215351402);
\draw[boundary] (0,0) ellipse[x radius=2.9,y radius=1.9];
\draw[caustic] (0,0) ellipse[x radius=2.384848003542,y radius=0.942072184071];
\draw[orbit] (P)--node[pos=.55,above=3pt,sloped] {$w_\lambda(0)$}(Q);
\draw[orbit] (Q)--(R);
\fill (1.961206896552,0.536006587489) circle[radius=1.1pt];
\fill (-2.158518133831,0.400574483287) circle[radius=1.1pt];
\foreach \point in {P,Q,R}{\fill (\point) circle[radius=1.7pt];}
\node[anchor=north west] at (2.96,-.14) {$P=r(0)$};
\node[anchor=south] at ($(Q)+(0,.10)$) {$Q$};
\node[anchor=east] at ($(R)+(-.10,0)$) {$R$};
\node at (-2.66,1.23) {$E$};
\node at (-.55,.22) {$C_\lambda$};
\draw[vector] (P)--++(0,1.10) node[above] {$\tau$};
\draw[vector] (P)--++(1.10,0) node[right] {$n$};
\draw[line width=.5pt] ($(P)+(0,.52)$)
  arc[start angle=90,end angle=150.275670852893,radius=.52];
\node at ($(P)+(-.71,.81)$) {$\alpha_\lambda$};

\draw[vector] ({2.9*cos(218)},{1.9*sin(218)})
  arc[start angle=218,end angle=243,x radius=2.9,y radius=1.9];
\node[align=center] at (0,-2.52) {(a) Billiard geometry};
\draw[vector] (4.72,.52)--(6.25,.52)
  node[midway,above=3pt] {$G_\lambda$}
  node[midway,below=3pt] {$\varphi\mapsto\theta$};
\begin{scope}[shift={(8.75,0)}]
\draw[boundary] (0,0) circle[radius=1.64];
\coordinate (A0) at (0:1.64);
\coordinate (A1) at (95.220679020038:1.64);
\coordinate (A2) at (190.441358040076:1.64);
\foreach \point in {A0,A1,A2}{
  \draw[densely dotted,line width=.4pt] (0,0)--(\point);
  \fill (\point) circle[radius=1.7pt];
}
\draw[vector] (0:1.27)
  arc[start angle=0,end angle=95.220679020038,radius=1.27];
\draw[vector] (95.220679020038:1.27)
  arc[start angle=95.220679020038,end angle=190.441358040076,radius=1.27];
\node at (47.610339510019: .91) {$2\pi\rho$};
\node at (142.831018530057: .91) {$2\pi\rho$};
\node[anchor=west,align=left] at ($(A0)+(.10,0)$) {$\theta_0=0$\\[-1pt]$u_0=K$};
\node[anchor=south] at ($(A1)+(0,.12)$) {$\theta_1=2\pi\rho$};
\node[anchor=north east] at ($(A2)+(-.06,-.10)$) {$\theta_2=4\pi\rho$};
\node[align=center] at (0,-2.52) {(b) Action coordinate};
\end{scope}
\end{tikzpicture}
\caption{A confocal billiard orbit and its action-angle representation.
Left: the first three impacts $P,Q,R$, starting at $P=r(0)=(a,0)$,
and two consecutive chords tangent to $C_\lambda$.
The vectors $\tau$ and $n$ give the positive tangent and outward normal
at $P$; $w_\lambda(0)$ is the outgoing direction and the marked angle
 is $\alpha_\lambda=\alpha_\lambda(0)$.
Right: the same impacts have action angles
$\theta_j=2\pi j\rho$ for $j=0,1,2$, understood modulo $2\pi$.
Here $\rho=\rho(\lambda)$ and $K=K(m_\lambda)$; the initial phase is
$\theta_0=0$, $u_0=K$. The diagram uses
$a=2.9$, $b=1.9$, and $\lambda=1.65$; periodicity is not assumed.}
\label{geo:fig:caustic-action-angle}
\end{figure}

We also record the incidence identity used below. If
$w_\lambda(\varphi)$ is the forward outgoing unit direction at
$(a\cos\varphi,b\sin\varphi)$, then
\cite[Lemma~4.1 and its proof]{PRR} gives
\[
-w_\lambda(\varphi)\cdot(b\cos\varphi,a\sin\varphi)=\lambda.
\]
The vector $(b\cos\varphi,a\sin\varphi)$ is the outward unit normal
multiplied by $\sqrt{D_E(\varphi)}$. The tangential component of
$w_\lambda$ is positive for small $\lambda$ and cannot vanish for
$0<\lambda<b$, since $\lambda^2<D_E(\varphi)$.
Thus the outgoing angle, measured from the positive tangent toward the
interior, satisfies
\begin{equation}\label{geo:eq:joachim}
0<\alpha_\lambda(\varphi)<\pi/2,\qquad
\sin\alpha_\lambda(\varphi)=\lambda/\sqrt{D_E(\varphi)},\qquad
\partial_\lambda\alpha_\lambda(\varphi)=\frac1{\sqrt{D_E(\varphi)-\lambda^2}}>0.
\end{equation}
Figure~\ref{geo:fig:caustic-action-angle} illustrates the incidence angle
and the action-angle convention.

We next justify using $\rho$ as the caustic parameter.
In the chosen principal-axis coordinates, the formulas above depend
smoothly on $(a,b,\lambda)$ for $a>b>\lambda>0$.
Let $f_\lambda$ denote the first-impact lift in the boundary angle
$\varphi$. To check the twist sign in this convention, consider a chord
from $r(\varphi)$ to $r(\psi)$ with $0<\psi-\varphi<2\pi$.
Write $L=|r(\psi)-r(\varphi)|$ and let $\alpha,\beta\in(0,\pi)$
be its incidence angles at the initial and terminal endpoints.
Holding $\varphi$ fixed, differentiation of the chord direction gives
\[
 \partial_\psi\alpha=\frac{|r'(\psi)|\sin\beta}{L}>0.
\]
On the confocal orbit $\psi=f_\lambda(\varphi)$, reflection and the
incidence identity at the terminal endpoint give
$|r'(\psi)|\sin\beta=\lambda$. Combining this with
\eqref{geo:eq:joachim} yields
\[
 \partial_\lambda f_\lambda(\varphi)
 =\frac{|r(f_\lambda(\varphi))-r(\varphi)|}
        {\lambda\sqrt{D_E(\varphi)-\lambda^2}}>0.
\]
Differentiate
$G_\lambda(f_\lambda(\varphi))=G_\lambda(\varphi)+2\pi\rho(\lambda)$
at fixed $\varphi$, then substitute $\varphi=\Phi_\lambda(\theta)$
and average over $\theta\in[0,2\pi]$.
The two terms involving
$(\partial_\lambda G_\lambda)\circ\Phi_\lambda$ have equal averages,
since this function is $2\pi$-periodic. Consequently,
\[
2\pi\rho'(\lambda)=\frac1{2\pi}\int_0^{2\pi}
\frac{\partial_\lambda f_\lambda(\Phi_\lambda(\theta))}
 {\Phi_\lambda'(\theta+2\pi\rho(\lambda))}\dd\theta>0.
\]
As $\lambda\downarrow0$, formula \eqref{geo:eq:rotation} gives
$\rho(\lambda)\to0$.
As $\lambda\uparrow b$ we have $K(m_\lambda)\to\infty$ and
\[
0\le K(m_\lambda)-F(\arcsin(\lambda/b)\mid m_\lambda)
\le \arccos(\lambda/b)\sqrt{\frac{a^2-\lambda^2}{b^2-\lambda^2}}=O(1),
\]
so $\rho(\lambda)\to1/2$.
Thus $\lambda\mapsto\rho(\lambda)$ is a smooth diffeomorphism from
$(0,b)$ onto $(0,1/2)$.
We write $\lambda_\rho$ for its inverse and
$\Phi_\rho=\Phi_{\lambda_\rho}$.
For a fixed compact interval $J\Subset(0,1/2)$, the implicit-function
theorem and compactness give uniform bounds for $\lambda_\rho$,
$\Phi_\rho-\operatorname{id}$, and $G_{\lambda_\rho}-\operatorname{id}$,
together with any prescribed finite number of angle and parameter
derivatives, for $\rho\in J$ and $(a,b)$ in a sufficiently small
compact neighborhood of the reference semiaxes within $\{a>b>0\}$.
The angle-dependent bounds are taken over one period; subtracting
$\operatorname{id}$ makes the degree-one lifts periodic.


\section{Exact normalization against the ellipse family}\label{geo:sec:gauge}
Write $\T=\R/(2\pi\mathbb Z)$ and fix
\[
E=\{(a\cos\varphi,b\sin\varphi):\varphi\in\T\},\qquad a>b>0.
\]
Unless stated otherwise, constants may depend on $E$ and on a specified
compact rotation interval $J\Subset(0,1/2)$.
We describe nearby boundaries in the elliptic graph chart
\begin{equation}\label{geo:eq:chart}
\mathcal X_E(s,\varphi)=
\big((a\cosh s+b\sinh s)\cos\varphi,
     (b\cosh s+a\sinh s)\sin\varphi\big).
\end{equation}
We restrict $s$ to a sufficiently small interval about zero, on which
$\mathcal X_E$ is a smooth coordinate chart onto a neighborhood of $E$.
This is the standard confocal elliptic coordinate shifted so that
$s=0$ represents $E$.
Put
\[
D_E(\varphi)=b^2\cos^2\varphi+a^2\sin^2\varphi.
\]
The transverse vector
$\partial_s\mathcal X_E(0,\varphi)=(b\cos\varphi,a\sin\varphi)$
equals the outward unit normal multiplied by $\sqrt{D_E(\varphi)}$.

\begin{figure}[htbp]
\centering
\begin{tikzpicture}[
  x=1cm,y=1cm,font=\small,
  every node/.style={inner sep=2pt},
  axis/.style={-{Stealth[length=1.6mm]},line width=.45pt},
  curve/.style={line width=.9pt},
  reference/.style={densely dashed,line width=.65pt},
  fibre/.style={densely dotted,line width=.5pt},
  direction/.style={-{Stealth[length=1.7mm]},line width=.6pt},
  declare function={
    mu(\p)=.12*cos(\p)+.085*sin(\p)+.035*cos(3*\p);
    aa(\s)=2.55*(exp(\s)+exp(-\s))/2+1.65*(exp(\s)-exp(-\s))/2;
    bb(\s)=1.65*(exp(\s)+exp(-\s))/2+2.55*(exp(\s)-exp(-\s))/2;
  }]
\node[anchor=west] at (0,2.36) {(a) Lifted graph coordinates};
\draw[axis] (-.17,0)--(5.91,0) node[right] {$\varphi$};
\draw[axis] (0,-1.92)--(0,1.98) node[above] {$s$};
\draw[fibre] (5.5,-1.77)--(5.5,1.77);
\draw[fibre] (0,1.72)--(5.5,1.72);
\draw[fibre] (0,-1.72)--(5.5,-1.72);
\node[anchor=east] at (-.08,1.72) {$s_+$};
\node[anchor=east] at (-.08,-1.72) {$s_-$};
\node[anchor=north east] at (0,-.06) {$0$};
\node[anchor=north] at (5.5,-.06) {$2\pi$};
\draw[curve] plot[domain=0:360,samples=241] ({5.5*\x/360},{8*mu(\x)});
\node[anchor=south] at (2.15,1.20) {$s=\mu(\varphi)$};
\draw[reference] (0,0)--(5.5,0);
\node[anchor=north] at (3.65,-.06) {$s=0$};
\coordinate (bcoord) at ({5.5/6},0);
\coordinate (pcoord) at ({5.5/6},{8*mu(60)});
\draw[line width=.65pt] (bcoord)--(pcoord);
\fill (bcoord) circle[radius=1.5pt];
\fill (pcoord) circle[radius=1.5pt];
\node[anchor=north] at ($(bcoord)+(0,-.09)$) {$\varphi_0$};
\node[font=\footnotesize,rotate=90] at ($(bcoord)+(-.27,{4*mu(60)})$) {$\mu(\varphi_0)$};
\node[align=center] at (2.75,-2.18) {The edges $\varphi=0$ and $2\pi$ are identified.};
\draw[direction] (6.12,.23)--(7.61,.23)
 node[midway,above=3pt] {$\mathcal X_E$};
\begin{scope}[shift={(10.57,0)},scale=.85]
\node[anchor=west] at (-3.04,2.78) {(b) The physical curve};
\foreach \s in {-.215,.215}{
  \draw[fibre] plot[domain=0:360,samples=181] ({aa(\s)*cos(\x)},{bb(\s)*sin(\x)});
}
\foreach \p in {0,60,120,180,240,300}{
  \draw[fibre] plot[domain=-.215:.215,samples=31]
    ({aa(\x)*cos(\p)},{bb(\x)*sin(\p)});
}
\draw[reference] (0,0) ellipse[x radius=2.55,y radius=1.65];
\draw[curve] plot[domain=0:360,samples=301]
  ({aa(mu(\x))*cos(\x)},{bb(mu(\x))*sin(\x)});
\coordinate (B) at ({2.55*cos(60)},{1.65*sin(60)});
\coordinate (P) at ({aa(mu(60))*cos(60)},{bb(mu(60))*sin(60)});
\draw[line width=.8pt] plot[domain=0:{mu(60)},samples=25]
  ({aa(\x)*cos(60)},{bb(\x)*sin(60)});
\fill (B) circle[radius=1.5pt];
\fill (P) circle[radius=1.5pt];
\node[anchor=north west] at ($(B)+(.10,-.08)$) {$B$};
\node[anchor=west] at ($(P)+(.13,.06)$) {$P$};
\draw[direction] (B)--++({.4*1.65*cos(60)},{.4*2.55*sin(60)})
 node[above] {$N_0$};
\draw[direction] (B)--++({-.4*2.55*sin(60)},{.4*1.65*cos(60)})
 node[below left=3pt] {$T_0$};
\node at (-1.38,-.60) {$E$};
\draw[line width=.4pt] (-1.55,-.79)--({2.55*cos(225)},{1.65*sin(225)});
\node at (.52,-2.45) {$\partial\Omega_\mu$};
\draw[line width=.4pt] (.52,-2.21)--({aa(mu(285))*cos(285)},{bb(mu(285))*sin(285)});
\end{scope}
\end{tikzpicture}
\caption{The elliptic graph chart. Left: a periodic displacement
$\mu(\varphi)$ in a cut of the coordinate cylinder. Right: its physical
image $r_\mu(\varphi)=\mathcal X_E(\mu(\varphi),\varphi)$; the dashed
ellipse is the zero graph and the light dotted curves are coordinate
lines. The marked points are $B=\mathcal X_E(0,\varphi_0)$ and
$P=r_\mu(\varphi_0)$. The arrows indicate the directions
$N_0=\partial_s\mathcal X_E(0,\varphi_0)$ and
$T_0=\partial_\varphi\mathcal X_E(0,\varphi_0)$, which are orthogonal.
The fibre joining $B$ to $P$ is a coordinate curve, not a straight
normal segment; $\mu(\varphi_0)$ is a coordinate displacement, not a
Euclidean distance.}
\label{geo:fig:elliptic-graph-chart}
\end{figure}
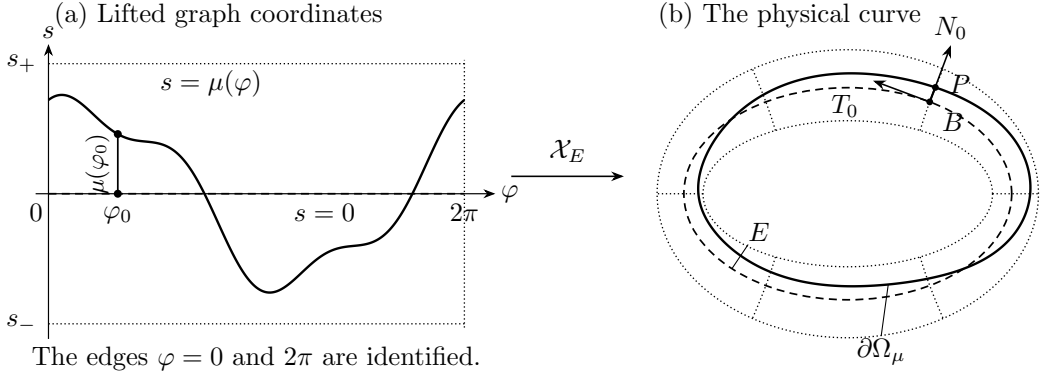

For $\mu\in C^2(\T)$ with sufficiently small $C^1$ norm, set
\[
r_\mu(\varphi)=\mathcal X_E(\mu(\varphi),\varphi),\qquad
\partial\Omega_\mu=r_\mu(\T).
\]
We assume that this boundary has positive curvature and hence is
strictly convex. Positive curvature is an additional hypothesis;
smallness in $C^1$ alone does not imply it.
The normalization theorem below is purely geometric and does not require
positive curvature or a caustic assumption; these conditions are imposed
only in its later billiard applications.
Figure~\ref{geo:fig:elliptic-graph-chart} illustrates the graph coordinates
and their physical image.

Fix a noncircular ellipse $E_*=E(a,b)$, $a>b>0$, and an anchor
rotation $\omega_0\in(0,1/2)$. We choose a nearby reference ellipse
so that the residual displacement has zero moments against its five
infinitesimal ellipse directions. The construction is geometric and
requires only a small $C^1$ norm. Its integer and Sobolev regularity
properties follow from the general coordinate estimates in
Appendix~\ref{app:sec:norms}.

\paragraph{Ellipse coordinates and tangent functions.}
Use the local parameters
\[
 t=(c_1,c_2,\alpha,\beta,\vartheta),\qquad
 a_t=a+\alpha,\qquad b_t=b+\beta,
\]
and, with $c=(c_1,c_2)$ and $R_\vartheta$ the counterclockwise rotation
through angle $\vartheta$, set
\[
 X_t(s,\psi)=c+R_\vartheta
 \big((a_t\cosh s+b_t\sinh s)\cos\psi,
      (b_t\cosh s+a_t\sinh s)\sin\psi\big).
\]
Then $E_t=X_t(0,\T)$ and $X_0=\mathcal X_{E_*}$. Restrict the parameter
neighborhood so that $a_t>b_t>0$ and these parameters form a
nonsingular local chart on the ellipse family. Put
\[
 N_t(\psi)=\partial_sX_t(0,\psi)
          =R_\vartheta(b_t\cos\psi,a_t\sin\psi),\qquad
 D_t(\psi)=|N_t(\psi)|^2
          =b_t^2\cos^2\psi+a_t^2\sin^2\psi.
\]
For $i=1,\ldots,5$, define
\[
 \tau_{t,i}(\psi)
 :=\frac{\langle\partial_{t_i}X_t(0,\psi),N_t(\psi)\rangle_{\R^2}}
          {D_t(\psi)},
\]
where $\partial_{t_i}$ is taken at fixed $s=0$ and $\psi$.
These are the infinitesimal graph displacements of
$E_{t+\epsilon e_i}$ over $E_t$, where $e_i$ is the $i$th coordinate
vector in parameter space. Indeed, $N_t$ and $\partial_\psi X_t(0,\psi)$
are orthogonal, and the unique decomposition
\[
 \partial_{t_i}X_t(0,\psi)
 =\tau_{t,i}(\psi)N_t(\psi)
   +\chi_{t,i}(\psi)\partial_\psi X_t(0,\psi)
\]
separates graph displacement from tangential reparametrization.
Explicitly,
\begin{align*}
 \tau_{t,1}
 &=\frac{b_t\cos\vartheta\cos\psi-a_t\sin\vartheta\sin\psi}{D_t},
 &\tau_{t,2}
 &=\frac{b_t\sin\vartheta\cos\psi+a_t\cos\vartheta\sin\psi}{D_t},\\
 \tau_{t,3}&=\frac{b_t\cos^2\psi}{D_t},
 &\tau_{t,4}&=\frac{a_t\sin^2\psi}{D_t},\\
 \tau_{t,5}&=\frac{(a_t^2-b_t^2)\sin\psi\cos\psi}{D_t}.
\end{align*}
The rotation factors in the first two formulas occur because $c_1,c_2$
are translations in fixed ambient directions. At $t=0$, write
$\tau_i=\tau_{0,i}$ and $D_E=D_0$; thus
\begin{align}\label{geo:eq:tangents}
\tau_1&=\frac{b\cos\psi}{D_E},&
\tau_2&=\frac{a\sin\psi}{D_E},&
\tau_3&=\frac{b\cos^2\psi}{D_E},\nonumber\\
\tau_4&=\frac{a\sin^2\psi}{D_E},&
\tau_5&=\frac{(a^2-b^2)\sin\psi\cos\psi}{D_E}.
\end{align}
Multiplication by $D_E$ and the double-angle identities show independence,
since $a>b$. In particular,
\begin{equation}\label{geo:eq:tangentspace}
 T_E=\frac1{D_E}\operatorname{span}_{\R}
 \{1,\cos\psi,\sin\psi,\cos2\psi,\sin2\psi\}.
\end{equation}
The functions $\tau_{t,i}$ depend smoothly on $t$ and remain independent
in a sufficiently small parameter neighborhood.

\paragraph{Action-angle moments.}
Let $\Phi_t$ be the action-angle parametrization of $E_t$ at $\omega_0$,
and let $G_t=\Phi_t^{-1}$. Section~\ref{geo:sec:convention} supplies their
smooth parameter dependence and locally uniform finite-order bounds.
For real periodic functions write
\[
 \langle f,g\rangle_0=\frac1{2\pi}\int_0^{2\pi}f(\xi)g(\xi)\dd\xi,
 \qquad \widetilde\tau_{t,i}=\tau_{t,i}\circ\Phi_t.
\]
Define the positive-definite Gram matrix and its dual weights by
\[
 (\mathcal G_t)_{ij}
 =\langle\widetilde\tau_{t,i},\widetilde\tau_{t,j}\rangle_0,
 \qquad
 \eta_{t,i}=\sum_{j=1}^5(\mathcal G_t^{-1})_{ij}\widetilde\tau_{t,j}.
\]
They vary smoothly with $t$ and satisfy
$\langle\widetilde\tau_{t,k},\eta_{t,i}\rangle_0=\delta_{ik}$.
In particular, the dual weights $\eta_{t,i}$ are distinct from the
tangent functions $\tau_{t,i}$.

\begin{theorem}[Exact ellipse normalization and regularity]\label{geo:thm:gauge}
There are a parameter neighborhood $\mathcal V$ of zero and constants
$\varepsilon,C>0$, depending on $E_*,\omega_0$, with the following
property. For every real $u\in C^2(\T)$ with
$\|u\|_{C^1}<\varepsilon$, there is a unique $t(u)\in\mathcal V$
such that the curve
\[
 \Gamma_u=\{X_0(u(\varphi),\varphi):\varphi\in\T\},
\]
written as an elliptic graph $v$ over $E_{t(u)}$, has action-coordinate
representative $U=v\circ\Phi_{t(u)}$ satisfying
\begin{equation}\label{geo:eq:gauge}
 \Pi_{E_{t(u)},i}U
 :=\frac1{2\pi}\int_0^{2\pi}U(\xi)\eta_{t(u),i}(\xi)\dd\xi=0,
 \qquad 1\le i\le5.
\end{equation}
The same selection has the following estimates.
\begin{enumerate}[label=(\roman*),ref=\roman*,leftmargin=*]
\item\label{geo:gauge:low} \emph{Parameter and low-norm bounds.} One has
\begin{equation}\label{geo:eq:gaugebounds}
 |t(u)|\le C\|u\|_{C^0},\qquad
 \|v\|_{C^1}+\|U\|_{C^1}\le C\|u\|_{C^1}.
\end{equation}
Moreover, $C^{-1}\|v\|_{C^1}\le\|U\|_{C^1}\le C\|v\|_{C^1}$.
\item\label{geo:gauge:integer} \emph{Integer regularity.}
For every integer $k\ge2$, if $u\in C^k(\T)$, then $v,U\in C^k(\T)$ and
\[
 \|v\|_{C^k}+\|U\|_{C^k}
 \le C_k\bigl(\|u\|_{C^k}+|t(u)|\bigr).
\]
The constant $C_k$ is independent of the high norm $\|u\|_{C^k}$.
\item\label{geo:gauge:sobolev} \emph{Sobolev regularity.}
For every real $s>11/2$ and every $M>0$, if $u\in H^s(\T)$ and
$\|u\|_{H^s}\le M$, then $v,U\in H^s(\T)$ and
\[
 \|v\|_{H^s}+\|U\|_{H^s}\le M_1,
\]
where $M_1$ may depend on $E_*,\omega_0,s,M$. No linear dependence
of $M_1$ on $M$ is asserted.
\end{enumerate}
The geometric neighborhood and the $C^1$ threshold can be chosen
independently of $k,s,M$. For fixed indices and norm bounds, all
constants are locally uniform when the reference ellipse varies in a
sufficiently small compact neighborhood of a fixed noncircular ellipse.
\end{theorem}

\begin{proof}
We construct the normalization from scalar moments on $C^0$ and then
apply the general norm estimates in Appendix~\ref{app:sec:norms}.

\smallskip\noindent
\emph{Step 1: represent the same curve over each nearby ellipse.}
The maps $X_t$ have smooth inverses $(S_t,\Psi_t)$ on a common
annular neighborhood, where $S_t$ is real valued and $\Psi_t$ takes
values in $\T$. For $z$ in a fixed small interval about zero, let
$A(t,\varphi,z)$ be the smooth real lift of
$\Psi_t(X_0(z,\varphi))$, chosen continuously from
$A(0,\varphi,z)=\varphi$, and set
\[
 B(t,\varphi,z)=S_t(X_0(z,\varphi)),\qquad
 A(t,\varphi+2\pi,z)=A(t,\varphi,z)+2\pi.
\]
The function $B$ is periodic in $\varphi$ and satisfies
$B(0,\varphi,z)=z$. For a $C^1$-small input define
\[
 h_{t,u}(\varphi)=A(t,\varphi,u(\varphi)),\qquad
 r_{t,u}(\varphi)=B(t,\varphi,u(\varphi)).
\]
Then
\[
 X_t(r_{t,u}(\varphi),h_{t,u}(\varphi))=X_0(u(\varphi),\varphi).
\]
Smoothness and the identities at $t=0$ give
\[
 \|h_{t,u}-\mathrm{id}\|_{C^1}\le C|t|(1+\|u\|_{C^1}),\qquad
 \|r_{t,u}\|_{C^1}\le C(|t|+\|u\|_{C^1}).
\]
Choose the parameter neighborhood and low-norm threshold so that
$h_{t,u}'\ge1/2$. This degree-one lift is an increasing circle
diffeomorphism, and the residual graph over $E_t$ is
\[
 v_{t,u}=r_{t,u}\circ h_{t,u}^{-1},\qquad
 \|v_{t,u}\|_{C^1}\le C(|t|+\|u\|_{C^1}).
\]
Indeed, $(v_{t,u}'\circ h_{t,u})=r_{t,u}'/h_{t,u}'$.
No high norm enters these choices.
Figure~\ref{geo:fig:exact-normalization} illustrates this representation
of the same physical curve; the parameter $t$ is selected in Step~3.

\begin{figure}[htbp]
\centering
\begin{tikzpicture}[
 x=1cm,y=1cm,font=\small,
 every node/.style={inner sep=2pt},
 boundary/.style={line width=.9pt},
 reference/.style={densely dashed,line width=.65pt},
 fibre/.style={line width=.65pt},
 arrow/.style={-{Stealth[length=1.7mm]},line width=.6pt}]
\def\physicalcurve{
  (2.765971101,0.128820761)
  (2.763409662,0.157470786)
  (2.760035952,0.186050397)
  (2.755849889,0.214550105)
  (2.750850900,0.242960167)
  (2.745037955,0.271270489)
  (2.738409605,0.299470532)
  (2.730964042,0.327549223)
  (2.722699153,0.355494875)
  (2.713612590,0.383295115)
  (2.703701843,0.410936826)
  (2.692964320,0.438406095)
  (2.681397429,0.465688175)
  (2.668998670,0.492767462)
  (2.655765725,0.519627480)
  (2.641696556,0.546250888)
  (2.626789498,0.572619492)
  (2.611043368,0.598714279)
  (2.594457557,0.624515470)
  (2.577032140,0.650002580)
  (2.558767971,0.675154504)
  (2.539666785,0.699949616)
  (2.519731293,0.724365882)
  (2.498965275,0.748381002)
  (2.477373664,0.771972551)
  (2.454962623,0.795118150)
  (2.431739615,0.817795644)
  (2.407713460,0.839983294)
  (2.382894380,0.861659985)
  (2.357294031,0.882805432)
  (2.330925518,0.903400407)
  (2.303803395,0.923426959)
  (2.275943648,0.942868637)
  (2.247363658,0.961710718)
  (2.218082151,0.979940424)
  (2.188119123,0.997547131)
  (2.157495751,1.014522568)
  (2.126234288,1.030861002)
  (2.094357941,1.046559407)
  (2.061890729,1.061617605)
  (2.028857338,1.076038393)
  (1.995282957,1.089827639)
  (1.961193109,1.102994350)
  (1.926613471,1.115550712)
  (1.891569699,1.127512094)
  (1.856087240,1.138897028)
  (1.820191155,1.149727147)
  (1.783905940,1.160027094)
  (1.747255359,1.169824398)
  (1.710262278,1.179149315)
  (1.672948518,1.188034638)
  (1.635334714,1.196515478)
  (1.597440197,1.204629012)
  (1.559282879,1.212414213)
  (1.520879169,1.219911542)
  (1.482243896,1.227162631)
  (1.443390255,1.234209940)
  (1.404329773,1.241096395)
  (1.365072283,1.247865023)
  (1.325625927,1.254558568)
  (1.285997169,1.261219100)
  (1.246190826,1.267887629)
  (1.206210112,1.274603710)
  (1.166056697,1.281405058)
  (1.125730780,1.288327168)
  (1.085231172,1.295402942)
  (1.044555389,1.302662336)
  (1.003699757,1.310132019)
  (0.962659526,1.317835055)
  (0.921428987,1.325790607)
  (0.880001609,1.334013669)
  (0.838370165,1.342514820)
  (0.796526884,1.351300023)
  (0.754463597,1.360370439)
  (0.712171895,1.369722300)
  (0.669643291,1.379346795)
  (0.626869390,1.389230019)
  (0.583842065,1.399352945)
  (0.540553633,1.409691453)
  (0.496997045,1.420216390)
  (0.453166073,1.430893685)
  (0.409055503,1.441684499)
  (0.364661323,1.452545430)
  (0.319980924,1.463428753)
  (0.275013283,1.474282705)
  (0.229759152,1.485051821)
  (0.184221237,1.495677301)
  (0.138404362,1.506097416)
  (0.092315628,1.516247960)
  (0.045964549,1.526062722)
  (-0.000636829,1.535473997)
  (-0.047473833,1.544413122)
  (-0.094529096,1.552811023)
  (-0.141782513,1.560598794)
  (-0.189211234,1.567708270)
  (-0.236789693,1.574072613)
  (-0.284489663,1.579626898)
  (-0.332280357,1.584308680)
  (-0.380128561,1.588058556)
  (-0.427998804,1.590820698)
  (-0.475853564,1.592543358)
  (-0.523653505,1.593179337)
  (-0.571357752,1.592686413)
  (-0.618924180,1.591027718)
  (-0.666309740,1.588172060)
  (-0.713470795,1.584094191)
  (-0.760363470,1.578775009)
  (-0.806944015,1.572201694)
  (-0.853169163,1.564367780)
  (-0.898996487,1.555273156)
  (-0.944384751,1.544923998)
  (-0.989294235,1.533332638)
  (-1.033687049,1.520517366)
  (-1.077527418,1.506502170)
  (-1.120781930,1.491316423)
  (-1.163419760,1.474994517)
  (-1.205412852,1.457575449)
  (-1.246736056,1.439102370)
  (-1.287367232,1.419622104)
  (-1.327287311,1.399184642)
  (-1.366480310,1.377842612)
  (-1.404933313,1.355650750)
  (-1.442636408,1.332665362)
  (-1.479582589,1.308943790)
  (-1.515767632,1.284543895)
  (-1.551189927,1.259523550)
  (-1.585850302,1.233940164)
  (-1.619751807,1.207850228)
  (-1.652899496,1.181308895)
  (-1.685300187,1.154369599)
  (-1.716962216,1.127083710)
  (-1.747895183,1.099500228)
  (-1.778109703,1.071665524)
  (-1.807617157,1.043623122)
  (-1.836429445,1.015413525)
  (-1.864558753,0.987074080)
  (-1.892017335,0.958638889)
  (-1.918817298,0.930138761)
  (-1.944970412,0.901601200)
  (-1.970487938,0.873050428)
  (-1.995380463,0.844507443)
  (-2.019657769,0.815990113)
  (-2.043328711,0.787513289)
  (-2.066401117,0.759088947)
  (-2.088881709,0.730726354)
  (-2.110776039,0.702432248)
  (-2.132088451,0.674211036)
  (-2.152822048,0.646065003)
  (-2.172978688,0.617994528)
  (-2.192558988,0.589998310)
  (-2.211562347,0.562073590)
  (-2.229986974,0.534216379)
  (-2.247829941,0.506421683)
  (-2.265087233,0.478683716)
  (-2.281753820,0.450996114)
  (-2.297823726,0.423352136)
  (-2.313290115,0.395744852)
  (-2.328145374,0.368167325)
  (-2.342381209,0.340612771)
  (-2.355988738,0.313074709)
  (-2.368958589,0.285547097)
  (-2.381281004,0.258024447)
  (-2.392945932,0.230501929)
  (-2.403943136,0.202975455)
  (-2.414262289,0.175441748)
  (-2.423893072,0.147898396)
  (-2.432825272,0.120343887)
  (-2.441048873,0.092777634)
  (-2.448554146,0.065199982)
  (-2.455331738,0.037612199)
  (-2.461372754,0.010016460)
  (-2.466668834,-0.017584187)
  (-2.471212233,-0.045185864)
  (-2.474995885,-0.072783920)
  (-2.478013470,-0.100372997)
  (-2.480259472,-0.127947102)
  (-2.481729227,-0.155499696)
  (-2.482418973,-0.183023781)
  (-2.482325879,-0.210512003)
  (-2.481448077,-0.237956751)
  (-2.479784674,-0.265350270)
  (-2.477335762,-0.292684764)
  (-2.474102416,-0.319952515)
  (-2.470086676,-0.347145986)
  (-2.465291527,-0.374257931)
  (-2.459720867,-0.401281501)
  (-2.453379464,-0.428210340)
  (-2.446272903,-0.455038679)
  (-2.438407535,-0.481761417)
  (-2.429790410,-0.508374200)
  (-2.420429210,-0.534873484)
  (-2.410332173,-0.561256595)
  (-2.399508019,-0.587521766)
  (-2.387965869,-0.613668178)
  (-2.375715161,-0.639695972)
  (-2.362765566,-0.665606262)
  (-2.349126899,-0.691401126)
  (-2.334809033,-0.717083584)
  (-2.319821811,-0.742657563)
  (-2.304174954,-0.768127848)
  (-2.287877973,-0.793500010)
  (-2.270940084,-0.818780326)
  (-2.253370116,-0.843975682)
  (-2.235176433,-0.869093451)
  (-2.216366849,-0.894141370)
  (-2.196948550,-0.919127385)
  (-2.176928027,-0.944059495)
  (-2.156311005,-0.968945571)
  (-2.135102383,-0.993793173)
  (-2.113306186,-1.018609345)
  (-2.090925518,-1.043400407)
  (-2.067962529,-1.068171737)
  (-2.044418391,-1.092927544)
  (-2.020293289,-1.117670647)
  (-1.995586418,-1.142402241)
  (-1.970295996,-1.167121672)
  (-1.944419297,-1.191826218)
  (-1.917952685,-1.216510872)
  (-1.890891677,-1.241168142)
  (-1.863231011,-1.265787861)
  (-1.834964736,-1.290357016)
  (-1.806086317,-1.314859601)
  (-1.776588752,-1.339276488)
  (-1.746464709,-1.363585330)
  (-1.715706675,-1.387760496)
  (-1.684307122,-1.411773031)
  (-1.652258681,-1.435590664)
  (-1.619554335,-1.459177843)
  (-1.586187614,-1.482495820)
  (-1.552152803,-1.505502768)
  (-1.517445156,-1.528153953)
  (-1.482061106,-1.550401937)
  (-1.445998481,-1.572196835)
  (-1.409256717,-1.593486608)
  (-1.371837059,-1.614217405)
  (-1.333742749,-1.634333936)
  (-1.294979214,-1.653779887)
  (-1.255554216,-1.672498373)
  (-1.215477993,-1.690432408)
  (-1.174763370,-1.707525415)
  (-1.133425842,-1.723721741)
  (-1.091483625,-1.738967195)
  (-1.048957669,-1.753209584)
  (-1.005871643,-1.766399260)
  (-0.962251872,-1.778489644)
  (-0.918127245,-1.789437746)
  (-0.873529076,-1.799204657)
  (-0.828490935,-1.807756011)
  (-0.783048442,-1.815062408)
  (-0.737239021,-1.821099800)
  (-0.691101631,-1.825849812)
  (-0.644676470,-1.829300023)
  (-0.598004653,-1.831444174)
  (-0.551127878,-1.832282322)
  (-0.504088079,-1.831820915)
  (-0.456927069,-1.830072815)
  (-0.409686194,-1.827057234)
  (-0.362405978,-1.822799617)
  (-0.315125797,-1.817331449)
  (-0.267883556,-1.810690004)
  (-0.220715396,-1.802918032)
  (-0.173655425,-1.794063390)
  (-0.126735482,-1.784178624)
  (-0.079984937,-1.773320509)
  (-0.033430517,-1.761549546)
  (0.012903811,-1.748929435)
  (0.058996940,-1.735526520)
  (0.104830661,-1.721409217)
  (0.150389683,-1.706647432)
  (0.195661599,-1.691311979)
  (0.240636829,-1.675473997)
  (0.285308520,-1.659204377)
  (0.329672426,-1.642573206)
  (0.373726755,-1.625649228)
  (0.417471993,-1.608499331)
  (0.460910716,-1.591188063)
  (0.504047374,-1.573777179)
  (0.546888083,-1.556325225)
  (0.589440386,-1.538887155)
  (0.631713029,-1.521513997)
  (0.673715727,-1.504252551)
  (0.715458929,-1.487145132)
  (0.756953596,-1.470229364)
  (0.798210978,-1.453538003)
  (0.839242402,-1.437098823)
  (0.880059068,-1.420934526)
  (0.920671855,-1.405062708)
  (0.961091142,-1.389495866)
  (1.001326637,-1.374241441)
  (1.041387219,-1.359301906)
  (1.081280790,-1.344674888)
  (1.121014145,-1.330353332)
  (1.160592846,-1.316325697)
  (1.200021113,-1.302576182)
  (1.239301724,-1.289084987)
  (1.278435932,-1.275828595)
  (1.317423383,-1.262780089)
  (1.356262056,-1.249909477)
  (1.394948211,-1.237184049)
  (1.433476346,-1.224568739)
  (1.471839174,-1.212026505)
  (1.510027606,-1.199518716)
  (1.548030751,-1.187005543)
  (1.585835931,-1.174446353)
  (1.623428707,-1.161800100)
  (1.660792926,-1.149025710)
  (1.697910778,-1.136082459)
  (1.734762869,-1.122930338)
  (1.771328315,-1.109530399)
  (1.807584847,-1.095845084)
  (1.843508927,-1.081838539)
  (1.879075887,-1.067476890)
  (1.914260078,-1.052728505)
  (1.949035022,-1.037564219)
  (1.983373587,-1.021957535)
  (2.017248166,-1.005884782)
  (2.050630852,-0.989325257)
  (2.083493634,-0.972261312)
  (2.115808580,-0.954678429)
  (2.147548027,-0.936565242)
  (2.178684765,-0.917913541)
  (2.209192210,-0.898718235)
  (2.239044582,-0.878977285)
  (2.268217052,-0.858691615)
  (2.296685897,-0.837864989)
  (2.324428624,-0.816503868)
  (2.351424086,-0.794617247)
  (2.377652580,-0.772216467)
  (2.403095922,-0.749315021)
  (2.427737506,-0.725928344)
  (2.451562347,-0.702073590)
  (2.474557095,-0.677769413)
  (2.496710044,-0.653035737)
  (2.518011108,-0.627893536)
  (2.538451794,-0.602364605)
  (2.558025150,-0.576471352)
  (2.576725705,-0.550236589)
  (2.594549388,-0.523683338)
  (2.611493452,-0.496834647)
  (2.627556368,-0.469713426)
  (2.642737733,-0.442342291)
  (2.657038159,-0.414743433)
  (2.670459161,-0.386938494)
  (2.683003049,-0.358948470)
  (2.694672811,-0.330793626)
  (2.705472003,-0.302493431)
  (2.715404639,-0.274066505)
  (2.724475081,-0.245530589)
  (2.732687937,-0.216902526)
  (2.740047961,-0.188198259)
  (2.746559958,-0.159432846)
  (2.752228692,-0.130620478)
  (2.757058806,-0.101774525)
  (2.761054743,-0.072907584)
  (2.764220674,-0.044031535)
  (2.766560432,-0.015157619)
  (2.768077464,0.013703484)
  (2.768774772,0.042541571)
  (2.768654887,0.071346815)
  (2.767719831,0.100109669)
  (2.765971101,0.128820761)
}
\def\oldfibre{
  (-1.064078572,-1.499479323)
  (-1.065695748,-1.504916576)
  (-1.067318781,-1.510362102)
  (-1.068947681,-1.515815929)
  (-1.070582456,-1.521278088)
  (-1.072223116,-1.526748609)
  (-1.073869671,-1.532227523)
  (-1.075522127,-1.537714859)
  (-1.077180496,-1.543210647)
  (-1.078844786,-1.548714918)
  (-1.080515006,-1.554227702)
  (-1.082191165,-1.559749029)
  (-1.083873273,-1.565278929)
  (-1.085561338,-1.570817434)
  (-1.087255371,-1.576364573)
  (-1.088955380,-1.581920377)
  (-1.090661375,-1.587484876)
  (-1.092373365,-1.593058101)
  (-1.094091359,-1.598640083)
  (-1.095815367,-1.604230853)
  (-1.097545399,-1.609830440)
  (-1.099281464,-1.615438876)
  (-1.101023571,-1.621056192)
  (-1.102771730,-1.626682418)
  (-1.104525951,-1.632317586)
  (-1.106286243,-1.637961727)
  (-1.108052616,-1.643614870)
  (-1.109825080,-1.649277049)
  (-1.111603644,-1.654948293)
  (-1.113388319,-1.660628634)
  (-1.115179113,-1.666318103)
  (-1.116976038,-1.672016731)
  (-1.118779102,-1.677724550)
  (-1.120588316,-1.683441591)
  (-1.122403689,-1.689167885)
  (-1.124225232,-1.694903465)
  (-1.126052955,-1.700648361)
  (-1.127886867,-1.706402605)
  (-1.129726979,-1.712166228)
  (-1.131573301,-1.717939263)
  (-1.133425842,-1.723721741)
}
\def\newfibre{
  (-1.086828738,-1.558913832)
  (-1.087721358,-1.562203539)
  (-1.088616380,-1.565496214)
  (-1.089513806,-1.568791864)
  (-1.090413638,-1.572090495)
  (-1.091315877,-1.575392113)
  (-1.092220526,-1.578696726)
  (-1.093127585,-1.582004339)
  (-1.094037058,-1.585314960)
  (-1.094948945,-1.588628594)
  (-1.095863249,-1.591945250)
  (-1.096779971,-1.595264932)
  (-1.097699113,-1.598587648)
  (-1.098620678,-1.601913404)
  (-1.099544666,-1.605242207)
  (-1.100471079,-1.608574064)
  (-1.101399920,-1.611908981)
  (-1.102331191,-1.615246965)
  (-1.103264892,-1.618588022)
  (-1.104201027,-1.621932159)
  (-1.105139597,-1.625279384)
  (-1.106080603,-1.628629701)
  (-1.107024048,-1.631983119)
  (-1.107969934,-1.635339643)
  (-1.108918261,-1.638699281)
  (-1.109869034,-1.642062039)
  (-1.110822252,-1.645427923)
  (-1.111777918,-1.648796941)
  (-1.112736035,-1.652169100)
  (-1.113696603,-1.655544405)
  (-1.114659625,-1.658922864)
  (-1.115625103,-1.662304483)
  (-1.116593039,-1.665689269)
  (-1.117563434,-1.669077229)
  (-1.118536290,-1.672468369)
  (-1.119511610,-1.675862697)
  (-1.120489395,-1.679260219)
  (-1.121469648,-1.682660941)
  (-1.122452370,-1.686064872)
  (-1.123437563,-1.689472016)
  (-1.124425229,-1.692882382)
  (-1.125415370,-1.696295975)
  (-1.126407988,-1.699712803)
  (-1.127403086,-1.703132872)
  (-1.128400664,-1.706556190)
  (-1.129400726,-1.709982763)
  (-1.130403273,-1.713412598)
  (-1.131408306,-1.716845701)
  (-1.132415829,-1.720282080)
  (-1.133425842,-1.723721741)
}
\begin{scope}[shift={(3.1,0)},scale=.90]
 \node[anchor=west] at (-3.12,2.36) {(a) Original reference};
 \draw[reference] (0,0) ellipse[x radius=2.55,y radius=1.65];
 \draw[boundary] plot coordinates {\physicalcurve} --cycle;
 \draw[fibre] plot coordinates {\oldfibre};
 \coordinate (Pold) at (-1.133425842,-1.723721741);
 \coordinate (Bold) at (-1.064078572,-1.499479323);
 \fill (Pold) circle[radius=1.6pt];
 \fill (Bold) circle[radius=1.3pt];
 \node[anchor=north east] at ($(Pold)+(-.10,-.08)$) {$P$};
 \node[anchor=south west] at ($(Bold)+(.08,.08)$) {$B_0$};
 \node at (.05,.18) {$E_*=X_0(0,\mathbb T)$};
 \draw[boundary] (1.05,1.91)--(1.45,1.91) node[right] {$\Gamma_u$};
\end{scope}
\begin{scope}[shift={(10.52,0)},scale=.90]
 \node[anchor=west] at (-3.12,2.36) {(b) Moment-normalized reference};
 \draw[reference,shift={(.12,-.07)},rotate=4.297183463481]
   (0,0) ellipse[x radius=2.63,y radius=1.61];
 \draw[boundary] plot coordinates {\physicalcurve} --cycle;
 \draw[fibre] plot coordinates {\newfibre};
 \coordinate (Pnew) at (-1.133425842,-1.723721741);
 \coordinate (Bnew) at (-1.086828738,-1.558913832);
 \fill (Pnew) circle[radius=1.6pt];
 \fill (Bnew) circle[radius=1.3pt];
 \node[anchor=north east] at ($(Pnew)+(-.10,-.08)$) {$P$};
 \node[anchor=south west] at ($(Bnew)+(.08,.08)$) {$B_t$};
 \node at (.05,.18) {$E_t=X_t(0,\mathbb T)$};
 \draw[boundary] (1.05,1.91)--(1.45,1.91) node[right] {$\Gamma_u$};
\end{scope}
\node (oldpair) at (3.10,-2.19) {$(u(\varphi),\varphi)$};
\node (newpair) at (10.52,-2.19) {$(v_{t,u}(\psi),\psi)$};
\draw[arrow] (oldpair)--node[midway,above=2pt]
 {$\psi=h_{t,u}(\varphi)$}(newpair);
\node (samepoint) at (6.81,-3.23) {$P\in\Gamma_u$};
\draw[arrow] (oldpair.south)--node[pos=.45,below left] {$X_0$}(samepoint.west);
\draw[arrow] (newpair.south)--node[pos=.45,below right] {$X_t$}(samepoint.east);
\node[align=center] at (6.81,-4.13)
 {$t=t(u):\qquad U=v_{t,u}\circ\Phi_t,\qquad
   \langle U,\eta_{t,i}\rangle_0=0\quad(1\le i\le5).$};
\end{tikzpicture}
\caption{The same physical boundary in two elliptic graph charts.
The solid curve $\Gamma_u$ and the marked physical point $P$ are
identical in both panels; the dashed reference ellipse changes from
$E_*$ to $E_t$. The fibre endpoints are
$B_0=X_0(0,\varphi)$ and $B_t=X_t(0,\psi)$, with
$\psi=h_{t,u}(\varphi)$. The diagram records equality of the two
physical representations, not equality of the angular parameters or
of the displacements. At $t=t(u)$, the residual has the five exact
zero moments in \eqref{geo:eq:gauge}; this is not a nearest-point or
least-distance projection onto the ellipse family.}
\label{geo:fig:exact-normalization}
\end{figure}
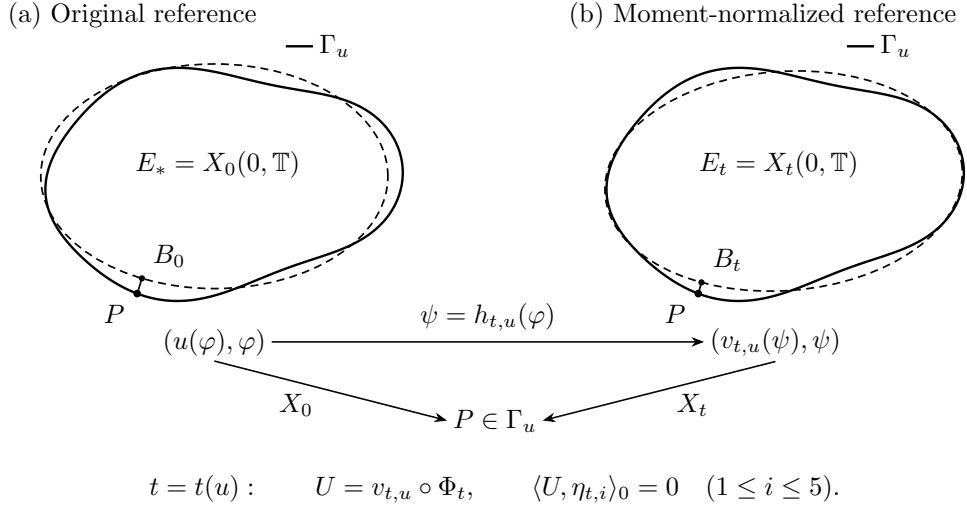

\smallskip\noindent
\emph{Step 2: define the moments and extend them smoothly to $C^0$.}
Initially for $C^1$-small $u$, set
\begin{equation}\label{geo:eq:Fmoment}
 F_i(t,u)=\frac1{2\pi}\int_0^{2\pi}
 v_{t,u}(\psi)\tau_{t,i}(\psi)G_t'(\psi)\dd\psi.
\end{equation}
With $U_{t,u}=v_{t,u}\circ\Phi_t$, changing variables
$\psi=\Phi_t(\xi)$ gives
$F_i(t,u)=\langle U_{t,u},\widetilde\tau_{t,i}\rangle_0$.
Thus
\[
 \bigl(\langle U_{t,u},\eta_{t,i}\rangle_0\bigr)_{i=1}^5
 =\mathcal G_t^{-1}F(t,u).
\]
The equations $F(t,u)=0$ are exactly the desired dual-moment equations.

Their integral form avoids assuming differentiability of the graph
reparametrization as a map from $C^1$ to $C^1$. Put
$w_{t,i}=\tau_{t,i}G_t'$ and change variables
$\psi=h_{t,u}(\varphi)$ in \eqref{geo:eq:Fmoment}. Periodicity permits
integration over any full angular period. The integrand becomes
\[
 B(t,\varphi,u)w_{t,i}(A(t,\varphi,u))
 \bigl[A_\varphi(t,\varphi,u)+A_z(t,\varphi,u)u'\bigr].
\]
Define the smooth periodic coefficients
\begin{align*}
 a_i(t,\varphi,z)
 &=B(t,\varphi,z)w_{t,i}(A(t,\varphi,z))A_\varphi(t,\varphi,z),\\
 b_i(t,\varphi,z)
 &=B(t,\varphi,z)w_{t,i}(A(t,\varphi,z))A_z(t,\varphi,z),
\end{align*}
and let $H_i(t,\varphi,z)=\int_0^z b_i(t,\varphi,\zeta)\dd\zeta$.
The chain rule and periodic integration give
\[
 \int_0^{2\pi}b_i(t,\varphi,u(\varphi))u'(\varphi)\dd\varphi
 =-\int_0^{2\pi}\partial_\varphi H_i(t,\varphi,u(\varphi))\dd\varphi,
\]
where $\partial_\varphi H_i$ holds $z$ fixed. Hence
\begin{equation}\label{geo:eq:C0extension}
 F_i(t,u)=\frac1{2\pi}\int_0^{2\pi}
 \mathcal C_i(t,\varphi,u(\varphi))\dd\varphi,
 \qquad \mathcal C_i=a_i-\partial_\varphi H_i.
\end{equation}
This formula contains no derivative of $u$ and defines a smooth map
on a product of a parameter neighborhood and a $C^0$ ball. Explicitly,
for $j\ge1$,
\[
 D_u^jF_i(t,u)[v_1,\ldots,v_j]
 =\frac1{2\pi}\int_0^{2\pi}
 \partial_z^j\mathcal C_i(t,\varphi,u(\varphi))
 \prod_{\ell=1}^j v_\ell(\varphi)\dd\varphi.
\]
Its absolute value is bounded by $C_j\prod_\ell\|v_\ell\|_{C^0}$
on a fixed small product neighborhood; mixed parameter derivatives
satisfy the same type of bound. This extends the scalar moments, not
the geometric graph construction to arbitrary $C^0$ inputs. For the
$C^1$-small inputs used here, both formulas for $F$ agree exactly.

\smallskip\noindent
\emph{Step 3: select the ellipse uniquely.}
Put $\mathcal G=\mathcal G_0$. At $(t,u)=(0,0)$ the residual vanishes.
To compute its parameter variation with the physical curve fixed,
let $k_t=h_{t,0}^{-1}$ and differentiate
\[
 X_t(v_{t,0}(\psi),\psi)=X_0(0,k_t(\psi)).
\]
At zero this yields
\[
 \partial_{t_j}X_t(0,\psi)\big|_0
 +N_0(\psi)\partial_{t_j}v_{t,0}(\psi)\big|_0
 =\partial_\psi X_0(0,\psi)\partial_{t_j}k_t(\psi)\big|_0.
\]
Projection onto $N_0$ gives $\partial_{t_j}v_{t,0}|_0=-\tau_j$.
Every derivative of a weight in \eqref{geo:eq:Fmoment} is multiplied
by $v_{0,0}=0$, so
\[
 D_tF(0,0)=-\mathcal G.
\]
For $z=(z_1,\ldots,z_5)\in\R^5$,
\[
 z^{\mathsf T}\mathcal Gz
 =\frac1{2\pi}\int_0^{2\pi}
 \left|\sum_{i=1}^5z_i\tau_i(\Phi_0(\xi))\right|^2\dd\xi>0
 \quad\text{if }z\ne0.
\]
Thus the Jacobian is invertible.

For a quantitative local selection define
$\mathscr T_u(t)=t+\mathcal G^{-1}F(t,u)$. Since
$D_t\mathscr T_0(0)=0$, choose a closed parameter ball
$\overline B_R(0)$ and a $C^0$ neighborhood on which
$\|D_t\mathscr T_u\|\le1/2$. At $t=0$,
\[
 F_i(0,u)=\frac1{2\pi}\int_0^{2\pi}
 u(\psi)\tau_i(\psi)G_0'(\psi)\dd\psi,
 \qquad |F(0,u)|\le C\|u\|_{C^0}.
\]
This identity also holds for continuous $u$ by
\eqref{geo:eq:C0extension}. Shrink its $C^0$ neighborhood so that
$|\mathscr T_u(0)|\le R/4$. Then $\mathscr T_u$ maps
$\overline B_R(0)$ into $B_{3R/4}(0)$ and is a contraction.
It has a unique fixed point $t(u)$ there, equivalently $F(t(u),u)=0$.
Set $\mathcal V=B_R(0)$. Moreover,
\[
 |t(u)|\le\tfrac12|t(u)|+C\|u\|_{C^0},
\]
which proves the parameter bound. The implicit-function theorem makes
$t(u)$ smooth on this $C^0$ neighborhood, with expansion
\[
 t(u)=\mathcal G^{-1}
 \bigl(\langle u\circ\Phi_0,\tau_i\circ\Phi_0\rangle_0\bigr)_{i=1}^5
 +O(\|u\|_{C^0}^2).
\]
This concerns only the parameter map, not differentiability of the
residual map in a fixed high-regularity space. Restricting to the
$C^1$-small inputs recovers the geometric residual from Step~1.
Step~2 gives the exact moments, while Step~1 and the parameter bound
give $\|v\|_{C^1}\le C\|u\|_{C^1}$.

\smallskip\noindent
\emph{Step 4: regularity and uniformity.}
The coordinate functions $A,B$ satisfy the periodicity and identity
hypotheses of Proposition~\ref{app:prop:graph-change}; their coefficient
bounds are uniform on a fixed small parameter and amplitude neighborhood.
Apply that proposition at the parameter $t=t(u)$ selected in Step~3.
Its integer estimate gives
\[
 \|v\|_{C^k}\le C_k(\|u\|_{C^k}+|t(u)|),\qquad k\ge2.
\]
If $u\in H^s$ with $s>11/2$, then $u\in C^2$ by the Fourier embedding
recorded in Appendix~\ref{app:sub:sobolev}. Thus exactly the same
selection applies. The Sobolev part of the proposition gives
$\|v\|_{H^s}\le C_{s,M}$ when $\|u\|_{H^s}\le M$; it requires no
smallness of this high norm.

Corollary~\ref{app:cor:smooth-composition}, applied to the uniformly
smooth family $\Phi_t$, gives the asserted $C^k$ and $H^s$ estimates
for $U=v\circ\Phi_t$. In particular,
\[
 \|U\|_{C^1}\le C\|v\|_{C^1},\qquad
 \|v\|_{C^1}\le C\|U\|_{C^1},
\]
since $v=U\circ G_t$ and both first derivatives are uniformly bounded.
This completes \eqref{geo:eq:gaugebounds} and all regularity claims.

The geometric choices in Steps~1--3 use only fixed coordinate bounds
and low norms, so the same neighborhood and threshold work for all
$k,s,M$. Finally, vary the reference ellipse in a sufficiently small
compact neighborhood, using smoothly chosen principal-axis frames.
The semiaxes remain positive and separated, the coordinate inverses
have uniform finite-order bounds, and the smallest eigenvalue of
$\mathcal G$ has a uniform positive lower bound. The action-coordinate
bounds at $\omega_0$ are uniform as well. The choices above and the
appendix estimates are therefore locally uniform as stated.
\end{proof}

\section{Calibration and uniform nonlinear caustic rows}\label{geo:sec:calibration}
We compare the reference and perturbed billiards in one fixed action
coordinate. Reference calibration supplies a quadratic lower bound for
excess action, while calibration of the perturbed caustic makes its orbit
an action minimizer. A global perturbation bound compares these actions
without assuming that their vertices are close. Only afterward do we
Taylor-expand in the boundary displacement at fixed reference vertices,
not in the orbit displacement. This distinction gives a smallness
threshold independent of the primitive denominator.

\subsection{Global calibration for a convex caustic}
Let $r:\R\to\R^2$ be a positively oriented, regular, strictly convex $C^2$ boundary parametrization, periodic with period $2\pi$. On the ordered edge strip
\[
\mathcal D=\{(x,y):0<y-x<2\pi\}
\]
define the negative-chord generating function
\[
h(x,y)=-|r(y)-r(x)|.
\]
Its mixed derivative is strictly negative. Indeed, if $\alpha,\beta\in(0,\pi)$ are the angles of the chord at its two ends, then
\begin{equation}\label{geo:eq:twist}
h_{12}(x,y)=-\frac{|r'(x)|\,|r'(y)|\sin\alpha\sin\beta}
 {|r(y)-r(x)|}<0.
\end{equation}
We use the canonical momenta $p=-h_1(x,y)$ at the initial point and $p'=h_2(x,y)$ at the terminal point.

\begin{lemma}[Global calibration]\label{geo:lem:calibration}
Suppose a caustic supplies a continuous periodic invariant momentum graph $p=\gamma(x)$, and the induced lift $f$ is an increasing degree-one homeomorphism with $x<f(x)<x+2\pi$. Let $H'=\gamma$. There is a constant $\kappa$ such that
\begin{equation}\label{geo:eq:calibration}
K(x,y):=h(x,y)-H(y)+H(x)-\kappa\ge0
\end{equation}
on $\mathcal D$, with equality precisely when $y=f(x)$. The inequality extends by continuity to its closed edge strip.
Set $\mathcal A=H(x+2\pi)-H(x)$. If integers $d\ge1$ and $p$ satisfy
$f^d(x)=x+2\pi p$ for every $x\in\R$, then every $d$-edge orbit on this
graph has negative perimeter $d\kappa+p\mathcal A$, independent of its
starting point, and minimizes negative perimeter among all $d$-edge
configurations with the same lifted endpoints and successive increments
in $[0,2\pi]$.
\end{lemma}

\begin{proof}
This is the invariant-graph calibration principle; compare
\cite[Appendix~D, Lemma~D.5]{KKZ}. We give the finite-strip argument to
include the grazing endpoints and the possibly nonperiodic primitive.
Invariance implies
\[
 h_1(x,f(x))=-\gamma(x),\qquad h_2(x,f(x))=\gamma(f(x)).
\]
Thus the $C^1$ function $\mathscr C(x,y)=h(x,y)-H(y)+H(x)$ has zero
gradient on the graph of $f$. On each compact interval that graph is
rectifiable: its length is at most the total increase of $x$ plus that
of $f(x)$. The chain rule along an arclength parametrization makes
$\mathscr C$ constant there. Overlapping intervals give a single
constant $\kappa$, so $K(x,f(x))=0$. This uses no differentiability of $f$.

Put $g=f^{-1}$. For fixed $y$, invariance at its predecessor gives
\[
 K_2(x,y)=h_2(x,y)-h_2(g(y),y)
        =\int_{g(y)}^x h_{12}(s,y)\,ds.
\]
Both endpoints of integration lie in $(y-2\pi,y)$. Since $h_{12}<0$
and $g$ is increasing, $K_2(x,y)<0$ for $x<y<f(x)$ and $K_2(x,y)>0$
for $f(x)<y<x+2\pi$. Hence $K(x,\cdot)$ has its unique minimum zero at
$f(x)$. Continuity extends the inequality to the closed strip.
Choosing an interior point on either side of $f(x)$ and using strict
monotonicity also gives
\[
 K(x,x)=-\kappa>0,\qquad
 K(x,x+2\pi)=-\kappa-\mathcal A>0,
 \qquad \mathcal A=\int_x^{x+2\pi}\gamma(s)\,ds.
\]
The last quantity is independent of $x$ by periodicity of $\gamma$.

For any competing configuration in the closed strip, telescoping yields
\[
 \sum_{j=0}^{d-1}h(z_j,z_{j+1})
 =d\kappa+H(z_d)-H(z_0)+\sum_{j=0}^{d-1}K(z_j,z_{j+1}).
\]
The last sum is nonnegative and is zero on a caustic orbit. If its
lifted endpoints differ by $2\pi p$, the primitive term is
$p\mathcal A$. This proves both the common negative perimeter
$d\kappa+p\mathcal A$ and its minimizing property, including degenerate
competing edges.
\end{proof}

\begin{remark}\label{geo:rem:caustic-graph}
The primitive $H$ is generally \emph{not periodic}: its drift contributes $p\mathcal A$ to the action of an orbit with winding number $p$. The lemma requires neither a quantitative bound on the perturbed caustic graph nor prior localization of the caustic.

The graph hypotheses also hold for nonsmooth convex caustics. Let $C$
be the compact convex body bounded by such a caustic strictly inside
the table. Since $C$ has nonempty interior, each boundary point $r(x)$
has exactly two supporting rays directed toward $C$. These rays vary
continuously: a limit of supporting rays still supports $C$, and the
side containing the interior of $C$ distinguishes the two choices.
Reflection interchanges the two supporting lines; see
\cite[p.~257 and Remark~6]{AB}. Thus the branch with $C$ on the left
of its directed ray is preserved by the billiard map and its inverse.
If $w(x)$ is its unit direction, its canonical momentum is
$\gamma(x)=-r'(x)\cdot w(x)$, which is continuous and periodic.
The billiard map restricts to a homeomorphism of this graph, so its
projection to the boundary is a circle homeomorphism. Its continuous
next-impact lift satisfies $x<f(x)<x+2\pi$ and
$f(x+2\pi)=f(x)+2\pi$; hence it has degree one and is increasing.
No differentiability or strict convexity of the caustic is required.
\end{remark}

\smallskip\noindent\emph{Common period.}\phantomsection\label{geo:par:common-period}
Let $F:\R\to\R$ be an increasing degree-one lift with reduced rotation
number $p/d$, and suppose every point of its circle map is periodic.
For each $x$, choose $q\ge1$ and $k\in\Z$ with
$F^q(x)=x+2\pi k$. Iterating this identity gives rotation number
$k/q=p/d$, so $q=d\ell$ and $k=p\ell$ for some $\ell\ge1$.
The increasing map $G=F^d-2\pi p$ then satisfies $G^\ell(x)=x$.
If $G(x)>x$ or $G(x)<x$, its iterates at $x$ would be strictly
monotone, a contradiction. Hence $F^d(x)=x+2\pi p$ for every $x$.

\subsection{Uniform comparison with the reference ellipse}
Use the action-angle coordinate of the reference caustic at $\rho\in J$, and put
\[
\alpha=2\pi\rho,\quad R_{0,\rho}(x)=r_0(\Phi_\rho(x)),\quad
R_{\mu,\rho}(x)=r_\mu(\Phi_\rho(x)),
\quad h_{\mu,\rho}(x,y)=-|R_{\mu,\rho}(y)-R_{\mu,\rho}(x)|.
\]
The reference caustic map is $x\mapsto x+\alpha$. Denote its calibrated excess by $K_{0,\rho}$.

\begin{lemma}[Uniform global coercivity]\label{geo:lem:coercivity}
There is $c=c(E,J)>0$ such that
\begin{equation}\label{geo:eq:coercivity}
K_{0,\rho}(x,y)\ge c(y-x-\alpha)^2,
\qquad 0\le y-x\le2\pi,\quad \rho\in J.
\end{equation}
\end{lemma}
\begin{proof}
Recall that $J\Subset(0,1/2)$ is compact, and write
$\alpha_\rho=2\pi\rho$ for the action-coordinate step denoted by
$\alpha$ in the statement. This is not an incidence angle.

\smallskip\noindent
\emph{Step 1: regularity of the reference calibration.}
For the reference ellipse the boundary map is
$x\mapsto x+\alpha_\rho$, so its invariant momentum is
\[
 \gamma_\rho(y)=h_{0,\rho,2}(y-\alpha_\rho,y).
\]
This identity makes $\gamma_\rho$ smooth, with smooth parameter
dependence: the reference edges remain separated from both diagonal
endpoints, and the action-angle maps have the uniform finite-order
bounds of Section~\ref{geo:sec:convention}.
Although Lemma~\ref{geo:lem:calibration} assumes only a continuous
graph, differentiating $\gamma_\rho$ here is justified by this
explicit reference identity.
For definiteness, take
\[
 H_\rho(x)=\int_0^x\gamma_\rho(s)\,ds,\qquad
 \kappa_\rho=h_{0,\rho}(0,\alpha_\rho)-H_\rho(\alpha_\rho).
\]
These choices define the same calibrated excess $K_{0,\rho}$ as
Lemma~\ref{geo:lem:calibration}, smoothly in the interior of the edge
strip and continuously on its closure, jointly with $\rho$.
The drift of $H_\rho$ cancels under simultaneous translation of both
endpoints, giving
\[
 K_{0,\rho}(x+2\pi,y+2\pi)=K_{0,\rho}(x,y).
\]

\smallskip\noindent
\emph{Step 2: a uniform quadratic bound near the reference edge.}
Holding $\rho$ fixed and differentiating the momentum identity yields
\[
 \gamma_\rho'(y)
 =h_{0,\rho,12}(y-\alpha_\rho,y)
    +h_{0,\rho,22}(y-\alpha_\rho,y).
\]
Since $\partial_{22}K_{0,\rho}=h_{0,\rho,22}-\gamma_\rho'$, we obtain
\[
 \partial_{22}K_{0,\rho}(x,x+\alpha_\rho)
 =-h_{0,\rho,12}(x,x+\alpha_\rho)>0.
\]
The function on the right is continuous and positive on the compact
set $\mathbb T\times J$. Hence
\[
 a_*:=\min_{x\in\mathbb T,\,\rho\in J}
       \big[-h_{0,\rho,12}(x,x+\alpha_\rho)\big]>0.
\]
Choose $\delta_*>0$ so small that the points
$(x,x+\alpha_\rho+s)$ with $|s|\le\delta_*$ remain in the open
edge strip for all $\rho\in J$, and, by uniform continuity,
\[
 \partial_{22}K_{0,\rho}(x,x+\alpha_\rho+s)\ge a_*/2.
\]
At the reference edge both $K_{0,\rho}$ and its second-coordinate
first derivative vanish. Taylor's formula with integral remainder
therefore gives, for $|s|\le\delta_*$,
\begin{align*}
 K_{0,\rho}(x,x+\alpha_\rho+s)
 &=s^2\int_0^1(1-\theta)
   \partial_{22}K_{0,\rho}(x,x+\alpha_\rho+\theta s)\,d\theta\\
 &\ge \frac{a_*}{4}s^2.
\end{align*}

\smallskip\noindent
\emph{Step 3: extend the estimate to the closed edge strip.}
Put $d=y-x$. On the compact set
\[
 \{(x,\rho,d)\in\mathbb T\times J\times[0,2\pi]:
                 |d-\alpha_\rho|\ge\delta_*\},
\]
the quotient
\[
 \frac{K_{0,\rho}(x,x+d)}{(d-\alpha_\rho)^2}
\]
is continuous and strictly positive. The only interior zero of the
numerator is $d=\alpha_\rho$, and the two endpoint values $d=0,2\pi$
are strictly positive by the endpoint argument in
Lemma~\ref{geo:lem:calibration}.
Let $c_{\mathrm{far}}>0$ be its minimum.
Taking $c=\min\{a_*/4,c_{\mathrm{far}}\}$ and combining the two regions
proves \eqref{geo:eq:coercivity}.
Only the reference edges require a uniform twist lower bound; no
such bound is assumed throughout the grazing part of the strip.
The same compactness argument is locally uniform when the reference
ellipse varies in a sufficiently small compact neighborhood.
\end{proof}

\begin{figure}[htbp]
\centering
\begin{tikzpicture}[
 x=1cm,y=1cm,font=\small,
 axis/.style={-{Stealth[length=1.6mm]},line width=.45pt},
 curve/.style={line width=.9pt},
 guide/.style={densely dashed,line width=.55pt},
 every node/.style={inner sep=2pt}]
\path[use as bounding box] (-.65,-.96) rectangle (15.0,4.34);
\node at (3.08,4.08) {(a) Ordered edge strip};
\node at (11.58,4.08) {(b) A calibration slice};
\fill[black!4] (.35,.40) rectangle (6.05,3.50);
\draw[line width=.55pt] (.35,.40) rectangle (6.05,3.50);
\draw[axis] (.35,.40)--(6.45,.40) node[right] {$x$};
\draw[axis] (.35,.40)--(.35,3.83) node[above] {$\zeta=y-x$};
\node[below left] at (.35,.40) {$0$};
\node[below] at (6.05,.40) {$2\pi$};
\node[left] at (.35,3.50) {$2\pi$};
\draw[densely dotted,line width=.65pt] (1.7750000000,.40)--(1.7750000000,3.50);
\node[below] at (1.7750000000,.40) {$x_0$};
\draw[curve] (.35,1.4333333333)--(6.05,1.4333333333);
\node[left] at (.35,1.4333333333) {$\alpha_\rho$};
\fill (1.7750000000,1.4333333333) circle[radius=1.7pt];
\node[fill=white] at (4.2,1.4333333333) {$K_{0,\rho}=0$};
\node at (3.55,2.73) {$K_{0,\rho}>0$};
\node at (3.55,.84) {$K_{0,\rho}>0$};
\node at (3.05,-.63) {$K_{0,\rho}(x,x+\zeta)\ge c(\zeta-\alpha_\rho)^2$};
\draw[axis] (8.85,.40)--(14.77,.40) node[right] {$\zeta$};
\draw[axis] (8.85,.40)--(8.85,3.80);
\node[anchor=west] at (8.93,3.77) {$K_{0,\rho}(x_0,x_0+\zeta)$};
\node[below left] at (8.85,.40) {$0$};
\node[below] at (14.35,.40) {$2\pi$};
\draw (10.6833333333,.35)--(10.6833333333,.45);
\node[below] at (10.6833333333,.35) {$\alpha_\rho$};
\draw[curve] plot coordinates {
(8.8500000000,0.6893164701)
(8.8729166667,0.6806003889)
(8.8958333333,0.6719021051)
(8.9187500000,0.6632392549)
(8.9416666667,0.6546291559)
(8.9645833333,0.6460886554)
(8.9875000000,0.6376339891)
(9.0104166667,0.6292806512)
(9.0333333333,0.6210432791)
(9.0562500000,0.6129355537)
(9.0791666667,0.6049701176)
(9.1020833333,0.5971585098)
(9.1250000000,0.5895111196)
(9.1479166667,0.5820371567)
(9.1708333333,0.5747446388)
(9.1937500000,0.5676403949)
(9.2166666667,0.5607300819)
(9.2395833333,0.5540182146)
(9.2625000000,0.5475082063)
(9.2854166667,0.5412024182)
(9.3083333333,0.5351022164)
(9.3312500000,0.5292080343)
(9.3541666667,0.5235194393)
(9.3770833333,0.5180352011)
(9.4000000000,0.5127533625)
(9.4229166667,0.5076713088)
(9.4458333333,0.5027858371)
(9.4687500000,0.4980932235)
(9.4916666667,0.4935892876)
(9.5145833333,0.4892694542)
(9.5375000000,0.4851288115)
(9.5604166667,0.4811621653)
(9.5833333333,0.4773640904)
(9.6062500000,0.4737289765)
(9.6291666667,0.4702510716)
(9.6520833333,0.4669245211)
(9.6750000000,0.4637434031)
(9.6979166667,0.4607017601)
(9.7208333333,0.4577936282)
(9.7437500000,0.4550130624)
(9.7666666667,0.4523541595)
(9.7895833333,0.4498110784)
(9.8125000000,0.4473780581)
(9.8354166667,0.4450494337)
(9.8583333333,0.4428196510)
(9.8812500000,0.4406832790)
(9.9041666667,0.4386350219)
(9.9270833333,0.4366697298)
(9.9500000000,0.4347824088)
(9.9729166667,0.4329682311)
(9.9958333333,0.4312225443)
(10.0187500000,0.4295408816)
(10.0416666667,0.4279189723)
(10.0645833333,0.4263527523)
(10.0875000000,0.4248383770)
(10.1104166667,0.4233722341)
(10.1333333333,0.4219509592)
(10.1562500000,0.4205714526)
(10.1791666667,0.4192308986)
(10.2020833333,0.4179267880)
(10.2250000000,0.4166569433)
(10.2479166667,0.4154195470)
(10.2708333333,0.4142131748)
(10.2937500000,0.4130368326)
(10.3166666667,0.4118899992)
(10.3395833333,0.4107726735)
(10.3625000000,0.4096854291)
(10.3854166667,0.4086294753)
(10.4083333333,0.4076067254)
(10.4312500000,0.4066198739)
(10.4541666667,0.4056724823)
(10.4770833333,0.4047690752)
(10.5000000000,0.4039152473)
(10.5229166667,0.4031177817)
(10.5458333333,0.4023847811)
(10.5687500000,0.4017258121)
(10.5916666667,0.4011520638)
(10.6145833333,0.4006765215)
(10.6375000000,0.4003141553)
(10.6604166667,0.4000821250)
(10.6833333333,0.4000000000)
(10.7062500000,0.4000899955)
(10.7291666667,0.4003772226)
(10.7520833333,0.4008899512)
(10.7750000000,0.4016598838)
(10.7979166667,0.4027224363)
(10.8208333333,0.4041170205)
(10.8437500000,0.4058873226)
(10.8666666667,0.4080815691)
(10.8895833333,0.4107527696)
(10.9125000000,0.4139589222)
(10.9354166667,0.4177631671)
(10.9583333333,0.4222338667)
(10.9812500000,0.4274445906)
(11.0041666667,0.4334739768)
(11.0270833333,0.4404054398)
(11.0500000000,0.4483266899)
(11.0729166667,0.4573290283)
(11.0958333333,0.4675063770)
(11.1187500000,0.4789540091)
(11.1416666667,0.4917669409)
(11.1645833333,0.5060379613)
(11.1875000000,0.5218552796)
(11.2104166667,0.5392997953)
(11.2333333333,0.5584420122)
(11.2562500000,0.5793386545)
(11.2791666667,0.6020290752)
(11.3020833333,0.6265315919)
(11.3250000000,0.6528399276)
(11.3479166667,0.6809199782)
(11.3708333333,0.7107071591)
(11.3937500000,0.7421046060)
(11.4166666667,0.7749824986)
(11.4395833333,0.8091787413)
(11.4625000000,0.8445011681)
(11.4854166667,0.8807313325)
(11.5083333333,0.9176298116)
(11.5312500000,0.9549428091)
(11.5541666667,0.9924096927)
(11.5770833333,1.0297709920)
(11.6000000000,1.0667763114)
(11.6229166667,1.1031916113)
(11.6458333333,1.1388053750)
(11.6687500000,1.1734333032)
(11.6916666667,1.2069213350)
(11.7145833333,1.2391469717)
(11.7375000000,1.2700190304)
(11.7604166667,1.2994760748)
(11.7833333333,1.3274838411)
(11.8062500000,1.3540319977)
(11.8291666667,1.3791305595)
(11.8520833333,1.4028062299)
(11.8750000000,1.4250988833)
(11.8979166667,1.4460583339)
(11.9208333333,1.4657414796)
(11.9437500000,1.4842098610)
(11.9666666667,1.5015276353)
(11.9895833333,1.5177599443)
(12.0125000000,1.5329716368)
(12.0354166667,1.5472262990)
(12.0583333333,1.5605855459)
(12.0812500000,1.5731085276)
(12.1041666667,1.5848516079)
(12.1270833333,1.5958681793)
(12.1500000000,1.6062085829)
(12.1729166667,1.6159201074)
(12.1958333333,1.6250470457)
(12.2187500000,1.6336307933)
(12.2416666667,1.6417099738)
(12.2645833333,1.6493205832)
(12.2875000000,1.6564961430)
(12.3104166667,1.6632678586)
(12.3333333333,1.6696647771)
(12.3562500000,1.6757139423)
(12.3791666667,1.6814405442)
(12.4020833333,1.6868680632)
(12.4250000000,1.6920184057)
(12.4479166667,1.6969120332)
(12.4708333333,1.7015680834)
(12.4937500000,1.7060044838)
(12.5166666667,1.7102380579)
(12.5395833333,1.7142846249)
(12.5625000000,1.7181590929)
(12.5854166667,1.7218755461)
(12.6083333333,1.7254473266)
(12.6312500000,1.7288871116)
(12.6541666667,1.7322069863)
(12.6770833333,1.7354185123)
(12.7000000000,1.7385327940)
(12.7229166667,1.7415605410)
(12.7458333333,1.7445121293)
(12.7687500000,1.7473976593)
(12.7916666667,1.7502270143)
(12.8145833333,1.7530099165)
(12.8375000000,1.7557559837)
(12.8604166667,1.7584747859)
(12.8833333333,1.7611759019)
(12.9062500000,1.7638689779)
(12.9291666667,1.7665637873)
(12.9520833333,1.7692702921)
(12.9750000000,1.7719987083)
(12.9979166667,1.7747595729)
(13.0208333333,1.7775638160)
(13.0437500000,1.7804228372)
(13.0666666667,1.7833485865)
(13.0895833333,1.7863536517)
(13.1125000000,1.7894513523)
(13.1354166667,1.7926558411)
(13.1583333333,1.7959822134)
(13.1812500000,1.7994466262)
(13.2041666667,1.8030664272)
(13.2270833333,1.8068602952)
(13.2500000000,1.8108483921)
(13.2729166667,1.8150525292)
(13.2958333333,1.8194963467)
(13.3187500000,1.8242055094)
(13.3416666667,1.8292079182)
(13.3645833333,1.8345339391)
(13.3875000000,1.8402166505)
(13.4104166667,1.8462921085)
(13.4333333333,1.8527996316)
(13.4562500000,1.8597821041)
(13.4791666667,1.8672862973)
(13.5020833333,1.8753632080)
(13.5250000000,1.8840684111)
(13.5479166667,1.8934624242)
(13.5708333333,1.9036110778)
(13.5937500000,1.9145858867)
(13.6166666667,1.9264644124)
(13.6395833333,1.9393306063)
(13.6625000000,1.9532751206)
(13.6854166667,1.9683955681)
(13.7083333333,1.9847967105)
(13.7312500000,2.0025905508)
(13.7541666667,2.0218962982)
(13.7770833333,2.0428401711)
(13.8000000000,2.0655549995)
(13.8229166667,2.0901795801)
(13.8458333333,2.1168577383)
(13.8687500000,2.1457370457)
(13.8916666667,2.1769671413)
(13.9145833333,2.2106976104)
(13.9375000000,2.2470753790)
(13.9604166667,2.2862415955)
(13.9833333333,2.3283279887)
(14.0062500000,2.3734527184)
(14.0291666667,2.4217157644)
(14.0520833333,2.4731939446)
(14.0750000000,2.5279356938)
(14.0979166667,2.5859557881)
(14.1208333333,2.6472302472)
(14.1437500000,2.7116916919)
(14.1666666667,2.7792254664)
(14.1895833333,2.8496668476)
(14.2125000000,2.9227996569)
(14.2354166667,2.9983565462)
(14.2583333333,3.0760211625)
(14.2812500000,3.1554322908)
(14.3041666667,3.2361899527)
(14.3270833333,3.3178632996)
(14.3500000000,3.4000000000)
};
\fill (8.85,0.6893164701) circle[radius=1.6pt];
\node[left=2pt] at (8.85,0.6893164701) {$-\kappa_\rho$};
\fill (14.35,3.40) circle[radius=1.6pt];
\node[anchor=south east] at (14.45,3.46) {$-\kappa_\rho-\mathcal A_\rho$};
\fill (10.6833333333,.40) circle[radius=1.8pt];
\node at (11.58,-.63) {$x_0=\pi/2,\quad \rho=1/3,\quad \alpha_\rho=2\pi/3$};
\end{tikzpicture}
\caption{Reference calibration on the closed ordered edge strip.
Left: one angular period in the coordinates $(x,\zeta)$, where
$\zeta=y-x$; the reference caustic is the zero-excess line
$\zeta=\alpha_\rho$. Right: a slice for the ellipse with semiaxes
$(a,b)=(2,1.2)$, at $\rho=1/3$ and $x_0=\pi/2$.
Both grazing values are positive. Here
$\mathcal A_\rho=H_\rho(2\pi)-H_\rho(0)$ is the retained primitive drift.
The plotted slice is illustrative; the uniform quadratic bound is proved
in Lemma~\ref{geo:lem:coercivity}.}
\label{geo:fig:reference-calibration}
\end{figure}
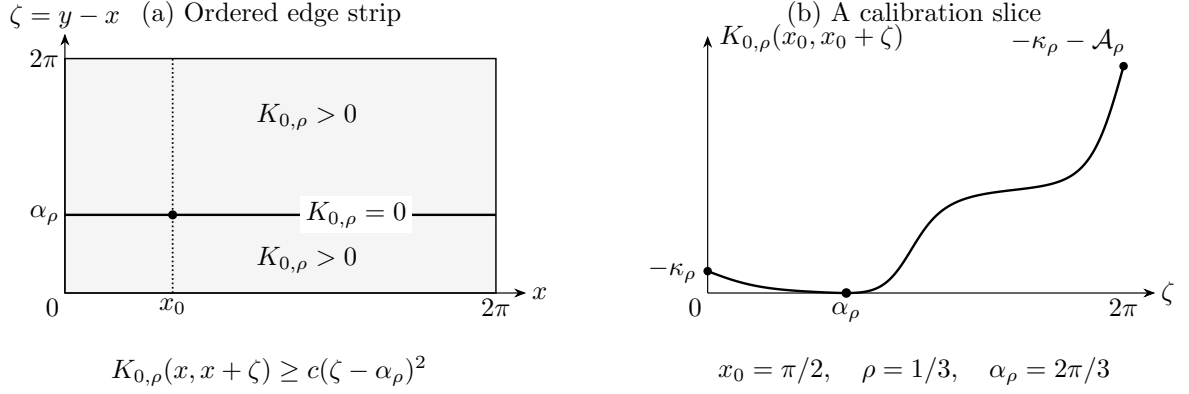

Figure~\ref{geo:fig:reference-calibration} illustrates the zero-excess line
and the positive grazing values of the reference calibration.

\begin{lemma}[A global perturbation bound, including grazing endpoints]\label{geo:lem:chord}
For $t=\norm\mu_{C^1}$ sufficiently small, the difference
$\mathcal R_\rho=h_{\mu,\rho}-h_{0,\rho}$ has supremum norm and first endpoint derivative norms at most $Ct$ throughout $\mathcal D$, uniformly for $\rho\in J$. It extends continuously to the closed strip and is $Ct$-Lipschitz there along line segments.
\end{lemma}
\begin{proof}
All constants below concern the reference ellipse, the compact
interval $J\Subset(0,1/2)$, and a fixed small $C^1$ neighborhood of
the zero displacement.

\smallskip\noindent
\emph{Step 1: control the boundary perturbation in the fixed coordinate.}
Put $\Delta R_\rho=R_{\mu,\rho}-R_{0,\rho}$.
The smooth elliptic chart gives
\[
 r_\mu'(\varphi)
 =\partial_s\mathcal X_E(\mu(\varphi),\varphi)\mu'(\varphi)
    +\partial_\varphi\mathcal X_E(\mu(\varphi),\varphi).
\]
Subtracting the formula at $\mu=0$ and using smoothness of the
coordinate coefficients shows that
$\|r_\mu-r_0\|_{C^1}\le Ct$.
Both curves are composed with the same reference map $\Phi_\rho$;
its uniform first-derivative bound consequently yields
\[
 \|\Delta R_\rho\|_{C^1}\le C_1t,
 \qquad \|R_{0,\rho}'\|_{C^0}\le C_0.
\]
No second derivative of $\mu$ is used.

\smallskip\noindent
\emph{Step 2: compare both chords on the same short-arc scale.}
For $(x,y)\in\mathcal D$, define
\[
 \sigma(x,y)=\min\{y-x,\,2\pi-(y-x)\}=d_{\mathbb T}(x,y).
\]
There is $c_*>0$, uniform in $\rho\in J$, such that
\[
 |R_{0,\rho}(y)-R_{0,\rho}(x)|\ge c_*\sigma(x,y).
\]
For an explicit justification, let
$\ell_J=\inf_{\rho\in J,\,x\in\mathbb T}\Phi_\rho'(x)>0$.
The inverse circle maps have Lipschitz constant at most
$\ell_J^{-1}$ for circle distance, so
\[
 d_{\mathbb T}(\Phi_\rho(x),\Phi_\rho(y))
 \ge\ell_J d_{\mathbb T}(x,y).
\]
For $r_0(\varphi)=(a\cos\varphi,b\sin\varphi)$ and
$\delta=d_{\mathbb T}(\varphi_1,\varphi_2)\in[0,\pi]$,
\[
 |r_0(\varphi_2)-r_0(\varphi_1)|
 \ge2b\sin(\delta/2)\ge\frac{2b}{\pi}\delta.
\]
Thus one may take $c_*=2b\ell_J/\pi$.

The perturbation $\Delta R_\rho$ is periodic. Integrating its
derivative from $x$ to $y$, or from $y$ to $x+2\pi$, according to
which arc is shorter, gives
\[
 |\Delta R_\rho(y)-\Delta R_\rho(x)|\le C_1t\sigma(x,y).
\]
Write the reference and perturbed chord vectors as
\[
 \mathbf c_0=R_{0,\rho}(y)-R_{0,\rho}(x),\qquad
 \mathbf c_\mu=R_{\mu,\rho}(y)-R_{\mu,\rho}(x).
\]
They satisfy
\[
 |\mathbf c_0|\ge c_*\sigma,\qquad
 |\mathbf c_\mu-\mathbf c_0|\le C_1t\sigma.
\]
After decreasing the smallness threshold so that $C_1t\le c_*/2$,
\[
 |\mathbf c_\mu|\ge(c_*-C_1t)\sigma\ge(c_*/2)\sigma>0.
\]
Thus both unit chord directions are well defined everywhere in
$\mathcal D$, however close the edge is to either grazing endpoint.

\smallskip\noindent
\emph{Step 3: normalize the chords and estimate endpoint derivatives.}
For nonzero vectors $v,w$, the identity
\[
 \frac{v}{|v|}-\frac{w}{|w|}
 =\frac{v-w}{|w|}
    +v\left(\frac1{|v|}-\frac1{|w|}\right)
\]
and the reverse triangle inequality imply
\[
 \left|\frac{v}{|v|}-\frac{w}{|w|}\right|
 \le\frac{2|v-w|}{|w|}.
\]
Consequently, with $e_\mu=\mathbf c_\mu/|\mathbf c_\mu|$ and
$e_0=\mathbf c_0/|\mathbf c_0|$,
\[
 |e_\mu-e_0|
 \le\frac{2C_1t\sigma}{c_*\sigma}
 =\frac{2C_1}{c_*}t.
\]
The same short-arc factor occurs in the numerator and denominator;
its cancellation is what makes this estimate uniform at grazing.
The reverse triangle inequality also gives the stronger value bound
\[
 |\mathcal R_\rho(x,y)|
 =\big||\mathbf c_0|-|\mathbf c_\mu|\big|
 \le C_1t\sigma(x,y)\le\pi C_1t.
\]

For any regular curve $R$ and its unit directed chord $e$, direct
differentiation of negative chord length gives
\[
 \partial_x\big(-|R(y)-R(x)|\big)=e\cdot R'(x),\qquad
 \partial_y\big(-|R(y)-R(x)|\big)=-e\cdot R'(y).
\]
Subtracting the reference and perturbed formulas yields
\begin{align*}
 \partial_x\mathcal R_\rho
 &=(e_\mu-e_0)\cdot R_{0,\rho}'(x)
                   +e_\mu\cdot\Delta R_\rho'(x),\\
 \partial_y\mathcal R_\rho
 &=-(e_\mu-e_0)\cdot R_{0,\rho}'(y)
                   -e_\mu\cdot\Delta R_\rho'(y).
\end{align*}
The bounds from Step~1 and $|e_\mu|=1$ now give
\[
 \sup_{\mathcal D}|\partial_x\mathcal R_\rho|
   +\sup_{\mathcal D}|\partial_y\mathcal R_\rho|\le Ct.
\]
Together with the value estimate, this proves all asserted interior
bounds.

\smallskip\noindent
\emph{Step 4: extend the Lipschitz estimate to the closed strip.}
The strip $\mathcal D$ is convex. For any two of its points,
integration of the derivative bounds along their connecting segment
gives
\[
 |\mathcal R_\rho(x,y)-\mathcal R_\rho(\widetilde x,\widetilde y)|
 \le Ct\big(|x-\widetilde x|+|y-\widetilde y|\big).
\]
The value estimate with the factor $\sigma$ shows that
$\mathcal R_\rho$ extends continuously, with value zero, at both
boundaries $y=x$ and $y=x+2\pi$. Approximating closed-strip endpoints
by interior points extends the displayed Lipschitz bound to the whole
closed strip. The lifted endpoints need not be close to one another;
only their connecting segment must remain in the strip, as convexity
ensures. No derivative bound across the strip boundary or estimate
for a perturbed caustic is required.
All constants used above are locally uniform for reference ellipses
in a sufficiently small compact neighborhood.
\end{proof}

\subsection{The uniform rational-caustic estimate}
\begin{theorem}[Interior orbit-sum remainder]\label{geo:thm:remainder}
Fix $E$ and $J\Subset(0,1/2)$. There are $\varepsilon,C>0$ such that the following holds. Let $\mu\in C^2(\T)$, $t=\norm\mu_{C^1}\le\varepsilon$, and suppose $\Omega_\mu$ is strictly convex with positive curvature. If it has an integrable rational caustic of reduced rotation $\rho=p/d\in J$, there is a constant $c_{\mu,\rho}$, independent of $\theta$, with
\begin{equation}\label{geo:eq:melnikov}
\left|2\lambda_\rho\sum_{j=0}^{d-1}
\mu\big(\Phi_\rho(\theta+2\pi j\rho)\big)-c_{\mu,\rho}\right|
\le C d^3 t^2\qquad\text{for all }\theta\in\R.
\end{equation}
No restriction of the form $d^\kappa t\ll1$ is imposed. The constants are uniform when the reference ellipse varies in a sufficiently small compact neighborhood. A bounded high norm is unnecessary for this theorem.
\end{theorem}

\begin{proof}
Fix a reduced rotation $\rho=p/d\in J$ and put $\alpha=2\pi\rho$.
The constants in the quantitative estimates depend only on the reference
ellipse, the compact interval $J$, and a fixed small graph-coordinate
neighborhood; they
are independent of $p,d$, the starting phase, and the shape of the
perturbed caustic.

\smallskip\noindent
\emph{Step 1: match the lifted endpoints in one fixed coordinate.}
Use the reference action-angle map $\Phi_\rho$ on both boundaries:
\[
 R_{0,\rho}=r_0\circ\Phi_\rho,\qquad
 R_{\mu,\rho}=r_\mu\circ\Phi_\rho.
\]
Write $G_\rho=\Phi_\rho^{-1}$. By
Remark~\ref{geo:rem:caustic-graph}, the perturbed caustic induces an
increasing degree-one boundary homeomorphism. Denote its lift in the
physical angular parameter by $f_\mu^{\mathrm{phys}}$; in the reference
action coordinate its lift is
\[
 \widetilde f_\mu
 =G_\rho\circ f_\mu^{\mathrm{phys}}\circ\Phi_\rho.
\]
Integrable rationality and the
\hyperref[geo:par:common-period]{common-period observation} give
$(f_\mu^{\mathrm{phys}})^d(\varphi)=\varphi+2\pi p$.
Since $G_\rho$ and $\Phi_\rho$ are increasing degree-one lifts,
\[
 \widetilde f_\mu^d(\theta)
 =G_\rho(\Phi_\rho(\theta)+2\pi p)=\theta+2\pi p,
 \qquad x<\widetilde f_\mu(x)<x+2\pi.
\]
For each $\theta\in\R$, set
\[
 x_j=\theta+j\alpha,\qquad y_j=\widetilde f_\mu^j(\theta),
 \qquad 0\le j\le d.
\]
Thus $x$ is the reference orbit and $y$ is the perturbed caustic
orbit, with
\[
 y_0=x_0=\theta,\qquad y_d=x_d=\theta+2\pi p,
 \qquad 0<y_{j+1}-y_j<2\pi.
\]
The entire periodic caustic is needed to obtain these identities for
every $\theta$; a single periodic orbit would not suffice.
The smooth increasing change of boundary coordinate preserves the
continuous invariant-graph hypotheses of
Lemma~\ref{geo:lem:calibration}. Both caustics can therefore be
calibrated in this same coordinate. We do not assume that the
perturbed map $\widetilde f_\mu$ is a translation.

\smallskip\noindent
\emph{Step 2: compare the two calibrated actions.}
Suppress $\rho$ in the generating functions and, for an ordered
configuration $z=(z_0,\ldots,z_d)$, write
\[
 A_\mu(z)=\sum_{j=0}^{d-1}h_\mu(z_j,z_{j+1}),\qquad
 A_0(z)=\sum_{j=0}^{d-1}h_0(z_j,z_{j+1}),\qquad
 \mathcal R=h_\mu-h_0.
\]
Define
\[
 \delta_j=y_j-x_j,\qquad \delta_0=\delta_d=0,\qquad
 \mathcal E=\sum_{j=0}^{d-1}(\delta_{j+1}-\delta_j)^2.
\]
For the reference calibration, Lemma~\ref{geo:lem:calibration} gives
\[
 A_0(z)=d\kappa_0+H_0(z_d)-H_0(z_0)
           +\sum_{j=0}^{d-1}K_0(z_j,z_{j+1}).
\]
Every reference edge has zero excess. Since $x$ and $y$ have the
same lifted endpoints, their primitive terms cancel, and
Lemma~\ref{geo:lem:coercivity} yields
\begin{align*}
 A_0(y)-A_0(x)
 &=\sum_{j=0}^{d-1}K_0(y_j,y_{j+1})\\
 &\ge c_*\sum_{j=0}^{d-1}(y_{j+1}-y_j-\alpha)^2
  =c_*\mathcal E,
\end{align*}
where $c_*=c(E,J)>0$.
On the perturbed boundary, $x$ is an admissible competing
configuration for the caustic orbit $y$, with the same endpoints.
Perturbed calibration gives $A_\mu(y)\le A_\mu(x)$.
Consequently, setting
\[
 B_\theta:=\sum_{j=0}^{d-1}
 \big[\mathcal R(x_j,x_{j+1})-\mathcal R(y_j,y_{j+1})\big],
\]
we have $A_0(y)-A_0(x)\le B_\theta$.

Let $L$ be the Lipschitz constant in Lemma~\ref{geo:lem:chord}
after factoring out $t=\|\mu\|_{C^1}$. The segment between each
pair of compared edges stays in the ordered strip: its increment is
\[
 (1-\eta)\alpha+\eta(y_{j+1}-y_j)\in(0,2\pi),
 \qquad 0\le\eta\le1.
\]
Hence each summand defining $B_\theta$ is at most
$Lt(|\delta_j|+|\delta_{j+1}|)$. The interior displacements
occur twice, and the endpoint contributions vanish, so
\begin{align}\label{geo:eq:energy}
 c_*\mathcal E
 &\le A_0(y)-A_0(x)\nonumber\\
 &\le B_\theta\le2Lt\sum_{j=1}^{d-1}|\delta_j|.
\end{align}
The comparison is global in the lifted coordinates; no closeness
of $x$ and $y$ has been used. It also applies when $p>1$, since
calibration requires ordered lifts, not a simple polygon.

\smallskip\noindent
\emph{Step 3: obtain the discrete energy and action-gap estimates.}
Telescoping from $\delta_0=0$ and applying Cauchy--Schwarz give
\[
 |\delta_j|
 =\left|\sum_{i=0}^{j-1}(\delta_{i+1}-\delta_i)\right|
 \le\sqrt j\,\sqrt{\mathcal E}
 \le\sqrt d\,\sqrt{\mathcal E},\qquad
 \sum_{j=1}^{d-1}|\delta_j|\le d^{3/2}\sqrt{\mathcal E}.
\]
If $\mathcal E=0$, all displacements vanish. Otherwise divide
\eqref{geo:eq:energy} by $c_*\sqrt{\mathcal E}$ to get
$\sqrt{\mathcal E}\le(2L/c_*)d^{3/2}t$. In either case,
\begin{equation}\label{geo:eq:localization}
 \mathcal E\le\frac{4L^2}{c_*^2}d^3t^2,\qquad
 \max_{0\le j\le d}|\delta_j|\le\frac{2L}{c_*}d^2t.
\end{equation}
Now let $\Delta_\theta=A_\mu(x)-A_\mu(y)\ge0$.
The exact relation between the two action comparisons is
\[
 \Delta_\theta=B_\theta-\big(A_0(y)-A_0(x)\big).
\]
The subtracted term is nonnegative. Using the Lipschitz estimate and
the energy bound therefore gives
\begin{equation}\label{geo:eq:gap}
 0\le\Delta_\theta
 \le2Lt\sum_{j=1}^{d-1}|\delta_j|
 \le\frac{4L^2}{c_*}d^3t^2.
\end{equation}
The power $d^3$ comes from the discrete factor $d^{3/2}$ and its
square in the energy estimate. The bound on the displacements in
\eqref{geo:eq:localization} need not be small; in particular, we
will not Taylor-expand the action in $y-x$.
Figure~\ref{geo:fig:action-comparison} summarizes the two calibrated
comparisons and their common endpoint condition.

\begin{figure}[htbp]
\centering
\begin{tikzpicture}[
 x=1cm,y=1cm,font=\small,
 axis/.style={-{Stealth[length=1.6mm]},line width=.45pt},
 arrow/.style={-{Stealth[length=1.8mm]},line width=.65pt},
 box/.style={draw,line width=.55pt,rounded corners=1.5pt,minimum width=1.25cm,
             minimum height=.7cm,fill=white},
 every node/.style={inner sep=2pt}]
\path[use as bounding box] (-.85,-1.80) rectangle (15.0,4.30);
\node at (3.10,4.05) {(a) Matched-endpoint displacements};
\node at (11.08,4.05) {(b) The two calibrated actions};
\draw[axis] (.42,.50)--(6.34,.50) node[right] {$j$};
\draw[axis] (.42,.50)--(.42,3.55);
\node[anchor=west] at (.54,3.46) {$\delta_j$};
\foreach \val/\height in {0.05/1.305555556,0.10/2.111111111,0.15/2.916666667}{
 \draw (.36,\height)--(.48,\height);
 \node[left] at (.34,\height) {$\val$};
}
\draw (0.4200000000,.45)--(0.4200000000,.55); \node[below] at (0.4200000000,.43) {$ 0 $};
\draw (1.2228571429,.45)--(1.2228571429,.55); \node[below] at (1.2228571429,.43) {$ 1 $};
\draw (2.0257142857,.45)--(2.0257142857,.55); \node[below] at (2.0257142857,.43) {$ 2 $};
\draw (2.8285714286,.45)--(2.8285714286,.55); \node[below] at (2.8285714286,.43) {$ 3 $};
\draw (3.6314285714,.45)--(3.6314285714,.55); \node[below] at (3.6314285714,.43) {$ 4 $};
\draw (4.4342857143,.45)--(4.4342857143,.55); \node[below] at (4.4342857143,.43) {$ 5 $};
\draw (5.2371428571,.45)--(5.2371428571,.55); \node[below] at (5.2371428571,.43) {$ 6 $};
\draw (6.0400000000,.45)--(6.0400000000,.55); \node[below] at (6.0400000000,.43) {$ 7 $};
\draw[line width=.8pt] plot coordinates {
(0.4200000000,0.5000000000)
(1.2228571429,1.5332496743)
(2.0257142857,2.4087294940)
(2.8285714286,1.3718629794)
(3.6314285714,1.3774747698)
(4.4342857143,2.9979643308)
(5.2371428571,1.9743150818)
(6.0400000000,0.5000000000)
};
\fill (0.4200000000,0.5000000000) circle[radius=1.7pt];
\fill (1.2228571429,1.5332496743) circle[radius=1.7pt];
\fill (2.0257142857,2.4087294940) circle[radius=1.7pt];
\fill (2.8285714286,1.3718629794) circle[radius=1.7pt];
\fill (3.6314285714,1.3774747698) circle[radius=1.7pt];
\fill (4.4342857143,2.9979643308) circle[radius=1.7pt];
\fill (5.2371428571,1.9743150818) circle[radius=1.7pt];
\fill (6.0400000000,0.5000000000) circle[radius=1.7pt];
\node at (3.13,-.04) {$\delta_0=\delta_d=0,\qquad d=7,\quad p=2$};
\node[box] (A0x) at (8.2,2.97) {$A_0(x)$};
\node[box] (A0y) at (13.9,2.97) {$A_0(y)$};
\node[box] (Amx) at (8.2,1.25) {$A_\mu(x)$};
\node[box] (Amy) at (13.9,1.25) {$A_\mu(y)$};
\draw[arrow] (A0x.east)--(A0y.west)
 node[midway,above=4pt] {$A_0(y)-A_0(x)$}
 node[midway,below=4pt] {$\ge c_*\mathcal E$};
\draw[arrow] (Amy.west)--(Amx.east)
 node[midway,above=4pt] {$\Delta_\theta=A_\mu(x)-A_\mu(y)$}
 node[midway,below=4pt] {$\ge0$};
\draw[densely dotted,line width=.45pt] (A0x.south)--(Amx.north);
\draw[densely dotted,line width=.45pt] (A0y.south)--(Amy.north);
\node at (11.05,-.04) {$B_\theta=[A_0(y)-A_0(x)]+\Delta_\theta$};
\draw[line width=.35pt] (0,-.48)--(14.75,-.48);
\node at (7.2,-1.17) {$\displaystyle
 \mathcal E=\sum_{j=0}^{d-1}(\delta_{j+1}-\delta_j)^2,
 \qquad c_*\mathcal E+\Delta_\theta\le B_\theta
       \le 2Lt\sum_{j=1}^{d-1}|\delta_j|$};
\end{tikzpicture}
\caption{Matched endpoints and the comparison of two calibrated actions.
Left: $\delta_j=y_j-x_j$ in one fixed reference action coordinate,
for the ellipse $(a,b)=(2,1.2)$ and its translate by $(0.06,-0.04)$,
with $\rho=2/7$ and starting phase $\theta=0.4$.
Segments join discrete samples; they are not billiard trajectories.
Right: the general comparison in Theorem~\ref{geo:thm:remainder}.
Horizontal arrows point from smaller to larger negative perimeter;
the vertical differences have no prescribed sign.
The argument does not require the bound $Cd^2t$ on the displacements
to be small. The later Taylor expansion varies the boundary at the fixed
reference vertices, not the displacements $\delta_j$.}
\label{geo:fig:action-comparison}
\end{figure}
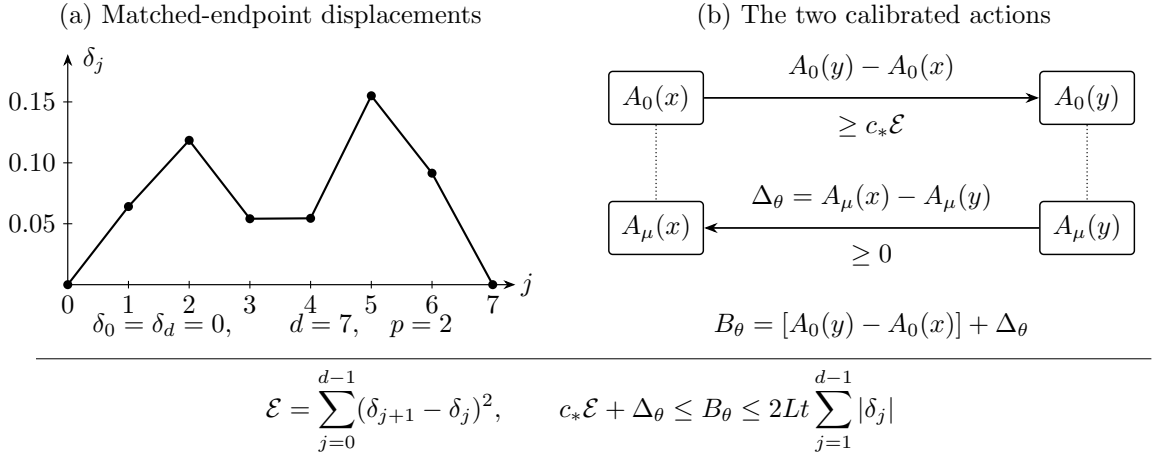

\smallskip\noindent
\emph{Step 4: first variation in the fixed elliptic graph chart.}
Put $\varphi_j=\Phi_\rho(x_j)$, $\mu_j=\mu(\varphi_j)$, and
$q_j(s)=\mathcal X_E(s\mu_j,\varphi_j)$. The elliptic Melnikov formula
\cite[Section~4, equation~(7)]{PRR} gives
\[
 \left.\frac{d}{ds}A_{s\mu}(x)\right|_{s=0}
 =-2\lambda_\rho\sum_{j=0}^{d-1}\mu_j.
\]
Here $n_j$ is the outward unit normal at $q_j(0)$ and $\beta_j$ is
the reference incidence angle; dots below denote $s$ derivatives.
Indeed, $\partial_1\mathcal X_E(0,\varphi_j)
=\sqrt{D_E(\varphi_j)}\,n_j$ and
$\sin\beta_j=\lambda_\rho/\sqrt{D_E(\varphi_j)}$;
the usual vertex contribution to positive perimeter is
$2\sin\beta_j\,\mu_j\sqrt{D_E(\varphi_j)}$.
The minus sign comes from our negative-chord convention.
The cyclic sum includes the common spatial endpoint:
$q_d(s)=q_0(s)$ and $\dot q_d(s)=\dot q_0(s)$, by periodicity and
$\varphi_d-\varphi_0=2\pi p$. Angular vertices, rather than spatial
endpoints, are fixed during this shape variation. The next step proves
the denominator-uniform remainder separately.

\smallskip\noindent
\emph{Step 5: bound the Taylor remainder uniformly by $Cdt^2$.}
Compactness of the reference caustic family gives
\[
 \ell_*:=\min_{\rho\in J,\,\xi\in\T}
 |R_{0,\rho}(\xi+2\pi\rho)-R_{0,\rho}(\xi)|>0.
\]
Indeed, $\Phi_\rho$ is an increasing degree-one diffeomorphism,
$0<\rho<1/2$, and the reference boundary is embedded. The positive
chord length varies continuously on $\T\times J$.
The step is $2\pi\rho$, not $2\pi/d$, so this lower bound does not
deteriorate as the primitive denominator grows within $J$.

The smooth chart on its fixed neighborhood gives, uniformly in
$j$, $\theta$, and $s\in[0,1]$,
\[
 |q_j(s)-q_j(0)|\le C_Xt,\qquad
 |\dot q_j(s)|\le C_Xt,\qquad
 |\ddot q_j(s)|\le C_Xt^2.
\]
For the derivative bounds, note explicitly that
\[
 \dot q_j(s)=\mu_j\,\partial_1\mathcal X_E(s\mu_j,\varphi_j),\qquad
 \ddot q_j(s)=\mu_j^2\,\partial_{11}\mathcal X_E(s\mu_j,\varphi_j),
\]
where $\partial_1$ differentiates the transverse coordinate of the
chart. No derivative of $\mu$ is taken in this shape interpolation.
Writing
\[
 \mathbf c_j(s)=q_{j+1}(s)-q_j(s),\qquad
 \ell_j(s)=|\mathbf c_j(s)|,
\]
and decreasing the fixed threshold so that $2C_Xt\le\ell_*/2$,
we obtain
\[
 \ell_j(s)\ge\ell_*/2,\qquad
 |\dot{\mathbf c}_j(s)|\le2C_Xt,\qquad
 |\ddot{\mathbf c}_j(s)|\le2C_Xt^2.
\]
For $e_j(s)=\mathbf c_j(s)/\ell_j(s)$, differentiating length twice
gives
\[
 \ell_j''(s)
 =\frac{|\dot{\mathbf c}_j(s)|^2
             -(e_j(s)\cdot\dot{\mathbf c}_j(s))^2}{\ell_j(s)}
       +e_j(s)\cdot\ddot{\mathbf c}_j(s).
\]
It follows that $|\ell_j''(s)|\le C_{\mathrm{edge}}t^2$ and hence
$|\frac{d^2}{ds^2}A_{s\mu}(x)|\le C_{\mathrm{edge}}dt^2$.
Taylor's formula with integral remainder now yields
\begin{equation}\label{geo:eq:firstvariation}
 \begin{aligned}
 A_\mu(x)
 &=A_0(x)-2\lambda_\rho\sum_{j=0}^{d-1}
             \mu(\Phi_\rho(\theta+j\alpha))+r_\theta,\\
 |r_\theta|&\le C_Tdt^2.
 \end{aligned}
\end{equation}
Only the actual perturbed boundary at $s=1$ needs the assumed
caustic and positive curvature. The intermediate boundaries are used
solely to differentiate these nonzero chords; no caustic persistence
or curvature assumption along the interpolation is required.

\smallskip\noindent
\emph{Step 6: identify the phase-independent constant and finish.}
Both caustic actions are independent of their starting point by
Lemma~\ref{geo:lem:calibration}. More explicitly, for the reference
and perturbed calibration primitives $H_0,H_\mu$ and their constants
$\kappa_0,\kappa_\mu$, put
\[
 \mathcal A_0=H_0(x+2\pi)-H_0(x),\qquad
 \mathcal A_\mu=H_\mu(x+2\pi)-H_\mu(x).
\]
These drifts are independent of $x$. Thus
\[
 S_0:=A_0(x(\theta))=d\kappa_0+p\mathcal A_0,\qquad
 S_\mu:=A_\mu(y(\theta))=d\kappa_\mu+p\mathcal A_\mu
\]
are constants in $\theta$. Define
\[
 c_{\mu,\rho}=S_0-S_\mu,\qquad
 \mathcal M_\rho(\theta)=2\lambda_\rho\sum_{j=0}^{d-1}
             \mu(\Phi_\rho(\theta+j\alpha)).
\]
The primitive drifts have not been discarded. From
\eqref{geo:eq:firstvariation} and the definition of $\Delta_\theta$,
\[
 \Delta_\theta=c_{\mu,\rho}-\mathcal M_\rho(\theta)+r_\theta,
 \qquad
 \mathcal M_\rho(\theta)-c_{\mu,\rho}=r_\theta-\Delta_\theta.
\]
Therefore \eqref{geo:eq:gap} gives
\[
 |\mathcal M_\rho(\theta)-c_{\mu,\rho}|
 \le C_Tdt^2+\frac{4L^2}{c_*}d^3t^2
 \le Cd^3t^2,
\]
which is \eqref{geo:eq:melnikov}, since $\alpha=2\pi\rho$.
The case $t=0$ is included: the energy and both errors then vanish.

All smallness requirements concern the fixed coordinate neighborhood,
the boundary $C^1$ perturbation estimate, and the reference chord
lower bound $\ell_*$. None involves $p$ or $d$, and the possibly
large bound $Cd^2t$ for the lifted displacements is never required
to be small. No high norm of $\mu$ or quantitative bound on the
perturbed caustic was used. The constants in
Lemmas~\ref{geo:lem:coercivity} and~\ref{geo:lem:chord}, the chart
bounds, and $\ell_*$ are locally uniform for reference ellipses in a
sufficiently small compact neighborhood. This proves all the stated
uniformity claims.
\end{proof}

\begin{corollary}[All detected Fourier harmonics]\label{geo:cor:fourier}
Under Theorem~\ref{geo:thm:remainder}, define
\[
g_\rho(\theta)=\mu(\Phi_\rho(\theta)),\qquad
\widehat g_\rho(n)=\frac1{2\pi}\int_0^{2\pi}g_\rho(\theta)e^{-in\theta}\dd\theta.
\]
For every nonzero integer $n$ divisible by the primitive denominator $d$,
\begin{equation}\label{geo:eq:Fourier}
|\widehat g_\rho(n)|\le C d^2 t^2.
\end{equation}
\end{corollary}
\begin{proof}
Fix the reduced rotation $\rho=p/d\in J$ and write
$t=\|\mu\|_{C^1}$, as in Theorem~\ref{geo:thm:remainder}.

\smallskip\noindent
\emph{Step 1: retain the uniform orbit-sum remainder.}
For a continuous periodic function $g$, define
\[
 S_\rho g(\theta)=\sum_{j=0}^{d-1}
 g(\theta+2\pi jp/d).
\]
The function $g_\rho=\mu\circ\Phi_\rho$ is periodic because $\mu$ is
periodic and $\Phi_\rho$ is a degree-one lift. By
\eqref{geo:eq:melnikov}, the error
\[
 \mathcal E_\rho(\theta)
 :=2\lambda_\rho S_\rho g_\rho(\theta)-c_{\mu,\rho}
\]
satisfies $\|\mathcal E_\rho\|_{L^\infty}\le C_0d^3t^2$.
For every integer $n$, our Fourier normalization gives
\[
 |\widehat{\mathcal E_\rho}(n)|
 \le\frac1{2\pi}\int_0^{2\pi}|\mathcal E_\rho(\theta)|\,d\theta
 \le C_0d^3t^2.
\]

\smallskip\noindent
\emph{Step 2: compute the exact Fourier multiplier.}
For any real shift $\delta$, a change of variable and periodicity give
\[
 \widehat{g(\cdot+\delta)}(n)
 =\frac{e^{in\delta}}{2\pi}
   \int_\delta^{2\pi+\delta}g(u)e^{-inu}\,du
 =e^{in\delta}\widehat g(n).
\]
Consequently, interchanging only a finite sum and an integral,
\[
 \widehat{S_\rho g}(n)
 =\widehat g(n)\sum_{j=0}^{d-1}e^{2\pi i njp/d}.
\]
If $d\mid n$, every summand is one. Otherwise coprimality of $p,d$
implies $z_n:=e^{2\pi i np/d}\ne1$, whereas $z_n^d=1$; the geometric
sum is then $(1-z_n^d)/(1-z_n)=0$. Thus
\[
 \widehat{S_\rho g}(n)
 =\begin{cases}
 d\widehat g(n),&d\mid n,\\
 0,&d\nmid n.
 \end{cases}
\]
This identity uses the primitive denominator $d$ and holds for positive
and negative frequencies alike.

\smallskip\noindent
\emph{Step 3: remove the constant and divide by the resonant multiplier.}
For $n\in d\Z\setminus\{0\}$, the constant $c_{\mu,\rho}$ has zero
$n$th Fourier coefficient. Hence
\[
 2\lambda_\rho d\,|\widehat g_\rho(n)|
 =|\widehat{\mathcal E_\rho}(n)|\le C_0d^3t^2.
\]
The positive continuous function $\rho\mapsto\lambda_\rho$ has a
positive minimum $\lambda_*:=\min_{\rho\in J}\lambda_\rho>0$, since
$J\Subset(0,1/2)$. Dividing gives
\[
 |\widehat g_\rho(n)|\le\frac{C_0}{2\lambda_*}d^2t^2,
\]
which proves \eqref{geo:eq:Fourier}. The constant is independent of the
nonzero resonant frequency and is locally uniform in the reference
ellipse. The zero mode is excluded because $c_{\mu,\rho}$ survives
there; at nonresonant modes the orbit-sum multiplier is zero, so this
argument supplies no bound on the corresponding coefficient of $g_\rho$.
\end{proof}

\begin{corollary}[Exact annihilation of ellipse tangents]\label{geo:cor:tangentkill}
For a fixed reduced rational $0<\rho=p/d<1/2$, every infinitesimal ellipse graph $\tau\in T_E$ satisfies
$\widehat{\tau\circ\Phi_\rho}(n)=0$ for every nonzero $n$ divisible by $d$.
\end{corollary}
\begin{proof}
Fix the reduced rotation $\rho=p/d$ and an index
$n\in d\Z\setminus\{0\}$. Choose a smooth family of ellipses $E_s$,
$E_0=E$, whose graph velocity is $\tau$, by varying the parameters in
\eqref{geo:eq:tangents}. In the fixed elliptic chart of $E$, its
displacement satisfies $\mu_s=s\tau+O(s^2)$ in $C^1$.
Choose a compact interval $J\Subset(0,1/2)$ containing $\rho$.
For sufficiently small $|s|$, the boundary $E_s$ is positively curved,
$\mu_s$ meets the fixed $C^1$ threshold, and its confocal caustic of
rotation $\rho$ exists. Corollary~\ref{geo:cor:fourier}, with the reference
ellipse $E$ fixed, therefore gives
\[
 \big|\widehat{\mu_s\circ\Phi_{E,\rho}}(n)\big|
 \le C d^2\|\mu_s\|_{C^1}^2=O(s^2).
\]
For $s\ne0$, divide by $|s|$ and let $s\to0$.
Since $\mu_s/s\to\tau$ in $C^1$ and $\Phi_{E,\rho}$ is fixed,
continuity of composition and of the Fourier coefficient gives
$\widehat{\tau\circ\Phi_{E,\rho}}(n)=0$.
The rotation, its primitive denominator, and the harmonic are held fixed
in this limit; no uniform-in-denominator differentiation is used.
\end{proof}

\section{The interior-band operator and analytic regularity}
\label{op:sec:operator}
Normalize the major semiaxis to one and write
\[
E_\eta:\quad x_1^2+\frac{x_2^2}{1-\eta}=1,\qquad
0<\eta=e^2<1,\qquad D_\eta(\varphi)=1-\eta\cos^2\varphi.
\]
Similarities preserve rotation numbers and the dimensionless elliptic
graph displacement. All parameter-uniform estimates in this section
are taken on compact intervals contained in $(0,1)$.

Common complex strips make the selected orbit-sum rows exponentially
localized in Fourier coordinates. Near a rational anchor their principal
amplitudes depend on the residue of the Fourier index; twist nevertheless
makes the principal operators invertible. A parametrix then gives
finite-dimensional analytic kernels. Identifying the solutions of all
rational trace equations in the band is a separate task for the later
global continuation argument.

\subsection{Uniform complex caustic coordinates}
\label{op:sec:coordinates}

We use elliptic integrals with parameter $m=k^2$:
\[
 F(\alpha\mid m)=\int_0^\alpha(1-m\sin^2t)^{-1/2}\dd t,
 \qquad K(m)=F(\pi/2\mid m).
\]
For $b=\sqrt{1-\eta}$ and $0<\lambda<b$, put
\[
 m=\frac{\eta}{1-\lambda^2},\qquad
 \rho(\eta,\lambda)
 =\frac{F(\arcsin(\lambda/b)\mid m)}{2K(m)}.
\]
Let $\lambda(\eta,\rho)$ denote the real inverse, and let
$m_{\eta,\rho}=\eta/(1-\lambda(\eta,\rho)^2)$. The coordinates
\eqref{geo:eq:G}--\eqref{geo:eq:Phi} take the form
\begin{equation}\label{op:eq:G}
 G_{\eta,\rho}(\varphi)
 =\frac{\pi}{2K(m_{\eta,\rho})}
   \int_0^\varphi\frac{\dd t}{\sqrt{1-m_{\eta,\rho}\cos^2t}},
 \qquad \Phi_{\eta,\rho}=G_{\eta,\rho}^{-1}.
\end{equation}
Both are degree-one lifts fixing zero.

\begin{lemma}[Common complex neighborhoods]\label{op:lem:strips}
For a compact real interval $I_\eta\subset(0,1)$ and a compact interval
$J\Subset(0,1/2)$, there are common complex parameter neighborhoods and
positive angle-strip widths on which $\lambda$, $G$, and $\Phi$ are
jointly holomorphic. The function $\lambda$ and the periodic parts
$G-\mathrm{id}$ and $\Phi-\mathrm{id}$ are uniformly bounded there
after shrinking the neighborhoods and strips. The coordinate compositions
used in this section satisfy the same bounds after subtracting their
degree-one linear parts. Any prescribed finite number of angle and parameter
derivatives is uniformly bounded on smaller neighborhoods.
\end{lemma}
\begin{proof}
Put $\mathcal P=I_\eta\times J$, and denote the horizontal strip of width
$\sigma>0$ by
\[
 S_\sigma=\{z\in\C:|\operatorname{Im}z|<\sigma\}.
\]
We distinguish the rotation function
$\mathscr R(\eta,\lambda):=\rho(\eta,\lambda)$ from the independent
rotation coordinate $\rho$. The source strip for $G$ and the target
strip on which its inverse is defined need not have the same width.

\smallskip\noindent
\emph{Step 1: strict positivity of the rotation derivative.}
Fix a real $\eta$ and suppress it in $f_\lambda$, $G_\lambda$, and
$\Phi_\lambda$. The incidence identity gives
\[
 \partial_\lambda\alpha_\lambda(\varphi)
 =\frac1{\sqrt{D_\eta(\varphi)-\lambda^2}}>0.
\]
For a chord with fixed initial point, the derivative of its direction
angle with respect to terminal arclength is the sine of the terminal
incidence angle divided by chord length. It is positive by strict
convexity. The positive billiard twist and the displayed identity
therefore give $\partial_\lambda f_\lambda(\varphi)>0$.
Differentiate the real conjugacy at fixed physical angle $\varphi$:
\begin{align*}
 &(\partial_\lambda G_\lambda)(f_\lambda(\varphi))
   +G_\lambda'(f_\lambda(\varphi))\partial_\lambda f_\lambda(\varphi)\\
 &\hspace{18mm}=(\partial_\lambda G_\lambda)(\varphi)
                   +2\pi\partial_\lambda\mathscr R(\eta,\lambda).
\end{align*}
Only after differentiating, substitute $\varphi=\Phi_\lambda(\theta)$.
The function
\[
 A_\lambda(\theta)
   =(\partial_\lambda G_\lambda)(\Phi_\lambda(\theta))
\]
is $2\pi$-periodic, since the degree-one identity for $G_\lambda$
becomes periodic after parameter differentiation. With
$\rho=\mathscr R(\eta,\lambda)$, the conjugacy and inverse derivative
formula give
\[
 A_\lambda(\theta+2\pi\rho)-A_\lambda(\theta)
 +\frac{\partial_\lambda f_\lambda(\Phi_\lambda(\theta))}
        {\Phi_\lambda'(\theta+2\pi\rho)}
 =2\pi\partial_\lambda\mathscr R(\eta,\lambda).
\]
The translated terms have equal averages, so
\begin{equation}\label{op:eq:rhoDerivative}
 2\pi\partial_\lambda\mathscr R(\eta,\lambda)
 =\frac1{2\pi}\int_0^{2\pi}
   \frac{\partial_\lambda f_\lambda(\Phi_\lambda(\theta))}
        {\Phi_\lambda'(\theta+2\pi\rho)}\dd\theta>0.
\end{equation}
No derivative of $\Phi_\lambda$ with respect to $\lambda$ is missing:
the substitution followed differentiation. The argument applies to
all real rotations, not only rational ones.

\smallskip\noindent
\emph{Step 2: one holomorphic inverse on a common parameter neighborhood.}
Section~\ref{geo:sec:convention} gives the real inverse
$\lambda(\eta,\rho)$ and its smooth parameter dependence. The set
\[
 \mathcal K=\{(\eta,\lambda(\eta,\rho)):(\eta,\rho)\in\mathcal P\}
\]
is compact. Consequently, for positive constants $\lambda_-$ and
$\delta_b$, and a constant $m_+<1$, the real family satisfies
\[
 \lambda\ge\lambda_-,\qquad
 b(\eta)-\lambda\ge\delta_b,\qquad
 m=\frac{\eta}{1-\lambda^2}\le m_+.
\]
In particular, $\lambda/b$ stays in a compact subinterval of $(0,1)$.
Choose the holomorphic branches of $b(\eta)=\sqrt{1-\eta}$ and
$\vartheta=\arcsin(\lambda/b)$ that agree with their real values near
$\mathcal K$. The representation
\[
 F(\vartheta\mid m)
 =\vartheta\int_0^1(1-m\sin^2(s\vartheta))^{-1/2}\,ds
\]
provides a common holomorphic branch of the numerator: on the real
compact family its radicand is bounded below by $1-m_+>0$, and
nearby complex radicands remain in the right half-plane. The analogous
fixed-interval integral for $K(m)$ is holomorphic and remains nonzero
after shrinking. Thus $\mathscr R$ is holomorphic on a neighborhood of
$\mathcal K$ with compatible branches throughout.

Consider the holomorphic map
\[
 \mathcal F(\eta,\lambda)=(\eta,\mathscr R(\eta,\lambda)).
\]
Its Jacobian determinant is $\partial_\lambda\mathscr R$, which is
positive on $\mathcal K$ by \eqref{op:eq:rhoDerivative}. It is therefore
locally biholomorphic near $\mathcal K$, and it is injective on
$\mathcal K$ by uniqueness of the real inverse. It is injective on
some whole complex neighborhood of $\mathcal K$ as well. Otherwise,
there would be distinct points $P_n,Q_n$, each at distance less than
$1/n$ from $\mathcal K$, with $\mathcal F(P_n)=\mathcal F(Q_n)$.
Compactness gives subsequential limits $P,Q\in\mathcal K$ with equal
images, so $P=Q$. Local injectivity near this common limit contradicts
the choice of $P_n,Q_n$ for large $n$.

The image of this biholomorphic restriction contains $\mathcal P$.
Choose a product neighborhood $\mathcal U=U_\eta\times U_\rho$ of
$\mathcal P$ with compact closure in that image. Its inverse has a
single holomorphic second component $\lambda(\eta,\rho)$, agreeing
with the real inverse. After shrinking $\mathcal U$, this function is
bounded and stays separated from $0$ and $b(\eta)$; also
\[
 m_{\eta,\rho}=\frac{\eta}{1-\lambda(\eta,\rho)^2}
 \quad\hbox{is holomorphic},\qquad
 |m_{\eta,\rho}|\le m_0<1,
\]
and $K(m_{\eta,\rho})$ stays away from zero. All subsequent shrinkings
retain $\mathcal P$ and leave positive margins inside the preceding
holomorphy domains.

\smallskip\noindent
\emph{Step 3: a common square-root branch and the degree-one identity.}
Choose $\sigma_0>0$ so small that $m_0\cosh^2\sigma_0<1$.
For $z=x+iy\in S_{\sigma_0}$,
\[
 |\cos z|^2=\cos^2x+\sinh^2y\le\cosh^2\sigma_0,
 \qquad |m_{\eta,\rho}\cos^2z|<1.
\]
The binomial series therefore defines
\[
 q_m(z)=(1-m\cos^2z)^{-1/2}
\]
as a single-valued jointly holomorphic function on this common product
domain, agreeing with the positive real branch. Integrating in the
simply connected strip defines $G$ by \eqref{op:eq:G}, jointly
holomorphically, with $G_{\eta,\rho}(0)=0$. The integrand is periodic
and, on the same branch, quadrant substitutions give
\[
 \int_0^{2\pi}q_m(t)\,dt=4K(m).
\]
It follows that
\[
 G_{\eta,\rho}(z+2\pi)=G_{\eta,\rho}(z)+2\pi.
\]
Thus $G_{\eta,\rho}-\mathrm{id}$ is periodic also for complex
parameters. On every smaller closed strip and relatively compact
parameter neighborhood it is uniformly bounded, by continuity on
one compact fundamental rectangle. This is a bound on the periodic
part, not on the unbounded lift $G$ itself.

\smallskip\noindent
\emph{Step 4: univalence on a common source strip.}
For real parameters and real $x$,
\[
 G_{\eta,\rho}'(x)
 =\frac{\pi}{2K(m_{\eta,\rho})
                 \sqrt{1-m_{\eta,\rho}\cos^2x}}>0.
\]
Its minimum on $\mathcal P\times\mathbb T$ is positive. By continuity
and periodicity, choose $0<\sigma<\sigma_0$ and shrink $\mathcal U$
so that, for some $a_*>0$,
\[
 \operatorname{Re}G_{\eta,\rho}'(z)\ge a_*
 \qquad ((\eta,\rho)\in\overline{\mathcal U},\
                         |\operatorname{Im}z|\le\sigma).
\]
The maps remain holomorphic on a neighborhood of this closed strip.
For $z_1\ne z_2$ in the convex strip $S_\sigma$,
\[
 \frac{G_{\eta,\rho}(z_2)-G_{\eta,\rho}(z_1)}{z_2-z_1}
 =\int_0^1G_{\eta,\rho}'(z_1+s(z_2-z_1))\,ds.
\]
The real part of the right side is at least $a_*$; in particular,
\[
 |G_{\eta,\rho}(z_2)-G_{\eta,\rho}(z_1)|
 \ge a_*|z_2-z_1|.
\]
Hence every $G_{\eta,\rho}$ is univalent on the whole strip.
Step~3 also supplies $B>0$ such that
\[
 |G_{\eta,\rho}(z)-z|\le B
 \qquad ((\eta,\rho)\in\overline{\mathcal U},\
                         |\operatorname{Im}z|\le\sigma).
\]

\smallskip\noindent
\emph{Step 5: a common target strip lies in every image.}
For a real parameter pair $\mathbf p\in\mathcal P$, the map $G_{\mathbf p}$
is real on the real axis. Vertical integration gives
\[
 \operatorname{Im}G_{\mathbf p}(x+i\sigma)
   =\int_0^\sigma\operatorname{Re}G_{\mathbf p}'(x+is)\,ds
   \ge a_*\sigma,
 \qquad
 \operatorname{Im}G_{\mathbf p}(x-i\sigma)\le-a_*\sigma.
\]
Uniform continuity on fundamental rectangles makes these signs persist
for all complex $\mathbf p=(\eta,\rho)$ in a smaller $\mathcal U$.
More precisely, with $\beta_*=a_*\sigma/2>0$, arrange
\[
 \operatorname{Im}G_{\mathbf p}(x+i\sigma)\ge\beta_*,\qquad
 \operatorname{Im}G_{\mathbf p}(x-i\sigma)\le-\beta_*
 \quad (x\in\R,\ \mathbf p\in\mathcal U).
\]
Fix $0<\tau<\beta_*$. We prove
$S_\tau\subset G_{\mathbf p}(S_\sigma)$ for every such parameter.
For fixed $\mathbf p$, the set
\[
 W=G_{\mathbf p}(S_\sigma)\cap S_\tau
\]
is relatively open in $S_\tau$ and contains zero. It is also relatively
closed. Indeed, if $w_n\in W$ and $w_n\to w\in S_\tau$, choose
$z_n\in S_\sigma$ with $G_{\mathbf p}(z_n)=w_n$. The bound on the
periodic part gives $|z_n|\le|w_n|+B$, so a subsequence converges to
$z\in\overline S_\sigma$. The boundary signs exclude
$\operatorname{Im}z=\sigma$ and $\operatorname{Im}z=-\sigma$, since
$|\operatorname{Im}w|<\tau<\beta_*$. Thus $z\in S_\sigma$ and
$G_{\mathbf p}(z)=w$, proving closedness. Connectedness of $S_\tau$
now gives $W=S_\tau$. The bound on $G-\mathrm{id}$ prevents horizontal
escape of preimages; the boundary signs prevent escape through either
edge of the strip. All choices are uniform in $\mathbf p$.

\smallskip\noindent
\emph{Step 6: a jointly holomorphic inverse with bounded periodic part.}
Univalence and Step~5 define a unique
$\Phi_{\eta,\rho}:S_\tau\to S_\sigma$ with
\[
 G_{\eta,\rho}(\Phi_{\eta,\rho}(w))=w.
\]
Since $G_{\eta,\rho}'$ never vanishes, the holomorphic implicit-function
theorem applied to $G_{\eta,\rho}(z)-w=0$ gives local solutions jointly
holomorphic in $(\eta,\rho,w)$. They agree on overlaps by uniqueness
of the preimage in $S_\sigma$. Hence $\Phi$ is jointly holomorphic on
the whole product domain and agrees with the real inverse.
The degree-one identity for $G$ and uniqueness imply
\[
 \Phi_{\eta,\rho}(0)=0,\qquad
 \Phi_{\eta,\rho}(w+2\pi)=\Phi_{\eta,\rho}(w)+2\pi.
\]
Moreover,
\[
 \Phi_{\eta,\rho}(w)-w
 =-(G_{\eta,\rho}-\mathrm{id})(\Phi_{\eta,\rho}(w)),\qquad
 |\Phi_{\eta,\rho}(w)-w|\le B.
\]
Thus its periodic part is uniformly bounded, and
\[
 \partial_w\Phi_{\eta,\rho}(w)
   =\frac1{G_{\eta,\rho}'(\Phi_{\eta,\rho}(w))},\qquad
 |\partial_w\Phi_{\eta,\rho}(w)|\le a_*^{-1}.
\]

\smallskip\noindent
\emph{Step 7: coordinate compositions and Cauchy estimates.}
For any fixed real anchor $\rho_*\in J$, the composition used below,
\[
 H_\eta(\rho,w)
   =G_{\eta,\rho}(\Phi_{\eta,\rho_*}(w)),
\]
is defined on the common domains above. It is jointly holomorphic,
a degree-one lift, and
\begin{align*}
 H_\eta(\rho,w)-w
 &= (G_{\eta,\rho}-\mathrm{id})(\Phi_{\eta,\rho_*}(w))
       +(\Phi_{\eta,\rho_*}(w)-w),
\end{align*}
so its periodic part is uniformly bounded. Any further coordinate
compositions occurring in this section can be accommodated by nested
strips and smaller parameter neighborhoods. To check the domain issue,
let $A_{\mathbf p}$ be an inner coordinate map and let $\mathbf p_0$ be a
nearby real parameter pair. On a smaller strip, its bounded first
angle derivative gives
\[
 |\operatorname{Im}A_{\mathbf p}(x+iy)|
 \le |A_{\mathbf p}(x)-A_{\mathbf p_0}(x)|+M_A|y|,
\]
since $A_{\mathbf p_0}(x)$ is real. Periodicity reduces the first term
to a compact real fundamental interval, where it tends uniformly to
zero as the complex parameters approach the real rectangle. Shrinking
the parameter neighborhood and the input strip therefore places the
inner image in the required outer strip. Repeating this for the
finite collection of compositions in question proves the asserted
common-domain and periodic-part bounds; no assertion about arbitrarily
long compositions on one fixed strip is needed.

Finally, let $Q$ be any of the bounded periodic parts just constructed.
On smaller parameter neighborhoods and angle strips, choose positive
Cauchy radii $\delta_\eta,\delta_\rho,\delta_w$ inside their preceding
domains. If $|Q|\le B_Q$ there, the iterated Cauchy formula gives
\[
 |\partial_\eta^a\partial_\rho^b\partial_w^j Q|
 \le B_Q\frac{a!\,b!\,j!}
   {\delta_\eta^a\delta_\rho^b\delta_w^j},
 \qquad a,b,j\in\mathbb N\cup\{0\}.
\]
The same argument applies to the bounded function $\lambda$, omitting
the angle variable. Periodicity makes these estimates uniform along
the real direction of each strip. The full coordinate maps differ from
their periodic parts only by the identity, whose positive-order
derivatives add at most a constant. Thus every prescribed finite set
of angle and parameter derivatives has the stated uniform bounds.
The neighborhoods depend on the fixed compact intervals; no uniformity
at the limiting degenerate parameter values is asserted.
\end{proof}

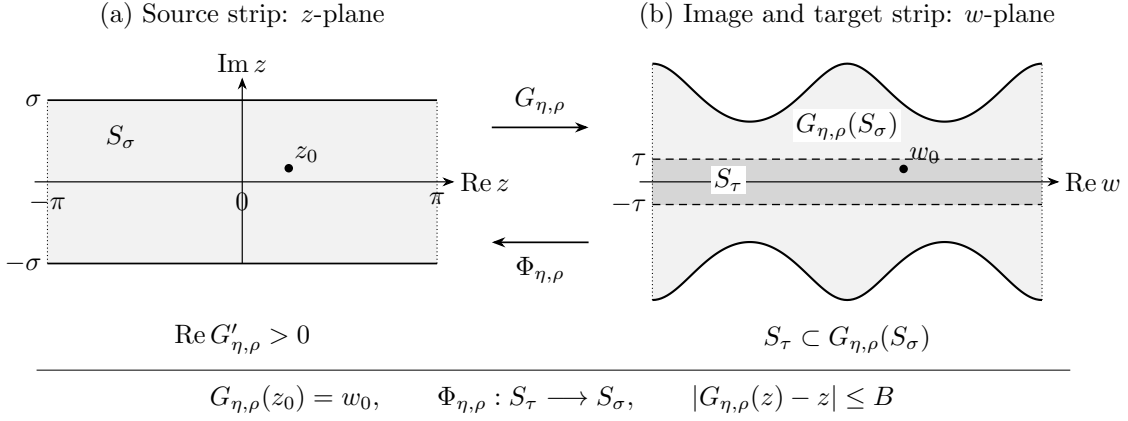
\begin{figure}[tbp]
\centering
\begin{tikzpicture}[
 x=1cm,y=1cm,font=\small,
 axis/.style={-{Stealth[length=1.6mm]},line width=.45pt},
 arr/.style={-{Stealth[length=1.8mm]},line width=.65pt},
 guide/.style={densely dashed,line width=.5pt},
 every node/.style={inner sep=2pt}]
\path[use as bounding box] (-.12,-.94) rectangle (15.2,4.72);
\node at (3.1,4.45) {(a) Source strip: $z$-plane};
\node at (11.1,4.45) {(b) Image and target strip: $w$-plane};
\fill[black!5] (0.5238940241,1.1700000000) rectangle (5.6761059759,3.3300000000);
\draw[line width=.65pt] (0.5238940241,1.1700000000)--(5.6761059759,1.1700000000);
\draw[line width=.65pt] (0.5238940241,3.3300000000)--(5.6761059759,3.3300000000);
\draw[densely dotted] (0.5238940241,1.1700000000)--(0.5238940241,3.3300000000);
\draw[densely dotted] (5.6761059759,1.1700000000)--(5.6761059759,3.3300000000);
\draw[axis] (0.3538940241,2.2500000000)--(5.9061059759,2.2500000000) node[right] {$\operatorname{Re}z$};
\draw[axis] (3.1000000000,1.1700000000)--(3.1000000000,3.6300000000) node[above] {$\operatorname{Im}z$};
\node[left] at (0.5238940241,3.3300000000) {$\sigma$};
\node[left] at (0.5238940241,1.1700000000) {$-\sigma$};
\node at (1.5,2.84) {$S_\sigma$};
\node at (3.1,.18) {$\operatorname{Re}G'_{\eta,\rho}>0$};
\node[below=2pt] at (0.5238940241,2.2500000000) {$-\pi$};
\node[below=2pt] at (3.1000000000,2.2500000000) {$0$};
\node[below=2pt] at (5.6761059759,2.2500000000) {$\pi$};
\fill[black!5] plot coordinates {
(8.5238940241,3.8128191161)
(8.5639280131,3.8102045624)
(8.6037899641,3.8024628285)
(8.6433170700,3.7898880882)
(8.6823635773,3.7729346423)
(8.7208062010,3.7521713369)
(8.7585467713,3.7282313839)
(8.7955123481,3.7017655402)
(8.8316534281,3.6734039622)
(8.8669409887,3.6437289951)
(8.9013630386,3.6132586906)
(8.9349211688,3.5824393510)
(8.9676274056,3.5516448195)
(8.9995015125,3.5211802800)
(9.0305687755,3.4912886934)
(9.0608582447,3.4621584688)
(9.0904013740,3.4339314077)
(9.1192309893,3.4067103196)
(9.1473805217,3.3805659716)
(9.1748834457,3.3555432207)
(9.2017728774,3.3316662914)
(9.2280812925,3.3089432337)
(9.2538403358,3.2873696332)
(9.2790806975,3.2669316609)
(9.3038320415,3.2476085551)
(9.3281229693,3.2293746213)
(9.3519810128,3.2122008318)
(9.3754326461,3.1960560928)
(9.3985033117,3.1809082402)
(9.4212174570,3.1667248154)
(9.4435985781,3.1534736635)
(9.4656692673,3.1411233883)
(9.4874512649,3.1296436950)
(9.5089655113,3.1190056420)
(9.5302322004,3.1091818233)
(9.5512708333,3.1001464951)
(9.5721002710,3.0918756610)
(9.5927387866,3.0843471243)
(9.6132041162,3.0775405161)
(9.6335135094,3.0714373064)
(9.6536837774,3.0660208017)
(9.6737313413,3.0612761349)
(9.6936722785,3.0571902493)
(9.7135223677,3.0537518804)
(9.7332971341,3.0509515367)
(9.7530118923,3.0487814809)
(9.7726817902,3.0472357138)
(9.7923218502,3.0463099606)
(9.8119470120,3.0460016605)
(9.8315721739,3.0463099606)
(9.8512122339,3.0472357138)
(9.8708821317,3.0487814809)
(9.8905968900,3.0509515367)
(9.9103716563,3.0537518804)
(9.9302217456,3.0571902493)
(9.9501626827,3.0612761349)
(9.9702102467,3.0660208017)
(9.9903805147,3.0714373064)
(10.0106899078,3.0775405161)
(10.0311552375,3.0843471243)
(10.0517937531,3.0918756610)
(10.0726231908,3.1001464951)
(10.0936618237,3.1091818233)
(10.1149285128,3.1190056420)
(10.1364427591,3.1296436950)
(10.1582247567,3.1411233883)
(10.1802954460,3.1534736635)
(10.2026765670,3.1667248154)
(10.2253907124,3.1809082402)
(10.2484613779,3.1960560928)
(10.2719130112,3.2122008318)
(10.2957710548,3.2293746213)
(10.3200619826,3.2476085551)
(10.3448133265,3.2669316609)
(10.3700536883,3.2873696332)
(10.3958127315,3.3089432337)
(10.4221211467,3.3316662914)
(10.4490105784,3.3555432207)
(10.4765135024,3.3805659716)
(10.5046630347,3.4067103196)
(10.5334926501,3.4339314077)
(10.5630357794,3.4621584688)
(10.5933252486,3.4912886934)
(10.6243925116,3.5211802800)
(10.6562666185,3.5516448195)
(10.6889728553,3.5824393510)
(10.7225309855,3.6132586906)
(10.7569530354,3.6437289951)
(10.7922405960,3.6734039622)
(10.8283816760,3.7017655402)
(10.8653472528,3.7282313839)
(10.9030878230,3.7521713369)
(10.9415304468,3.7729346423)
(10.9805769541,3.7898880882)
(11.0201040600,3.8024628285)
(11.0599660109,3.8102045624)
(11.1000000000,3.8128191161)
(11.1400339891,3.8102045624)
(11.1798959400,3.8024628285)
(11.2194230459,3.7898880882)
(11.2584695532,3.7729346423)
(11.2969121770,3.7521713369)
(11.3346527472,3.7282313839)
(11.3716183240,3.7017655402)
(11.4077594040,3.6734039622)
(11.4430469646,3.6437289951)
(11.4774690145,3.6132586906)
(11.5110271447,3.5824393510)
(11.5437333815,3.5516448195)
(11.5756074884,3.5211802800)
(11.6066747514,3.4912886934)
(11.6369642206,3.4621584688)
(11.6665073499,3.4339314077)
(11.6953369653,3.4067103196)
(11.7234864976,3.3805659716)
(11.7509894216,3.3555432207)
(11.7778788533,3.3316662914)
(11.8041872685,3.3089432337)
(11.8299463117,3.2873696332)
(11.8551866735,3.2669316609)
(11.8799380174,3.2476085551)
(11.9042289452,3.2293746213)
(11.9280869888,3.2122008318)
(11.9515386221,3.1960560928)
(11.9746092876,3.1809082402)
(11.9973234330,3.1667248154)
(12.0197045540,3.1534736635)
(12.0417752433,3.1411233883)
(12.0635572409,3.1296436950)
(12.0850714872,3.1190056420)
(12.1063381763,3.1091818233)
(12.1273768092,3.1001464951)
(12.1482062469,3.0918756610)
(12.1688447625,3.0843471243)
(12.1893100922,3.0775405161)
(12.2096194853,3.0714373064)
(12.2297897533,3.0660208017)
(12.2498373173,3.0612761349)
(12.2697782544,3.0571902493)
(12.2896283437,3.0537518804)
(12.3094031100,3.0509515367)
(12.3291178683,3.0487814809)
(12.3487877661,3.0472357138)
(12.3684278261,3.0463099606)
(12.3880529880,3.0460016605)
(12.4076781498,3.0463099606)
(12.4273182098,3.0472357138)
(12.4469881077,3.0487814809)
(12.4667028659,3.0509515367)
(12.4864776323,3.0537518804)
(12.5063277215,3.0571902493)
(12.5262686587,3.0612761349)
(12.5463162226,3.0660208017)
(12.5664864906,3.0714373064)
(12.5867958838,3.0775405161)
(12.6072612134,3.0843471243)
(12.6278997290,3.0918756610)
(12.6487291667,3.1001464951)
(12.6697677996,3.1091818233)
(12.6910344887,3.1190056420)
(12.7125487351,3.1296436950)
(12.7343307327,3.1411233883)
(12.7564014219,3.1534736635)
(12.7787825430,3.1667248154)
(12.8014966883,3.1809082402)
(12.8245673539,3.1960560928)
(12.8480189872,3.2122008318)
(12.8718770307,3.2293746213)
(12.8961679585,3.2476085551)
(12.9209193025,3.2669316609)
(12.9461596642,3.2873696332)
(12.9719187075,3.3089432337)
(12.9982271226,3.3316662914)
(13.0251165543,3.3555432207)
(13.0526194783,3.3805659716)
(13.0807690107,3.4067103196)
(13.1095986260,3.4339314077)
(13.1391417553,3.4621584688)
(13.1694312245,3.4912886934)
(13.2004984875,3.5211802800)
(13.2323725944,3.5516448195)
(13.2650788312,3.5824393510)
(13.2986369614,3.6132586906)
(13.3330590113,3.6437289951)
(13.3683465719,3.6734039622)
(13.4044876519,3.7017655402)
(13.4414532287,3.7282313839)
(13.4791937990,3.7521713369)
(13.5176364227,3.7729346423)
(13.5566829300,3.7898880882)
(13.5962100359,3.8024628285)
(13.6360719869,3.8102045624)
(13.6761059759,3.8128191161)
} -- plot coordinates {
(13.6761059759,0.6871808839)
(13.6360719869,0.6897954376)
(13.5962100359,0.6975371715)
(13.5566829300,0.7101119118)
(13.5176364227,0.7270653577)
(13.4791937990,0.7478286631)
(13.4414532287,0.7717686161)
(13.4044876519,0.7982344598)
(13.3683465719,0.8265960378)
(13.3330590113,0.8562710049)
(13.2986369614,0.8867413094)
(13.2650788312,0.9175606490)
(13.2323725944,0.9483551805)
(13.2004984875,0.9788197200)
(13.1694312245,1.0087113066)
(13.1391417553,1.0378415312)
(13.1095986260,1.0660685923)
(13.0807690107,1.0932896804)
(13.0526194783,1.1194340284)
(13.0251165543,1.1444567793)
(12.9982271226,1.1683337086)
(12.9719187075,1.1910567663)
(12.9461596642,1.2126303668)
(12.9209193025,1.2330683391)
(12.8961679585,1.2523914449)
(12.8718770307,1.2706253787)
(12.8480189872,1.2877991682)
(12.8245673539,1.3039439072)
(12.8014966883,1.3190917598)
(12.7787825430,1.3332751846)
(12.7564014219,1.3465263365)
(12.7343307327,1.3588766117)
(12.7125487351,1.3703563050)
(12.6910344887,1.3809943580)
(12.6697677996,1.3908181767)
(12.6487291667,1.3998535049)
(12.6278997290,1.4081243390)
(12.6072612134,1.4156528757)
(12.5867958838,1.4224594839)
(12.5664864906,1.4285626936)
(12.5463162226,1.4339791983)
(12.5262686587,1.4387238651)
(12.5063277215,1.4428097507)
(12.4864776323,1.4462481196)
(12.4667028659,1.4490484633)
(12.4469881077,1.4512185191)
(12.4273182098,1.4527642862)
(12.4076781498,1.4536900394)
(12.3880529880,1.4539983395)
(12.3684278261,1.4536900394)
(12.3487877661,1.4527642862)
(12.3291178683,1.4512185191)
(12.3094031100,1.4490484633)
(12.2896283437,1.4462481196)
(12.2697782544,1.4428097507)
(12.2498373173,1.4387238651)
(12.2297897533,1.4339791983)
(12.2096194853,1.4285626936)
(12.1893100922,1.4224594839)
(12.1688447625,1.4156528757)
(12.1482062469,1.4081243390)
(12.1273768092,1.3998535049)
(12.1063381763,1.3908181767)
(12.0850714872,1.3809943580)
(12.0635572409,1.3703563050)
(12.0417752433,1.3588766117)
(12.0197045540,1.3465263365)
(11.9973234330,1.3332751846)
(11.9746092876,1.3190917598)
(11.9515386221,1.3039439072)
(11.9280869888,1.2877991682)
(11.9042289452,1.2706253787)
(11.8799380174,1.2523914449)
(11.8551866735,1.2330683391)
(11.8299463117,1.2126303668)
(11.8041872685,1.1910567663)
(11.7778788533,1.1683337086)
(11.7509894216,1.1444567793)
(11.7234864976,1.1194340284)
(11.6953369653,1.0932896804)
(11.6665073499,1.0660685923)
(11.6369642206,1.0378415312)
(11.6066747514,1.0087113066)
(11.5756074884,0.9788197200)
(11.5437333815,0.9483551805)
(11.5110271447,0.9175606490)
(11.4774690145,0.8867413094)
(11.4430469646,0.8562710049)
(11.4077594040,0.8265960378)
(11.3716183240,0.7982344598)
(11.3346527472,0.7717686161)
(11.2969121770,0.7478286631)
(11.2584695532,0.7270653577)
(11.2194230459,0.7101119118)
(11.1798959400,0.6975371715)
(11.1400339891,0.6897954376)
(11.1000000000,0.6871808839)
(11.0599660109,0.6897954376)
(11.0201040600,0.6975371715)
(10.9805769541,0.7101119118)
(10.9415304468,0.7270653577)
(10.9030878230,0.7478286631)
(10.8653472528,0.7717686161)
(10.8283816760,0.7982344598)
(10.7922405960,0.8265960378)
(10.7569530354,0.8562710049)
(10.7225309855,0.8867413094)
(10.6889728553,0.9175606490)
(10.6562666185,0.9483551805)
(10.6243925116,0.9788197200)
(10.5933252486,1.0087113066)
(10.5630357794,1.0378415312)
(10.5334926501,1.0660685923)
(10.5046630347,1.0932896804)
(10.4765135024,1.1194340284)
(10.4490105784,1.1444567793)
(10.4221211467,1.1683337086)
(10.3958127315,1.1910567663)
(10.3700536883,1.2126303668)
(10.3448133265,1.2330683391)
(10.3200619826,1.2523914449)
(10.2957710548,1.2706253787)
(10.2719130112,1.2877991682)
(10.2484613779,1.3039439072)
(10.2253907124,1.3190917598)
(10.2026765670,1.3332751846)
(10.1802954460,1.3465263365)
(10.1582247567,1.3588766117)
(10.1364427591,1.3703563050)
(10.1149285128,1.3809943580)
(10.0936618237,1.3908181767)
(10.0726231908,1.3998535049)
(10.0517937531,1.4081243390)
(10.0311552375,1.4156528757)
(10.0106899078,1.4224594839)
(9.9903805147,1.4285626936)
(9.9702102467,1.4339791983)
(9.9501626827,1.4387238651)
(9.9302217456,1.4428097507)
(9.9103716563,1.4462481196)
(9.8905968900,1.4490484633)
(9.8708821317,1.4512185191)
(9.8512122339,1.4527642862)
(9.8315721739,1.4536900394)
(9.8119470120,1.4539983395)
(9.7923218502,1.4536900394)
(9.7726817902,1.4527642862)
(9.7530118923,1.4512185191)
(9.7332971341,1.4490484633)
(9.7135223677,1.4462481196)
(9.6936722785,1.4428097507)
(9.6737313413,1.4387238651)
(9.6536837774,1.4339791983)
(9.6335135094,1.4285626936)
(9.6132041162,1.4224594839)
(9.5927387866,1.4156528757)
(9.5721002710,1.4081243390)
(9.5512708333,1.3998535049)
(9.5302322004,1.3908181767)
(9.5089655113,1.3809943580)
(9.4874512649,1.3703563050)
(9.4656692673,1.3588766117)
(9.4435985781,1.3465263365)
(9.4212174570,1.3332751846)
(9.3985033117,1.3190917598)
(9.3754326461,1.3039439072)
(9.3519810128,1.2877991682)
(9.3281229693,1.2706253787)
(9.3038320415,1.2523914449)
(9.2790806975,1.2330683391)
(9.2538403358,1.2126303668)
(9.2280812925,1.1910567663)
(9.2017728774,1.1683337086)
(9.1748834457,1.1444567793)
(9.1473805217,1.1194340284)
(9.1192309893,1.0932896804)
(9.0904013740,1.0660685923)
(9.0608582447,1.0378415312)
(9.0305687755,1.0087113066)
(8.9995015125,0.9788197200)
(8.9676274056,0.9483551805)
(8.9349211688,0.9175606490)
(8.9013630386,0.8867413094)
(8.8669409887,0.8562710049)
(8.8316534281,0.8265960378)
(8.7955123481,0.7982344598)
(8.7585467713,0.7717686161)
(8.7208062010,0.7478286631)
(8.6823635773,0.7270653577)
(8.6433170700,0.7101119118)
(8.6037899641,0.6975371715)
(8.5639280131,0.6897954376)
(8.5238940241,0.6871808839)
} -- cycle;
\fill[black!16] (8.5238940241,1.9500000000) rectangle (13.6761059759,2.5500000000);
\draw[line width=.8pt] plot coordinates {
(8.5238940241,3.8128191161)
(8.5639280131,3.8102045624)
(8.6037899641,3.8024628285)
(8.6433170700,3.7898880882)
(8.6823635773,3.7729346423)
(8.7208062010,3.7521713369)
(8.7585467713,3.7282313839)
(8.7955123481,3.7017655402)
(8.8316534281,3.6734039622)
(8.8669409887,3.6437289951)
(8.9013630386,3.6132586906)
(8.9349211688,3.5824393510)
(8.9676274056,3.5516448195)
(8.9995015125,3.5211802800)
(9.0305687755,3.4912886934)
(9.0608582447,3.4621584688)
(9.0904013740,3.4339314077)
(9.1192309893,3.4067103196)
(9.1473805217,3.3805659716)
(9.1748834457,3.3555432207)
(9.2017728774,3.3316662914)
(9.2280812925,3.3089432337)
(9.2538403358,3.2873696332)
(9.2790806975,3.2669316609)
(9.3038320415,3.2476085551)
(9.3281229693,3.2293746213)
(9.3519810128,3.2122008318)
(9.3754326461,3.1960560928)
(9.3985033117,3.1809082402)
(9.4212174570,3.1667248154)
(9.4435985781,3.1534736635)
(9.4656692673,3.1411233883)
(9.4874512649,3.1296436950)
(9.5089655113,3.1190056420)
(9.5302322004,3.1091818233)
(9.5512708333,3.1001464951)
(9.5721002710,3.0918756610)
(9.5927387866,3.0843471243)
(9.6132041162,3.0775405161)
(9.6335135094,3.0714373064)
(9.6536837774,3.0660208017)
(9.6737313413,3.0612761349)
(9.6936722785,3.0571902493)
(9.7135223677,3.0537518804)
(9.7332971341,3.0509515367)
(9.7530118923,3.0487814809)
(9.7726817902,3.0472357138)
(9.7923218502,3.0463099606)
(9.8119470120,3.0460016605)
(9.8315721739,3.0463099606)
(9.8512122339,3.0472357138)
(9.8708821317,3.0487814809)
(9.8905968900,3.0509515367)
(9.9103716563,3.0537518804)
(9.9302217456,3.0571902493)
(9.9501626827,3.0612761349)
(9.9702102467,3.0660208017)
(9.9903805147,3.0714373064)
(10.0106899078,3.0775405161)
(10.0311552375,3.0843471243)
(10.0517937531,3.0918756610)
(10.0726231908,3.1001464951)
(10.0936618237,3.1091818233)
(10.1149285128,3.1190056420)
(10.1364427591,3.1296436950)
(10.1582247567,3.1411233883)
(10.1802954460,3.1534736635)
(10.2026765670,3.1667248154)
(10.2253907124,3.1809082402)
(10.2484613779,3.1960560928)
(10.2719130112,3.2122008318)
(10.2957710548,3.2293746213)
(10.3200619826,3.2476085551)
(10.3448133265,3.2669316609)
(10.3700536883,3.2873696332)
(10.3958127315,3.3089432337)
(10.4221211467,3.3316662914)
(10.4490105784,3.3555432207)
(10.4765135024,3.3805659716)
(10.5046630347,3.4067103196)
(10.5334926501,3.4339314077)
(10.5630357794,3.4621584688)
(10.5933252486,3.4912886934)
(10.6243925116,3.5211802800)
(10.6562666185,3.5516448195)
(10.6889728553,3.5824393510)
(10.7225309855,3.6132586906)
(10.7569530354,3.6437289951)
(10.7922405960,3.6734039622)
(10.8283816760,3.7017655402)
(10.8653472528,3.7282313839)
(10.9030878230,3.7521713369)
(10.9415304468,3.7729346423)
(10.9805769541,3.7898880882)
(11.0201040600,3.8024628285)
(11.0599660109,3.8102045624)
(11.1000000000,3.8128191161)
(11.1400339891,3.8102045624)
(11.1798959400,3.8024628285)
(11.2194230459,3.7898880882)
(11.2584695532,3.7729346423)
(11.2969121770,3.7521713369)
(11.3346527472,3.7282313839)
(11.3716183240,3.7017655402)
(11.4077594040,3.6734039622)
(11.4430469646,3.6437289951)
(11.4774690145,3.6132586906)
(11.5110271447,3.5824393510)
(11.5437333815,3.5516448195)
(11.5756074884,3.5211802800)
(11.6066747514,3.4912886934)
(11.6369642206,3.4621584688)
(11.6665073499,3.4339314077)
(11.6953369653,3.4067103196)
(11.7234864976,3.3805659716)
(11.7509894216,3.3555432207)
(11.7778788533,3.3316662914)
(11.8041872685,3.3089432337)
(11.8299463117,3.2873696332)
(11.8551866735,3.2669316609)
(11.8799380174,3.2476085551)
(11.9042289452,3.2293746213)
(11.9280869888,3.2122008318)
(11.9515386221,3.1960560928)
(11.9746092876,3.1809082402)
(11.9973234330,3.1667248154)
(12.0197045540,3.1534736635)
(12.0417752433,3.1411233883)
(12.0635572409,3.1296436950)
(12.0850714872,3.1190056420)
(12.1063381763,3.1091818233)
(12.1273768092,3.1001464951)
(12.1482062469,3.0918756610)
(12.1688447625,3.0843471243)
(12.1893100922,3.0775405161)
(12.2096194853,3.0714373064)
(12.2297897533,3.0660208017)
(12.2498373173,3.0612761349)
(12.2697782544,3.0571902493)
(12.2896283437,3.0537518804)
(12.3094031100,3.0509515367)
(12.3291178683,3.0487814809)
(12.3487877661,3.0472357138)
(12.3684278261,3.0463099606)
(12.3880529880,3.0460016605)
(12.4076781498,3.0463099606)
(12.4273182098,3.0472357138)
(12.4469881077,3.0487814809)
(12.4667028659,3.0509515367)
(12.4864776323,3.0537518804)
(12.5063277215,3.0571902493)
(12.5262686587,3.0612761349)
(12.5463162226,3.0660208017)
(12.5664864906,3.0714373064)
(12.5867958838,3.0775405161)
(12.6072612134,3.0843471243)
(12.6278997290,3.0918756610)
(12.6487291667,3.1001464951)
(12.6697677996,3.1091818233)
(12.6910344887,3.1190056420)
(12.7125487351,3.1296436950)
(12.7343307327,3.1411233883)
(12.7564014219,3.1534736635)
(12.7787825430,3.1667248154)
(12.8014966883,3.1809082402)
(12.8245673539,3.1960560928)
(12.8480189872,3.2122008318)
(12.8718770307,3.2293746213)
(12.8961679585,3.2476085551)
(12.9209193025,3.2669316609)
(12.9461596642,3.2873696332)
(12.9719187075,3.3089432337)
(12.9982271226,3.3316662914)
(13.0251165543,3.3555432207)
(13.0526194783,3.3805659716)
(13.0807690107,3.4067103196)
(13.1095986260,3.4339314077)
(13.1391417553,3.4621584688)
(13.1694312245,3.4912886934)
(13.2004984875,3.5211802800)
(13.2323725944,3.5516448195)
(13.2650788312,3.5824393510)
(13.2986369614,3.6132586906)
(13.3330590113,3.6437289951)
(13.3683465719,3.6734039622)
(13.4044876519,3.7017655402)
(13.4414532287,3.7282313839)
(13.4791937990,3.7521713369)
(13.5176364227,3.7729346423)
(13.5566829300,3.7898880882)
(13.5962100359,3.8024628285)
(13.6360719869,3.8102045624)
(13.6761059759,3.8128191161)
};
\draw[line width=.8pt] plot coordinates {
(8.5238940241,0.6871808839)
(8.5639280131,0.6897954376)
(8.6037899641,0.6975371715)
(8.6433170700,0.7101119118)
(8.6823635773,0.7270653577)
(8.7208062010,0.7478286631)
(8.7585467713,0.7717686161)
(8.7955123481,0.7982344598)
(8.8316534281,0.8265960378)
(8.8669409887,0.8562710049)
(8.9013630386,0.8867413094)
(8.9349211688,0.9175606490)
(8.9676274056,0.9483551805)
(8.9995015125,0.9788197200)
(9.0305687755,1.0087113066)
(9.0608582447,1.0378415312)
(9.0904013740,1.0660685923)
(9.1192309893,1.0932896804)
(9.1473805217,1.1194340284)
(9.1748834457,1.1444567793)
(9.2017728774,1.1683337086)
(9.2280812925,1.1910567663)
(9.2538403358,1.2126303668)
(9.2790806975,1.2330683391)
(9.3038320415,1.2523914449)
(9.3281229693,1.2706253787)
(9.3519810128,1.2877991682)
(9.3754326461,1.3039439072)
(9.3985033117,1.3190917598)
(9.4212174570,1.3332751846)
(9.4435985781,1.3465263365)
(9.4656692673,1.3588766117)
(9.4874512649,1.3703563050)
(9.5089655113,1.3809943580)
(9.5302322004,1.3908181767)
(9.5512708333,1.3998535049)
(9.5721002710,1.4081243390)
(9.5927387866,1.4156528757)
(9.6132041162,1.4224594839)
(9.6335135094,1.4285626936)
(9.6536837774,1.4339791983)
(9.6737313413,1.4387238651)
(9.6936722785,1.4428097507)
(9.7135223677,1.4462481196)
(9.7332971341,1.4490484633)
(9.7530118923,1.4512185191)
(9.7726817902,1.4527642862)
(9.7923218502,1.4536900394)
(9.8119470120,1.4539983395)
(9.8315721739,1.4536900394)
(9.8512122339,1.4527642862)
(9.8708821317,1.4512185191)
(9.8905968900,1.4490484633)
(9.9103716563,1.4462481196)
(9.9302217456,1.4428097507)
(9.9501626827,1.4387238651)
(9.9702102467,1.4339791983)
(9.9903805147,1.4285626936)
(10.0106899078,1.4224594839)
(10.0311552375,1.4156528757)
(10.0517937531,1.4081243390)
(10.0726231908,1.3998535049)
(10.0936618237,1.3908181767)
(10.1149285128,1.3809943580)
(10.1364427591,1.3703563050)
(10.1582247567,1.3588766117)
(10.1802954460,1.3465263365)
(10.2026765670,1.3332751846)
(10.2253907124,1.3190917598)
(10.2484613779,1.3039439072)
(10.2719130112,1.2877991682)
(10.2957710548,1.2706253787)
(10.3200619826,1.2523914449)
(10.3448133265,1.2330683391)
(10.3700536883,1.2126303668)
(10.3958127315,1.1910567663)
(10.4221211467,1.1683337086)
(10.4490105784,1.1444567793)
(10.4765135024,1.1194340284)
(10.5046630347,1.0932896804)
(10.5334926501,1.0660685923)
(10.5630357794,1.0378415312)
(10.5933252486,1.0087113066)
(10.6243925116,0.9788197200)
(10.6562666185,0.9483551805)
(10.6889728553,0.9175606490)
(10.7225309855,0.8867413094)
(10.7569530354,0.8562710049)
(10.7922405960,0.8265960378)
(10.8283816760,0.7982344598)
(10.8653472528,0.7717686161)
(10.9030878230,0.7478286631)
(10.9415304468,0.7270653577)
(10.9805769541,0.7101119118)
(11.0201040600,0.6975371715)
(11.0599660109,0.6897954376)
(11.1000000000,0.6871808839)
(11.1400339891,0.6897954376)
(11.1798959400,0.6975371715)
(11.2194230459,0.7101119118)
(11.2584695532,0.7270653577)
(11.2969121770,0.7478286631)
(11.3346527472,0.7717686161)
(11.3716183240,0.7982344598)
(11.4077594040,0.8265960378)
(11.4430469646,0.8562710049)
(11.4774690145,0.8867413094)
(11.5110271447,0.9175606490)
(11.5437333815,0.9483551805)
(11.5756074884,0.9788197200)
(11.6066747514,1.0087113066)
(11.6369642206,1.0378415312)
(11.6665073499,1.0660685923)
(11.6953369653,1.0932896804)
(11.7234864976,1.1194340284)
(11.7509894216,1.1444567793)
(11.7778788533,1.1683337086)
(11.8041872685,1.1910567663)
(11.8299463117,1.2126303668)
(11.8551866735,1.2330683391)
(11.8799380174,1.2523914449)
(11.9042289452,1.2706253787)
(11.9280869888,1.2877991682)
(11.9515386221,1.3039439072)
(11.9746092876,1.3190917598)
(11.9973234330,1.3332751846)
(12.0197045540,1.3465263365)
(12.0417752433,1.3588766117)
(12.0635572409,1.3703563050)
(12.0850714872,1.3809943580)
(12.1063381763,1.3908181767)
(12.1273768092,1.3998535049)
(12.1482062469,1.4081243390)
(12.1688447625,1.4156528757)
(12.1893100922,1.4224594839)
(12.2096194853,1.4285626936)
(12.2297897533,1.4339791983)
(12.2498373173,1.4387238651)
(12.2697782544,1.4428097507)
(12.2896283437,1.4462481196)
(12.3094031100,1.4490484633)
(12.3291178683,1.4512185191)
(12.3487877661,1.4527642862)
(12.3684278261,1.4536900394)
(12.3880529880,1.4539983395)
(12.4076781498,1.4536900394)
(12.4273182098,1.4527642862)
(12.4469881077,1.4512185191)
(12.4667028659,1.4490484633)
(12.4864776323,1.4462481196)
(12.5063277215,1.4428097507)
(12.5262686587,1.4387238651)
(12.5463162226,1.4339791983)
(12.5664864906,1.4285626936)
(12.5867958838,1.4224594839)
(12.6072612134,1.4156528757)
(12.6278997290,1.4081243390)
(12.6487291667,1.3998535049)
(12.6697677996,1.3908181767)
(12.6910344887,1.3809943580)
(12.7125487351,1.3703563050)
(12.7343307327,1.3588766117)
(12.7564014219,1.3465263365)
(12.7787825430,1.3332751846)
(12.8014966883,1.3190917598)
(12.8245673539,1.3039439072)
(12.8480189872,1.2877991682)
(12.8718770307,1.2706253787)
(12.8961679585,1.2523914449)
(12.9209193025,1.2330683391)
(12.9461596642,1.2126303668)
(12.9719187075,1.1910567663)
(12.9982271226,1.1683337086)
(13.0251165543,1.1444567793)
(13.0526194783,1.1194340284)
(13.0807690107,1.0932896804)
(13.1095986260,1.0660685923)
(13.1391417553,1.0378415312)
(13.1694312245,1.0087113066)
(13.2004984875,0.9788197200)
(13.2323725944,0.9483551805)
(13.2650788312,0.9175606490)
(13.2986369614,0.8867413094)
(13.3330590113,0.8562710049)
(13.3683465719,0.8265960378)
(13.4044876519,0.7982344598)
(13.4414532287,0.7717686161)
(13.4791937990,0.7478286631)
(13.5176364227,0.7270653577)
(13.5566829300,0.7101119118)
(13.5962100359,0.6975371715)
(13.6360719869,0.6897954376)
(13.6761059759,0.6871808839)
};
\draw[densely dotted] (8.5238940241,0.6871808839)--(8.5238940241,3.8128191161);
\draw[densely dotted] (13.6761059759,0.6871808839)--(13.6761059759,3.8128191161);
\draw[guide] (8.5238940241,1.9500000000)--(13.6761059759,1.9500000000);
\draw[guide] (8.5238940241,2.5500000000)--(13.6761059759,2.5500000000);
\draw[axis] (8.3538940241,2.2500000000)--(13.9061059759,2.2500000000) node[right] {$\operatorname{Re}w$};
\node[left] at (8.5238940241,2.5500000000) {$\tau$};
\node[left] at (8.5238940241,1.9500000000) {$-\tau$};
\node[fill=white,inner sep=1pt] at (9.53,2.28) {$S_\tau$};
\node[fill=white,inner sep=1pt] at (11.1,2.99) {$G_{\eta,\rho}(S_\sigma)$};
\node at (11.1,.18) {$S_\tau\subset G_{\eta,\rho}(S_\sigma)$};
\draw[arr] (6.4,2.97)--(7.68,2.97) node[midway,above=3pt] {$G_{\eta,\rho}$};
\draw[arr] (7.68,1.45)--(6.4,1.45) node[midway,below=3pt] {$\Phi_{\eta,\rho}$};
\fill (3.7150000000,2.4300000000) circle[radius=1.6pt];\node[above right] at (3.7150000000,2.4300000000) {$z_0$};
\fill (11.8458314960,2.4206545972) circle[radius=1.6pt];\node[above right] at (11.8458314960,2.4206545972) {$w_0$};
\draw[line width=.3pt] (.38,-.24)--(14.05,-.24);
\node at (7.2,-.62) {$G_{\eta,\rho}(z_0)=w_0,\qquad
 \Phi_{\eta,\rho}:S_\tau\longrightarrow S_\sigma,\qquad
 |G_{\eta,\rho}(z)-z|\le B$};
\end{tikzpicture}
\caption{Common source and target strips in Lemma~\ref{op:lem:strips}.
The image boundaries in the right panel are a numerical illustration for
$\eta=0.4$, $\rho=2/7$, $\sigma=0.18$, and $\tau=0.05$;
the same magnification is used for both imaginary axes. Dotted sides mark
one displayed period, not boundaries of the infinite strips. The darker
horizontal band is $S_\tau$. The proof separately establishes univalence
on $S_\sigma$ and inclusion of a common $S_\tau$ in every image, including
for nearby complex parameters; the plotted widths are not certified
uniform constants for that family.}
\label{op:fig:common-strips}
\end{figure}

Figure~\ref{op:fig:common-strips} illustrates the distinct roles of the
source strip and the smaller target strip in this construction.

\subsection{Selected rows and exact moments}
\label{op:sec:arbitrary-band}

Fix a nonempty open interval $I\subset(0,1/2)$ and a reduced rational
$\rho_*=p/d\in I$. Choose an integer $N\ge3$ and a compact interval
$J\Subset I$ such that, for every $q\ge N$,
\begin{equation}\label{band:eq:selection}
 \rho_q=\frac{\lceil q p/d\rceil}{q}\in J.
\end{equation}
Such choices exist because $0\le\rho_q-\rho_*<1/q$. The reduced
denominator $d_q$ of $\rho_q$ divides $q$. If $q\equiv s\pmod d$, put
$r_s\in\{0,\ldots,d-1\}$ for the residue of $-ps$ modulo $d$; then
\[
 \rho_q=\rho_*+\frac{r_s}{dq}.
\]
Set
\[
 H_\eta(\rho,x)=G_{\eta,\rho}(\Phi_{\eta,\rho_*}(x)),
 \qquad Y_\eta(x)=\partial_\rho H_\eta(\rho_*,x).
\]
The common-strip result, Lemma~\ref{op:lem:strips}, gives
\begin{equation}\label{band:eq:amplitude}
 \begin{aligned}
 H_\eta(\rho_q,x)
 &=x+\frac{r_s}{dq}Y_\eta(x)+O(q^{-2}),\\
 H_\eta'(\rho_q,x)e^{\mp iq(H_\eta(\rho_q,x)-x)}
 &=e^{\mp ir_sY_\eta(x)/d}+O(q^{-1}),
 \end{aligned}
\end{equation}
uniformly on common smaller strips and compact eccentricity intervals.
Here the prime is the derivative in $x$.

For $\gamma\ge0$, use the Fourier sequence space
\[
 \mathcal X_\gamma
 =\{(u_n)_{n\in\mathbb Z}:
       \sup_n(1+|n|)^\gamma|u_n|<\infty\}.
\]
For $\gamma>2$ the Fourier series and its first derivative converge
absolutely, giving a continuous embedding into $C^1$.
For $\gamma=0$ this is a sequence space; every continuous function
supplies an element of it. Finite Fourier sums need not be norm dense,
so all matrix products below are justified by absolute summation.
Define $\widetilde A_\eta$ by setting its output coefficients to zero
for $|n|<N$ and, for $|n|=q\ge N$, setting
\begin{equation}\label{band:eq:operator}
 (\widetilde A_\eta U)_n
 =\frac1{2\pi}\int_0^{2\pi}
 U(x)e^{-inH_\eta(\rho_q,x)}
       \partial_xH_\eta(\rho_q,x)\,\mathrm dx.
\end{equation}
For general input sequences this formula means the corresponding
matrix row. That row has an analytic input kernel, so its Fourier sum
converges absolutely for bounded input coefficients. For continuous
inputs the row equals the displayed integral: approximate its analytic
kernel uniformly by its absolutely convergent Fourier series.

Let $I_0=\{-2,-1,0,1,2\}$ and define
\[
 t_{\eta,j}(x)
 =\frac{e^{ij\Phi_{\eta,\rho_*}(x)}}
        {D_\eta(\Phi_{\eta,\rho_*}(x))},\qquad j\in I_0.
\]
By \eqref{geo:eq:tangentspace}, these functions span the complexified
ellipse tangent space in action coordinates. Set
\begin{equation}\label{op:eq:moments}
 \mathsf G_{ij}(\eta)
 =\frac1{2\pi}\int_0^{2\pi}t_{\eta,j}t_{\eta,-i}\dd x,
 \qquad
 \ell_{\eta,i}(U)
 =\sum_{k\in I_0}(\mathsf G^{-1})_{ik}
       \frac1{2\pi}\int_0^{2\pi}U t_{\eta,-k}\dd x.
\end{equation}
For real parameters, $\mathsf G$ is a positive-definite Hermitian Gram
matrix. For complex parameters, the displayed formula defines its
holomorphic continuation without conjugating the parameter. It stays
invertible near every compact real interval. The definition gives $\ell_{\eta,i}(t_{\eta,j})=\delta_{ij}$.
For real $\eta$ and real $U$, vanishing of these five moments is
equivalent to the exact real tangent-moment normalization in
Section~\ref{geo:sec:gauge}.

\subsection{Parametrix and analytic kernels}

\smallskip\noindent\textbf{Localized-matrix calculus.}\phantomsection\label{op:par:localized-calculus}
We use the following elementary consequences of exponential localization.
Write $\langle n\rangle=1+|n|$. If
\[
 |M_{nj}|\le Ce^{-a|n-j|},\qquad
 |R_{nj}|\le C\langle n\rangle^{-1}e^{-a|n-j|},
 \qquad a>0,
\]
then $M$ is bounded on every $H^s$, $s\in\R$, and every
$\mathcal X_\gamma$, $\gamma\ge0$; $R$ gains one derivative on the
Sobolev scale and is compact on each of these spaces. Products preserve
exponential localization after decreasing $a$; a product containing $R$
retains its inverse output-frequency factor. These realizations have the
same absolutely convergent rows on polynomially growing Fourier
sequences. Locally uniform holomorphic entry bounds imply operator-norm
holomorphy on each fixed space.

For completeness, the weight inequalities
\[
 \frac{\langle n\rangle^s}{\langle j\rangle^s}
 \le\langle n-j\rangle^{|s|},\qquad
 \frac1{\langle j\rangle}\le
       \frac{\langle n-j\rangle}{\langle n\rangle}
\]
reduce the assertions to convolution with
$e^{-a|\ell|}\langle\ell\rangle^{|s|}$ in $\ell^2$ or weighted
$\ell^\infty$. The same sums justify matrix multiplication on
polynomially growing inputs. If $F_K$ truncates output frequencies,
\[
 \|(I-F_K)R\|_{H^s\to H^s}
 +\|(I-F_K)R\|_{\mathcal X_\gamma\to\mathcal X_\gamma}
 \le C_{s,\gamma}(1+K)^{-1}.
\]
Thus compactness is an operator-norm truncation statement, not a density
assertion in $\mathcal X_\gamma$. Parameter Cauchy estimates in the
same weighted sums prove norm holomorphy. We will also use Jaffard's
inverse-localization theorem \cite[Proposition~2]{Jaffard}: an
exponentially localized matrix invertible on $\ell^2(\Z)$ has an
exponentially localized inverse, with a possibly smaller exponent.

\begin{proposition}[The operator for an arbitrary interior band]
\label{band:prop:operator}
For each $\gamma\ge0$, $\widetilde A_\eta$ is a bounded,
locally norm-holomorphic operator on $\mathcal X_\gamma$. It has a
two-sided parametrix $Q_\eta$ whose errors are compact and satisfy
\begin{equation}\label{band:eq:error}
 |(Q_\eta\widetilde A_\eta-I)_{nj}|
 +|(\widetilde A_\eta Q_\eta-I)_{nj}|
 \le \frac{C}{1+|n|}e^{-\sigma|n-j|}.
\end{equation}
Its kernel consists of exponentially decaying Fourier sequences, with
a uniform exponential bound on each compact real eccentricity interval.
The same assertions hold if finitely many zero low rows are replaced by
the exact tangent-dual moment rows at the anchor $\rho_*$.
\end{proposition}

\begin{proof}
Fix the anchor $p/d$, the tail cutoff, and a compact real eccentricity
interval $I_\eta\Subset(0,1)$. All strip widths and bounds below are
uniform on this interval. The common complex neighborhoods are supplied
by Lemma~\ref{op:lem:strips}.

\smallskip\noindent
\emph{Principal amplitudes and geometric invertibility.}
For $|n|=q\ge N$, the selected row is
\[
 (\widetilde A_\eta U)_n=\widehat{b_{\eta,n}U}(n),\qquad
 b_{\eta,n}=\partial_xH_\eta(\rho_q,x)
             e^{-in(H_\eta(\rho_q,x)-x)}.
\]
Set $a_{\eta,s}^{\pm}=e^{\mp ir_sY_\eta/d}$, $0\le s<d$.
Equation~\eqref{band:eq:amplitude} gives
$b_{\eta,\pm q}=a_{\eta,s}^{\pm}+O(q^{-1})$ on a common strip when
$q\equiv s\pmod d$. These amplitudes need not be close to constants.
Let $\mathsf P_s$ project onto frequencies congruent to $s$ modulo $d$.
Multiplication precedes projection in
\begin{equation}\label{band:eq:principal}
 T_\eta^+=\sum_{s=0}^{d-1}\mathsf P_sM_{a_{\eta,s}^+},\qquad
 T_\eta^-=\sum_{s=0}^{d-1}\mathsf P_{-s}M_{a_{\eta,s}^-}.
\end{equation}
Contour shifting gives exponential localization of these matrices and
of $\widetilde A_\eta$, and for $|n|\ge N$,
\[
 |(\widetilde A_\eta-T_\eta^{\operatorname{sgn}n})_{nj}|
 \le C\langle n\rangle^{-1}e^{-a|n-j|}.
\]
Here $T_\eta^{\operatorname{sgn}n}$ means $T_\eta^+$ for $n>0$ and
$T_\eta^-$ for $n<0$. Any finite initial part outside the Taylor regime
is included by increasing $C$, without changing the prescribed cutoff.

Put $h=2\pi p/d$, $x_j=x+jh$, and let $s(r)$ invert the permutation
$r_s\equiv-ps\pmod d$. The root-of-unity formula for $\mathsf P_s$
gives, for $g=T_\eta^+U$ and $g_r=\mathsf P_{s(r)}g$,
\[
 g_r(x)=\frac1d\sum_{j=0}^{d-1}\xi_j(x)^rU(x_j),\qquad
 \xi_j(x)=\exp\!\left[\frac{i}{d}(2\pi j-Y_\eta(x_j))\right].
\]
For the physical caustic map $f_\rho$, differentiation of the conjugacy
at the fixed physical angle $\Phi_{\eta,\rho_*}(x_j)$ gives
\[
 \left.\partial_\rho f_\rho(\Phi_{\eta,\rho_*}(x_j))\right|_{\rho_*}
 =\Phi_{\eta,\rho_*}'(x_{j+1})
      [2\pi+Y_\eta(x_j)-Y_\eta(x_{j+1})]>0.
\]
The positivity uses the twist and $\partial_\rho\lambda>0$ from
Lemma~\ref{op:lem:strips}. Thus the lifted arguments
$\vartheta_j=(2\pi j-Y_\eta(x_j))/d$ strictly increase, and
$\vartheta_d-\vartheta_0=2\pi$. The nodes are distinct, so
\[
 \mathsf V_{rj}=\xi_j^r,\qquad
 \det\mathsf V=\prod_{j<k}(\xi_k-\xi_j)\ne0.
\]
Compactness and holomorphy make $\mathsf V^{-1}$ uniformly bounded on
smaller common complex parameter neighborhoods and angle strips.

The inverse on smooth functions is
\[
 (T_\eta^+)^{-1}g(x)
 =d\sum_{r=0}^{d-1}(\mathsf V(x)^{-1})_{0r}\mathsf P_{s(r)}g(x).
\]
To check that this pointwise fiber solution represents a single periodic
function, note
\[
 \mathsf V_{rj}(x+h)=e^{-2\pi ir/d}\mathsf V_{r,j+1}(x),\qquad
 g_r(x+h)=e^{-2\pi ir/d}g_r(x),
\]
with cyclic columns. Uniqueness of $v=d\mathsf V^{-1}(g_r)_r$ implies
$v_j(x+h)=v_{j+1}(x)$; also $v$ is $2\pi$-periodic. Hence
$U(x)=v_0(x)$ has exactly these fiber values and $T_\eta^+U=g$.
The reverse identity follows by inverting the fiber equation of an
existing $U$. Replacing $\xi_j^r$ by $\xi_j^{-r}$ and $\mathsf P_{s(r)}$
by $\mathsf P_{-s(r)}$ proves the analogous assertions for $T_\eta^-$.
Both inverses are finite sums of analytic multipliers followed by
residue projections, and therefore have exponentially localized
matrices. This argument also covers $p>1$.
Figure~\ref{op:fig:vandermonde-nodes} illustrates the rational fiber and
the separated auxiliary nodes used in this inversion.

\begin{figure}[tbp]
\centering
\begin{tikzpicture}[
 x=1cm,y=1cm,font=\small,
 axis/.style={-{Stealth[length=1.6mm]},line width=.45pt},
 arr/.style={-{Stealth[length=1.8mm]},line width=.65pt},
 guide/.style={densely dashed,line width=.5pt},
 every node/.style={inner sep=2pt}]
\path[use as bounding box] (0,-2.05) rectangle (15.2,5.15);
\node at (3.3,4.90) {(a) Rational action-coordinate fiber};
\node at (11.25,4.90) {(b) Twist-separated auxiliary nodes};
\draw[line width=.7pt] (3.3000000000,2.5200000000) circle[radius=1.57];
\draw[line width=.7pt] (11.2500000000,2.5200000000) circle[radius=1.57];
\draw[black!38,line width=.45pt,postaction={decorate},decoration={markings,mark=at position .55 with {\arrow{Stealth[length=1.25mm]}}}] (4.8211925021,2.9084242160)--(2.5828172137,3.9166205108);
\draw[black!38,line width=.45pt,postaction={decorate},decoration={markings,mark=at position .55 with {\arrow{Stealth[length=1.25mm]}}}] (2.5828172137,3.9166205108)--(2.0979838647,1.5100211831);
\draw[black!38,line width=.45pt,postaction={decorate},decoration={markings,mark=at position .55 with {\arrow{Stealth[length=1.25mm]}}}] (2.0979838647,1.5100211831)--(4.5521302924,1.5728623485);
\draw[black!38,line width=.45pt,postaction={decorate},decoration={markings,mark=at position .55 with {\arrow{Stealth[length=1.25mm]}}}] (4.5521302924,1.5728623485)--(3.9447657311,3.9514947265);
\draw[black!38,line width=.45pt,postaction={decorate},decoration={markings,mark=at position .55 with {\arrow{Stealth[length=1.25mm]}}}] (3.9447657311,3.9514947265)--(1.7609219623,2.8300625645);
\draw[black!38,line width=.45pt,postaction={decorate},decoration={markings,mark=at position .55 with {\arrow{Stealth[length=1.25mm]}}}] (1.7609219623,2.8300625645)--(3.3401884337,0.9505144506);
\draw[black!38,line width=.45pt,postaction={decorate},decoration={markings,mark=at position .55 with {\arrow{Stealth[length=1.25mm]}}}] (3.3401884337,0.9505144506)--(4.8211925021,2.9084242160);
\fill (4.8211925021,2.9084242160) circle[radius=1.8pt];\node[fill=white,inner sep=1pt] at (5.1409336013,2.9900675226) {$x_0$};
\fill (12.8123394541,2.3650954161) circle[radius=1.8pt];\node[fill=white,inner sep=1pt] at (13.1407292757,2.3325358539) {$\xi_0$};
\fill (2.5828172137,3.9166205108) circle[radius=1.8pt];\node[fill=white,inner sep=1pt] at (2.4320717873,4.2101776882) {$x_1$};
\fill (12.0592775126,3.8653512209) circle[radius=1.8pt];\node[fill=white,inner sep=1pt] at (12.2293804293,4.1481320508) {$\xi_1$};
\fill (2.0979838647,1.5100211831) circle[radius=1.8pt];\node[fill=white,inner sep=1pt] at (1.8453307917,1.2977326420) {$x_2$};
\fill (11.1794000199,4.0884118218) circle[radius=1.8pt];\node[fill=white,inner sep=1pt] at (11.1645605337,4.4180780009) {$\xi_2$};
\fill (4.5521302924,1.5728623485) circle[radius=1.8pt];\node[fill=white,inner sep=1pt] at (4.8153169143,1.3737824599) {$x_3$};
\fill (9.7371305188,2.9396736027) circle[radius=1.8pt];\node[fill=white,inner sep=1pt] at (9.4191388444,3.0278852517) {$\xi_3$};
\fill (3.9447657311,3.9514947265) circle[radius=1.8pt];\node[fill=white,inner sep=1pt] at (4.0802897382,4.2523821531) {$x_4$};
\fill (9.7641506193,2.0128988090) circle[radius=1.8pt];\node[fill=white,inner sep=1pt] at (9.4518383291,1.9063106606) {$\xi_4$};
\fill (1.7609219623,2.8300625645) circle[radius=1.8pt];\node[fill=white,inner sep=1pt] at (1.4374214831,2.8952349507) {$x_5$};
\fill (11.0248324955,0.9662305207) circle[radius=1.8pt];\node[fill=white,inner sep=1pt] at (10.9775042939,0.6396420314) {$\xi_5$};
\fill (3.3401884337,0.9505144506) circle[radius=1.8pt];\node[fill=white,inner sep=1pt] at (3.3486356841,0.6206225835) {$x_6$};
\fill (12.2384410542,1.3002113780) circle[radius=1.8pt];\node[fill=white,inner sep=1pt] at (12.4462025497,1.0438226868) {$\xi_6$};
\draw[black!65,-{Stealth[length=1.3mm]},line width=.5pt] (12.6595998015,2.4864083387) arc[start angle=-1.3651365462,end angle=53.2419888124,radius=1.4100000000000001];
\draw[black!65,-{Stealth[length=1.3mm]},line width=.5pt] (11.8842267233,3.7793079304) arc[start angle=63.2687502272,end angle=86.8477758175,radius=1.4100000000000001];
\draw[black!65,-{Stealth[length=1.3mm]},line width=.5pt] (11.0812291525,3.9198629937) arc[start angle=96.8745372323,end angle=158.7662806730,radius=1.4100000000000001];
\draw[black!65,-{Stealth[length=1.3mm]},line width=.5pt] (9.8668864996,2.7940384004) arc[start angle=168.7930420878,end angle=193.1145319056,radius=1.4100000000000001];
\draw[black!65,-{Stealth[length=1.3mm]},line width=.5pt] (9.9534507159,1.9658700930) arc[start angle=203.1412933204,end angle=256.0246963683,radius=1.4100000000000001];
\draw[black!65,-{Stealth[length=1.3mm]},line width=.5pt] (11.1529066587,1.1133469216) arc[start angle=266.0514577831,end angle=303.2895858752,radius=1.4100000000000001];
\draw[black!65,-{Stealth[length=1.3mm]},line width=.5pt] (12.2172966160,1.4941163532) arc[start angle=313.3163472900,end angle=348.6081020390,radius=1.4100000000000001];
\draw[densely dotted,black!65,line width=.45pt] (11.2500000000,2.5200000000)--(12.0592775126,3.8653512209);
\draw[densely dotted,black!65,line width=.45pt] (11.2500000000,2.5200000000)--(11.1794000199,4.0884118218);
\node[align=center,fill=white,inner sep=3pt] at (11.25,2.45) {$\xi_j=e^{i\vartheta_j}$\\$\vartheta_{j+1}>\vartheta_j$};
\draw[arr] (6.15,2.55)--(8.35,2.55) node[midway,above=4pt,align=center] {principal rows};
\node at (3.3,.28) {$x_j=x+2\pi pj/d\pmod{2\pi}$};
\node at (11.25,.28) {$\vartheta_d-\vartheta_0=2\pi$};
\draw[line width=.3pt] (.45,-.14)--(14.48,-.14);
\node at (7.45,-.76) {$\displaystyle g_r(x)=\frac1d\sum_{j=0}^{d-1}\xi_j(x)^r U(x_j),
 \qquad \mathsf V_{rj}(x)=\xi_j(x)^r$};
\node at (7.45,-1.62) {$\det\mathsf V(x)=\prod_{j<k}(\xi_k(x)-\xi_j(x))\ne0$};
\end{tikzpicture}
\caption{Rational fibers and the Vandermonde nodes for $\eta=0.4$,
$\rho_*=p/d=2/7$, and $x=0.25$. In the left panel the points are plotted
on the action-coordinate circle, with arrows following increasing $j$;
the connecting segments are not physical billiard chords. The right panel
shows the corresponding computed auxiliary nodes $\xi_j$. Their unequal
angular gaps distinguish this fiber from the equally spaced Fourier
matrix. Twist makes every lifted gap positive and their sum $2\pi$,
so the nodes are distinct. The proof also checks cyclic compatibility of
the inverse fibers; a single plotted fiber is not a uniform invertibility
certificate.}
\label{op:fig:vandermonde-nodes}
\end{figure}
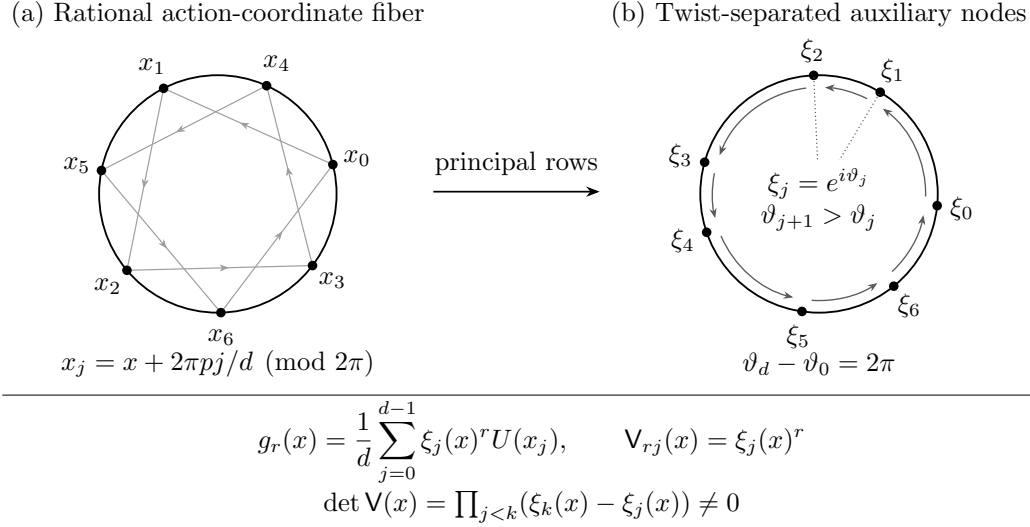

\smallskip\noindent
\emph{Parametrix and finite-row insertions.}
Let $H_\pm$ project onto strictly positive and strictly negative
frequencies, and $P_0$ onto the zero frequency. Set
\[
 S_\eta=P_0+H_+T_\eta^+H_++H_-T_\eta^-H_-,\qquad
 Q_\eta=P_0+H_+(T_\eta^+)^{-1}H_++H_-(T_\eta^-)^{-1}H_-.
\]
For $T=T_\eta^+$ and $H=H_+$,
\[
 HTHT^{-1}H-H=-HT(I-H)T^{-1}H,\qquad
 HT^{-1}HTH-H=-HT^{-1}(I-H)TH.
\]
For $n,j>0$ the intermediate indices satisfy $k\le0$, giving the bound
\[
 C\sum_{k\le0}e^{-a|n-k|}e^{-a|k-j|}\le C'e^{-a(n+j)}.
\]
The negative block is identical with signs reversed.
Thus the two compression errors decay like
$e^{-a(|n|+|j|)}$. On the other hand,
$R_\eta=\widetilde A_\eta-S_\eta$ satisfies
\[
 |(R_\eta)_{nj}|\le C\langle n\rangle^{-1}e^{-a_1|n-j|}
\]
for some $a_1>0$: same-sign high rows use the amplitude estimate,
whereas opposite-sign and finitely many low rows have the stronger
absolute-frequency decay. The identities
\[
 Q_\eta\widetilde A_\eta-I=(Q_\eta S_\eta-I)+Q_\eta R_\eta,\qquad
 \widetilde A_\eta Q_\eta-I=(S_\eta Q_\eta-I)+R_\eta Q_\eta
\]
and the \hyperref[op:par:localized-calculus]{localized-matrix calculus}
prove \eqref{band:eq:error}, boundedness, norm holomorphy, and compactness.
All products are absolutely convergent on the asserted sequence spaces.

The moment kernels in \eqref{op:eq:moments} are analytic: their
geometric denominator is separated from zero on a smaller strip and
the Gram matrix has a uniformly holomorphic inverse. An insertion
$\mathcal F_\eta$ in finitely many zero rows satisfies
$|(\mathcal F_\eta)_{nj}|\le Ce^{-a_2(|n|+|j|)}$.
The same is true, at a smaller rate, of $Q_\eta\mathcal F_\eta$ and
$\mathcal F_\eta Q_\eta$, so the same parametrix applies to
$\mathscr A_\eta=\widetilde A_\eta+\mathcal F_\eta$.
For real parameters, function conjugation is
$(\mathcal CU)_n=\overline{u_{-n}}$; opposite selected rows are conjugate
and $\ell_{\eta,-i}(\mathcal CU)=\overline{\ell_{\eta,i}(U)}$.
Thus inserting $\ell_{\eta,i}$ at index $i$, $|i|\le2$, preserves
reality. No such assertion is needed for an arbitrary unpaired insertion.

\smallskip\noindent
\emph{Analyticity of the kernel.}
Write $B_\eta=Q_\eta\mathscr A_\eta-I$. On the closed kernel,
$B_\eta=-I$; compactness makes that kernel finite dimensional.
To obtain analyticity, fix $0<\tau<\sigma$, where $\sigma$ is an
error-decay rate, and use
$\|u\|_{\mathcal E_\tau}=\sup_n e^{\tau|n|}|u_n|$.
With $H_K=I-F_K$,
\[
 \|H_KB_\eta H_K\|_{\mathcal E_\tau\to\mathcal E_\tau}
 \le\frac{C}{1+K}\sum_{\ell\in\Z}e^{-(\sigma-\tau)|\ell|}.
\]
Its norm on $\mathcal X_0$ is likewise $O(K^{-1})$. Choose one $K$
so both norms are at most $1/2$, uniformly on $I_\eta$. A kernel vector
satisfies
\[
 (I+H_KB_\eta H_K)H_KU=-H_KB_\eta F_KU,
 \qquad
 \|H_KB_\eta F_KU\|_{\mathcal E_\tau}
 \le C_K\|U\|_{\mathcal X_0}.
\]
Solve first in $\mathcal E_\tau$ by a Neumann series. The solution is
in $\mathcal X_0$, where uniqueness identifies it with the original
tail. Including the finite low block gives
\[
 \sup_n e^{\tau|n|}|u_n|
 \le C_{\rm ker}\|U\|_{\mathcal X_0},\qquad \mathscr A_\eta U=0.
\]
Its Fourier series converges normally on every smaller strip. The
constants are uniform on $I_\eta$, with the finite insertion fixed;
no injectivity or classification of the kernel has been used.
\end{proof}

\begin{corollary}[Sobolev parametrices and analytic distributional kernels]
\label{band:cor:sobolev-kernels}
For every real $s$, the operators in Proposition~\ref{band:prop:operator}
act boundedly on the periodic Sobolev space $H^s(\mathbb T)$, and their
parametrix errors map $H^s$ boundedly into $H^{s+1}$. The errors are compact as operators on $H^s$. The bounds and
operator-norm continuity are locally uniform in the eccentricity.
Every distributional vector in the kernel of $\widetilde A_\eta$
is represented by an analytic periodic function. These statements also hold after inserting the five
tangent-dual moment rows.
\end{corollary}
\begin{proof}
Let $\mathscr A_\eta$ denote either allowed operator and put
$B_\eta^L=Q_\eta\mathscr A_\eta-I$ and
$B_\eta^R=\mathscr A_\eta Q_\eta-I$.
Proposition~\ref{band:prop:operator} gives common exponential bounds
for $\mathscr A_\eta,Q_\eta$, and the additional factor
$\langle n\rangle^{-1}$ for both errors, locally uniformly on complex
parameter neighborhoods. The
\hyperref[op:par:localized-calculus]{localized-matrix calculus}
therefore gives boundedness on every $H^s$, the one-derivative gain,
compactness, and operator-norm continuity. The parametrix identities
hold first on Fourier modes and then on every $H^s$; their row sums
are absolutely convergent and agree on intersections of these spaces.

If $U\in\mathcal D'(\T)$, its Fourier coefficients have some polynomial
bound $|u_n|\le C_U\langle n\rangle^m$; see, for example,
\cite[Chapters~2 and 12]{Dyatlov}. Thus $U\in H^{-m-1}$.
The kernel equation and the compatible left parametrix give
$U=-B_\eta^LU\in H^{-m}$, and $m+1$ iterations place the same
kernel vector in $L^2$. Its Fourier sequence then belongs to
$\mathcal X_0$, with $\|U\|_{\mathcal X_0}\le\|U\|_{L^2}$.
The exponential kernel estimate in Proposition~\ref{band:prop:operator}
makes $\sum_nu_ne^{inx}$ normally convergent on a strip. Fourier
uniqueness identifies its boundary restriction with the original
distribution. This proves analyticity, not merely smoothness, for both
the unmodified and moment-modified operator.
\end{proof}

\section{Jacobi geometry and the arithmetic obstruction}\label{arith:section}
Fix a noncircular ellipse with major semiaxis $a=1$, and write
$\eta=1-b^2\in(0,1)$. Throughout this section, $\eta$ is fixed.

In the continuation argument, the Jacobi covering relates the action
and boundary-parameter planes, and the focal fibers identify the possible
one-step companions. The finite-grid obstruction and Picard--Fuchs source
identities are used together later to exclude persistent nonfocal
collisions. The two-sheeted geometry concerns the coordinate map, not the
possible monodromy of the unknown displacement.

The caustic parameter $m\in(\eta,1)$ and rotation number $\rho(m)$
are related by
\[
J(m)=F\left(\arcsin\sqrt{\frac{m-\eta}{m(1-\eta)}}\,\middle|\,m\right)
=2K(m)\rho(m),\qquad h(m)=2J(m).
\]
The physical Jacobi map is
$\mathcal P_m(u)=\sn(u\mid m)-i\cn(u\mid m)$. We abbreviate $K(m)$ by $K$
and use the following conventions throughout the continuation proof.
\begin{center}
\begin{tabular}{@{}ll@{}}
\toprule
Notation & Meaning or relation\\
\midrule
$\varphi$ & Physical ellipse angle\\
$x$ & Caustic action angle, $\varphi=\Phi_{E,\rho}(x)$\\
$z$ & Physical complex coordinate, $z=e^{i\varphi}=\mathcal P_m(u)$\\
$\zeta$ & Action complex coordinate, $\zeta=e^{ix}$\\
$u$ & Jacobi coordinate, $u=K+2Kx/\pi$\\
$m\in(\eta,1)$ & Elliptic parameter $m=k^2$; $k$ is the modulus\\
$K'=K(1-m)$ & Complementary complete elliptic integral\\
$K_m=dK/dm$ & Derivative of the complete elliptic integral\\
\bottomrule
\end{tabular}
\end{center}

For an analytic displacement $v$ and a chosen rational base $\rho_*$,
the periodic function
$U(x)=v(\Phi_{E,\rho_*}(x))$ also defines a holomorphic germ in
$\zeta=e^{ix}$, denoted by the same letter. The physical germ $V$
satisfies $V(e^{i\varphi})=v(\varphi)$. Initially these two germs are
defined on annuli about their respective unit circles.

\begin{lemma}[Jacobi covering geometry]\label{arith:jacobi}
The map $\mathcal P_m$ has degree two on the torus with lattice
$4K\mathbb Z+4iK'\mathbb Z$. Its deck involution is
$u\mapsto2K+2iK'-u$. Its lower cylinder
$-2K'<\operatorname{Im}u<0$, modulo $4K\mathbb Z$, maps properly with
degree two onto $\{|z|>1\}\cup\{\infty\}$.
\end{lemma}
\begin{proof}
Fix $m\in(\eta,1)$ and suppress it in the Jacobi functions. All square
roots of positive real parameters below are positive. Put
$\Lambda=4K\Z+4iK'\Z$, where $K'=K(1-m)$ is the complementary
integral, not a parameter derivative.

\smallskip\noindent
\emph{Step 1: identify the full period lattice.}
The identity $\sn^2u+\cn^2u=1$ and the standard translations and parity
formulas give, as meromorphic identities,
\[
 (\sn u-i\cn u)(\sn u+i\cn u)=1,
\]
\[
 \mathcal P_m(u+2K)=-\mathcal P_m(u),\qquad
 \mathcal P_m(u+2iK')=\mathcal P_m(u)^{-1},\qquad
 \mathcal P_m(-u)=-\mathcal P_m(u)^{-1}.
\]
Here we use the modulus convention $k=\sqrt m$ in
\cite[Table~22.4.3 and Eq.~22.6.1]{DLMF}.
In particular, $4K$ and $4iK'$ are periods. Conversely,
\[
 \sn u=\frac{\mathcal P_m(u)+\mathcal P_m(u)^{-1}}2,\qquad
 \cn u=\frac{\mathcal P_m(u)^{-1}-\mathcal P_m(u)}{2i},
\]
so every period of $\mathcal P_m$ is a common period of $\sn$ and $\cn$.
Their full period lattices are given by
\cite[Table~22.4.1]{DLMF}:
\[
 \Lambda_{\mathrm{sn}}=4K\Z+2iK'\Z,\qquad
 \Lambda_{\mathrm{cn}}=4K\Z+(2K+2iK')\Z.
\]
If $4Kr+2iK's=4Kt+(2K+2iK')v$ with integer coefficients, comparison
of real and imaginary parts gives $s=v=2(r-t)$. Thus
$\Lambda_{\mathrm{sn}}\cap\Lambda_{\mathrm{cn}}=\Lambda$.
This proves that $\Lambda$ is the full period lattice of $\mathcal P_m$.
It therefore defines a meromorphic map on the torus
$\mathcal T_m=\C/\Lambda$.

\smallskip\noindent
\emph{Step 2: count all zeros and poles, including the cancellations.}
The half-period formulas in \cite[Table~22.4.3]{DLMF} give
\[
 \sn(t+iK')=\frac1{\sqrt m\,\sn t},\qquad
 \cn(t+iK')=-\frac{i\dn t}{\sqrt m\,\sn t}.
\]
Together with the translation by $2iK'$, they imply
\[
 \mathcal P_m(t+iK')=\frac{1-\dn t}{\sqrt m\,\sn t},\qquad
 \mathcal P_m(t-iK')=\frac{1+\dn t}{\sqrt m\,\sn t}.
\]
Using $\sn t=t+O(t^3)$ and $\dn t=1-mt^2/2+O(t^4)$ as $t\to0$,
we obtain
\[
 \mathcal P_m(t+iK')=\frac{\sqrt m}{2}t+O(t^3),\qquad
 \mathcal P_m(t-iK')=\frac{2}{\sqrt m\,t}+O(t).
\]
Thus the common Jacobi poles cancel at $iK'$, leaving a simple zero,
whereas $-iK'$ is a simple pole. The shift by $2K$ gives the same
orders at $2K+iK'$ and $2K-iK'$, respectively.

These four points exhaust the possible singularities in a fundamental
cell: the common poles of $\sn$ and $\cn$ are
$2rK+(2s+1)iK'$, $r,s\in\Z$
\cite[Table~22.4.1]{DLMF}. No pole of their linear combination can
occur elsewhere. Nor can $\mathcal P_m$ vanish at a point where both
Jacobi functions are finite, by the product identity in Step~1.
Consequently, modulo $\Lambda$, the complete lists are
\[
 \begin{array}{c|c}
 \text{simple zeros}&iK',\ 2K+iK'\\
 \text{simple poles}&-iK',\ 2K-iK'.
 \end{array}
\]
A nonconstant meromorphic map from a compact connected Riemann surface
to the sphere is onto, and its degree is the total pole order.
Hence $\mathcal P_m:\mathcal T_m\to\widehat{\C}$ has degree two.

\smallskip\noindent
\emph{Step 3: identify the involution and the ramification.}
The affine map $\iota(u)=2K+2iK'-u$ descends to a nontrivial
holomorphic involution of $\mathcal T_m$. Step~1 gives
\[
 \mathcal P_m(\iota(u))
 =-\mathcal P_m(2iK'-u)
 =-\mathcal P_m(-u)^{-1}
 =\mathcal P_m(u).
\]
At a fixed point of $\iota$, a local coordinate $t$ centered there
makes $\iota$ act by $t\mapsto-t$. In a local target coordinate the
map $\mathcal P_m$ is therefore even and nonconstant, so its local
degree is at least two, and hence exactly two. Conversely, a
ramification point distinct from its involution partner would give a
fiber of total multiplicity at least three, contradicting degree two.
Thus the fixed points are precisely the ramification points, and
$\iota$ exchanges the two points of every regular fiber. The degree-two
map is a branched covering; it is an ordinary two-sheeted covering
after its branch values and their preimages are removed.

\smallskip\noindent
\emph{Step 4: identify the full preimage of the spherical exterior.}
On $\mathcal T_m$, consider the two cylinders
\[
 \mathcal C_+=\{0<\operatorname{Im}u<2K'\}/(4K\Z),\qquad
 \mathcal C_-=\{-2K'<\operatorname{Im}u<0\}/(4K\Z).
\]
For real $u$, the Jacobi functions are real and
$|\mathcal P_m(u)|^2=\sn^2u+\cn^2u=1$. The reciprocal translation
identity gives modulus one also on the circles
$\operatorname{Im}u=\pm2K'$.
The function $\mathcal P_m$ has no poles on the closed upper cylinder,
and its zeros lie in its interior. On the compact cylinder obtained
by identifying the real-period sides, the maximum-modulus principle
and its strict form therefore give
\[
 |\mathcal P_m(u)|<1\qquad(u\in\mathcal C_+).
\]
On the closed lower cylinder, $1/\mathcal P_m$ is holomorphic, with
zeros at the two poles of $\mathcal P_m$ and boundary modulus one.
The same argument gives
\[
 |\mathcal P_m(u)|>1\qquad(u\in\mathcal C_-),
\]
with the poles interpreted as taking the value $\infty$.
The two cylinders and their boundary circles exhaust the torus.
For the spherical exterior
$\mathcal D_{\mathrm{ext}}=\{|z|>1\}\cup\{\infty\}$, we have therefore
proved the full-preimage identity
\[
 \mathcal P_m^{-1}(\mathcal D_{\mathrm{ext}})=\mathcal C_-.
\]

\smallskip\noindent
\emph{Step 5: restrict the proper map and retain its degree.}
For every compact set $L\subset\mathcal D_{\mathrm{ext}}$, its
preimage under the torus map is compact and lies entirely in
$\mathcal C_-$. Thus the restriction to $\mathcal C_-$ is proper.
It is onto and contains the entire fiber of every exterior value,
so its degree is still two, counted with multiplicity. This proves
the asserted exterior mapping property. When $\infty$ is omitted
from the target, its two preimages $-iK'$ and $2K-iK'$ must likewise
be omitted from the source.
\end{proof}

In the physical coordinate $z=e^{i\varphi}$, we write
$D_E(z)=1-\eta(z+z^{-1})^2/4$ for the Laurent extension of
$D_E(\varphi)$.

\begin{lemma}[Focal fibers]\label{arith:focal}
With $a=1$, write $b=\sqrt{1-\eta}$ and $c=\sqrt\eta>0$.
The exterior zeros of $D_E$ are $\pm(a+b)/c$; we call them the
exterior focal values. Their lower-cylinder preimages are, modulo
the period lattice, the pairs
\[
K\pm J-iK',\qquad 3K\pm J-iK'.
\]
Each pair is the full fiber of one exterior focal value.
Its two points are unramified for $\eta<m<1$ and differ by the real
shift $2J=h$.
\end{lemma}
\begin{proof}
\smallskip\noindent
\emph{Step 1: locate the exterior zeros in the boundary-parameter plane.}
Since $a=1$, the equation $D_E(z)=0$ is
\[
 z+z^{-1}=\pm\frac2c,\qquad c=\sqrt\eta.
\]
Its four roots are
\[
 \pm\frac{1+b}{c},\qquad \pm\frac{1-b}{c},
 \qquad b=\sqrt{1-\eta}.
\]
The positive roots satisfy
$((1+b)/c)((1-b)/c)=1$, with $(1+b)/c>1$.
Hence precisely
\[
 z_+=\frac{1+b}{c},\qquad z_-=-\frac{1+b}{c}
\]
lie in the exterior. These are values of the complexified boundary
coordinate $z=e^{i\varphi}$, not Cartesian coordinates of the foci.

\smallskip\noindent
\emph{Step 2: evaluate the two proposed focal pairs.}
Set
\[
 \vartheta=\arcsin\sqrt{\frac{m-\eta}{m(1-\eta)}}.
\]
For $\eta<m<1$, the expression under the square root lies strictly
between zero and one, so $0<\vartheta<\pi/2$ and
$0<J=F(\vartheta\mid m)<K$. The inverse-integral definition and the
Jacobi square identities give
\[
 \sn^2J=\frac{m-\eta}{m(1-\eta)},\qquad
 \cn^2J=\frac{\eta(1-m)}{m(1-\eta)},\qquad
 \dn^2J=\frac{1-m}{1-\eta}.
\]
All three Jacobi functions are positive at $J\in(0,K)$, fixing the
signs of their square roots. The real quarter-period formulas
\cite[Table~22.4.3]{DLMF} therefore yield
\[
 \sn(K\pm J)=\frac{\cn J}{\dn J}=\sqrt{\frac\eta m},\qquad
 \dn(K\pm J)=\frac{\sqrt{1-m}}{\dn J}=b.
\]
The lower imaginary shift established in
Lemma~\ref{arith:jacobi} now evaluates the physical values directly:
\[
 \mathcal P_m(K\pm J-iK')
 =\frac{1+\dn(K\pm J)}{\sqrt m\,\sn(K\pm J)}
 =\frac{1+b}{\sqrt\eta}=z_+.
\]
The shift by $2K$ changes the sign of $\mathcal P_m$, and hence
\[
 \mathcal P_m(3K\pm J-iK')=z_-.
\]
All four points lie in the lower cylinder, since their imaginary
part is $-K'$.

\smallskip\noindent
\emph{Step 3: verify the full fibers, their shift, and unramifiedness.}
The two points in each pair differ by $2J$, where $0<2J<2K<4K$.
They are therefore distinct modulo the period lattice. Each pair
already supplies two distinct preimages of its value, so the
degree-two conclusion of Lemma~\ref{arith:jacobi} excludes any
additional preimage. In the displayed order, the minus member is
sent to the plus member by the real shift $2J=h$.

For an explicit local regularity check, the Jacobi derivative
formulas \cite[Table~22.13.1]{DLMF} give
\[
 \mathcal P_m'(u)
 =\cn u\,\dn u+i\sn u\,\dn u
 =i\dn u\,\mathcal P_m(u).
\]
The imaginary-shift formula also gives
$\sn(K\pm J-iK')=1/\sqrt\eta$; at the negative pair its sign is
reversed. Thus at every displayed preimage,
\[
 \sn^2u=\frac1\eta,\qquad
 \dn^2u=1-m\sn^2u=1-\frac m\eta\ne0.
\]
Since $\mathcal P_m(u)=z_\pm$ is finite and nonzero, we obtain
$\mathcal P_m'(u)\ne0$. Equivalently, a local logarithm
$\varphi=-i\log\mathcal P_m$ satisfies
\[
 (\partial_u\varphi)^2=\dn^2u
 =1-m\cos^2\varphi=1-\frac m\eta\ne0.
\]
This proves unramifiedness in the ordinary physical base coordinate.
The strict inequality $m>\eta$ is used here; no conclusion at the
limiting value $m=\eta$ is asserted.
\end{proof}

\begin{figure}[tbp]
\centering
\begin{tikzpicture}[
 x=1cm,y=1cm,font=\small,
 axis/.style={-{Stealth[length=1.6mm]},line width=.45pt},
 arr/.style={-{Stealth[length=1.8mm]},line width=.65pt},
 focal/.style={circle,draw,fill=white,inner sep=0pt,minimum size=4.3pt,line width=.7pt},
 every node/.style={inner sep=2pt}]
\path[use as bounding box] (-.45,-1.55) rectangle (15.50,5.80);
\def\ssJJ{0.5405375566328798}
\node at (3.88,5.55) {(a) Lower Jacobi cylinder};
\node at (11.78,5.55) {(b) Boundary-parameter plane};
\fill[black!5] (.80,1.70) rectangle (7.00,4.70);
\draw[densely dashed,line width=.55pt] (.80,4.70)--(7.00,4.70);
\draw[densely dashed,line width=.55pt] (.80,1.70)--(7.00,1.70);
\draw[densely dotted,line width=.6pt] (.80,1.70)--(.80,4.70);
\draw[densely dotted,line width=.6pt] (7.00,1.70)--(7.00,4.70);
\draw[arr] (.80,4.03)--(.80,4.39);
\draw[arr] (7.00,4.03)--(7.00,4.39);
\node[left] at (.70,4.70) {$0$};
\node[left] at (.70,3.20) {$-iK'$};
\node[left] at (.70,1.70) {$-2iK'$};
\foreach \xx/\lab in {.80/0,2.35/K,3.90/2K,5.45/3K,7.00/4K}{
 \node[above=2pt] at (\xx,4.70) {$\lab$};
}
\node at (3.90,4.29) {$\mathcal C_-:\quad |\mathcal P_m|>1$};
\draw[black!40,line width=.4pt] (.80,3.20)--(7.00,3.20);
\coordinate (fpA) at ({.80+1.55*(1-\ssJJ)},3.20);
\coordinate (fpB) at ({.80+1.55*(1+\ssJJ)},3.20);
\coordinate (fnA) at ({.80+1.55*(3-\ssJJ)},3.20);
\coordinate (fnB) at ({.80+1.55*(3+\ssJJ)},3.20);
\foreach \aa/\bb in {fpA/fpB,fnA/fnB}{
 \draw[black!40,densely dotted] (\aa)--($(\aa)+(0,.55)$);
 \draw[black!40,densely dotted] (\bb)--($(\bb)+(0,.55)$);
 \draw[arr] ($(\aa)+(0,.55)$)--node[above] {$h=2J$} ($(\bb)+(0,.55)$);
}
\foreach \pt/\lab in {fpA/u_+^-,fpB/u_+^+,fnA/u_-^-,fnB/u_-^+}{
 \node[focal] at (\pt) {};
 \node[below=6pt] at (\pt) {$\lab$};
}
\foreach \xx in {2.35,5.45}{
 \fill (\xx,3.20) ++(0,.075)--++(.075,-.075)--++(-.075,-.075)--++(-.075,.075)--cycle;
}
\foreach \xx in {.80,3.90,7.00}{
 \draw[line width=.8pt] (\xx-.055,3.145)--(\xx+.055,3.255);
 \draw[line width=.8pt] (\xx-.055,3.255)--(\xx+.055,3.145);
}
\node at (3.90,2.20) {vertical sides identified modulo $4K$};
\node at (3.90,1.35) {both horizontal boundaries map to $|z|=1$};
\draw[arr] (7.25,3.20)--(8.25,3.20);
\node[above] at (7.75,3.38) {$\mathcal P_m$};
\node[below,align=center,font=\footnotesize] at (7.75,3.02) {degree $2$};
\fill[black!5] (8.55,1.70) rectangle (15.10,4.70);
\fill[white] (11.78,3.20) circle[radius=.84];
\draw[line width=.7pt] (11.78,3.20) circle[radius=.84];
\draw[axis,black!55] (8.65,3.20)--(15.18,3.20);
\node[below left] at (15.10,3.13) {$\operatorname{Re}z$};
\node at (11.78,3.50) {$|z|<1$};
\node[below] at (11.78,3.15) {$0$};
\node[above] at (11.78,4.27) {$\mathcal D_{\mathrm{ext}}=\{|z|>1\}\cup\{\infty\}$};
\foreach \xx/\lab in {9.26/z_-,14.30/z_+}{
 \node[focal] at (\xx,3.20) {};
 \node[above=6pt] at (\xx,3.20) {$\lab$};
}
\foreach \xx/\lab in {10.10/\beta_-,13.46/\beta_+}{
 \fill (\xx,3.20) ++(0,.075)--++(.075,-.075)--++(-.075,-.075)--++(-.075,.075)--cycle;
 \node[below=6pt] at (\xx,3.20) {$\lab$};
}
\node at (11.78,2.03) {$\beta_\pm=\pm2,\qquad z_\pm=\pm3$};
\node at (11.78,1.35) {$\eta=9/25,\qquad m=16/25$};
\node[focal] at (1.00,.66) {};
\node[anchor=west] at (1.16,.66) {regular focal points};
\fill (5.75,.66) ++(0,.075)--++(.075,-.075)--++(-.075,-.075)--++(-.075,.075)--cycle;
\node[anchor=west] at (5.92,.66) {ramification / branch points};
\draw[line width=.8pt] (12.10,.605)--(12.21,.715);
\draw[line width=.8pt] (12.10,.715)--(12.21,.605);
\node[anchor=west] at (12.31,.66) {poles over $\infty$};
\node at (7.78,.01) {$u_+^\pm=K\pm J-iK',\qquad
 u_-^\pm=3K\pm J-iK'$};
\node at (7.78,-.58) {$\mathcal P_m^{-1}(\mathcal D_{\mathrm{ext}})=\mathcal C_-
 \qquad\text{(degree counted with multiplicity)}$};
\node[font=\footnotesize] at (7.78,-1.14) {The repeated pole on the two vertical sides is one point of the cylinder.};
\end{tikzpicture}
\caption{The degree-two Jacobi map and its regular focal fibers, illustrated
at $\eta=9/25$ and $m=16/25$. The source coordinates are scaled by $K$ and
$K'$; the vertical sides are identified. Open circles mark the two
preimages of each exterior focal value $z_\pm$, and the arrows show the
real shift $h=2J$. Diamonds mark the ramification points $K-iK'$ and
$3K-iK'$ and their branch values
$\beta_\pm=\pm(1+\sqrt{1-m})/\sqrt m$; these are distinct from $z_\pm$.
Crosses mark the two poles over $\infty$, with the boundary pole repeated
only for display. The target shading denotes the spherical exterior;
its rectangular display window is not an additional boundary.
The covering is branched, whereas both focal fibers are unramified.}
\label{arith:fig:covering-focal}
\end{figure}
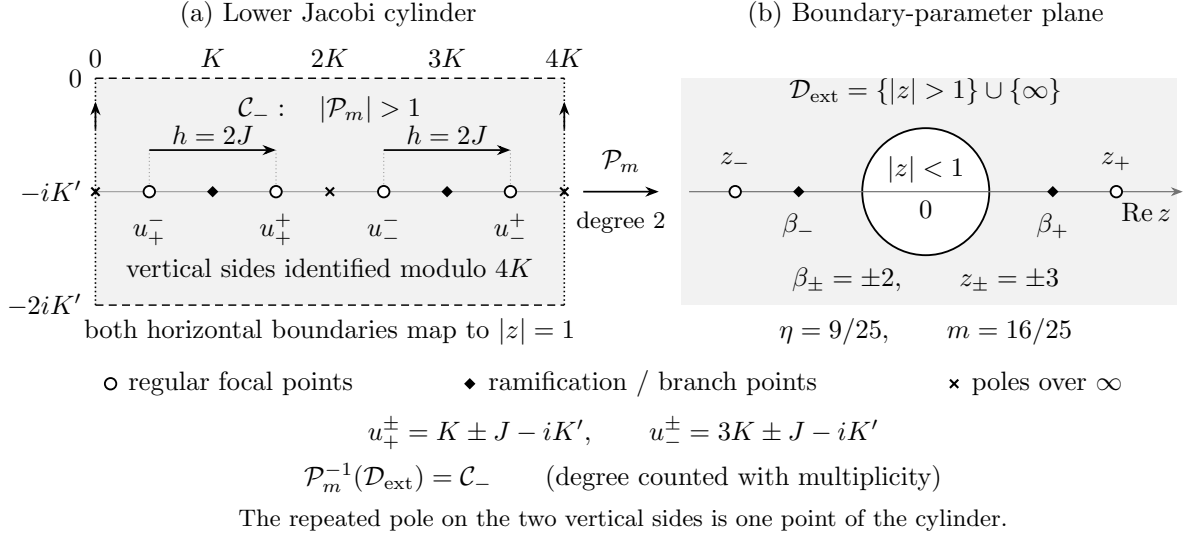

For $\eta<m<1$, these exterior focal values are regular values of the
Jacobi covering, not its branch values.
Figure~\ref{arith:fig:covering-focal} depicts the distinction on the lower
cylinder and in the boundary-parameter plane.

\begin{lemma}[Finite-grid obstruction]\label{arith:grid}
Let $f_1,\ldots,f_S$ be complex analytic on a neighborhood of a real open
interval $I$. Suppose that for every sufficiently large odd prime $q$ and
every integer $p$ with $p/q\in I$, at least one index satisfies
$qf_s(p/q)\in\mathbb Z$. Then at least one $f_s$ is affine on $I$.
The same conclusion holds for any fixed subfamily that still supplies
a lattice value at every grid point in the hypothesis.
\end{lemma}
\begin{proof}
\smallskip\noindent
\emph{Step 1: choose three lattice values of one function in a short block.}
Fix $x_0\in I$ and a compact interval $I_{x_0}\Subset I$ whose
interior contains $x_0$. For each sufficiently large odd prime $q$,
put $p_0(q)=\lfloor qx_0\rfloor-S$ and consider the $2S+1$
consecutive numerators
\[
 p_0(q),\ p_0(q)+1,\ldots,\ p_0(q)+2S.
\]
Their quotients by $q$ lie in $I_{x_0}$ for large $q$ and all tend to
$x_0$. Assign to each quotient one index supplying its lattice value.
If every index occurred at most twice, there would be at most $2S$
points. Hence one index $s$ occurs at least three times. Write three
of its numerators in increasing order as
\[
 p,\quad p+a,\quad p+a+b,
 \qquad a,b\ge1,\quad a+b\le2S.
\]
The index $s$, the starting numerator $p$, and the gaps $a,b$ may
depend on $q$. At the three selected points, however, the same
function has integer values after multiplication by $q$.

\smallskip\noindent
\emph{Step 2: form an integer second difference with uniform remainder.}
Define
\[
 Q_q=bqf_s(p/q)-(a+b)qf_s((p+a)/q)
                   +aqf_s((p+a+b)/q)\in\Z.
\]
Set $x_q=p/q$. All functions and their first three derivatives are
bounded on the fixed compact interval $I_{x_0}$. Taylor expansion
along this real interval gives, uniformly in the finite choices of
$s,a,b$,
\[
 f_s(x_q+t)=f_s(x_q)+t f_s'(x_q)
                  +\frac{t^2}{2}f_s''(x_q)+O(|t|^3).
\]
The weights in $Q_q$ cancel the constant and linear terms because
\[
 b-(a+b)+a=0,\qquad -(a+b)a+a(a+b)=0.
\]
The quadratic coefficient is determined by
\[
 -(a+b)a^2+a(a+b)^2=ab(a+b).
\]
Substitution at $t=0,a/q,(a+b)/q$ therefore yields
\[
 Q_q=\frac{ab(a+b)}{2q}f_s''(x_q)+O(q^{-2}).
\]
The error is uniform: the unscaled Taylor errors are $O(q^{-3})$,
the gaps and integer weights are bounded in terms of $S$, and the
outer factor is $q$. The second derivatives are uniformly bounded
on $I_{x_0}$, so $|Q_q|\le C_{x_0}/q$. For large $q$ this is less
than one, and integrality forces $Q_q=0$. Although the functions are
complex-valued, the hypothesis makes each selected value $qf_s$
an ordinary integer, so this discreteness argument is unchanged.

\smallskip\noindent
\emph{Step 3: obtain a pointwise zero of one second derivative.}
There are only finitely many triples $(s,a,b)$ with
$1\le s\le S$, $a,b\ge1$, and $a+b\le2S$. Along an infinite
subsequence of the odd primes, this triple is fixed. Multiply the
last expansion by $q$, use $Q_q=0$, and pass to the limit
$x_q\to x_0$. We obtain
\[
 0=\frac{ab(a+b)}2 f_s''(x_0).
\]
Since $a,b>0$, this gives $f_s''(x_0)=0$. The point $x_0$ was
arbitrary. Thus for each $x\in I$ at least one $f_s''(x)$ vanishes;
the index may still depend on $x$.

\smallskip\noindent
\emph{Step 4: use analyticity to select one affine member on all of $I$.}
Suppose no $f_s''$ vanishes identically on $I$. The identity theorem
then makes the zeros of each $f_s''$ isolated on the connected
component of the common analytic neighborhood containing $I$.
In particular, each has only finitely many zeros on any fixed
nondegenerate compact subinterval $I_0\Subset I$. Their finite union
cannot cover $I_0$, contradicting Step~3. Consequently there is one
index $s_*$ for which $f_{s_*}''\equiv0$ on $I$. Since $I$ is an
interval, integration gives
\[
 f_{s_*}(x)=Ax+B\qquad(x\in I)
\]
for constants $A,B\in\C$. This is a single index valid on the whole
interval, not an index fixed only along one grid subsequence.
Finally, any fixed subfamily that still covers every required grid
point satisfies the same hypotheses. Applying the argument with
its cardinality in place of $S$ proves the subfamily assertion.
\end{proof}

Figure~\ref{arith:fig:finite-grid} summarizes the local selection and the
separate analytic step in Lemma~\ref{arith:grid}.

\begin{figure}[tbp]
\centering
\begin{tikzpicture}[
 x=1cm,y=1cm,font=\small,
 arr/.style={-{Stealth[length=1.6mm]},line width=.6pt},
 both/.style={<->,>=Stealth,line width=.5pt},
 every node/.style={inner sep=2pt}]
\path[use as bounding box] (-.12,-.66) rectangle (15.52,5.98);
\node at (3.40,5.68) {(a) One block: $2S+1=7$ points};
\node at (11.25,5.68) {(b) From local blocks to one affine member};
\node at (3.40,5.13) {one possible supplier assignment for $S=3$};
\draw[arr,black!60] (.30,3.77)--(6.46,3.77);
\foreach \xx/\lab in {.60/1,2.40/1,5.10/1}{
 \fill (\xx,3.77) circle[radius=2.0pt];
 \node[above=5pt] at (\xx,3.77) {$\lab$};
}
\foreach \xx/\lab in {1.50/2,4.20/2}{
 \draw[fill=white,line width=.6pt] (\xx-.055,3.715) rectangle (\xx+.055,3.825);
 \node[above=5pt] at (\xx,3.77) {$\lab$};
}
\foreach \xx/\lab in {3.30/3,6.00/3}{
 \draw[fill=white,line width=.6pt] (\xx,3.85)--(\xx-.065,3.72)--(\xx+.065,3.72)--cycle;
 \node[above=5pt] at (\xx,3.77) {$\lab$};
}
\foreach \xx in {.60,2.40,5.10}{
 \draw[black!45,densely dotted] (\xx,3.85)--(\xx,4.52);
}
\draw[both] (.60,4.52)--node[above] {$a/q=2/q$} (2.40,4.52);
\draw[both] (2.40,4.52)--node[above] {$b/q=3/q$} (5.10,4.52);
\node[below=5pt] at (.60,3.77) {$p/q$};
\node[below=5pt] at (2.40,3.77) {$(p+2)/q$};
\node[below=5pt] at (5.10,3.77) {$(p+5)/q$};
\node[align=center] at (3.32,2.72) {$A=qf_1(p/q),\quad B=qf_1((p+2)/q),$\\[3pt]
 $C=qf_1((p+5)/q)\quad\text{are integers.}$};
\node[draw,rounded corners=2pt,fill=black!3,inner sep=8pt] at (3.32,1.76)
 {$Q_q=3A-5B+2C\in\mathbb Z$};
\node at (3.32,.91) {$3-5+2=0,\qquad -5\cdot2+2\cdot5=0$};
\node at (3.32,.27) {$Q_q=\dfrac{15}{q}f_1''(p/q)+O(q^{-2})$};
\node[draw,rounded corners=2pt,fill=black!3,minimum width=6.65cm,
 minimum height=.62cm] (qzero) at (11.25,4.75)
 {$Q_q\in\mathbb Z,\quad Q_q=O(q^{-1})\quad\Longrightarrow\quad Q_q=0$};
\node[align=center,font=\footnotesize] at (11.25,4.09)
 {for all sufficiently large allowed $q$};
\draw[arr] (11.25,3.85)--(11.25,3.41);
\node[align=center,font=\footnotesize] at (11.25,3.10)
 {Fix $x_0$; pass to a subsequence with fixed $(s,a,b)$.};
\node[draw,rounded corners=2pt,fill=black!3,minimum width=6.65cm,
 minimum height=.62cm] (pointzero) at (11.25,2.46)
 {$\forall x_0\in I\quad\exists s:\quad f_s''(x_0)=0$};
\draw[arr] (11.25,2.06)--(11.25,1.62);
\node[align=center,font=\footnotesize] at (11.25,1.29)
 {Finiteness of the family and the identity theorem};
\node[draw,rounded corners=2pt,fill=black!3,minimum width=6.65cm,
 minimum height=.69cm] at (11.25,.55)
 {$\exists s_*:\quad f_{s_*}''\equiv0\ \text{on }I\quad\Longrightarrow\quad f_{s_*}\ \text{affine}$};
\node[font=\footnotesize] at (7.65,-.40)
 {The supplying index may change with the grid point; one finite block alone does not imply affinity.};
\end{tikzpicture}
\caption{The finite-grid argument, with a schematic supplier assignment
for $S=3$. Among seven grid points one supplier occurs at least three
times. In the displayed choice the gaps are $a=2$, $b=3$, and the weights
$3,-5,2$ cancel the constant and linear Taylor terms; $A,B,C$ denote the
three scaled values of that same function, not ordinate coordinates of
the drawing. The right panel separates the pointwise subsequence argument
from the use of analyticity to obtain one affine member on all of $I$.
Neither a single grid block nor a single supplying index fixed in advance
is assumed to establish the conclusion.}
\label{arith:fig:finite-grid}
\end{figure}

\begin{lemma}[Picard--Fuchs source of a marked point]\label{arith:pf}
All parameter derivatives below hold $\eta$ fixed, and the derivatives
of $I_\varphi$ also hold the physical angle $\varphi$ fixed. Let
\[
 \mathcal L=m(1-m)\partial_m^2+(1-2m)\partial_m-\tfrac14,
 \qquad I_\varphi(m)=\int_0^\varphi(1-m\cos^2t)^{-1/2}\,dt.
\]
The following identities hold for the real branches on $\eta<m<1$
and their compatible analytic continuations, with all square roots
continued together with the corresponding elliptic integrals:
\[
 4\mathcal L I_\varphi
 =R_\varphi(m):=\frac{\sin\varphi\cos\varphi}
 {(1-m\cos^2\varphi)^{3/2}},\qquad
 \mathcal L K=\mathcal L K'=0,
 \qquad \mathcal L J=-\frac{\sqrt{\eta(1-\eta)}}{4(m-\eta)^{3/2}}.
\]
\end{lemma}
\begin{proof}
The fixed-amplitude equation is the classical inhomogeneous
Picard--Fuchs equation \cite[Eq.~19.4.8]{DLMF}. Under $m=k^2$,
\[
 D_k=2k\partial_m,\qquad
 k(1-k^2)D_k^2+(1-3k^2)D_k-k=4k\mathcal L,
\]
so it reads
\[
 4\mathcal L_mF(\vartheta\mid m)
 =-\frac{\sin\vartheta\cos\vartheta}
         {(1-m\sin^2\vartheta)^{3/2}}
\]
with $\vartheta$ fixed. The endpoint $\vartheta=\pi/2$ gives
$\mathcal LK=0$; invariance of the homogeneous equation under
$m\mapsto1-m$ gives $\mathcal LK'=0$, where $K'=K(1-m)$.
Since $I_\varphi=K-F(\pi/2-\varphi\mid m)$, the stated positive
source $R_\varphi$ follows with the physical endpoint $\varphi$ fixed.

For the moving endpoint in $J$, put
\[
 \delta=m-\eta,\qquad e=1-m,\qquad C_\eta=\sqrt{\eta(1-\eta)},
 \qquad\vartheta(m)=\arcsin\sqrt{\frac{m-\eta}{m(1-\eta)}}.
\]
For real $\eta<m<1$, use positive square roots. Then
\[
 \vartheta'=\frac{\sqrt\eta}{2m\sqrt\delta\sqrt e},\qquad
 \vartheta''=\vartheta'
     \left(-\frac1m-\frac1{2\delta}+\frac1{2e}\right).
\]
Writing $f=F_\vartheta$ and evaluating partial derivatives at
$(\vartheta(m),m)$ gives
\[
 f=\frac{\sqrt{1-\eta}}{\sqrt e},\quad
 f_m=\frac{\delta\sqrt{1-\eta}}{2me^{3/2}},\quad
 f_\vartheta=\frac{C_\eta\sqrt\delta}{e},\quad
 (\mathcal L_mF)(\vartheta(m),m)=-\frac{C_\eta\sqrt\delta}{4me}.
\]
The total chain rule, rather than the fixed-amplitude equation alone,
therefore gives
\begin{align*}
 \mathcal LJ={}&(\mathcal L_mF)(\vartheta(m),m)
 +me\bigl(2f_m\vartheta'+f_\vartheta(\vartheta')^2+f\vartheta''\bigr)
 +(1-2m)f\vartheta'\\
 ={}&\frac{C_\eta}{4\delta^{3/2}}
 \left[-1+\frac{2\delta}{me}(m-e+1-2m)\right]
 =-\frac{C_\eta}{4\delta^{3/2}},
\end{align*}
since $e=1-m$. All derivatives hold $\eta$ fixed.
Finally the identities continue as identities of holomorphic germs.
For complex fixed $\varphi$, use a locally continued integration path;
continue each square root together with its elliptic integral and do
not reset it to an independent principal branch. This proves the
assertions on every compatible continuation.
\end{proof}

\section{An exact Borel--Volterra continuation mechanism}
\label{vol:sec:borel-volterra}

We first continue the known symbol coefficients, then convert the
vanishing selected positive rows into an exact Volterra equation whose
forcing is determined by the original source-circle data. The equation
is restarted along regular paths with the full integration history,
including loops, retained. A proper Jacobi covering then transfers the
continued germs to the physical coordinate outside one fixed closed
countable set with finite relative derived set. The conclusion is
all-path continuation, not single-valuedness, finite monodromy, or
meromorphicity.

Use $\zeta=e^{ix}$ for the complex action coordinate and
$z=e^{i\varphi}$ for the physical coordinate.
Fix a noncircular reference ellipse and an arbitrary reduced rational
base $\rho_0=p_0/P\in I$. For each $s\in\{0,\ldots,P-1\}$, let
$r_s$ be the representative of $-p_0s\pmod P$ in $\{0,\ldots,P-1\}$. For all sufficiently large positive integers $n\equiv s
\pmod P$, the permitted selection is
\[
\rho_n=\rho_0+\frac{r_s}{Pn}.
\]
The selected positive row applied to $U$ is exactly the $n$th Fourier
coefficient of $b_nU$, where
\[
b_n(x)=H_n'(x)\exp[-in(H_n(x)-x)],\qquad
H_n=G_{\rho_n}\circ\Phi_{\rho_0}.
\]
Recall that $Y(x)=\left.\partial_\rho
(G_\rho\circ\Phi_{\rho_0})(x)\right|_{\rho=\rho_0}$.

\begin{lemma}[A convergent symbol expansion]
\label{vol:lem:convergent-symbol}
For each residue $s$ there is a function $b_s(\varepsilon,\zeta)$, holomorphic
near $\varepsilon=0$ and the real action circle, such that
$b_n=b_s(1/n,\cdot)$ for all sufficiently large $n$ in that residue. It has a convergent
expansion
\[
b_s(\varepsilon,\zeta)=\sum_{k\ge0}a_{s,k}(\zeta)\varepsilon^k,
\qquad
a_{s,0}(\zeta)=\exp[-ir_sY(\zeta)/P].
\]
There is a fixed closed discrete Jacobi set $S\subset\mathbb C^*$ such that all
$a_{s,k}$ are single-valued holomorphic on $\mathbb C^*\setminus S$.
For every compact subset $\mathcal K$ of that domain there are
$M_{\mathcal K},R_{\mathcal K}>0$, independent of $s$ and $k$, such
that
\[
\sup_{\mathcal K}|a_{s,k}|\le M_{\mathcal K}R_{\mathcal K}^{-k}.
\]
Consequently
\[
B_s(t,\zeta)=\sum_{k\ge1}a_{s,k}(\zeta)\frac{t^{k-1}}{(k-1)!}
\]
is entire in $t$, holomorphic in $\zeta\notin S$, and satisfies
$\sup_{\mathcal K}|B_s(t,\cdot)|\le(M_{\mathcal K}/R_{\mathcal K})e^{|t|/R_{\mathcal K}}$.
For every sufficiently large $n\equiv s\pmod P$,
\begin{equation}\label{vol:eq:borel-laplace-symbol}
b_n(\zeta)=a_{s,0}(\zeta)+\int_0^\infty e^{-nt}B_s(t,\zeta)\,dt
\end{equation}
uniformly on a fixed neighborhood of the real action circle.
\end{lemma}
\begin{proof}
The reference ellipse and the rational base $\rho_0=p_0/P$ are fixed.
Put $c_s=r_s/P$ for $0\le s<P$. We first construct the symbol near the
real action circle, then continue its Taylor coefficients. The parameter
radius in the latter construction will depend on a compact subset of the
coefficient domain; no global positive radius is asserted.

\smallskip\noindent
\emph{Step 1: the initial convergent symbol and its periodic descent.}
By Lemma~\ref{op:lem:strips}, the functions
\[
 H_s(\varepsilon,x)
 =G_{\rho_0+c_s\varepsilon}(\Phi_{\rho_0}(x))
\]
are jointly holomorphic for $\varepsilon$ near zero and $x$ in a common
strip about the real axis. Since $H_s(0,x)=x$, the quotient
\[
 Q_s(\varepsilon,x)
 =\frac{H_s(\varepsilon,x)-x}{\varepsilon}
 =\int_0^1\partial_\varepsilon H_s(t\varepsilon,x)\,dt
\]
has a removable singularity at $\varepsilon=0$. Thus
\[
 b_s(\varepsilon,x)
 =\partial_xH_s(\varepsilon,x)\exp[-iQ_s(\varepsilon,x)]
\]
is jointly holomorphic and has a convergent Taylor expansion in
$\varepsilon$. The degree-one identity for $H_s$ makes $H_s-x$,
$Q_s$, and $\partial_xH_s$ periodic in $x$. Consequently $b_s$ descends
to a holomorphic function of $\zeta=e^{ix}$ in an annulus about the
unit circle. Moreover,
\[
 Q_s(0,x)=c_sY(x),\qquad \partial_xH_s(0,x)=1,
 \qquad a_{s,0}(\zeta)=e^{-ic_sY(\zeta)}.
\]
For sufficiently large $n\equiv s\pmod P$, substitution of
$\varepsilon=1/n$ gives exactly the selected amplitude $b_n$.
For $r_s=0$ the symbol is identically one, so all its positive-order
coefficients vanish.

\smallskip\noindent
\emph{Step 2: compute the first parameter derivative at a marked point.}
Set
\[
 m_0=m(\rho_0),\qquad K_0=K(m_0),\qquad K_{\mathrm c,0}=K(1-m_0),
 \qquad u=K_0+\frac{2K_0}{\pi}x.
\]
All Jacobi functions below have the frozen parameter $m_0$ unless
another parameter is displayed. For derivatives at a fixed physical
angle, write
\[
 \mathcal G(m,\varphi)
 =\frac{\pi}{2K(m)}I_\varphi(m),\qquad
 I_\varphi(m)=\int_0^\varphi(1-m\cos^2t)^{-1/2}\,dt.
\]
Thus $\mathcal G(m(\rho),\varphi)=G_\rho(\varphi)$ locally. With
$\alpha=\varphi+\pi/2$, the coordinate convention gives
\[
 \mathcal G(m,\varphi)
 =\frac{\pi}{2K(m)}F(\alpha\mid m)-\frac\pi2.
\]
Let $E(\alpha\mid m)=\int_0^\alpha(1-m\sin^2t)^{1/2}\,dt$ and
$E(m)=E(\pi/2\mid m)$ denote the second-kind elliptic integrals.
Converting \cite[(19.4.1) and (19.4.5)]{DLMF} from $k$ to $m=k^2$
gives, with $\alpha$ fixed and $K_m=dK/dm$,
\begin{align*}
 \partial_mF(\alpha\mid m)
 &=\frac{E(\alpha\mid m)}{2m(1-m)}-\frac{F(\alpha\mid m)}{2m}
   -\frac{\sin\alpha\cos\alpha}
          {2(1-m)\sqrt{1-m\sin^2\alpha}},\\
 K_m(m)&=\frac{E(m)}{2m(1-m)}-\frac{K(m)}{2m}.
\end{align*}
The terms proportional to $F/(2m)$ cancel upon differentiating the
normalization. Therefore
\[
 \partial_m\mathcal G(m,\varphi)
 =\frac{\pi}{4K(m)m(1-m)}
  \left[E(\alpha\mid m)-\frac{E(m)}{K(m)}F(\alpha\mid m)
       -m\frac{\sin\alpha\cos\alpha}{\sqrt{1-m\sin^2\alpha}}\right].
\]
At $m=m_0$ and $\varphi=\Phi_{\rho_0}(x)$, we have
$F(\alpha\mid m_0)=u$, $\sin\alpha=\sn u$, $\cos\alpha=\cn u$,
and the compatible square root is $\dn u$.
Initially near the real axis, put
\[
 \mathcal E(u)=\int_0^u\dn^2t\,dt,\qquad
 Z(u)=\mathcal E(u)-\frac{E(m_0)}{K_0}u.
\]
With the continuously lifted amplitude, $E(\alpha\mid m_0)=\mathcal E(u)$:
as functions of $u$, both sides have derivative $\dn^2u$ and agree
at $u=0$. This is \cite[(22.16.31)]{DLMF} on its initial real
interval, extended with the compatible amplitude; the normalization
of $Z$ is \cite[(22.16.32)]{DLMF}. We obtain
\begin{equation}\label{vol:eq:marked-first-derivative}
 \left.\partial_m\mathcal G(m,\Phi_{\rho_0}(x))\right|_{m=m_0}
 =\frac{\pi}{4K_0m_0(1-m_0)}
   \left[Z(u)-m_0\frac{\sn u\cn u}{\dn u}\right].
\end{equation}
The theta-function representation \cite[(22.16.30)]{DLMF} supplies a
single-valued meromorphic continuation of $Z$, and
\cite[(22.16.34)]{DLMF} gives $Z(u+2K_0)=Z(u)$. The Jacobi quotient
in the displayed formula is meromorphic and $2K_0$-periodic as well.
In particular, the right side is $2\pi$-periodic in $x$, and hence
defines a single-valued meromorphic function $D_0(\zeta)$ on
$\mathbb C^*$. On the original annulus the chain rule gives
\[
 Y(\zeta)=\left.\frac{dm}{d\rho}\right|_{\rho_0}D_0(\zeta).
\]
All differentiation above held $\varphi$ fixed before substituting
the base action coordinate.

\smallskip\noindent
\emph{Step 3: normalize the Picard--Fuchs equation.}
Use $\mathcal L=m(1-m)\partial_m^2+(1-2m)\partial_m-1/4$ from
Lemma~\ref{arith:pf}. Since $I_\varphi=(2/\pi)K\mathcal G$, the
product rule gives, at fixed $\varphi$,
\begin{align*}
 \mathcal L(K\mathcal G)
 &=K m(1-m)\partial_m^2\mathcal G
   +\big[2m(1-m)K_m+(1-2m)K\big]\partial_m\mathcal G
   +(\mathcal LK)\mathcal G.
\end{align*}
By that lemma, $\mathcal LK=0$ and
$\mathcal LI_\varphi=R_\varphi/4$, where
$R_\varphi=\sin\varphi\cos\varphi(1-m\cos^2\varphi)^{-3/2}$.
Dividing the preceding identity therefore yields
\begin{equation}\label{vol:eq:normalized-pf-symbol}
 \partial_m^2\mathcal G+
 \left(\frac{2K_m(m)}{K(m)}+\frac{1-2m}{m(1-m)}\right)
       \partial_m\mathcal G
 =\frac{\pi R_\varphi(m)}{8K(m)m(1-m)}.
\end{equation}
In particular there is no zeroth-order term. Initially on the action
annulus define
\[
 W(m,x)=\mathcal G(m,\Phi_{\rho_0}(x))-x.
\]
The derivatives in the parameter equation keep $x$, and therefore
$\Phi_{\rho_0}(x)$, fixed. Subtracting $x$ leaves the same equation,
with initial data $W(m_0,x)=0$ and
$\partial_mW(m_0,x)=D_0(e^{ix})$. A zeroth-order term would introduce
an additional source involving $x$, and hence a choice of $\log\zeta$;
its cancellation is what permits the construction below.

\smallskip\noindent
\emph{Step 4: express the source on a fixed punctured Jacobi domain.}
Set
\[
 \chi(\zeta)=\frac{\sn^2u}{\dn^2u},\qquad
 R_0(\zeta)=-\frac{\sn u\cn u}{\dn^3u}.
\]
These are single-valued meromorphic functions of $\zeta$: the
expressions on the right have real period $4K_0$ in $u$.
With $\varphi=\Phi_{\rho_0}(x)$ fixed we have
$\cos\varphi=\sn u$, $\sin\varphi=-\cn u$, and
\[
 1-m\cos^2\varphi
 =\dn^2u-(m-m_0)\sn^2u.
\]
Consequently the source in \eqref{vol:eq:normalized-pf-symbol} is
obtained from the germ
\[
 R(m,\zeta)
 =R_0(\zeta)\big[1-(m-m_0)\chi(\zeta)\big]^{-3/2}.
\]
The binomial factor is the analytic germ equal to one at $m=m_0$;
the factor $\dn^{-3}u$ carries the compatible initial square root.
Neither factor is independently reset to a principal branch during
continuation.

The Jacobi pole and zero tables
\cite[Tables~22.4.1--22.4.2]{DLMF} place all poles of $\sn$, $\cn$,
and $\dn$, and all zeros of $\dn$, in the containing set
\[
 u=jK_0+(2\ell+1)iK_{\mathrm c,0},\qquad j,\ell\in\mathbb Z.
\]
Also $Z'(u)=\dn^2u-E(m_0)/K_0$; since $Z$ is meromorphic, each pole
of $Z$ must be a pole of $\dn$. Put
\[
 \mathfrak q=e^{-\pi K_{\mathrm c,0}/K_0}\in(0,1),\qquad
 S=\{i^j\mathfrak q^{\,\ell+1/2}:
                  j=0,1,2,3,\ \ell\in\mathbb Z\}.
\]
Indeed, the change of coordinate is
\[
 \zeta=\exp\!\left(\frac{i\pi(u-K_0)}{2K_0}\right),
\]
and it sends $u=jK_0+(2\ell+1)iK_{\mathrm c,0}$ to
$i^{j-1}\mathfrak q^{\ell+1/2}$, which is the same set after relabeling
the four phases. Its only possible spherical accumulation points
are zero and infinity. Thus $S$ is closed and discrete in
$\mathbb C^*$ and avoids the unit circle. In the lower cylinder its
representatives are $u=jK_0-iK_{\mathrm c,0}$ modulo $4K_0$.
The set may include removable points. All possible poles of $D_0$,
$\chi$, and $R_0$ lie in $S$.

\smallskip\noindent
\emph{Step 5: solve the parameter equation with a compact-uniform radius.}
Let $\mathcal D=\mathbb C^*\setminus S$ and fix
$\mathcal K\Subset\mathcal D$. Choose relatively compact open
neighborhoods
\[
 \mathcal K\subset\Omega_0\Subset\Omega_1\Subset\mathcal D.
\]
The functions $D_0$, $\chi$, and $R_0$ are bounded on
$\overline{\Omega_1}$. Choose $r>0$ such that the disk
$|m-m_0|<r$ avoids $0,1$ and zeros of $K(m)$, on the local branch of
$K$ at $m_0$, and such that
\[
 r\sup_{\overline{\Omega_1}}|\chi|<\tfrac12.
\]
The binomial series for $R$ then converges normally on this product
domain. Define
\[
 A(m)=\frac{2K_m(m)}{K(m)}+\frac{1-2m}{m(1-m)},\qquad
 \mathcal F(m,\zeta)=\frac{\pi R(m,\zeta)}{8K(m)m(1-m)}.
\]
We solve
\[
 W_{mm}+A(m)W_m=\mathcal F(m,\zeta),\qquad
 W(m_0,\zeta)=0,\quad W_m(m_0,\zeta)=D_0(\zeta).
\]
An explicit integrating factor is
\[
 \mathfrak a(m)=m(1-m)K(m)^2,\qquad
 \frac{\mathfrak a'(m)}{\mathfrak a(m)}=A(m).
\]
Thus
\[
 \partial_m(\mathfrak a W_m)=\frac\pi8K(m)R(m,\zeta),
\]
and the solution is
\begin{align*}
 W_m(m,\zeta)
 &=\frac{1}{\mathfrak a(m)}
    \left[\mathfrak a(m_0)D_0(\zeta)
          +\frac\pi8\int_{m_0}^{m}K(v)R(v,\zeta)\,dv\right],\\
 W(m,\zeta)&=\int_{m_0}^{m}W_m(v,\zeta)\,dv.
\end{align*}
The integrals are taken in the parameter disk; their integrands are
holomorphic there, so they are path independent. They define a
jointly holomorphic, single-valued solution. No integration or
logarithm in $\zeta$ is involved. On a smaller parameter disk the
formulas give uniform bounds
\[
 |W_m(m,\zeta)|\le C_{\mathcal K},\qquad
 |W(m,\zeta)|\le C_{\mathcal K}|m-m_0|
 \quad (\zeta\in\Omega_1).
\]
Near the original annulus this solution equals the function from
Step~3, by uniqueness of the parameter initial-value problem.
Constructions for different compact neighborhoods agree as germs
at $m_0$ on their overlaps, again by uniqueness.

For explicit coefficient descent, write
\[
 W=\sum_{j\ge0}w_j(\zeta)(m-m_0)^j,\qquad
 A=\sum_{j\ge0}A_j(m-m_0)^j,\qquad
 \mathcal F=\sum_{j\ge0}F_j(\zeta)(m-m_0)^j.
\]
Here $w_0=0$, $w_1=D_0$, and each $F_j$ is meromorphic on
$\mathbb C^*$ with possible poles only in $S$, by the binomial
source expansion. Equating coefficients gives
\[
 (j+2)(j+1)w_{j+2}
 =F_j-\sum_{\ell=0}^{j}A_\ell(j-\ell+1)w_{j-\ell+1},
 \qquad j\ge0.
\]
Induction proves the same single-valued meromorphic property for
every $w_j$. The integral construction proves actual convergence
on each compact product neighborhood, not just a formal recursion.

\smallskip\noindent
\emph{Step 6: return to the symbol and bound all Taylor coefficients.}
The change $m=m(\rho)$ is holomorphic near $\rho_0$ by
Lemma~\ref{op:lem:strips}. Put
\[
 m_s(\varepsilon)=m(\rho_0+c_s\varepsilon),\qquad
 \mathcal W_s(\varepsilon,\zeta)=W(m_s(\varepsilon),\zeta).
\]
For the finitely many residues, choose a common sufficiently small
$\varepsilon$-disk with $|m_s(\varepsilon)-m_0|<r/2$. On
$\Omega_1$ the preceding bounds and the analytic parameter change
give $|\mathcal W_s|\le C_{\mathcal K}|\varepsilon|$.
The quotient $\mathcal W_s/\varepsilon$ therefore extends
holomorphically at zero and is uniformly bounded. Cauchy's estimate
in $\zeta$, applied on $\Omega_1$ and restricted to $\Omega_0$, also
gives
\[
 |\zeta\partial_\zeta\mathcal W_s|
 \le C_{\mathcal K}|\varepsilon|.
\]
This is why a larger spatial neighborhood was chosen before
differentiating. Since $\partial_x=i\zeta\partial_\zeta$, the
continued symbol germ is
\[
 b_s(\varepsilon,\zeta)
 =\bigl(1+i\zeta\partial_\zeta\mathcal W_s(\varepsilon,\zeta)\bigr)
   \exp\!\left[-i\frac{\mathcal W_s(\varepsilon,\zeta)}{\varepsilon}\right].
\]
It agrees with Step~1 on the original annulus. Its value at zero is
$\exp[-ic_s m'(\rho_0)D_0]=\exp[-ic_sY]$, as required.
Each Taylor coefficient is $\exp[-ic_sY]$ times a finite polynomial
expression in meromorphic coefficients of $\mathcal W_s$ and their
spatial derivatives. Consequently all $a_{s,k}$ are single-valued
holomorphic on $\mathcal D$; they may have essential singularities
at $S$, since an exponential of a meromorphic function occurs.
Uniqueness on overlapping neighborhoods identifies them with the
initial Taylor coefficients throughout this domain.

Choose a circle $|\varepsilon|=R_{\mathcal K}$ strictly inside the
common disk just constructed. The symbol is uniformly bounded there
for $\zeta\in\mathcal K$ and all residues, say by $M_{\mathcal K}$.
Cauchy's coefficient formula now gives
\[
 \sup_{\zeta\in\mathcal K}|a_{s,k}(\zeta)|
 \le M_{\mathcal K}R_{\mathcal K}^{-k},\qquad k\ge0.
\]
Both constants may depend on $\mathcal K$, but neither depends on
$s$ or $k$. In particular, the globally continued coefficients do
not imply convergence on all of $\mathcal D$ for one fixed nonzero
$\varepsilon$.

\smallskip\noindent
\emph{Step 7: prove the Borel bound and exact Laplace reconstruction.}
The coefficient estimate yields, for every complex $t$,
\begin{align*}
 \sup_{\zeta\in\mathcal K}|B_s(t,\zeta)|
 &\le\sum_{k\ge1}M_{\mathcal K}R_{\mathcal K}^{-k}
                         \frac{|t|^{k-1}}{(k-1)!}\\
 &=\frac{M_{\mathcal K}}{R_{\mathcal K}}
                         e^{|t|/R_{\mathcal K}}.
\end{align*}
Normal convergence on compact products proves joint holomorphy on
$\mathbb C\times\mathcal D$, in particular entire dependence on $t$.
For the exact identity, choose a fixed closed annulus
$\mathcal K_0$ containing the unit circle in its interior and lying
inside the initial symbol neighborhood and $\mathcal D$.
Write $M=M_{\mathcal K_0}$ and $R=R_{\mathcal K_0}$, decreasing $R$
if necessary to remain within the initial symbol domain.
For $n>R^{-1}$,
\[
 e^{-nt}|B_s(t,\zeta)|\le\frac MR e^{-(n-R^{-1})t},\qquad t\ge0,
\]
uniformly on $\mathcal K_0$. More explicitly,
\begin{align*}
 \sum_{k\ge1}\sup_{\mathcal K_0}|a_{s,k}|
       \int_0^\infty e^{-nt}\frac{t^{k-1}}{(k-1)!}\,dt
 &\le M\sum_{k\ge1}(nR)^{-k}
 =\frac{M}{nR-1}<\infty.
\end{align*}
Thus termwise integration is justified uniformly. The elementary
Laplace moment is
\[
 \int_0^\infty e^{-nt}\frac{t^{k-1}}{(k-1)!}\,dt=n^{-k}.
\]
For all sufficiently large $n\equiv s\pmod P$, it follows that
\[
 \int_0^\infty e^{-nt}B_s(t,\zeta)\,dt
 =\sum_{k\ge1}a_{s,k}(\zeta)n^{-k}
 =b_s(1/n,\zeta)-a_{s,0}(\zeta)
 =b_n(\zeta)-a_{s,0}(\zeta).
\]
This proves \eqref{vol:eq:borel-laplace-symbol} on a fixed neighborhood
of the real action circle. The equality is exact, because the
$\varepsilon$-series converges; it is not only an asymptotic formula.
The exponential type and the parameter radius remain compact-dependent.
No assertion that one fixed $n$ gives a Laplace integral on the whole
of $\mathbb C^*\setminus S$ is needed. The construction concerns
only the known coordinate symbols and uses no continuation or
vanishing-row property of an unknown displacement.
\end{proof}

\begin{proposition}[Exact local Volterra equation]
\label{vol:prop:exact-volterra}
Suppose that $U$ is analytic in an annulus about $|\zeta|=1$ and all its
sufficiently large selected positive rows vanish. Choose a source circle
$|\zeta|=R_c$ inside that annulus, sufficiently close to one that the
matrix $M(w)$ defined below is invertible whenever $R_ce^w$ lies on that
circle. Put $\omega=e^{2\pi i/P}$,
$\zeta_l=R_c\omega^l$, and $u_l(w)=U(\zeta_le^w)$.
Then, initially for small $w$ with $\operatorname{Re}w>0$,
\begin{equation}\label{vol:eq:volterra-system}
M(w)u(w)+\int_0^w\mathcal B(w,\sigma)u(\sigma)\,d\sigma=g(w),
\end{equation}
where
\[
M_{sl}(w)=\frac1P\omega^{-sl}a_{s,0}(\zeta_le^w),\qquad
\mathcal B_{sl}(w,\sigma)=\frac1P\omega^{-sl}
 B_s(w-\sigma,\zeta_le^{\sigma}).
\]
The vector $g$ is known from the original data on the fixed source
circle and is holomorphic on the whole right half-plane.
\end{proposition}
\begin{proof}
The reference ellipse and rational base are fixed. The only values of the
unknown function used to construct the forcing will be its values on one
source circle. We derive the equation first on the inner side of that
circle, then continue the identity locally across it.

\smallskip\noindent
\emph{Step 1: choose a common source collar and a fixed cutoff.}
The source-circle condition is attainable. For $\zeta=e^{ix}$, the matrix
with entries $P^{-1}\omega^{-sl}a_{s,0}(\omega^l\zeta)$ is, after
reordering the rows by $r_s$ and the columns by $l\equiv p_0j\pmod P$,
the Vandermonde matrix of Proposition~\ref{band:prop:operator} divided
by $P$. Indeed its reordered entries are
\[
 \frac1P\exp\!\left[\frac{ir_s}{P}
   \left(2\pi j-Y\!\left(x+\frac{2\pi p_0j}{P}\right)\right)\right].
\]
Its determinant is nonzero on the unit circle. The coefficients are
holomorphic on a common annular neighborhood by
Lemma~\ref{vol:lem:convergent-symbol}, so compactness gives a smaller
annulus on which this determinant remains nonzero. Choose $R_c$
sufficiently close to one within that annulus and the original data
annulus. This justifies the stated invertibility along $|\zeta|=R_c$;
it is not an assertion of invertibility throughout the exterior.

Choose radii $0<R_-<R_c<R_+$ so that the closed collar
\[
 \mathcal A=\{\xi\in\mathbb C:R_-\le|\xi|\le R_+\}
\]
has a neighborhood contained in the original holomorphy annulus of $U$
and in the fixed neighborhood where
\eqref{vol:eq:borel-laplace-symbol} holds. In particular it avoids the
coefficient singularities. Lemma~\ref{vol:lem:convergent-symbol} and
boundedness of $U$ on this collar give constants $C,C_0>0$ such that
\[
 \sup_{\xi\in\mathcal A}|B_s(\tau,\xi)U(\xi)|
 \le Ce^{C_0|\tau|},\qquad \tau\in\mathbb C,\quad 0\le s<P.
\]
The functions $a_{s,0}U$ are bounded on the same collar. Increase one
integer cutoff $N\ge1$ so that all selected rows with $n\ge N$ vanish,
the exact symbol formula applies to them on this collar, and
\[
 N>C_0.
\]
This $N$ is fixed for the remainder of the proof and does not depend on
$w$.

\smallskip\noindent
\emph{Step 2: separate high and low Laurent modes by Cauchy transforms.}
For $f$ holomorphic near $\mathcal A$, define, with counterclockwise
orientation,
\[
 C_Nf(\zeta)=\frac{\zeta^N}{2\pi i}
  \int_{|\xi|=R_c}\frac{f(\xi)\xi^{-N}}{\xi-\zeta}\,d\xi.
\]
Write $C_N^{\rm in}f$ and $C_N^{\rm out}f$ for the functions initially
defined on $|\zeta|<R_c$ and $|\zeta|>R_c$, respectively. If
$f(\xi)=\sum_{n\in\mathbb Z}f_n\xi^n$ in the source collar, the two
geometric expansions of $(\xi-\zeta)^{-1}$ give
\[
 C_N^{\rm in}f(\zeta)=\sum_{n\ge N}f_n\zeta^n,
 \qquad
 C_N^{\rm out}f(\zeta)=-\sum_{n<N}f_n\zeta^n.
\]
The second expression includes all negative modes and all nonnegative
modes below $N$, with the displayed minus sign. The Laurent expansions
on the common collar therefore give the jump identity
\[
 C_N^{\rm in}f=f+C_N^{\rm out}f.
\]
These are identities of holomorphic germs across the source circle,
not principal-value boundary formulas.

We will use two concrete representations of those germs. Moving the
inside contour outward to $|\xi|=R_+$ extends $C_N^{\rm in}f$ to
$|\zeta|<R_+$. Moving the outside contour inward to $|\xi|=R_-$
extends $C_N^{\rm out}f$ to $|\zeta|>R_-$. Cauchy's theorem identifies
these with the original transforms on their initial domains; their
extensions overlap on the entire collar and satisfy the same jump.
Only the analyticity of $f$ in the collar is needed, not a holomorphic
extension of $f$ to zero or infinity.

Set
\[
 F_s^{\rm in/out}(\tau,\zeta)
   =C_N^{\rm in/out}(B_s(\tau,\cdot)U)(\zeta).
\]
The fixed-contour representations make these transforms jointly
holomorphic in their indicated variables and entire in $\tau$. For
$R=R_c$ or $R_+$ and $|\zeta|<R$, the elementary contour estimate is
\[
 |F_s^{\rm in}(\tau,\zeta)|
 \le Ce^{C_0|\tau|}
       \frac{(|\zeta|/R)^N}{1-|\zeta|/R},
\]
using the original inside function for $R=R_c$ and its continuation for
$R=R_+$. The constants may depend on the fixed circles and $N$.

\smallskip\noindent
\emph{Step 3: sum the vanishing rows inside the source circle.}
For a function $f$ on the source collar let $[f]_n$ denote its Laurent
coefficient at index $n$. These are the Fourier coefficients on the
unit action circle whenever the same germ is defined there. The
selected-row condition and \eqref{vol:eq:borel-laplace-symbol} give
\[
 0=[b_nU]_n=[a_{s,0}U]_n+
       \int_0^\infty e^{-nt}[B_s(t,\cdot)U]_n\,dt,
 \qquad n\ge N,\quad n\equiv s\pmod P.
\]
Taking the coefficient under the integral is justified by the uniform
Laplace bound on the source circle. The finite identity
\[
 \frac1P\sum_{l=0}^{P-1}\omega^{l(n-s)}
 =\begin{cases}1,&n\equiv s\pmod P,\\0,&n\not\equiv s\pmod P\end{cases}
\]
selects exactly the required residue. Multiplying the row identities
by $\zeta^n$ and summing therefore yields, for $|\zeta|<R_c$,
\[
 \frac1P\sum_{l=0}^{P-1}\omega^{-sl}
 \left\{C_N^{\rm in}(a_{s,0}U)(\omega^l\zeta)
   +\int_0^\infty
      F_s^{\rm in}(t,e^{-t}\omega^l\zeta)\,dt\right\}=0.
\]
To justify this summation, write $\chi=|\zeta|/R_c<1$. The source-circle
coefficient estimate gives
$|[B_s(t,\cdot)U]_n|\le Ce^{C_0t}R_c^{-n}$ for $t\ge0$. Hence
\[
 \sum_{n\ge N}e^{-nt}|[B_s(t,\cdot)U]_n|\,|\zeta|^n
 \le C\frac{\chi^N e^{-(N-C_0)t}}{1-\chi e^{-t}}.
\]
The right side is integrable on $[0,\infty)$ because $\chi<1$ and
$N>C_0$. All interchanges of the Laurent sum, the Laplace integral,
and the finite residue average are thus absolutely convergent.
The identity concerns transforms of source data; it does not evaluate
$U$ throughout the inside disk.

\smallskip\noindent
\emph{Step 4: deform the Laplace contour on the inner side.}
Choose $\delta>0$ with
\[
 R_ce^{-\delta}>R_-,\qquad R_ce^\delta<R_+,
\]
and first take $|w|<\delta$ with $\operatorname{Re}w<0$. Put
$\zeta=R_ce^w$. For fixed $s,l$, define
\[
 \mathscr I_{sl}(\tau;w)
   =F_s^{\rm in}(\tau,\zeta_l e^{w-\tau}).
\]
For $T>0$, the region between the segment $[0,T]$ and the broken path
$[0,w]\cup[w,w+T]\cup[w+T,T]$ consists of points
$\tau=vw+aT$, $0\le v,a\le1$. On that region,
\[
 \operatorname{Re}(w-\tau)
   =(1-v)\operatorname{Re}w-aT\le0,
 \qquad |\zeta_l e^{w-\tau}|\le R_c.
\]
The outward-contour representation on $|\xi|=R_+$ therefore makes
$\mathscr I_{sl}(\tau;w)$ holomorphic on a neighborhood of the entire
finite deformation region, including $\tau=w$. At that endpoint its
spatial argument is $\zeta_l$; no singular Cauchy boundary value is used.
Cauchy's theorem gives the finite identity
\[
 \int_0^T\mathscr I_{sl}(\tau;w)\,d\tau
 =\int_0^w\mathscr I_{sl}(\tau;w)\,d\tau
   +\int_w^{w+T}\mathscr I_{sl}(\tau;w)\,d\tau
   +\int_{w+T}^T\mathscr I_{sl}(\tau;w)\,d\tau.
\]
On the last segment, $\tau=T+vw$ and
$|\zeta_l e^{w-\tau}|\le R_ce^{-T}$. The inside-transform estimate
from Step~2 and $|\tau|\le T+|w|$ imply
\[
 \left|\int_{w+T}^T\mathscr I_{sl}(\tau;w)\,d\tau\right|
 \le C_w e^{-(N-C_0)T}\longrightarrow0.
\]
The original and translated rays converge by the same estimates.
Parametrizing the translated ray by $\tau=w+t$ and letting
$T\to\infty$ yields
\begin{align*}
 &\int_0^\infty F_s^{\rm in}(t,\zeta_l e^{w-t})\,dt\\
 &\quad=\int_0^w
     F_s^{\rm in}(\tau,\zeta_l e^{w-\tau})\,d\tau
    +\int_0^\infty
     F_s^{\rm in}(w+t,\zeta_l e^{-t})\,dt.
\end{align*}
Only the Cauchy transforms, not a continuation of $U$, were evaluated
along the deformation.

\smallskip\noindent
\emph{Step 5: extract the Cauchy jump and identify the forcing.}
For $|w|<\delta$ and $\tau\in[0,w]$, the point
$\zeta_l e^{w-\tau}$ lies in the original collar. The jump identity
there gives
\[
 F_s^{\rm in}(\tau,z)
   =B_s(\tau,z)U(z)+F_s^{\rm out}(\tau,z).
\]
In its term containing $U$, substitute $\sigma=w-\tau$. The reversal
of the endpoints and $d\tau=-d\sigma$ cancel, leaving
\begin{align*}
 &\int_0^w
   B_s(\tau,\zeta_l e^{w-\tau})U(\zeta_l e^{w-\tau})\,d\tau\\
 &\quad=\int_0^w
   B_s(w-\sigma,\zeta_l e^\sigma)U(\zeta_l e^\sigma)\,d\sigma.
\end{align*}
For the nonintegral term, use
\[
 C_N^{\rm in}(a_{s,0}U)(\zeta_l e^w)
 =a_{s,0}(\zeta_l e^w)U(\zeta_l e^w)
     +C_N^{\rm out}(a_{s,0}U)(\zeta_l e^w).
\]
Substituting these formulas and the contour identity into Step~3,
the terms containing $U$ directly are
\[
 \sum_{l=0}^{P-1}M_{sl}(w)u_l(w)
   +\sum_{l=0}^{P-1}\int_0^w
                 \mathcal B_{sl}(w,\sigma)u_l(\sigma)\,d\sigma,
\]
with exactly the matrix and kernel in the statement. The remaining
terms, moved to the right side, define
\begin{align*}
 g_s(w)=-\frac1P\sum_{l=0}^{P-1}\omega^{-sl}
 \bigg\{&C_N^{\rm out}(a_{s,0}U)(\zeta_l e^w)\\
 &+\int_0^w F_s^{\rm out}(w-\sigma,\zeta_l e^\sigma)\,d\sigma\\
 &+\int_0^\infty F_s^{\rm in}(w+t,\zeta_l e^{-t})\,dt\bigg\}.
\end{align*}
The sign of the Volterra term is positive, whereas all three forcing
terms carry the outer minus sign. We have proved
\eqref{vol:eq:volterra-system} so far for small
$\operatorname{Re}w<0$. We next justify its continuation across zero
and the larger domain claimed for $g$.

\smallskip\noindent
\emph{Step 6: prove holomorphy of the three forcing terms.}
The initial outside-transform term is holomorphic for
$\operatorname{Re}w>0$, since $|\zeta_l e^w|>R_c$. Its representation
on $|\xi|=R_-$ also makes it holomorphic for $w$ near zero on both
sides of the imaginary axis.

For the finite integral, parametrize its straight segment by
$\sigma=vw$:
\[
 \int_0^w F_s^{\rm out}(w-\sigma,\zeta_l e^\sigma)\,d\sigma
 =w\int_0^1F_s^{\rm out}((1-v)w,\zeta_l e^{vw})\,dv.
\]
For $\operatorname{Re}w>0$ and $0\le v\le1$, the evaluation point
has modulus at least $R_c$. In the inward-contour representation,
\[
 |\zeta_l e^{vw}-\xi|\ge R_c-R_->0,
 \qquad |\xi|=R_-.
\]
For each $v$ the integrand is holomorphic in $w$, and the fixed gap and
Borel bound give uniform bounds on every compact $w$-set in that
half-plane, uniformly in $v$. Integration over $[0,1]$ is therefore
holomorphic. This also treats $v=0$ without a principal value. For
$|w|<\delta$, the lower bound $|\zeta_l e^{vw}|>R_-$ gives the same
holomorphic representation across $w=0$.

For the infinite term, the evaluation point $\zeta_l e^{-t}$ is
independent of $w$ and lies in $|\zeta|\le R_c$ for every $t\ge0$.
Set $q_+=R_c/R_+<1$. The outward-contour estimate from Step~2 gives,
for every $w\in\mathbb C$ and $t\ge0$,
\begin{align*}
 |F_s^{\rm in}(w+t,\zeta_l e^{-t})|
 &\le Ce^{C_0|w+t|}
           \frac{q_+^N e^{-Nt}}{1-q_+e^{-t}}\\
 &\le C' e^{C_0|w|}e^{-(N-C_0)t}.
\end{align*}
In particular, the bound remains uniform at $t=0$ because
$1-q_+>0$. For each $t$ the integrand is entire in $w$, and this
integrable majorant is uniform on every compact subset of the
$w$-plane. The infinite integral is consequently entire in $w$.
No value of $U$ at $\zeta_l e^{-t}$ is used: that point is the
argument of its inside Cauchy transform.

Combining the three terms proves that $g$ is holomorphic on the whole
right half-plane and also on a neighborhood of zero. The same fixed
$N$ suffices; the bound is compact-uniform in $w$, not a claim that
$g$ is bounded on the entire half-plane.

\smallskip\noindent
\emph{Step 7: continue the local identity and verify its data dependence.}
For $w$ in a sufficiently small disk about zero, all arguments of $U$
in $u_l(w)$ and in the Volterra integral remain in the original collar.
The left side of \eqref{vol:eq:volterra-system} is holomorphic there,
as is the forcing by Step~6. Their difference vanishes on the part of
the disk with $\operatorname{Re}w<0$. The identity theorem therefore
proves the equation throughout the disk, in particular for the small
$\operatorname{Re}w>0$ required by the proposition. At zero it gives
$M(0)u(0)=g(0)$, since the finite integral vanishes.

Every term defining $g$ is a Cauchy transform of $a_{s,0}U$ or
$B_s(\tau,\cdot)U$ on the original source circle. The full trace of $U$
on that circle, together with the known symbols, therefore determines
the forcing. The inward and outward contours merely represent analytic
extensions of these same transforms; they introduce neither new unknown
data nor an assumed continuation of $U$. Low and negative Laurent modes
are retained in the outside transforms, and no negative-row hypothesis
has been used.

The source-circle invertibility from Step~1 gives an invertible leading
matrix locally, but the derivation of the identity did not divide by
$M$. The conclusion is an exact local equation for the original germ
with globally holomorphic forcing in the stated half-plane, not a global
single-valued solution or global invertibility of its principal matrix.
\end{proof}

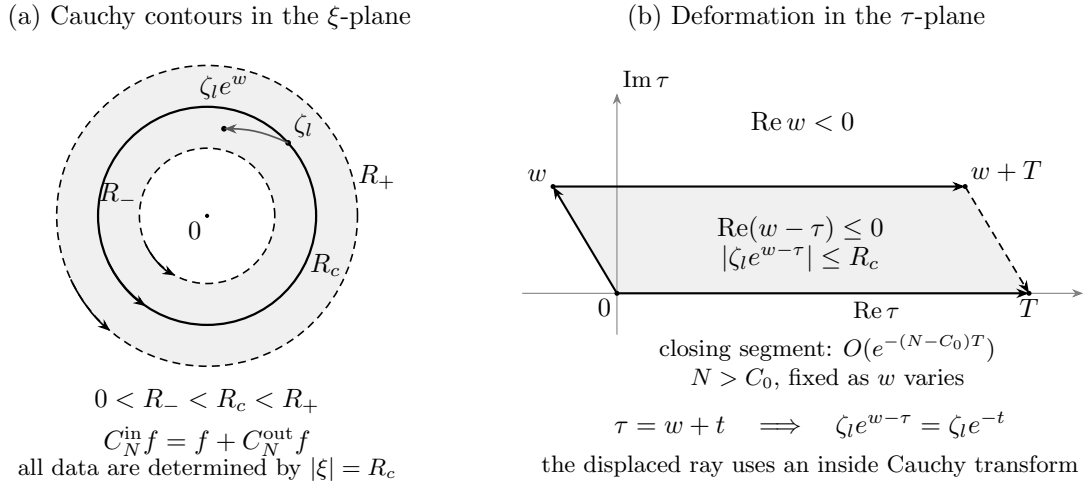
\begin{figure}[htbp]
\centering
\resizebox{.99\linewidth}{!}{%
\begin{tikzpicture}[
 x=1cm,y=1cm,font=\small,
 arr/.style={-{Stealth[length=1.7mm]},line width=.65pt},
 axis/.style={-{Stealth[length=1.5mm]},line width=.4pt,black!50},
 every node/.style={inner sep=2pt}]
\path[use as bounding box] (0,-.68) rectangle (15.80,5.77);
\node at (3.78,5.49) {(a) Cauchy contours in the $\xi$-plane};
\node at (11.72,5.49) {(b) Deformation in the $\tau$-plane};
\coordinate (O) at (3.75,2.85);
\fill[black!6] (O) circle[radius=2.00];
\fill[white] (O) circle[radius=.90];
\draw[densely dashed,line width=.6pt] (O) circle[radius=2.00];
\draw[line width=.85pt] (O) circle[radius=1.45];
\draw[densely dashed,line width=.6pt] (O) circle[radius=.90];
\draw[arr] ($(O)+(205:2.00)$) arc[start angle=205,end angle=230,radius=2.00];
\draw[arr] ($(O)+(205:1.45)$) arc[start angle=205,end angle=236,radius=1.45];
\draw[arr] ($(O)+(205:.90)$) arc[start angle=205,end angle=242,radius=.90];
\fill (O) circle[radius=.028];
\node[below left] at (O) {$0$};
\node[anchor=west] at ($(O)+(15:2.00)$) {$R_+$};
\node[anchor=west] at ($(O)+(-27:1.45)$) {$R_c$};
\node[anchor=east] at ($(O)+(162:.90)$) {$R_-$};
\coordinate (zl) at ($(O)+(42:1.45)$);
\coordinate (zin) at ($(O)+(79:1.18)$);
\fill (zl) circle[radius=.038];
\fill (zin) circle[radius=.038];
\node[anchor=south west] at (zl) {$\zeta_l$};
\node[above=9pt] at (zin) {$\zeta_l e^w$};
\draw[arr,black!65] (zl) to[bend right=12] (zin);
\node at (3.75,.39) {$0<R_-<R_c<R_+$};
\node at (3.75,-.15) {$C_N^{\rm in}f=f+C_N^{\rm out}f$};
\node[font=\footnotesize] at (3.75,-.52) {all data are determined by $|\xi|=R_c$};
\coordinate (zero) at (9.20,1.82);
\coordinate (w) at (8.35,3.24);
\coordinate (T) at (14.68,1.82);
\coordinate (wT) at (13.83,3.24);
\fill[black!6] (zero)--(w)--(wT)--(T)--cycle;
\draw[axis] (7.94,1.82)--(15.40,1.82);
\draw[axis] (9.20,1.27)--(9.20,4.48);
\node[above right,font=\footnotesize] at (9.20,4.48) {$\operatorname{Im}\tau$};
\node[below,font=\footnotesize] at (12.65,1.78) {$\operatorname{Re}\tau$};
\draw[arr,line width=.8pt] (zero)--(T);
\draw[arr,line width=.8pt] (zero)--(w);
\draw[arr,line width=.8pt] (w)--(wT);
\draw[arr,densely dashed] (wT)--(T);
\foreach \pt in {zero,w,T,wT}{\fill (\pt) circle[radius=.033];}
\node[below left] at (zero) {$0$};
\node[above left] at (w) {$w$};
\node[above right] at (wT) {$w+T$};
\node[below] at (T) {$T$};
\node at (11.67,4.11) {$\operatorname{Re}w<0$};
\node[align=center] at (11.64,2.48)
 {$\operatorname{Re}(w-\tau)\le0$\\[-1pt]
  $|\zeta_l e^{w-\tau}|\le R_c$};
\node[align=center,font=\footnotesize] at (12.00,.91)
 {closing segment: $O(e^{-(N-C_0)T})$\\$N>C_0$, fixed as $w$ varies};
\node at (11.77,.07) {$\tau=w+t\quad\Longrightarrow\quad
 \zeta_l e^{w-\tau}=\zeta_l e^{-t}$};
\node[font=\footnotesize] at (11.77,-.48)
 {the displaced ray uses an inside Cauchy transform};
\end{tikzpicture}
}
\caption{Source contours and the exact Laplace deformation in
Proposition~\ref{vol:prop:exact-volterra}. Left: all three circles lie
in a collar where $U$ is initially holomorphic. Moving the inside Cauchy
contour to $R_+$ and the outside contour to $R_-$ gives overlapping
holomorphic representations; neither a value of $U$ at zero nor a
principal-value boundary integral is used. Right: for small $w$ with
$\operatorname{Re}w<0$, the original segment $[0,T]$ is replaced by
$[0,w]$, $[w,w+T]$, and $[w+T,T]$. The closing contribution tends to
zero for the fixed cutoff $N>C_0$. The translated ray has spatial
argument $\zeta_l e^{-t}$ only inside a Cauchy transform. Both panels
are schematic; they do not prescribe uniform radii or continue $U$
outside the known collar. The resulting local identity extends across
$w=0$ by holomorphy.}
\label{vol:fig:source-contours}
\end{figure}

The source-contour representations and the Laplace-contour
deformation in the preceding proof are illustrated in
Figure~\ref{vol:fig:source-contours}.

\Needspace{18\baselineskip}
\begin{corollary}[Continuation in the exterior and interior sectors]
\label{vol:cor:sector-continuation}
Put $\mathcal S=\bigcup_{l=0}^{P-1}\omega^{-l}S$ and define
the single-valued \emph{principal matrix}
\[
\widehat M_{sl}(\zeta)=\frac1P\omega^{-sl}a_{s,0}(\omega^l\zeta),
\qquad 0\le s,l<P.
\]
Thus $M(w)=\widehat M(R_ce^w)$. Let
$\mathcal Z=\{\zeta\in\mathbb C^*\setminus\mathcal S:
                  \det\widehat M(\zeta)=0\}$. Then $\mathcal S\cup\mathcal Z$ is
closed and countable in $\mathbb C^*$, and its derived set is contained
in the discrete set $\mathcal S$.
The germ of $U$ continues along every path, starting in the original
annulus, contained in
\[
\{\,|\zeta|>R_c\,\}\setminus(\mathcal S\cup\mathcal Z).
\]
If the sufficiently large selected negative rows also vanish, an
analogous statement holds for an interior region by applying those rows
and replacing $\zeta$ by $1/\zeta$.
Neither statement asserts single-valuedness or meromorphicity of the
continued solution.
\end{corollary}
\begin{proof}
We use the hypotheses, source circle, and notation of
Proposition~\ref{vol:prop:exact-volterra}. In particular, $U$ is initially
holomorphic in an annulus, its sufficiently large selected positive rows
vanish, and the forcing $g$ is holomorphic on the right half-plane and
on a neighborhood of zero. The local equation
\eqref{vol:eq:volterra-system} holds across $w=0$. Put
$\mathfrak E=\mathcal S\cup\mathcal Z$.

\smallskip\noindent
\emph{Step 1: locate the possible accumulation points of the exceptional set.}
The finite union $\mathcal S=\bigcup_l\omega^{-l}S$ is closed and discrete
in $\mathbb C^*$, and
\[
 \zeta\notin\mathcal S\quad\Longleftrightarrow\quad
 \omega^l\zeta\notin S\quad\hbox{for every }0\le l<P.
\]
Thus every entry of $\widehat M$ is single-valued and holomorphic on
$\mathbb C^*\setminus\mathcal S$. This domain is connected: a polygonal
path in $\mathbb C^*$ meets only finitely many points of the closed
discrete set on its compact image, and small local detours avoid them.

For $\zeta=e^{ix}$, reorder the columns by $l\equiv p_0j\pmod P$ and the
rows by the residues $r_s\equiv-p_0s\pmod P$. The entries become
\[
 \frac1P\exp\!\left[\frac{ir_s}{P}
  \left(2\pi j-Y\!\left(x+\frac{2\pi p_0j}{P}\right)\right)\right].
\]
This is the Vandermonde matrix of
Proposition~\ref{band:prop:operator}, divided by $P$. Its nodes are
distinct by that proposition's twist argument, so
$\det\widehat M(e^{ix})\ne0$ for every real $x$. The determinant is
therefore not identically zero on its connected holomorphy domain.
Its zero set $\mathcal Z$ is discrete relative to that domain and is
countable, for example by a compact exhaustion.

Every point of $\mathbb C^*\setminus\mathcal S$ has a neighborhood where
$\mathfrak E$ is just the zero set of this nonzero holomorphic function.
In that neighborhood the set is closed and has no accumulation point.
All points of $\mathcal S$ are already included in $\mathfrak E$.
Consequently $\mathfrak E$ is closed and countable in $\mathbb C^*$, and
\[
 \mathfrak E'\subset\mathcal S.
\]
The union need not be discrete: determinant zeros may accumulate at
coefficient singularities. All closedness and accumulation assertions
here are relative to $\mathbb C^*$.

\smallskip\noindent
\emph{Step 2: lift the path and identify the regular kernel domain.}
Define
\[
 \Omega_\zeta=\{\zeta:|\zeta|>R_c\}\setminus\mathfrak E,
 \qquad
 \Omega_w=\{w:\operatorname{Re}w>0,\ R_ce^w\notin\mathfrak E\}.
\]
A continuous path $\gamma:[0,1]\to\Omega_\zeta$ has, after choosing its
initial logarithm, a continuous lift $w_\gamma:[0,1]\to\Omega_w$ with
$R_ce^{w_\gamma(t)}=\gamma(t)$. This is a pathwise lift, not a choice of
a single-valued logarithm on $\Omega_\zeta$.

On $\Omega_w$, the matrix $M(w)^{-1}$ is holomorphic. Moreover,
\[
 \mathcal B_{sl}(w,\sigma)
  =\frac1P\omega^{-sl}B_s(w-\sigma,\zeta_l e^\sigma)
\]
is jointly holomorphic whenever $\sigma$ is a regular spatial history
point: if $R_ce^\sigma\notin\mathcal S$, then every
$\zeta_l e^\sigma\notin S$. Lemma~\ref{vol:lem:convergent-symbol} makes
$B_s$ entire in its first variable. Hence, for each such $\sigma$,
$\mathcal B(w,\sigma)$ is entire in $w$; the difference $w-\sigma$ need
not be small. Joint bounds hold on compact sets of these pairs. The
same facts hold for $\sigma$ in the initial regular disk through zero.

\smallskip\noindent
\emph{Step 3: restart the equation with its entire integration history.}
Suppose that the vector germ has been continued along a finite piecewise
smooth history $\Gamma$ from zero to $p\in\Omega_w$. Record the branch
on each traversal, denoted by $u_\Gamma(\sigma)$, and define
\[
 h_\Gamma(w)=g(w)-\int_\Gamma
                    \mathcal B(w,\sigma)u_\Gamma(\sigma)\,d\sigma.
\]
The history integral is the sum of integrals on its finitely many
pieces; repeated visits to a point retain their respective branches.
On each piece the chosen branch is bounded along the compact path.
The kernel is entire in $w$ with compact-uniform bounds, so integration
shows that the history term is entire in $w$. Thus $h_\Gamma$ is
holomorphic on the right half-plane. No earlier part of the history,
including a loop, is removed.

Choose a disk $D(p,r)$ whose closure has a regular neighborhood in
$\Omega_w$. On this disk the restarted equation is
\[
 M(w)v(w)+\int_p^w\mathcal B(w,\sigma)v(\sigma)\,d\sigma
      =h_\Gamma(w),
\]
where the new integral follows the straight segment. On the Banach space
$H^\infty(D(p,r);\mathbb C^P)$, define
\[
 \begin{aligned}
 F_\Gamma(w)&=M(w)^{-1}h_\Gamma(w),\\
 (\mathcal V_pv)(w)&=\int_p^w
          M(w)^{-1}\mathcal B(w,\sigma)v(\sigma)\,d\sigma.
 \end{aligned}
\]
The forcing is bounded and holomorphic on the disk. Parametrizing the
segment by $\sigma=p+t(w-p)$ shows that $\mathcal V_p$ preserves
holomorphy. For any fixed vector norm and its induced matrix norm,
\[
 \|\mathcal V_p\|\le rL_p,\qquad
 L_p=\sup_{w,\sigma\in\overline{D(p,r)}}
                \|M(w)^{-1}\mathcal B(w,\sigma)\|.
\]
If $rL_p\le1/2$, the Neumann series gives the unique bounded holomorphic
solution
\[
 v=\sum_{k\ge0}(-\mathcal V_p)^kF_\Gamma,
 \qquad \|v\|_\infty\le2\|F_\Gamma\|_\infty.
\]
The permitted radius depends on the coefficient bounds, not on the size
of the accumulated history or of the solution. Only the forcing changes
when the history changes.

\smallskip\noindent
\emph{Step 4: verify agreement of consecutive germs.}
Suppose a previously constructed branch $u_{\rm old}$ is holomorphic on
a disk centered at $a$, with history $\Gamma_{\rm old}$ ending at $a$.
Choose the next center $p$ inside that disk and append the segment
$[a,p]$, using this branch, to form $\Gamma$. For $w$ sufficiently close
to $p$, the triangle with vertices $a,p,w$ stays in the old disk.
Holding $w$ fixed, Cauchy's theorem in $\sigma$ gives
\[
 \begin{aligned}
 \int_a^w\mathcal B(w,\sigma)u_{\rm old}(\sigma)\,d\sigma
  &=\int_a^p\mathcal B(w,\sigma)u_{\rm old}(\sigma)\,d\sigma\\
  &\quad+\int_p^w\mathcal B(w,\sigma)u_{\rm old}(\sigma)\,d\sigma.
 \end{aligned}
\]
The old branch therefore satisfies the restarted equation with the
updated forcing $h_\Gamma$ on a neighborhood of $p$. On a smaller disk
centered at $p$, both old and new branches are bounded and solve the
same contraction equation. They agree there, and the identity theorem
gives agreement on the connected overlap of the two disks.

This proves equality of germs, not merely equality of their values at
$p$. Only the terminal segment inside the old disk has been deformed;
all earlier loops and their branch determinations remain in $\Gamma$.
Agreement is asserted for consecutive germs along this history, not for
all branches encountered at the same point after different histories.

\smallskip\noindent
\emph{Step 5: start the history and follow every continuous path.}
Proposition~\ref{vol:prop:exact-volterra} supplies the original vector
$u_l(w)=U(\zeta_l e^w)$ and its equation on a disk crossing $w=0$. Shrink that disk to avoid
coefficient singularities and zeros of $\det M$.
Choose a small positive real point $p_*$ in that disk and use the known
branch along $[0,p_*]$ as the initial history. The preceding restart
argument applies at $p_*$, now strictly inside the right half-plane.
This separates the boundary seed from all subsequent continuation disks.

Write the original annulus of $U$ as $r_0<|\zeta|<R_0$, with
$r_0<R_c<R_0$. The part of its logarithmic preimage to the right of the
source circle is the convex strip
\[
 \mathcal A_w=\{w:0<\operatorname{Re}w<\log(R_0/R_c)\}.
\]
If the requested path starts in that annulus, then $w_\gamma(0)$ lies
in this strip. Connect $p_*$ to $w_\gamma(0)$ there by a polygonal path
avoiding the lifted exceptional set. Such a connector exists: the
lifted set is closed and countable, and an intermediate point outside
the countably many lines through an endpoint and an excluded point
gives a two-segment path in the convex strip. Both endpoints are regular.
The original vector is holomorphic throughout $\mathcal A_w$, so the
overlap argument identifies continuation along the connector with its
original germ at the required starting point.

The connector and $w_\gamma([0,1])$ have a compact union in $\Omega_w$.
Choose a compact regular neighborhood $\mathcal K$ of this union, and put
\[
 L=\sup_{w,\sigma\in\mathcal K}
             \|M(w)^{-1}\mathcal B(w,\sigma)\|<\infty.
\]
Choose one $0<r<\pi$ so that disks of radius $r$ centered on this union
have their closures in the interior of $\mathcal K$ and $rL\le1/2$. This radius
may depend on the chosen path, but is independent of all branch values.
By uniform continuity, a finite subdivision of the concatenated
connector and lifted path has each subpath contained in the disk of
radius $r/2$ centered at its entry point.

Apply Steps~3--4 successively on the disks of radius $r$. Straight
segments between successive entry points update the integration history;
each is compactly contained in the preceding solution disk. The
holomorphic solution on that disk continues its germ along the entire
original subpath, whether or not that subpath is rectifiable. Thus no
integral is taken along a possibly infinite-length original path.
There are only finitely many steps. Each updated history is finite and
piecewise smooth, and its size does not affect the contraction bound,
so the process reaches the endpoint. It retains any loops of the
polygonal history rather than resetting the integral after a return.

Since $r<\pi$, the map $w\mapsto R_ce^w$ is injective on each disk.
Its local inverse and the zeroth component $u_0(w)$ give a holomorphic
$\zeta$-germ along the original path. This proves the exterior conclusion
for every continuous path with compact parameter interval whose image,
including its endpoint, lies in $\Omega_\zeta$.

\smallskip\noindent
\emph{Step 6: convert negative rows into the interior continuation problem.}
To distinguish angular and complex variables, write
$U_{\rm ang}(x)=U(e^{ix})$.
For a positive integer $n$, the selected row at index $-n$ is
\[
 \frac1{2\pi}\int_0^{2\pi}
             U_{\rm ang}(x)e^{inH_n(x)}H_n'(x)\,dx.
\]
Set $y=-x$, $\widetilde U_{\rm ang}(y)=U_{\rm ang}(-y)$, and
$\widetilde H_n(y)=-H_n(-y)$. Since
$\widetilde H_n'(y)=H_n'(-y)$, substitution and periodicity turn this
row into
\[
 \frac1{2\pi}\int_0^{2\pi}
     \widetilde U_{\rm ang}(y)e^{-in\widetilde H_n(y)}
                     \widetilde H_n'(y)\,dy.
\]
It is a positive selected row in the reciprocal complex coordinate
$\xi=e^{iy}=1/\zeta$, with $\widetilde U(\xi)=U(1/\xi)$.
No conjugation of the unknown $U$ occurs.

The transformed symbols have the same properties as in
Lemma~\ref{vol:lem:convergent-symbol}. Indeed, put
$\widetilde H_s(\varepsilon,y)=-H_s(\varepsilon,-y)$. The real
reference phase maps are real-valued for real $\varepsilon$ and $x$.
Their amplitudes consequently satisfy, initially for real $\varepsilon,y$,
\[
 \widetilde b_s(\varepsilon,e^{iy})
       =\overline{b_s(\varepsilon,e^{-iy})}.
\]
Taking Taylor coefficients and then continuing holomorphically gives
\[
 \widetilde a_{s,k}(\xi)=\overline{a_{s,k}(\overline\xi)},\qquad
 \widetilde B_s(t,\xi)=\overline{B_s(\overline t,\overline\xi)}.
\]
These are holomorphic reflected functions, not antiholomorphic ones.
Their coefficient set is $\widetilde S=\{\overline z:z\in S\}$; they
inherit single-valuedness, entire dependence on $t$, and compact-uniform
exponential bounds. Their convergent Taylor series and the same Laplace
moments give the exact symbol formula for the transformed positive rows.

Define the transformed principal matrix by the same formula as
$\widehat M$, using $\widetilde a_{s,0}$, and let $\mathsf J_{\rm rev}$
be the column permutation $l\mapsto-l\pmod P$. On the unit circle,
\[
 \widetilde{\widehat M}_{sl}(e^{iy})
     =\overline{\widehat M_{s,-l}(e^{-iy})},\qquad
 \widetilde{\widehat M}(e^{iy})
     =\overline{\widehat M(e^{-iy})}\,\mathsf J_{\rm rev}.
\]
Thus the transformed principal matrix is invertible there. The
construction of Proposition~\ref{vol:prop:exact-volterra}, which uses
only these symbol properties and the vanishing positive rows, applies
to $\widetilde U$. Steps~1--5 then give its exterior continuation off
the transformed closed countable exceptional set
$\widetilde{\mathfrak E}$. If its source radius is
$\widetilde R_c$, inversion gives continuation of $U$ along paths in
\[
 \{0<|\zeta|<\widetilde R_c^{-1}\}
       \setminus\{\xi^{-1}:\xi\in\widetilde{\mathfrak E}\}.
\]
Inversion is a biholomorphism of $\mathbb C^*$, so it preserves the
relative closedness, countability, and stated accumulation property.
This proves the interior assertion without a reality assumption on $U$;
no regularity at zero is asserted.

In both directions, local uniqueness holds for a fixed accumulated
history. Different histories ending at the same point can give different
branches. The argument therefore proves neither single-valuedness nor
finite monodromy, meromorphicity, or removability at the exceptional set.
\end{proof}

\begin{remark}[Limit of the action-coordinate argument]
This continuation excludes a full circular natural boundary on the exterior side and, under the additional negative-row hypothesis, on the interior side. It does not prove finite monodromy at determinant zeros, pole-type singularities, or removability of their clusters at $\mathcal S$. The next lemma transfers this exterior continuation directly to the
physical coordinate.
\end{remark}

\paragraph{Exceptional sets and their ambient spaces.}
The table records the ambient space for each derived-set assertion.
Its last row anticipates Lemma~\ref{vol:lem:physical-exterior-continuation}
and uses $D=\{z\in\C:|z|>1\}$. All five sets are independent of the
auxiliary source radius $R_c$.
\begin{center}
\begingroup
\small
\renewcommand{\arraystretch}{1.18}
\begin{tabular}{@{}p{.25\linewidth}p{.23\linewidth}p{.46\linewidth}@{}}
\toprule
Set & Ambient space & Property used in the continuation argument\\
\midrule
$S$ & $\C^*$ & Closed and discrete; it contains the possible
coefficient singularities.\\
$\mathcal S=\bigcup_{l=0}^{P-1}\omega^{-l}S$
 & $\C^*$ & Closed and discrete.\\
$\mathcal Z=\{\det\widehat M=0\}$
 & $\C^*\setminus\mathcal S$ & Relatively closed and discrete;
its accumulation points in $\C^*$ can lie in $\mathcal S$.\\
$\mathfrak E=\mathcal S\cup\mathcal Z$
 & $\C^*$ & Closed and countable, with
$\mathfrak E'\subset\mathcal S$; not necessarily discrete.\\
$B$ & $D$ & Closed and countable, with a finite relative derived
set $B'$; not necessarily finite.\\
\bottomrule
\end{tabular}
\endgroup
\end{center}

\begin{lemma}[Physical exterior continuation with a finite accumulation set]
\label{vol:lem:physical-exterior-continuation}
For a fixed noncircular reference ellipse and the fixed rational base used
in the Borel--Volterra construction, let $V(z)$ be the physical germ of an analytic displacement whose selected positive rows vanish. There is a closed countable set
$B\subset\{|z|>1\}$ with finite derived set such that the physical
boundary germ $V(z)$ continues along every path in
$\{|z|>1\}\setminus B$ starting in its original exterior collar.
Closedness and derived sets here are relative to the finite exterior
domain $\{|z|>1\}\subset\mathbb C$.
Here $z=e^{i\varphi}$ is the physical variable. Continuations may have
monodromy; no assertion of meromorphicity is made.
\end{lemma}
\begin{proof}
Fix the base parameter $m_0=m(\rho_0)$ and put
$K=K(m_0)$, $K'=K(1-m_0)$, and
$\mathcal P=\mathcal P_{m_0}$. As throughout this section, $K'$ is the
complementary integral, not a derivative. Write
\[
 D=\{z\in\C:|z|>1\},\qquad
 D_{\mathrm{sph}}=D\cup\{\infty\}.
\]
The variable $z=e^{i\varphi}$ is the complexified boundary parameter.
We keep the finite target $D$ separate from its spherical version and
specify the ambient space whenever a derived set is used.

\smallskip\noindent
\emph{Step 1: put the proper Jacobi map in action coordinates.}
By Lemma~\ref{arith:jacobi}, the lower cylinder
\[
 \mathcal C_-=\{-2K'<\operatorname{Im}u<0\}/(4K\Z)
\]
is the full inverse image of $D_{\mathrm{sph}}$ under the compact-torus
map, and $\mathcal P:\mathcal C_-\to D_{\mathrm{sph}}$ is proper of
degree two, counted with multiplicity. The map
\[
 \chi(u)=\exp\!\left(\frac{i\pi(u-K)}{2K}\right),\qquad
 |\chi(u)|=\exp\!\left(-\frac{\pi\operatorname{Im}u}{2K}\right),
\]
identifies this cylinder biholomorphically with
\[
 A=\{1<|\zeta|<R_*\},\qquad R_*=e^{\pi K'/K}.
\]
Define $\Psi=\mathcal P\circ\chi^{-1}$ on $A$. This is single-valued:
a local logarithm used to recover $u$ changes by a multiple of $4K$,
which is a period of $\mathcal P$.

The two simple poles $-iK'$ and $2K-iK'$ correspond to
\[
 \Pi=\{-i e^{\pi K'/(2K)},\ i e^{\pi K'/(2K)}\},\qquad
 X=A\setminus\Pi.
\]
Thus $\Psi:X\to D$ is holomorphic, onto, and still proper of degree two.
Indeed, for any compact $L\subset D$, its full preimage under the torus
map is compact and lies in $\mathcal C_-$; it contains neither pole,
since $\infty\notin L$. Its action-coordinate image is therefore compact
in $X$. Removing the poles together with infinity is essential to this
finite-target properness statement.

The Jacobi derivative identities used in Lemma~\ref{arith:focal} give
\[
 \mathcal P'(u)=i\operatorname{dn}(u)\mathcal P(u).
\]
On the pole-deleted lower cylinder, $\mathcal P$ is finite and nonzero.
Its two critical points are consequently $K-iK'$ and $3K-iK'$, the
zeros of $\operatorname{dn}$ there. The identity
\[
 \mathcal P(t-iK')=
 \frac{1+\operatorname{dn}t}{\sqrt{m_0}\operatorname{sn}t}
\]
from Lemma~\ref{arith:jacobi}, evaluated at $t=K$ and shifted by $2K$,
gives their branch-value set
\[
 C_{\mathrm{br}}=
 \left\{\frac{1+\sqrt{1-m_0}}{\sqrt{m_0}},
          -\frac{1+\sqrt{1-m_0}}{\sqrt{m_0}}\right\}\subset D.
\]
These are branch values of the base Jacobi map, not the focal values of
the reference ellipse.

\smallskip\noindent
\emph{Step 2: restrict the action exceptional set to the finite annulus.}
Use the notation of Corollary~\ref{vol:cor:sector-continuation} and put
\[
 \mathfrak E=\mathcal S\cup\mathcal Z,\qquad
 \Sigma=\mathfrak E\cap A,\qquad \Sigma_X=\Sigma\cap X.
\]
That corollary gives closedness and countability of $\mathfrak E$ in
$\C^*$ and $\mathfrak E'\subset\mathcal S$, where $\mathcal S$ is
closed and discrete in $\C^*$. Hence $\Sigma_X$ is closed and countable
relative to $X$, and its derived set taken in $X$ satisfies
\[
 \Sigma_X'\subset\mathcal S\cap X.
\]
The closure $\overline A$ is compact in $\C^*$, so
$\mathcal S\cap\overline A$ is finite. Thus $\Sigma_X'$ is contained
in a fixed finite set, although $\Sigma$ itself need not be finite.

For completeness, the containing set can be read explicitly from
Lemma~\ref{vol:lem:convergent-symbol}. With
$\mathfrak q=e^{-\pi K'/K}$, we have $R_*=\mathfrak q^{-1}$ and
\[
 S=\{i^j\mathfrak q^{\ell+1/2}:
                 j=0,1,2,3,\ \ell\in\Z\},\qquad
 \mathcal S=\bigcup_{l=0}^{P-1}\omega^{-l}S.
\]
The inequalities $1<\mathfrak q^{\ell+1/2}<\mathfrak q^{-1}$ force
$\ell=-1$. Therefore
\[
 \mathcal S\cap A=
 \{\omega^{-l}i^j\mathfrak q^{-1/2}:
                  0\le l<P,\ 0\le j<4\},
\]
with repetitions removed. It has at most $4P$ points; the poles in
$\Pi$ are omitted when intersecting with $X$. This counts possible
accumulation locations, not all determinant zeros.

\smallskip\noindent
\emph{Step 3: project the set and control finite accumulation.}
Define the physical exceptional set by
\[
 B=\Psi(\Sigma_X)\cup C_{\mathrm{br}}\subset D.
\]
It is countable as the image of a countable set together with a finite
set. To prove relative closedness, suppose
$z_n\in\Psi(\Sigma_X)$ and $z_n\to z_\infty\in D$. Choose
$\zeta_n\in\Sigma_X$ with $\Psi(\zeta_n)=z_n$.
The set $L=\{z_\infty\}\cup\{z_n:n\ge1\}$ is compact in $D$.
Properness makes $\Psi^{-1}(L)$ compact in $X$, so a subsequence
converges to $\zeta_\infty\in X$. Relative closedness of $\Sigma_X$
gives $\zeta_\infty\in\Sigma_X$, and continuity gives
$\Psi(\zeta_\infty)=z_\infty$. Thus $\Psi(\Sigma_X)$, and hence $B$,
is closed in $D$. This conclusion uses properness, not merely
continuity of the projection.

Next take $z_\infty\in B'$, with the derived set now taken in $D$, and
choose pairwise distinct $z_n\in B$ converging to $z_\infty$.
Discard the finitely many terms belonging to $C_{\mathrm{br}}$, then
choose $\zeta_n\in\Sigma_X$ mapping to $z_n$. These preimages are
pairwise distinct because their images are. The same properness
argument gives a subsequence converging in $X$, and its limit belongs
to $\Sigma_X'$. Consequently,
\[
 B'\subset\Psi(\Sigma_X')
       \subset\Psi(\mathcal S\cap X).
\]
The last set is finite by Step~2. In particular, preimages of a sequence
converging to a finite exterior value cannot escape to a pole or a
boundary component of the cylinder. The finite bound consists of
physical images of rotated action-coordinate coefficient points; it
need not reduce to the two branch values or to the focal values.

\smallskip\noindent
\emph{Step 4: remove the full exceptional fibers and choose the initial sheet.}
Set
\[
 D_B=D\setminus B,\qquad X_B=\Psi^{-1}(D_B).
\]
The restriction $\Psi:X_B\to D_B$ is proper and has no critical point,
since all branch values have been deleted. Every fiber has two regular
points. Local holomorphic inverses at these two points, together with
properness to exclude further preimages over a sufficiently small
target disk, show that this is an unbranched two-sheeted covering.
We use the full inverse image of $D_B$, not merely $X\setminus\Sigma_X$:
if one member of a fiber is exceptional, the entire physical value has
been removed. In particular $X_B\cap\Sigma_X=\varnothing$.

The initial sheet is determined by the original coordinate relation.
On the real action circle,
\[
 \Psi(e^{ix})=e^{i\Phi_{\rho_0}(x)},\qquad
 U(e^{ix})=V(\Psi(e^{ix})).
\]
This boundary map is an analytic diffeomorphism. Its derivative is
nonzero there because $\Phi_{\rho_0}'>0$. The holomorphic inverse
function theorem and compactness of the circle give a biholomorphic
restriction near it: otherwise collisions arbitrarily near the circle
would limit to the same point, contradicting local injectivity there.
Shrink the representatives of the two original germs accordingly.
There is then a physical exterior collar $\mathcal C_0$ with an inverse
branch $\psi_0$ whose image lies in the known annulus of $U$ adjacent
to $|\zeta|=1$, and
\[
 V=U\circ\psi_0\quad\hbox{on }\mathcal C_0.
\]
The image is on the $|\zeta|>1$ side by the strict cylinder inequalities
of Lemma~\ref{arith:jacobi}. This does not require the other Jacobi
sheet, near $|\zeta|=R_*$, to lie in the initially known annulus.

Let $\gamma:[0,1]\to D_B$ be a continuous path starting in
$\mathcal C_0\setminus B$. Lift it through the covering with initial
point $\psi_0(\gamma(0))$. The resulting continuous path
\[
 \widetilde\gamma:[0,1]\to X_B,\qquad
 \Psi\circ\widetilde\gamma=\gamma,
\]
avoids $\Sigma$. The path-lifting property applies to continuous paths,
not only piecewise smooth ones: local inverse disks and a finite
subdivision of the compact parameter interval give the lift.

\smallskip\noindent
\emph{Step 5: choose an admissible source circle without changing $B$.}
Compactness of the lifted path in $A$ gives
\[
 r_\gamma:=\min_{0\le t\le1}|\widetilde\gamma(t)|>1.
\]
Choose $1<R_c<r_\gamma$ sufficiently close to one that its circle lies
in the original annulus of $U$ and meets the source-circle and symbol
conditions of Proposition~\ref{vol:prop:exact-volterra}. Such a choice
is possible: the principal determinant is nonzero on the unit circle,
and hence on a fixed smaller annular neighborhood, while $U$ is
holomorphic around that circle. The whole lifted path now lies in
\[
 \{|\zeta|>R_c\}\setminus\mathfrak E
\]
and starts in the original action annulus. Corollary~\ref{vol:cor:sector-continuation}
therefore continues $U$ along this entire lift, with all branch history
retained. Only the selected positive rows are used.

The geometric lift was constructed before choosing $R_c$, without
using any continued values of $U$. Moreover, the set $B$ in Step~3 is
fixed before the path is chosen. Indeed, $S$, $\mathcal S$, and
\[
 \mathcal Z=\{\zeta\notin\mathcal S:\det\widehat M(\zeta)=0\},\qquad
 \widehat M_{sl}(\zeta)=P^{-1}\omega^{-sl}a_{s,0}(\omega^l\zeta),
\]
are defined in the actual action variable and do not depend on $R_c$.
Changing $R_c$ changes the logarithmic representation
$M(w)=\widehat M(R_ce^w)$ and the forcing, not these exceptional sets
or their physical images. Thus one fixed $B$ works for every path,
although the auxiliary source radius may depend on that path.

\smallskip\noindent
\emph{Step 6: compose local germs and cover every original starting point.}
Along $\gamma$, choose successive small physical disks on which the
inverse branches $\psi_j$ follow $\widetilde\gamma$. Let $U_j$ be the
action germ reached by continuation on the corresponding lifted
neighborhood, shrinking the disks as necessary. Define
\[
 V_j=U_j\circ\psi_j.
\]
These are holomorphic germs. The first equals the original $V$ by the
choice of $\psi_0$. Consecutive inverse branches agree on the overlap
belonging to the chosen lift, and consecutive action germs agree there
by analytic continuation. Hence the $V_j$ agree as consecutive germs
and continue $V$ along the entire physical path, including its endpoint.
No derivative or finite-length hypothesis on the original path is used.

This also treats a path starting farther out in any originally chosen
exterior collar where $V$ is known. Choose a regular point in the
smaller collar $\mathcal C_0$ and connect it to the specified starting
point inside the known physical collar, avoiding $B$. Such a polygonal
connector exists because removing a closed countable subset from an
open annulus does not destroy path connectivity: cover an annular path
by finitely many convex disks, choose regular joining points, and in
each disk choose an intermediate point outside the countably many lines
through a joining point and an excluded point. The resulting segments
avoid $B$. Apply Steps~4--6 to the concatenation of this connector and
the requested path. Along the connector, uniqueness of analytic
continuation identifies the constructed germs with the already known
function $V$. Thus the germ at the requested starting point is the
original one, not a separately selected branch.

Only consecutive germs on the chosen lift have been identified. There
is no assumption that a continued action germ is invariant under the
Jacobi deck involution. A physical loop may exchange the two sheets,
and the action germ itself may acquire further monodromy. The degree-two
geometric covering therefore imposes no two-branch bound on the
continued function. The result is precisely all-path continuation
outside the fixed closed countable set $B$, with the finite derived-set
bound proved in Step~3. No regularity at infinity, single-valuedness,
finite monodromy, or meromorphicity is asserted.
\end{proof}

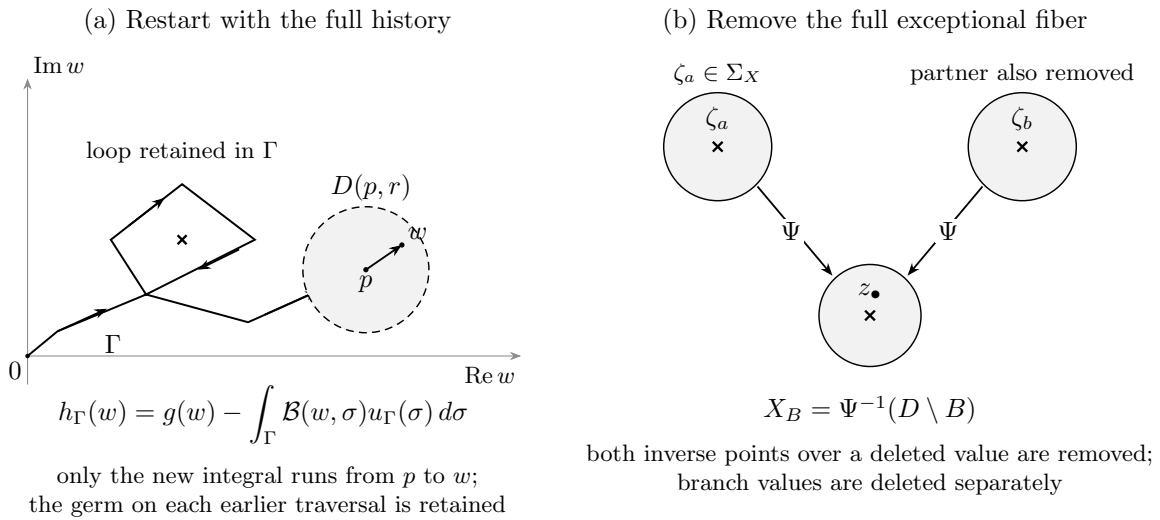
\begin{figure}[htbp]
\centering
\resizebox{.99\linewidth}{!}{%
\begin{tikzpicture}[
 x=1cm,y=1cm,font=\small,
 arr/.style={-{Stealth[length=1.7mm]},line width=.7pt},
 axis/.style={-{Stealth[length=1.5mm]},line width=.4pt,black!45},
 every node/.style={inner sep=2pt}]
\path[use as bounding box] (-.04,-1.44) rectangle (15.80,5.68);
\node at (3.74,5.36) {(a) Restart with the full history};
\node at (11.84,5.36) {(b) Remove the full exceptional fiber};
\draw[axis] (.54,.90)--(7.12,.90);
\draw[axis] (.54,.52)--(.54,4.60);
\node[below left,font=\footnotesize] at (7.12,.83) {$\operatorname{Re}w$};
\node[above right,font=\footnotesize] at (.54,4.60) {$\operatorname{Im}w$};
\coordinate (seed) at (.54,.90);
\coordinate (entry) at (2.12,1.72);
\coordinate (exception) at (2.60,2.45);
\coordinate (p) at (5.05,2.05);
\coordinate (neww) at (5.53,2.38);
\draw[line width=.75pt] (seed)--(.94,1.23)--(entry)
 --(1.65,2.45)--(2.60,3.19)--(3.57,2.45)--(entry)
 --(3.48,1.35)--(p);
\draw[arr] (.94,1.23)--(1.60,1.53);
\draw[arr] (1.88,2.63)--(2.37,3.01);
\draw[arr] (3.36,2.32)--(2.81,2.05);
\draw[arr] (3.70,1.45)--(4.41,1.77);
\draw[line width=.9pt] ($(exception)+(-.06,-.06)$)--($(exception)+(.06,.06)$);
\draw[line width=.9pt] ($(exception)+(-.06,.06)$)--($(exception)+(.06,-.06)$);
\node[font=\footnotesize] at (2.60,3.62) {loop retained in $\Gamma$};
\node at (1.68,1.06) {$\Gamma$};
\fill[black!5] (p) circle[radius=.84];
\draw[densely dashed,line width=.55pt] (p) circle[radius=.84];
\draw[arr] (p)--(neww);
\fill (seed) circle[radius=.03];
\fill (p) circle[radius=.035];
\fill (neww) circle[radius=.035];
\node[below left] at (seed) {$0$};
\node[below] at (p) {$p$};
\node[above right] at (neww) {$w$};
\node at (5.12,3.15) {$D(p,r)$};
\node[align=center] at (3.68,.12)
 {$h_\Gamma(w)=g(w)-\displaystyle\int_\Gamma
   \mathcal B(w,\sigma)u_\Gamma(\sigma)\,d\sigma$};
\node[font=\footnotesize,align=center] at (3.74,-.91)
 {only the new integral runs from $p$ to $w$;\\
  the germ on each earlier traversal is retained};
\coordinate (a) at (9.74,3.69);
\coordinate (b) at (13.80,3.69);
\coordinate (zbad) at (11.77,1.44);
\foreach \pt in {a,b}{
 \fill[black!5] (\pt) circle[radius=.72];
 \draw[line width=.55pt] (\pt) circle[radius=.72];
}
\fill[black!5] (zbad) circle[radius=.68];
\draw[line width=.55pt] (zbad) circle[radius=.68];
\foreach \pt in {a,b,zbad}{
 \draw[line width=.9pt] ($(\pt)+(-.065,-.065)$)--($(\pt)+(.065,.065)$);
 \draw[line width=.9pt] ($(\pt)+(-.065,.065)$)--($(\pt)+(.065,-.065)$);
}
\node[above=4pt] at (a) {$\zeta_a$};
\node[above=4pt] at (b) {$\zeta_b$};
\node[font=\footnotesize] at (9.74,4.66) {$\zeta_a\in\Sigma_X$};
\node[font=\footnotesize] at (13.80,4.66) {partner also removed};
\node[above=4pt] at (zbad) {$z_\bullet$};
\draw[arr] (10.27,3.16)--(11.27,1.97);
\draw[arr] (13.27,3.16)--(12.27,1.97);
\node[fill=white,inner sep=1pt] at (10.72,2.56) {$\Psi$};
\node[fill=white,inner sep=1pt] at (12.81,2.56) {$\Psi$};
\node at (11.77,.18) {$X_B=\Psi^{-1}(D\setminus B)$};
\node[font=\footnotesize,align=center] at (11.77,-.61)
 {both inverse points over a deleted value are removed;\\
  branch values are deleted separately};
\end{tikzpicture}
}
\caption{The two continuation requirements in
Corollary~\ref{vol:cor:sector-continuation} and
Lemma~\ref{vol:lem:physical-exterior-continuation}. Left: a schematic
regular history $\Gamma$ in the logarithmic plane retains its loop
around an excluded point (cross) and the branch on each traversal. Only the terminal integral is
restarted on $D(p,r)$; the entire preceding integral remains in
$h_\Gamma$. This does not assert that different histories give the
same germ. Right: two local inverse sheets of the proper Jacobi map
$\Psi:X\to D$ over an ordinary value $z_\bullet\in B$ are shown. If
one preimage belongs to $\Sigma_X$, deleting its physical image
removes both preimages, whether or not the other belongs to
$\Sigma_X$. The regular source is the full inverse image of
$D\setminus B$; all branch values are removed as well. The marked
exceptional points are schematic, not a computation of the actual
exceptional sets. Degree two of $\Psi$ does not bound the number of
branches of the continued unknown.}
\label{vol:fig:history-fibers}
\end{figure}

Figure~\ref{vol:fig:history-fibers} summarizes why continuation
retains both the complete integration history and the complete
inverse image of the regular physical domain.

\Needspace{8\baselineskip}
\begin{remark}
For every rational caustic, real translations in its own Jacobi
parameter preserve the lower cylinder. Hence their physical images
remain in the spherical exterior $\{|z|>1\}\cup\{\infty\}$. To obtain
paths in the finite exterior, one must also avoid the finitely many
shifted Jacobi poles; the trace path constructions below do so.
Thus the exterior trace constructions can be carried out entirely
in the physical exterior. It does not itself prove trivial
local monodromy, remove the finite accumulation set, or classify the
isolated singularities. Reality provides the corresponding interior
continuation statement by reflection across the unit circle.
\end{remark}

\section{Multiple rational bases and removal of finite accumulation}
\label{cont:sec:multiple-bases}

The construction of Section~\ref{vol:sec:borel-volterra} leaves a
closed countable exceptional set for each rational base, with only finitely
many possible accumulation locations. We first extend each construction
through a collar of the unit circle on both geometric inverse sheets
(Lemmas~\ref{cont:lem:arbitrary-rational-base}--\ref{cont:lem:physical-collar-domain}).
Lemma~\ref{cont:lem:countable-gluing} then allows us to switch constructions
while preserving the germ already reached. After choosing two bases with
disjoint accumulation bounds in Lemma~\ref{cont:lem:disjoint-base-bounds},
we intersect their exceptional sets. The final collar separation excludes
accumulation anywhere in the finite plane. The remaining set is locally
finite, not necessarily finite, and monodromy is not eliminated here.

\begin{lemma}[Continuation from an arbitrary interior rational base]
\label{cont:lem:arbitrary-rational-base}
Fix a noncircular reference ellipse, and let $V$ be the physical germ
of a real analytic displacement whose rational traces are constant
throughout an open band $I\subset(0,1/2)$.
Let $\rho_*=p/d\in I$ be reduced. For $n\equiv s\pmod d$, let
$r_s\in\{0,\ldots,d-1\}$ be the representative of $-ps\pmod d$. For all sufficiently large $n$,
the rotations
\[
\rho_n=\frac{\lceil np/d\rceil}{n}
       =\rho_*+\frac{r_s}{dn}
\]
are permitted and have reduced denominator dividing $n$.
The exact construction of Section~\ref{vol:sec:borel-volterra} applies with this base.
In particular, it gives physical exterior continuation outside a
relatively closed countable set whose derived set is contained in the
finite set obtained by shifting the Jacobi coefficient points in the
lower cylinder through the $d$ action fibers and taking their finite
physical images.
\end{lemma}
\begin{proof}
Let $v(\varphi)$ denote the given real analytic displacement, so that its
physical germ satisfies $V(e^{i\varphi})=v(\varphi)$. Fix the reduced
anchor $\rho_*=p/d$ throughout the proof.

\smallskip\noindent
\emph{Step 1: distinguish the selected index from the primitive denominator.}
For each positive integer $n$, put $a_n=\lceil np/d\rceil$. Then
\[
 0\le da_n-np<d,\qquad da_n-np\equiv-np\pmod d.
\]
If $n\equiv s\pmod d$, uniqueness of the representative in
$\{0,\ldots,d-1\}$ gives $da_n-np=r_s$, and hence
\[
 \rho_n=\frac{a_n}{n}=\rho_*+\frac{r_s}{dn},\qquad
 0\le\rho_n-\rho_*<\frac1n.
\]
Since $I$ is open and contains $\rho_*$, choose a compact interval
$J\Subset I$ containing all these rotations for sufficiently large $n$.
Write $\rho_n=p_n/q_n$ in lowest terms. Its primitive denominator is
\[
 q_n=\frac{n}{\gcd(a_n,n)},\qquad q_n\mid n.
\]
We do not assume $q_n=n$. If $d\mid n$, the selected rotation is the
anchor itself; these repetitions cause no difficulty.

\smallskip\noindent
\emph{Step 2: turn trace constancy into exact vanishing of the selected rows.}
Put $g_n(y)=v(\Phi_{E,\rho_n}(y))$. The trace at the reduced rotation
$p_n/q_n$ is constant by hypothesis. Its Fourier coefficient at frequency
$k$ is
\[
 \widehat g_n(k)\sum_{j=0}^{q_n-1}e^{2\pi i k j p_n/q_n}
 =\begin{cases}
     q_n\widehat g_n(k),&q_n\mid k,\\
     0,&q_n\nmid k.
   \end{cases}
\]
Taking the nonzero frequency $k=n$ and using $q_n\mid n$ gives
$\widehat g_n(n)=0$.

Now express the displacement in the fixed anchor coordinate:
\[
 U(x)=v(\Phi_{E,\rho_*}(x)),\qquad
 H_n(x)=G_{E,\rho_n}(\Phi_{E,\rho_*}(x)).
\]
Here $G_{E,\rho_n}$ is the inverse action-coordinate map. We have
$g_n(H_n(x))=U(x)$. The map $H_n$ is an increasing degree-one lift, so
changing variables $y=H_n(x)$ over one period yields
\[
 0=\widehat g_n(n)
   =\frac1{2\pi}\int_0^{2\pi}
       U(x)e^{-inH_n(x)}H_n'(x)\,dx
   =(\widetilde A U)_n.
\]
The integrand in the $y$ variable is periodic, so a shifted integration
interval of length $2\pi$ gives the same coefficient. Thus the
Jacobian is exactly $H_n'$, and every sufficiently large selected
positive row vanishes. No nonlinear caustic estimate or later kernel
classification is used here.

\smallskip\noindent
\emph{Step 3: apply the construction at this base and retain its accumulation bound.}
Analyticity of $v$ and of the anchor coordinate makes $U$ a holomorphic
germ in $\zeta=e^{ix}$ on an annulus about the unit circle. All hypotheses
of the construction in Section~\ref{vol:sec:borel-volterra} therefore
hold with $p_0=p$ and $P=d$, after increasing its cutoff if necessary.
Proposition~\ref{band:prop:operator} supplies the required principal
invertibility for this arbitrary reduced anchor; it does not require
$p=1$. Lemma~\ref{vol:lem:convergent-symbol} and
Corollary~\ref{vol:cor:sector-continuation} supply action-coordinate
continuation, and Lemma~\ref{vol:lem:physical-exterior-continuation}
transfers it to the physical exterior.

For clarity, the finite containing set for accumulation is the one in
that last lemma. Put $m_0=m(\rho_*)$, $K=K(m_0)$, $K'=K(1-m_0)$,
$\mathfrak q=e^{-\pi K'/K}$, and $\omega=e^{2\pi i/d}$. In the action
annulus $A=\{1<|\zeta|<\mathfrak q^{-1}\}$, let $\Pi$ denote the two
Jacobi poles and let $\Psi:A\setminus\Pi\to\{|z|>1\}$ be the proper
Jacobi map from Lemma~\ref{vol:lem:physical-exterior-continuation}.
Then the exceptional set $B$ furnished there satisfies
\[
 B'\subset F_{\rho_*},\qquad
 F_{\rho_*}:=
 \Psi\!\left(
   \{\omega^{-l}i^j\mathfrak q^{-1/2}:0\le l<d,\ 0\le j<4\}
      \setminus\Pi\right).
\]
This is a finite set of finite physical images, with repetitions
removed. Rotating $\zeta$ by $\omega^{-l}$ translates the Jacobi
coordinate by $-4Kl/d$; coprimality of $p,d$ makes these the same
$d$ fibers as the base orbit translations, up to ordering and sign.
Thus it is precisely the accumulation bound described in the statement.
The set $B$ itself need not be finite, and no uniformity over all
rational anchors or restriction on monodromy is claimed.
\end{proof}

\begin{lemma}[A common continuation domain containing the unit circle]
\label{cont:lem:physical-collar-domain}
For each rational base in Lemma~\ref{cont:lem:arbitrary-rational-base}, its
physical continuation statement can be made on
$\Omega_*=\{|z|>r_*\}$ for some $0<r_*<1$.
The projected exceptional set is relatively closed and countable in
$\Omega_*$, has the same finite derived-set bound as before, and is
disjoint from a neighborhood of the unit circle.
Every branch reachable from the original collar can therefore be
continued through that collar.
\end{lemma}
\begin{proof}
Fix one reduced base $\rho_*=p/d$ from
Lemma~\ref{cont:lem:arbitrary-rational-base}. Put
\[
 m_0=m(\rho_*),\quad K=K(m_0),\quad K'=K(1-m_0),\quad
 L=\frac{\pi K'}K,\quad R_{\mathrm a}=e^L,\quad
 \omega=e^{2\pi i/d}.
\]
All Jacobi functions below have the fixed parameter $m_0$; $K'$ is the
complementary integral, not a derivative. Let
\[
 \mathcal T=\C/(4K\Z+4iK'\Z),\qquad
 \mathcal P(u)=\sn u-i\cn u,
 \qquad
 \chi(u)=\exp\!\left(\frac{i\pi(u-K)}{2K}\right).
\]
The map $\chi$ is a coordinate on a suitable real-period cylinder,
not on the whole torus. We will choose a cylinder that contains the
full preimage of the enlarged physical domain.

\smallskip\noindent
\emph{Step 1: identify both regular preimages of the unit circle.}
By Lemma~\ref{arith:jacobi}, the full torus preimage of $|z|=1$ is
\[
 \Gamma_0=\{\operatorname{Im}u=0\}/(4K\Z),\qquad
 \Gamma_1=\{\operatorname{Im}u=-2K'\}/(4K\Z).
\]
The intervening lower cylinder maps to the spherical exterior; the
upper cylinder maps to $|z|<1$. Neither boundary circle contains a
zero or pole of $\mathcal P$. The derivative identity
$\mathcal P'=i\dn\mathcal P$ gives, for $u=t+iy$ near these circles,
\[
 \partial_y\log|\mathcal P(t+iy)|=-\operatorname{Re}\dn(t+iy).
\]
For real $t$, $\dn t>0$, whereas $\dn(t-2iK')=-\dn t$.
Thus this normal derivative is strictly negative on $\Gamma_0$ and
strictly positive on $\Gamma_1$, uniformly over one real period.
In particular both components are regular, and nearby modulus level
sets form collars on their upper-cylinder sides.

\smallskip\noindent
\emph{Step 2: verify principal invertibility on the second component.}
Use the coefficient set and principal matrix of
Corollary~\ref{vol:cor:sector-continuation}, now with $P=d$:
\[
 \mathfrak E=\mathcal S\cup\mathcal Z,\qquad
 \widehat M_{sl}(\zeta)=\frac1d\omega^{-sl}a_{s,0}(\omega^l\zeta),
 \qquad a_{s,0}(\zeta)=e^{-ir_sY(\zeta)/d}.
\]
That corollary gives $\mathfrak E'\subset\mathcal S$, with
$\mathfrak E$ closed and countable and $\mathcal S$ closed and discrete
in $\C^*$. It also gives $\det\widehat M(e^{ix})\ne0$ for real $x$.

To treat $\Gamma_1$, keep the base fixed and use
\eqref{vol:eq:marked-first-derivative} together with
$Y=m'(\rho_*)D_0$. In the notation $u=K+2Kx/\pi$, this reads
\[
 Y(x)=C_{\mathrm b}
   \left[Z(u)-m_0\frac{\sn u\cn u}{\dn u}\right],\qquad
 C_{\mathrm b}=\frac{m'(\rho_*)\pi}{4Km_0(1-m_0)}.
\]
Here $Z$ is the single-valued meromorphic Jacobi zeta function at
parameter $m_0$. The translation identities
\cite[Table~22.4.3]{DLMF} give
\[
 \sn(u-2iK')=\sn u,\qquad
 \cn(u-2iK')=-\cn u,\qquad
 \dn(u-2iK')=-\dn u.
\]
Hence $\sn u\cn u/\dn u$ is unchanged by this translation. Also
$Z'(u)=\dn^2u-E(m_0)/K$, so
\[
 \frac{d}{du}\big[Z(u-2iK')-Z(u)\big]=0.
\]
The two terms are already single-valued meromorphic functions, by
\cite[Eq.~22.16.30]{DLMF}; their difference therefore is a constant
$\kappa_0$. Its possible singularities are removable because its
meromorphic derivative is zero. Consequently
\[
 Y(x-iL)=Y(x)+c_*,\qquad c_*=C_{\mathrm b}\kappa_0.
\]
This is evaluation at two spatial points of the fixed-base function
$Y$; neither the shift nor $K,K'$ is being parameter-differentiated.

The shift $x\mapsto x-iL$ sends $\zeta$ to $R_{\mathrm a}\zeta$.
Applying it at every fiber gives the matrix identity
\[
 \widehat M(R_{\mathrm a}\zeta)=\mathsf D_*\widehat M(\zeta),\qquad
 \mathsf D_*=
 \operatorname{diag}_{0\le s<d}\big(e^{-ir_sc_*/d}\big),
\]
where both sides are interpreted away from their coefficient sets.
Every row factor is nonzero. Since the $r_s$ permute
$0,\ldots,d-1$, $\det\mathsf D_*=e^{-ic_*(d-1)/2}\ne0$.
It follows that the principal determinant is nonzero on
$|\zeta|=R_{\mathrm a}$ as well as on $|\zeta|=1$.

Both circles avoid $\mathcal S$. Indeed its Jacobi representatives
have imaginary parts $(2\ell+1)K'$, while rotation of the action
coordinate changes only their real parts. Compactness gives regular
neighborhoods of both $\Gamma_0$ and $\Gamma_1$ on which the
coefficient set and determinant zeros are absent. Checking both
components is necessary: a continued germ can return on either sheet.

\smallskip\noindent
\emph{Step 3: construct the enlarged annular preimage.}
Choose $0<\varepsilon<K'$ small enough that the normal derivative
signs from Step~1 persist in the two boundary collars, and that
\[
 r_{\mathrm a}:=e^{-\pi\varepsilon/(2K)}
\]
lies inside the original action annulus and the open neighborhood of
admissible source radii in Proposition~\ref{vol:prop:exact-volterra}.
These choices are possible because $r_{\mathrm a}\to1$ as
$\varepsilon\to0$.

The upper-cylinder core
$\{\varepsilon\le\operatorname{Im}u\le2K'-\varepsilon\}/(4K\Z)$
is compact and has $|\mathcal P|<1$. Let $\nu_\varepsilon<1$ be its
maximum. The two inner branch values have modulus
\[
 \beta_{\mathrm{in}}=
 \frac{1-\sqrt{1-m_0}}{\sqrt{m_0}}
 =\frac{\sqrt{m_0}}{1+\sqrt{1-m_0}}<1.
\]
They are reciprocal to the exterior branch values of
Lemma~\ref{vol:lem:physical-exterior-continuation}, by the reciprocal
Jacobi translation identity. Choose $r_*<1$ sufficiently close to one
that $r_*>\max\{\nu_\varepsilon,\beta_{\mathrm{in}}\}$ and the level
$|\mathcal P|=r_*$ lies in the two transverse collars of Step~1.

More explicitly, the implicit function theorem and the uniform normal
derivative signs give periodic boundary graphs
$y_+(t)\in(0,\varepsilon)$ and
$y_-(t)\in(-2K'-\varepsilon,-2K')$ for this level. In the two collars
the superlevel set is on the lower-cylinder side of these graphs;
in the upper core it is empty. Thus the full preimage
\[
 Y_*:=\mathcal P^{-1}
          \big(\{|z|>r_*\}\cup\{\infty\}\big)
\]
has the annular representative $y_-(t)<y<y_+(t)$. In particular it
contains the whole lower cylinder and is contained in
\[
 \mathcal N_\varepsilon=
 \{-2K'-\varepsilon<\operatorname{Im}u<\varepsilon\}/(4K\Z).
\]
The height of this cylinder is less than $4K'$, so its projection to
$\mathcal T$ is injective. We henceforth use that representative of
$Y_*$. Its action-coordinate image $\mathscr A_*:=\chi(Y_*)$ is
annular and satisfies
\[
 r_{\mathrm a}<|\zeta|
       <R_{\mathrm a}e^{\pi\varepsilon/(2K)}.
\]
No other component of the torus preimage has been omitted.

\smallskip\noindent
\emph{Step 4: obtain a proper finite-target map and one source circle.}
Let $\Pi=\{-i e^{\pi K'/(2K)},\ i e^{\pi K'/(2K)}\}$ be the two
pole points in action coordinates and put
\[
 X_*:=\mathscr A_*\setminus\Pi,\qquad
 \Omega_*:=\{z\in\C:|z|>r_*\},\qquad
 \Psi_*:=\mathcal P\circ\chi^{-1}:X_*\longrightarrow\Omega_*.
\]
This map is holomorphic, onto, and proper of degree two. For any compact
$Q\subset\Omega_*$, its full torus preimage is compact, lies in $Y_*$,
and contains neither pole, since the target values are finite. Its
image under $\chi$ is therefore compact in $X_*$. This proves the
finite-target properness without treating infinity as an ordinary point.
By the choice of $r_*$, the only branch values in $\Omega_*$ are the
same exterior pair
\[
 C_{\mathrm{br}}=
 \left\{\frac{1+\sqrt{1-m_0}}{\sqrt{m_0}},
          -\frac{1+\sqrt{1-m_0}}{\sqrt{m_0}}\right\}.
\]

Choose $R_c$ slightly below $r_{\mathrm a}$ while remaining in the
same open interval of admissible source radii. Then
\[
 R_c<r_{\mathrm a}<1,
\]
and its circle lies in the known action annulus, with the required
principal invertibility and symbol bounds. Proposition~\ref{vol:prop:exact-volterra}
does not require $R_c>1$. Every point of $X_*$ has $|\zeta|>R_c$,
so this single source circle works for all lifted paths in the
enlarged domain. The sets $\mathcal S,\mathcal Z$ are defined in
$\zeta$ and do not change with this choice of $R_c$.

\smallskip\noindent
\emph{Step 5: project exceptions without adding accumulation locations.}
Set
\[
 \Sigma_*:=\mathfrak E\cap X_*,\qquad
 B_*:=\Psi_*(\Sigma_*)\cup C_{\mathrm{br}}.
\]
The set $\Sigma_*$ is closed and countable relative to $X_*$. If
$z_n\in\Psi_*(\Sigma_*)$ converges to a point $z\in\Omega_*$, choose
preimages $\zeta_n\in\Sigma_*$. Properness over the compact set
$\{z\}\cup\{z_n\}$ gives a subsequence converging in $X_*$. Its
limit lies in $\Sigma_*$ by relative closedness and maps to $z$.
Thus $B_*$ is relatively closed and countable in $\Omega_*$.

For a convergent sequence of distinct points of $B_*$, discard the
finitely many terms in $C_{\mathrm{br}}$ and choose preimages in
$\Sigma_*$. They are distinct because their images are. The same
properness argument gives a limit in the derived set of $\Sigma_*$,
taken in $X_*$. Since $\mathfrak E'\subset\mathcal S$, this proves
\[
 B_*'\subset\Psi_*(\Sigma_*')
           \subset\Psi_*(\mathcal S\cap X_*),
\]
where $B_*'$ is taken relative to the finite domain $\Omega_*$.

To identify this containing set, put $\mathfrak q=e^{-\pi K'/K}$ and
$A=\{1<|\zeta|<\mathfrak q^{-1}\}$. The set $S$ in
Lemma~\ref{vol:lem:convergent-symbol} has representatives
$u=jK+(2\ell+1)iK'$. Its rotated copies have the same imaginary
heights. Neither added collar contains such a point, because
$\varepsilon<K'$. Therefore
\[
 \mathcal S\cap X_*
   =\mathcal S\cap(A\setminus\Pi)
   =\{\omega^{-l}i^j\mathfrak q^{-1/2}:
          0\le l<d,\ 0\le j<4\}\setminus\Pi.
\]
On this set $\Psi_*$ agrees with the original map $\Psi$ from
Lemma~\ref{vol:lem:physical-exterior-continuation}. Hence
\[
 B_*'\subset F_{\rho_*}:=
    \Psi\big(\mathcal S\cap(A\setminus\Pi)\big),
\]
the same finite containing set as in the preceding lemma. The
exceptional set may contain additional isolated determinant zeros;
equality of the actual derived sets is not asserted. The exterior
branch values are already images of coefficient points, since zeros
of $\dn$ are included in $S$, but their explicit addition is harmless.

\smallskip\noindent
\emph{Step 6: make the full preimage of a physical collar exception-free.}
By Step~2 choose torus neighborhoods $\mathcal U_0,\mathcal U_1$ of
$\Gamma_0,\Gamma_1$, with closures in the chart containing $Y_*$,
small enough that on them the action coefficient set and determinant
zeros are absent and $\mathcal P$ has no critical point. They may be
chosen inside $Y_*$ since the two circles map to modulus one.

The full preimage of $|z|=1$ is contained in
$\mathcal U_0\cup\mathcal U_1$. Compactness of $\mathcal T$ gives
$\delta>0$ such that the collar
\[
 \mathcal C=\{1-\delta<|z|<1+\delta\}\subset\Omega_*
\]
has its entire torus preimage in this union. Indeed, failure would give
points $z_n$ with $|z_n|\to1$ and preimages outside the union. A
convergent torus subsequence would limit to a preimage of the unit
circle outside that open union, a contradiction. Thus no point of
$\mathcal C$ is a projected action exception or a branch value, and
\[
 B_*\cap\mathcal C=\varnothing.
\]
This uses regularity on both source components, not merely regularity
of the initial physical branch. Shrink $\delta$ further so that the
collar lies in the original holomorphy annulus of $V$ and its inverse
sheet adjacent to $\Gamma_0$ lies in the known annulus of $U$. On
that sheet the original germs satisfy $V=U\circ\psi_0$.

\smallskip\noindent
\emph{Step 7: continue along enlarged-domain paths and through the collar.}
Restrict $\Psi_*$ to the full inverse image of $\Omega_*\setminus B_*$.
As in Lemma~\ref{vol:lem:physical-exterior-continuation}, properness and
the exclusion of all branch values make it an unbranched two-sheeted
covering. Every continuous physical path with compact parameter
interval in this domain lifts, starting on the original inverse
sheet, to a path in $X_*$ avoiding $\mathfrak E$. Its action radius
is larger than the fixed $R_c$ by Step~4. Corollary~\ref{vol:cor:sector-continuation}
therefore continues $U$ along the whole lift. Composing consecutive
action germs with the local inverse branches of $\Psi_*$ gives
consecutive physical germs, with agreement on their overlaps. This
proves the enlarged all-path continuation statement. No rectifiability
of the given path or prior knowledge of $U$ near the second source
component is required. Other starting points in an originally chosen
collar are handled by a connector in the known physical annulus,
avoiding the closed countable set, exactly as in the final step of
Lemma~\ref{vol:lem:physical-exterior-continuation}.

Finally, let a germ at $z_1\in\mathcal C$ be reachable from the original
collar along a path $\Gamma$ in $\Omega_*\setminus B_*$. For any
further path $\gamma$ in $\mathcal C$ starting at $z_1$, the
concatenation $\Gamma*\gamma$ remains in the regular domain. All-path
continuation along it starts on $\gamma$ with the already reached
germ, by uniqueness along $\Gamma$. Thus every germ reachable in the
stated regular domain can be continued along every further collar
path, on either Jacobi sheet. The history is not reset to the
original branch. This asserts neither a single-valued function on
the whole collar nor trivial or finite monodromy. No additional
negative-row hypothesis or reflection argument is used, and no
uniform collar over all rational bases is claimed.
\end{proof}

\begin{figure}[tbp]
\centering
\begin{tikzpicture}[x=1cm,y=1cm,font=\small,
 arr/.style={-{Stealth[length=1.8mm]},line width=.7pt},
 every node/.style={inner sep=2pt}]
\path[use as bounding box] (-.1,-.40) rectangle (15.75,5.42);
\node at (3.75,5.15) {(a) Both geometric preimage components};
\node at (11.84,5.15) {(b) A collar of the physical unit circle};
\fill[black!3] plot coordinates {(0.650000,1.073029) (0.727500,1.081725) (0.805000,1.089859) (0.882500,1.097289) (0.960000,1.103904) (1.037500,1.109616) (1.115000,1.114361) (1.192500,1.118094) (1.270000,1.120780) (1.347500,1.122400) (1.425000,1.122942) (1.502500,1.122400) (1.580000,1.120780) (1.657500,1.118094) (1.735000,1.114361) (1.812500,1.109616) (1.890000,1.103904) (1.967500,1.097289) (2.045000,1.089859) (2.122500,1.081725) (2.200000,1.073029) (2.277500,1.063947) (2.355000,1.054691) (2.432500,1.045507) (2.510000,1.036669) (2.587500,1.028476) (2.665000,1.021232) (2.742500,1.015231) (2.820000,1.010733) (2.897500,1.007946) (2.975000,1.007001) (3.052500,1.007946) (3.130000,1.010733) (3.207500,1.015231) (3.285000,1.021232) (3.362500,1.028476) (3.440000,1.036669) (3.517500,1.045507) (3.595000,1.054691) (3.672500,1.063947) (3.750000,1.073029) (3.827500,1.081725) (3.905000,1.089859) (3.982500,1.097289) (4.060000,1.103904) (4.137500,1.109616) (4.215000,1.114361) (4.292500,1.118094) (4.370000,1.120780) (4.447500,1.122400) (4.525000,1.122942) (4.602500,1.122400) (4.680000,1.120780) (4.757500,1.118094) (4.835000,1.114361) (4.912500,1.109616) (4.990000,1.103904) (5.067500,1.097289) (5.145000,1.089859) (5.222500,1.081725) (5.300000,1.073029) (5.377500,1.063947) (5.455000,1.054691) (5.532500,1.045507) (5.610000,1.036669) (5.687500,1.028476) (5.765000,1.021232) (5.842500,1.015231) (5.920000,1.010733) (5.997500,1.007946) (6.075000,1.007001) (6.152500,1.007946) (6.230000,1.010733) (6.307500,1.015231) (6.385000,1.021232) (6.462500,1.028476) (6.540000,1.036669) (6.617500,1.045507) (6.695000,1.054691) (6.772500,1.063947) (6.850000,1.073029) (6.850000,4.466971) (6.772500,4.476053) (6.695000,4.485309) (6.617500,4.494493) (6.540000,4.503331) (6.462500,4.511524) (6.385000,4.518768) (6.307500,4.524769) (6.230000,4.529267) (6.152500,4.532054) (6.075000,4.532999) (5.997500,4.532054) (5.920000,4.529267) (5.842500,4.524769) (5.765000,4.518768) (5.687500,4.511524) (5.610000,4.503331) (5.532500,4.494493) (5.455000,4.485309) (5.377500,4.476053) (5.300000,4.466971) (5.222500,4.458275) (5.145000,4.450141) (5.067500,4.442711) (4.990000,4.436096) (4.912500,4.430384) (4.835000,4.425639) (4.757500,4.421906) (4.680000,4.419220) (4.602500,4.417600) (4.525000,4.417058) (4.447500,4.417600) (4.370000,4.419220) (4.292500,4.421906) (4.215000,4.425639) (4.137500,4.430384) (4.060000,4.436096) (3.982500,4.442711) (3.905000,4.450141) (3.827500,4.458275) (3.750000,4.466971) (3.672500,4.476053) (3.595000,4.485309) (3.517500,4.494493) (3.440000,4.503331) (3.362500,4.511524) (3.285000,4.518768) (3.207500,4.524769) (3.130000,4.529267) (3.052500,4.532054) (2.975000,4.532999) (2.897500,4.532054) (2.820000,4.529267) (2.742500,4.524769) (2.665000,4.518768) (2.587500,4.511524) (2.510000,4.503331) (2.432500,4.494493) (2.355000,4.485309) (2.277500,4.476053) (2.200000,4.466971) (2.122500,4.458275) (2.045000,4.450141) (1.967500,4.442711) (1.890000,4.436096) (1.812500,4.430384) (1.735000,4.425639) (1.657500,4.421906) (1.580000,4.419220) (1.502500,4.417600) (1.425000,4.417058) (1.347500,4.417600) (1.270000,4.419220) (1.192500,4.421906) (1.115000,4.425639) (1.037500,4.430384) (0.960000,4.436096) (0.882500,4.442711) (0.805000,4.450141) (0.727500,4.458275) (0.650000,4.466971)} -- cycle;
\fill[black!16] plot coordinates {(0.650000,4.126248) (0.727500,4.131115) (0.805000,4.135677) (0.882500,4.139854) (0.960000,4.143578) (1.037500,4.146798) (1.115000,4.149477) (1.192500,4.151587) (1.270000,4.153106) (1.347500,4.154022) (1.425000,4.154329) (1.502500,4.154022) (1.580000,4.153106) (1.657500,4.151587) (1.735000,4.149477) (1.812500,4.146798) (1.890000,4.143578) (1.967500,4.139854) (2.045000,4.135677) (2.122500,4.131115) (2.200000,4.126248) (2.277500,4.121178) (2.355000,4.116025) (2.432500,4.110925) (2.510000,4.106033) (2.587500,4.101509) (2.665000,4.097521) (2.742500,4.094224) (2.820000,4.091757) (2.897500,4.090231) (2.975000,4.089714) (3.052500,4.090231) (3.130000,4.091757) (3.207500,4.094224) (3.285000,4.097521) (3.362500,4.101509) (3.440000,4.106033) (3.517500,4.110925) (3.595000,4.116025) (3.672500,4.121178) (3.750000,4.126248) (3.827500,4.131115) (3.905000,4.135677) (3.982500,4.139854) (4.060000,4.143578) (4.137500,4.146798) (4.215000,4.149477) (4.292500,4.151587) (4.370000,4.153106) (4.447500,4.154022) (4.525000,4.154329) (4.602500,4.154022) (4.680000,4.153106) (4.757500,4.151587) (4.835000,4.149477) (4.912500,4.146798) (4.990000,4.143578) (5.067500,4.139854) (5.145000,4.135677) (5.222500,4.131115) (5.300000,4.126248) (5.377500,4.121178) (5.455000,4.116025) (5.532500,4.110925) (5.610000,4.106033) (5.687500,4.101509) (5.765000,4.097521) (5.842500,4.094224) (5.920000,4.091757) (5.997500,4.090231) (6.075000,4.089714) (6.152500,4.090231) (6.230000,4.091757) (6.307500,4.094224) (6.385000,4.097521) (6.462500,4.101509) (6.540000,4.106033) (6.617500,4.110925) (6.695000,4.116025) (6.772500,4.121178) (6.850000,4.126248) (6.850000,4.389613) (6.772500,4.395349) (6.695000,4.401181) (6.617500,4.406954) (6.540000,4.412496) (6.462500,4.417620) (6.385000,4.422141) (6.307500,4.425879) (6.230000,4.428676) (6.152500,4.430407) (6.075000,4.430993) (5.997500,4.430407) (5.920000,4.428676) (5.842500,4.425879) (5.765000,4.422141) (5.687500,4.417620) (5.610000,4.412496) (5.532500,4.406954) (5.455000,4.401181) (5.377500,4.395349) (5.300000,4.389613) (5.222500,4.384109) (5.145000,4.378950) (5.067500,4.374229) (4.990000,4.370021) (4.912500,4.366382) (4.835000,4.363355) (4.757500,4.360973) (4.680000,4.359257) (4.602500,4.358222) (4.525000,4.357876) (4.447500,4.358222) (4.370000,4.359257) (4.292500,4.360973) (4.215000,4.363355) (4.137500,4.366382) (4.060000,4.370021) (3.982500,4.374229) (3.905000,4.378950) (3.827500,4.384109) (3.750000,4.389613) (3.672500,4.395349) (3.595000,4.401181) (3.517500,4.406954) (3.440000,4.412496) (3.362500,4.417620) (3.285000,4.422141) (3.207500,4.425879) (3.130000,4.428676) (3.052500,4.430407) (2.975000,4.430993) (2.897500,4.430407) (2.820000,4.428676) (2.742500,4.425879) (2.665000,4.422141) (2.587500,4.417620) (2.510000,4.412496) (2.432500,4.406954) (2.355000,4.401181) (2.277500,4.395349) (2.200000,4.389613) (2.122500,4.384109) (2.045000,4.378950) (1.967500,4.374229) (1.890000,4.370021) (1.812500,4.366382) (1.735000,4.363355) (1.657500,4.360973) (1.580000,4.359257) (1.502500,4.358222) (1.425000,4.357876) (1.347500,4.358222) (1.270000,4.359257) (1.192500,4.360973) (1.115000,4.363355) (1.037500,4.366382) (0.960000,4.370021) (0.882500,4.374229) (0.805000,4.378950) (0.727500,4.384109) (0.650000,4.389613)} -- cycle;
\fill[black!16] plot coordinates {(0.650000,1.150387) (0.727500,1.155891) (0.805000,1.161050) (0.882500,1.165771) (0.960000,1.169979) (1.037500,1.173618) (1.115000,1.176645) (1.192500,1.179027) (1.270000,1.180743) (1.347500,1.181778) (1.425000,1.182124) (1.502500,1.181778) (1.580000,1.180743) (1.657500,1.179027) (1.735000,1.176645) (1.812500,1.173618) (1.890000,1.169979) (1.967500,1.165771) (2.045000,1.161050) (2.122500,1.155891) (2.200000,1.150387) (2.277500,1.144651) (2.355000,1.138819) (2.432500,1.133046) (2.510000,1.127504) (2.587500,1.122380) (2.665000,1.117859) (2.742500,1.114121) (2.820000,1.111324) (2.897500,1.109593) (2.975000,1.109007) (3.052500,1.109593) (3.130000,1.111324) (3.207500,1.114121) (3.285000,1.117859) (3.362500,1.122380) (3.440000,1.127504) (3.517500,1.133046) (3.595000,1.138819) (3.672500,1.144651) (3.750000,1.150387) (3.827500,1.155891) (3.905000,1.161050) (3.982500,1.165771) (4.060000,1.169979) (4.137500,1.173618) (4.215000,1.176645) (4.292500,1.179027) (4.370000,1.180743) (4.447500,1.181778) (4.525000,1.182124) (4.602500,1.181778) (4.680000,1.180743) (4.757500,1.179027) (4.835000,1.176645) (4.912500,1.173618) (4.990000,1.169979) (5.067500,1.165771) (5.145000,1.161050) (5.222500,1.155891) (5.300000,1.150387) (5.377500,1.144651) (5.455000,1.138819) (5.532500,1.133046) (5.610000,1.127504) (5.687500,1.122380) (5.765000,1.117859) (5.842500,1.114121) (5.920000,1.111324) (5.997500,1.109593) (6.075000,1.109007) (6.152500,1.109593) (6.230000,1.111324) (6.307500,1.114121) (6.385000,1.117859) (6.462500,1.122380) (6.540000,1.127504) (6.617500,1.133046) (6.695000,1.138819) (6.772500,1.144651) (6.850000,1.150387) (6.850000,1.413752) (6.772500,1.418822) (6.695000,1.423975) (6.617500,1.429075) (6.540000,1.433967) (6.462500,1.438491) (6.385000,1.442479) (6.307500,1.445776) (6.230000,1.448243) (6.152500,1.449769) (6.075000,1.450286) (5.997500,1.449769) (5.920000,1.448243) (5.842500,1.445776) (5.765000,1.442479) (5.687500,1.438491) (5.610000,1.433967) (5.532500,1.429075) (5.455000,1.423975) (5.377500,1.418822) (5.300000,1.413752) (5.222500,1.408885) (5.145000,1.404323) (5.067500,1.400146) (4.990000,1.396422) (4.912500,1.393202) (4.835000,1.390523) (4.757500,1.388413) (4.680000,1.386894) (4.602500,1.385978) (4.525000,1.385671) (4.447500,1.385978) (4.370000,1.386894) (4.292500,1.388413) (4.215000,1.390523) (4.137500,1.393202) (4.060000,1.396422) (3.982500,1.400146) (3.905000,1.404323) (3.827500,1.408885) (3.750000,1.413752) (3.672500,1.418822) (3.595000,1.423975) (3.517500,1.429075) (3.440000,1.433967) (3.362500,1.438491) (3.285000,1.442479) (3.207500,1.445776) (3.130000,1.448243) (3.052500,1.449769) (2.975000,1.450286) (2.897500,1.449769) (2.820000,1.448243) (2.742500,1.445776) (2.665000,1.442479) (2.587500,1.438491) (2.510000,1.433967) (2.432500,1.429075) (2.355000,1.423975) (2.277500,1.418822) (2.200000,1.413752) (2.122500,1.408885) (2.045000,1.404323) (1.967500,1.400146) (1.890000,1.396422) (1.812500,1.393202) (1.735000,1.390523) (1.657500,1.388413) (1.580000,1.386894) (1.502500,1.385978) (1.425000,1.385671) (1.347500,1.385978) (1.270000,1.386894) (1.192500,1.388413) (1.115000,1.390523) (1.037500,1.393202) (0.960000,1.396422) (0.882500,1.400146) (0.805000,1.404323) (0.727500,1.408885) (0.650000,1.413752)} -- cycle;
\draw[densely dashed,line width=.65pt] plot coordinates {(0.650000,4.466971) (0.727500,4.458275) (0.805000,4.450141) (0.882500,4.442711) (0.960000,4.436096) (1.037500,4.430384) (1.115000,4.425639) (1.192500,4.421906) (1.270000,4.419220) (1.347500,4.417600) (1.425000,4.417058) (1.502500,4.417600) (1.580000,4.419220) (1.657500,4.421906) (1.735000,4.425639) (1.812500,4.430384) (1.890000,4.436096) (1.967500,4.442711) (2.045000,4.450141) (2.122500,4.458275) (2.200000,4.466971) (2.277500,4.476053) (2.355000,4.485309) (2.432500,4.494493) (2.510000,4.503331) (2.587500,4.511524) (2.665000,4.518768) (2.742500,4.524769) (2.820000,4.529267) (2.897500,4.532054) (2.975000,4.532999) (3.052500,4.532054) (3.130000,4.529267) (3.207500,4.524769) (3.285000,4.518768) (3.362500,4.511524) (3.440000,4.503331) (3.517500,4.494493) (3.595000,4.485309) (3.672500,4.476053) (3.750000,4.466971) (3.827500,4.458275) (3.905000,4.450141) (3.982500,4.442711) (4.060000,4.436096) (4.137500,4.430384) (4.215000,4.425639) (4.292500,4.421906) (4.370000,4.419220) (4.447500,4.417600) (4.525000,4.417058) (4.602500,4.417600) (4.680000,4.419220) (4.757500,4.421906) (4.835000,4.425639) (4.912500,4.430384) (4.990000,4.436096) (5.067500,4.442711) (5.145000,4.450141) (5.222500,4.458275) (5.300000,4.466971) (5.377500,4.476053) (5.455000,4.485309) (5.532500,4.494493) (5.610000,4.503331) (5.687500,4.511524) (5.765000,4.518768) (5.842500,4.524769) (5.920000,4.529267) (5.997500,4.532054) (6.075000,4.532999) (6.152500,4.532054) (6.230000,4.529267) (6.307500,4.524769) (6.385000,4.518768) (6.462500,4.511524) (6.540000,4.503331) (6.617500,4.494493) (6.695000,4.485309) (6.772500,4.476053) (6.850000,4.466971)};
\draw[densely dashed,line width=.65pt] plot coordinates {(0.650000,1.073029) (0.727500,1.081725) (0.805000,1.089859) (0.882500,1.097289) (0.960000,1.103904) (1.037500,1.109616) (1.115000,1.114361) (1.192500,1.118094) (1.270000,1.120780) (1.347500,1.122400) (1.425000,1.122942) (1.502500,1.122400) (1.580000,1.120780) (1.657500,1.118094) (1.735000,1.114361) (1.812500,1.109616) (1.890000,1.103904) (1.967500,1.097289) (2.045000,1.089859) (2.122500,1.081725) (2.200000,1.073029) (2.277500,1.063947) (2.355000,1.054691) (2.432500,1.045507) (2.510000,1.036669) (2.587500,1.028476) (2.665000,1.021232) (2.742500,1.015231) (2.820000,1.010733) (2.897500,1.007946) (2.975000,1.007001) (3.052500,1.007946) (3.130000,1.010733) (3.207500,1.015231) (3.285000,1.021232) (3.362500,1.028476) (3.440000,1.036669) (3.517500,1.045507) (3.595000,1.054691) (3.672500,1.063947) (3.750000,1.073029) (3.827500,1.081725) (3.905000,1.089859) (3.982500,1.097289) (4.060000,1.103904) (4.137500,1.109616) (4.215000,1.114361) (4.292500,1.118094) (4.370000,1.120780) (4.447500,1.122400) (4.525000,1.122942) (4.602500,1.122400) (4.680000,1.120780) (4.757500,1.118094) (4.835000,1.114361) (4.912500,1.109616) (4.990000,1.103904) (5.067500,1.097289) (5.145000,1.089859) (5.222500,1.081725) (5.300000,1.073029) (5.377500,1.063947) (5.455000,1.054691) (5.532500,1.045507) (5.610000,1.036669) (5.687500,1.028476) (5.765000,1.021232) (5.842500,1.015231) (5.920000,1.010733) (5.997500,1.007946) (6.075000,1.007001) (6.152500,1.007946) (6.230000,1.010733) (6.307500,1.015231) (6.385000,1.021232) (6.462500,1.028476) (6.540000,1.036669) (6.617500,1.045507) (6.695000,1.054691) (6.772500,1.063947) (6.850000,1.073029)};
\draw[black!40,line width=.25pt] plot coordinates {(0.650000,4.389613) (0.727500,4.384109) (0.805000,4.378950) (0.882500,4.374229) (0.960000,4.370021) (1.037500,4.366382) (1.115000,4.363355) (1.192500,4.360973) (1.270000,4.359257) (1.347500,4.358222) (1.425000,4.357876) (1.502500,4.358222) (1.580000,4.359257) (1.657500,4.360973) (1.735000,4.363355) (1.812500,4.366382) (1.890000,4.370021) (1.967500,4.374229) (2.045000,4.378950) (2.122500,4.384109) (2.200000,4.389613) (2.277500,4.395349) (2.355000,4.401181) (2.432500,4.406954) (2.510000,4.412496) (2.587500,4.417620) (2.665000,4.422141) (2.742500,4.425879) (2.820000,4.428676) (2.897500,4.430407) (2.975000,4.430993) (3.052500,4.430407) (3.130000,4.428676) (3.207500,4.425879) (3.285000,4.422141) (3.362500,4.417620) (3.440000,4.412496) (3.517500,4.406954) (3.595000,4.401181) (3.672500,4.395349) (3.750000,4.389613) (3.827500,4.384109) (3.905000,4.378950) (3.982500,4.374229) (4.060000,4.370021) (4.137500,4.366382) (4.215000,4.363355) (4.292500,4.360973) (4.370000,4.359257) (4.447500,4.358222) (4.525000,4.357876) (4.602500,4.358222) (4.680000,4.359257) (4.757500,4.360973) (4.835000,4.363355) (4.912500,4.366382) (4.990000,4.370021) (5.067500,4.374229) (5.145000,4.378950) (5.222500,4.384109) (5.300000,4.389613) (5.377500,4.395349) (5.455000,4.401181) (5.532500,4.406954) (5.610000,4.412496) (5.687500,4.417620) (5.765000,4.422141) (5.842500,4.425879) (5.920000,4.428676) (5.997500,4.430407) (6.075000,4.430993) (6.152500,4.430407) (6.230000,4.428676) (6.307500,4.425879) (6.385000,4.422141) (6.462500,4.417620) (6.540000,4.412496) (6.617500,4.406954) (6.695000,4.401181) (6.772500,4.395349) (6.850000,4.389613)};
\draw[black!40,line width=.25pt] plot coordinates {(0.650000,4.126248) (0.727500,4.131115) (0.805000,4.135677) (0.882500,4.139854) (0.960000,4.143578) (1.037500,4.146798) (1.115000,4.149477) (1.192500,4.151587) (1.270000,4.153106) (1.347500,4.154022) (1.425000,4.154329) (1.502500,4.154022) (1.580000,4.153106) (1.657500,4.151587) (1.735000,4.149477) (1.812500,4.146798) (1.890000,4.143578) (1.967500,4.139854) (2.045000,4.135677) (2.122500,4.131115) (2.200000,4.126248) (2.277500,4.121178) (2.355000,4.116025) (2.432500,4.110925) (2.510000,4.106033) (2.587500,4.101509) (2.665000,4.097521) (2.742500,4.094224) (2.820000,4.091757) (2.897500,4.090231) (2.975000,4.089714) (3.052500,4.090231) (3.130000,4.091757) (3.207500,4.094224) (3.285000,4.097521) (3.362500,4.101509) (3.440000,4.106033) (3.517500,4.110925) (3.595000,4.116025) (3.672500,4.121178) (3.750000,4.126248) (3.827500,4.131115) (3.905000,4.135677) (3.982500,4.139854) (4.060000,4.143578) (4.137500,4.146798) (4.215000,4.149477) (4.292500,4.151587) (4.370000,4.153106) (4.447500,4.154022) (4.525000,4.154329) (4.602500,4.154022) (4.680000,4.153106) (4.757500,4.151587) (4.835000,4.149477) (4.912500,4.146798) (4.990000,4.143578) (5.067500,4.139854) (5.145000,4.135677) (5.222500,4.131115) (5.300000,4.126248) (5.377500,4.121178) (5.455000,4.116025) (5.532500,4.110925) (5.610000,4.106033) (5.687500,4.101509) (5.765000,4.097521) (5.842500,4.094224) (5.920000,4.091757) (5.997500,4.090231) (6.075000,4.089714) (6.152500,4.090231) (6.230000,4.091757) (6.307500,4.094224) (6.385000,4.097521) (6.462500,4.101509) (6.540000,4.106033) (6.617500,4.110925) (6.695000,4.116025) (6.772500,4.121178) (6.850000,4.126248)};
\draw[black!40,line width=.25pt] plot coordinates {(0.650000,1.150387) (0.727500,1.155891) (0.805000,1.161050) (0.882500,1.165771) (0.960000,1.169979) (1.037500,1.173618) (1.115000,1.176645) (1.192500,1.179027) (1.270000,1.180743) (1.347500,1.181778) (1.425000,1.182124) (1.502500,1.181778) (1.580000,1.180743) (1.657500,1.179027) (1.735000,1.176645) (1.812500,1.173618) (1.890000,1.169979) (1.967500,1.165771) (2.045000,1.161050) (2.122500,1.155891) (2.200000,1.150387) (2.277500,1.144651) (2.355000,1.138819) (2.432500,1.133046) (2.510000,1.127504) (2.587500,1.122380) (2.665000,1.117859) (2.742500,1.114121) (2.820000,1.111324) (2.897500,1.109593) (2.975000,1.109007) (3.052500,1.109593) (3.130000,1.111324) (3.207500,1.114121) (3.285000,1.117859) (3.362500,1.122380) (3.440000,1.127504) (3.517500,1.133046) (3.595000,1.138819) (3.672500,1.144651) (3.750000,1.150387) (3.827500,1.155891) (3.905000,1.161050) (3.982500,1.165771) (4.060000,1.169979) (4.137500,1.173618) (4.215000,1.176645) (4.292500,1.179027) (4.370000,1.180743) (4.447500,1.181778) (4.525000,1.182124) (4.602500,1.181778) (4.680000,1.180743) (4.757500,1.179027) (4.835000,1.176645) (4.912500,1.173618) (4.990000,1.169979) (5.067500,1.165771) (5.145000,1.161050) (5.222500,1.155891) (5.300000,1.150387) (5.377500,1.144651) (5.455000,1.138819) (5.532500,1.133046) (5.610000,1.127504) (5.687500,1.122380) (5.765000,1.117859) (5.842500,1.114121) (5.920000,1.111324) (5.997500,1.109593) (6.075000,1.109007) (6.152500,1.109593) (6.230000,1.111324) (6.307500,1.114121) (6.385000,1.117859) (6.462500,1.122380) (6.540000,1.127504) (6.617500,1.133046) (6.695000,1.138819) (6.772500,1.144651) (6.850000,1.150387)};
\draw[black!40,line width=.25pt] plot coordinates {(0.650000,1.413752) (0.727500,1.408885) (0.805000,1.404323) (0.882500,1.400146) (0.960000,1.396422) (1.037500,1.393202) (1.115000,1.390523) (1.192500,1.388413) (1.270000,1.386894) (1.347500,1.385978) (1.425000,1.385671) (1.502500,1.385978) (1.580000,1.386894) (1.657500,1.388413) (1.735000,1.390523) (1.812500,1.393202) (1.890000,1.396422) (1.967500,1.400146) (2.045000,1.404323) (2.122500,1.408885) (2.200000,1.413752) (2.277500,1.418822) (2.355000,1.423975) (2.432500,1.429075) (2.510000,1.433967) (2.587500,1.438491) (2.665000,1.442479) (2.742500,1.445776) (2.820000,1.448243) (2.897500,1.449769) (2.975000,1.450286) (3.052500,1.449769) (3.130000,1.448243) (3.207500,1.445776) (3.285000,1.442479) (3.362500,1.438491) (3.440000,1.433967) (3.517500,1.429075) (3.595000,1.423975) (3.672500,1.418822) (3.750000,1.413752) (3.827500,1.408885) (3.905000,1.404323) (3.982500,1.400146) (4.060000,1.396422) (4.137500,1.393202) (4.215000,1.390523) (4.292500,1.388413) (4.370000,1.386894) (4.447500,1.385978) (4.525000,1.385671) (4.602500,1.385978) (4.680000,1.386894) (4.757500,1.388413) (4.835000,1.390523) (4.912500,1.393202) (4.990000,1.396422) (5.067500,1.400146) (5.145000,1.404323) (5.222500,1.408885) (5.300000,1.413752) (5.377500,1.418822) (5.455000,1.423975) (5.532500,1.429075) (5.610000,1.433967) (5.687500,1.438491) (5.765000,1.442479) (5.842500,1.445776) (5.920000,1.448243) (5.997500,1.449769) (6.075000,1.450286) (6.152500,1.449769) (6.230000,1.448243) (6.307500,1.445776) (6.385000,1.442479) (6.462500,1.438491) (6.540000,1.433967) (6.617500,1.429075) (6.695000,1.423975) (6.772500,1.418822) (6.850000,1.413752)};
\draw[line width=.7pt] (.65,4.250000)--(6.85,4.250000);
\node[anchor=east] at (.61,4.250000) {$\Gamma_0$};
\draw[line width=.7pt] (.65,1.290000)--(6.85,1.290000);
\node[anchor=east] at (.61,1.290000) {$\Gamma_1$};
\draw[black!55,densely dotted] (0.650000,1.073029)--(0.650000,4.466971);
\draw[black!55,densely dotted] (6.850000,1.073029)--(6.850000,4.466971);
\node[font=\footnotesize] at (3.75,3.59) {$|\mathcal P_m(u)|>1$ in the lower cylinder};
\node[font=\footnotesize] at (3.75,2.97) {the two poles are omitted in $X_*$};
\node[font=\footnotesize] at (3.75,.42) {vertical sides identified by $u\mapsto u+4K$};
\node[font=\footnotesize] at (3.75,-.05) {dashed curves: $|\mathcal P_m(u)|=r_*$};
\draw[line width=.8pt] (1.370000,2.715000)--(1.480000,2.825000);
\draw[line width=.8pt] (1.370000,2.825000)--(1.480000,2.715000);
\draw[line width=.8pt] (4.470000,2.715000)--(4.580000,2.825000);
\draw[line width=.8pt] (4.470000,2.825000)--(4.580000,2.715000);
\fill (3.517500,4.250000) circle[radius=.04];
\node[below=3pt,fill=white] at (3.517500,4.250000) {$u^{(0)}$};
\fill (2.432500,1.290000) circle[radius=.04];
\node[above=3pt,fill=white] at (2.432500,1.290000) {$u^{(1)}$};
\draw[arr] (7.28,2.85)--(8.73,2.85);
\node at (8.00,3.20) {$\Psi_*$};
\node[font=\footnotesize,align=center] at (8.00,2.35) {proper\\degree two};
\begin{scope}[shift={(11.95,2.95)},scale=1.45]
 \fill[black!3] (-1.65,-1.45) rectangle (1.65,1.45);
 \fill[white] (0,0) circle[radius=.82];
 \fill[even odd rule,black!16] (0,0) circle[radius=1.12] (0,0) circle[radius=.88];
 \draw[densely dashed,line width=.65pt] (0,0) circle[radius=.82];
 \draw[black!40,line width=.3pt] (0,0) circle[radius=.88];
 \draw[black!40,line width=.3pt] (0,0) circle[radius=1.12];
 \draw[line width=.7pt] (0,0) circle[radius=1];
 \fill (0.904607,0.426247) circle[radius=.025];
 \node[above right=1pt,fill=white] at (0.904607,0.426247) {$z_0$};
 \node[font=\footnotesize] at (0,0) {$|z|\le r_*$};
 \node[font=\footnotesize,fill=white] at (-.96,-.84) {$|z|=1$};
\end{scope}
\node[font=\footnotesize,align=center] at (11.95,.64)
 {one physical point has two geometric preimages;\\regularity must be checked at both};
\node[font=\footnotesize] at (11.95,-.05)
 {$\widehat M(R_{\mathrm a}\zeta)=\mathsf D_*\widehat M(\zeta),\quad\det\mathsf D_*\ne0$};
\end{tikzpicture}

\caption{Both geometric sheets of the collar construction. The left panel
shows one real period of the enlarged Jacobi cylinder; vertical sides are
identified, the crosses are the two omitted poles, and $u^{(0)},u^{(1)}$
map to the same point $z_0$ on the physical unit circle. The level curves
use $m_0=16/25$ at anchor $\rho_*=1/4$, with illustrated inner boundary
$r_*=0.82$ and shaded physical collar $0.88<|z|<1.12$. This is a
single-parameter geometric illustration, not a certified choice of
exception-free widths. The proof separately checks the principal matrix
on both components and shrinks the collars to avoid coefficient and
determinant exceptions; no initial knowledge of the unknown on the second
sheet is assumed.}
\label{cont:fig:two-sheet-collar}
\end{figure}
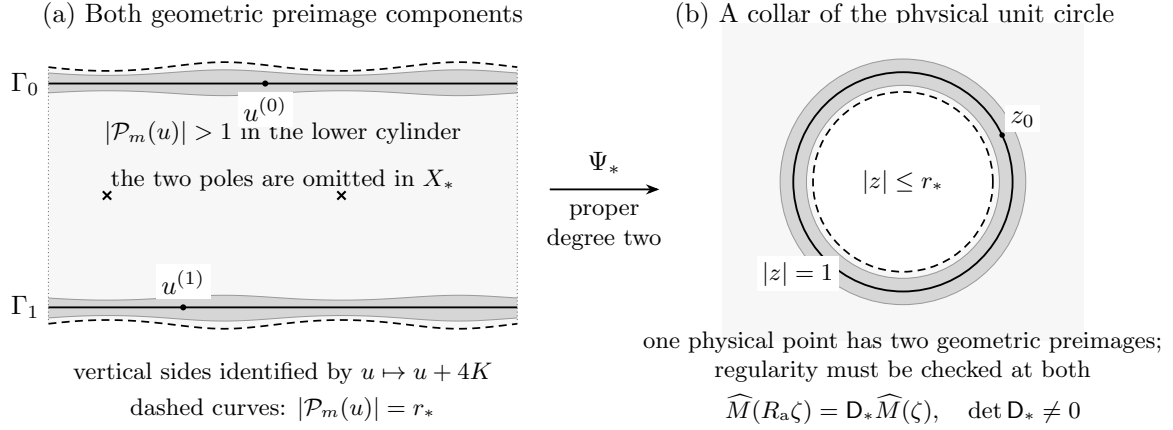

Figure~\ref{cont:fig:two-sheet-collar} illustrates why the preceding
collar construction checks both preimage components of the unit circle.
\Needspace{6\baselineskip}

\begin{lemma}[Gluing all-path continuation across countable exceptional sets]
\label{cont:lem:countable-gluing}
Let $\Omega\subset\mathbb C$ be a domain, and let $B_1,B_2\subset\Omega$
be relatively closed countable sets. Fix a holomorphic germ at
$z_0\in\Omega\setminus(B_1\cup B_2)$.
Suppose that it admits analytic continuation along every path in each
$D_i=\Omega\setminus B_i$ starting at $z_0$.
Then it admits analytic continuation along every path in
\[
D_1\cup D_2=\Omega\setminus(B_1\cap B_2)
\]
starting at $z_0$. Different paths may lead to different germs;
no restriction on monodromy is imposed.
\end{lemma}
\begin{proof}
Put $D_{12}=D_1\cap D_2=\Omega\setminus(B_1\cup B_2)$, and denote the
initial germ at $z_0$ by $\mathfrak f_0$. All paths below are continuous
maps from $[0,1]$; no rectifiability is assumed. We first show that the
two domains supply the same \emph{sets of reachable germs} on their
overlap, rather than selecting a single branch there.

\smallskip\noindent
\emph{Step 1: avoid the additional countable set without changing the homotopy class.}
Fix $i\in\{1,2\}$ and a path $\gamma:[0,1]\to D_i$ whose endpoints lie
in $D_{12}$. Relative closedness of $B_i$ makes $D_i$ open. Compactness
of the path image therefore permits a number $\delta>0$ such that
\[
 \{z\in\C:\operatorname{dist}(z,\gamma([0,1]))\le\delta\}
 \subset D_i.
\]
Choose a partition $0=t_0<\cdots<t_N=1$ so fine that
\[
 |\gamma(t)-\gamma(t_j)|<\delta/8
 \quad(t_j\le t\le t_{j+1}).
\]
For $0<j<N$, choose a point $v_j\in D_{12}$ with
$|v_j-\gamma(t_j)|<\delta/16$, and keep
$v_0=\gamma(0)$ and $v_N=\gamma(1)$ fixed. Such choices are possible
because the small disks lie in $D_i$ and a countable set contains no
open disk.

Both $v_j$ and $v_{j+1}$ lie in the disk
$D(\gamma(t_j),\delta/4)$. In that disk choose a point $w_j$ outside
all lines through $v_j$ and a point of $B_{3-i}$, and outside all lines
through $v_{j+1}$ and a point of $B_{3-i}$. The endpoints are not in
$B_{3-i}$, so these are well-defined lines. There are only countably
many of them; their union has planar measure zero and cannot fill the
disk. The two segments $[v_j,w_j]$ and $[w_j,v_{j+1}]$ then avoid
$B_{3-i}$ and stay inside the disk, which itself lies in $D_i$.

Parametrize these two segments on the two halves of
$[t_j,t_{j+1}]$. The resulting polygonal path $\widetilde\gamma$ lies
in $D_{12}$, has the prescribed endpoints, and satisfies
\[
 \sup_{0\le t\le1}|\widetilde\gamma(t)-\gamma(t)|<\delta.
\]
Consequently
\[
 H(r,t)=(1-r)\gamma(t)+r\widetilde\gamma(t),\qquad 0\le r,t\le1,
\]
is a homotopy in $D_i$ with endpoints fixed: each $H(r,t)$ remains
within distance $\delta$ of $\gamma(t)$. The homotopy need not avoid
$B_{3-i}$; only its final path must do so. This proves the required
approximation even when $B_{3-i}$ has accumulation points.

\smallskip\noindent
\emph{Step 2: identify all reachable germs on the overlap.}
For $z\in D_{12}$, let $\mathscr G_i(z)$ be the set of germs obtained
from $\mathfrak f_0$ by continuation along paths from $z_0$ to $z$ in
$D_i$. Define $\mathscr G_{12}(z)$ in the same way using paths in
$D_{12}$. The inclusion
$\mathscr G_{12}(z)\subset\mathscr G_i(z)$ is immediate.

For the reverse inclusion, choose a representing path in $D_i$ for a
germ in $\mathscr G_i(z)$, and apply Step~1. Its replacement path lies
in $D_{12}$ and is homotopic to it in $D_i$ with endpoints fixed.
Homotopy invariance of analytic continuation applies here because the
hypothesis supplies continuation from $\mathfrak f_0$ along every path
of that homotopy. More explicitly, for any one path of the homotopy,
its continuation is represented by a finite chain of holomorphic
elements on disks, with agreement near the successive joining points.
Compactness permits positive margins inside these disks and their
joining neighborhoods. Uniformly nearby paths use the same chain and
have the same terminal germ. The terminal germ is thus locally
constant in the homotopy parameter, and hence constant on $[0,1]$.
This argument does not require $D_i$ to be simply connected.

The replacement path consequently reaches the same germ at $z$. Thus
$\mathscr G_i(z)\subset\mathscr G_{12}(z)$, and applying the argument
for both indices yields
\[
 \boxed{\mathscr G_1(z)=\mathscr G_{12}(z)=\mathscr G_2(z),
        \qquad z\in D_{12}.}
\]
This is equality of sets of germs; their common set need not consist
of a single germ. In addition, any germ in $\mathscr G_i(z)$ continues
along every further path in $D_i$ starting at $z$: concatenate that
path with a path representing the germ, and use the original
all-path hypothesis. Uniqueness along the representing path ensures
that the continuation starts with the specified germ at $z$.

\smallskip\noindent
\emph{Step 3: switch domains while preserving the germ along the given path.}
Let $\gamma:[0,1]\to D_1\cup D_2$ start at $z_0$. The inverse images
$\gamma^{-1}(D_1)$ and $\gamma^{-1}(D_2)$ form an open cover of
$[0,1]$. A Lebesgue number for this cover gives a finite subdivision
\[
 0=t_0<t_1<\cdots<t_N=1
\]
for which every closed subpath $\gamma([t_{j-1},t_j])$ lies in one
$D_{i_j}$. Continue $\mathfrak f_0$ along the first subpath. Inductively,
the germ at the end of a subpath is reachable from $z_0$ through its
current domain. If the next domain is the same, the final observation
of Step~2 extends that germ along the next subpath.

If the index changes, the joining point lies in $D_{12}$, because it
belongs to both adjacent closed subpaths. The displayed equality of
reachable-germ sets represents the current germ by a path from $z_0$
in the next domain. Concatenation with the next actual subpath then
continues exactly this germ and shows that the new endpoint germ is
reachable through the next domain. A different representing path is
used only to invoke that domain's continuation hypothesis; it does
not replace the germ already reached along the original path.

There are finitely many subpaths, so this constructs continuation
along all of $\gamma$, including its endpoint. At no stage are all
branches at a point identified, and no loop is required to have
trivial monodromy. The conclusion is all-path continuation on
$D_1\cup D_2=\Omega\setminus(B_1\cap B_2)$ with arbitrary monodromy
still allowed.
\end{proof}

\begin{lemma}[Two bases with disjoint accumulation bounds]
\label{cont:lem:disjoint-base-bounds}
For a fixed noncircular ellipse and any nonempty open interior rotation
interval, one can choose two reduced rational bases $\rho_i=p_i/d_i$
in that interval with $4\mid d_i$ whose finite physical exterior
accumulation bounds in Lemma~\ref{cont:lem:arbitrary-rational-base} are disjoint.
\end{lemma}
\begin{proof}
Let $I\subset(0,1/2)$ be the given nonempty open interval. Normalize the
major semiaxis to $a=1$ and put $\eta=1-b^2\in(0,1)$. The ellipse and
$\eta$ remain fixed throughout. For a reduced base $\rho=p/d\in I$,
write $m=m(\rho)$, $K=K(m)$, $K'=K(1-m)$, and
\[
 \mathcal P_m(u)=\sn(u\mid m)-i\cn(u\mid m).
\]
As before, $K'$ denotes the complementary integral, not
differentiation. We separate the geometric accumulation bound from the
arithmetic condition that would force every second bound to meet it.

\smallskip\noindent
\emph{Step 1: identify the finite bound for denominators divisible by four.}
By Lemma~\ref{cont:lem:arbitrary-rational-base}, the possible
accumulation points are contained in the finite physical images of
\[
 jK-\frac{4K\ell}{d}-iK',\qquad
 0\le j<4,\quad 0\le\ell<d,
\]
modulo the real period $4K$, with the poles omitted. These are the
lower-cylinder Jacobi coefficient points shifted through all $d$
action fibers. If $4\mid d$, each $jK$ is an integer multiple of
$4K/d$, and $j=0$ with $\ell$ varying already supplies the entire
$d$-grid. Hence a valid finite containing set is
\[
 F_{p/d}=\left\{
   \mathcal P_m\!\left(\frac{4Kj}{d}-iK'\right):0\le j<d,
   \quad\mathcal P_m\!\left(\frac{4Kj}{d}-iK'\right)\ne\infty
              \right\}.
\]
All its points lie in the finite exterior $\{|z|>1\}$ by
Lemma~\ref{arith:jacobi}. Repetitions are removed, and retaining
removable coefficient points only enlarges a permissible bound. The
set describes possible accumulation locations, not the whole
exceptional set of the continuation construction.

\smallskip\noindent
\emph{Step 2: exclude the focal values and fix a first base.}
For a reduced base $p/d$ with $4\mid d$, its numerator $p$ is odd.
Lemma~\ref{arith:focal} gives the complete lower-cylinder fibers of
the two exterior focal values as
\[
 K\pm J-iK',\qquad 3K\pm J-iK',\qquad J=\frac{2Kp}{d}.
\]
Their prospective indices in the preceding grid are
\[
 \frac{d}{4K}(K\pm J)=\frac d4\pm\frac p2,
 \qquad
 \frac{d}{4K}(3K\pm J)=\frac{3d}{4}\pm\frac p2.
\]
These are nonintegral half-integers. Adding a real period changes an
index by an integer multiple of $d$ and cannot make it integral.
Since these are the full fibers, neither exterior focal value belongs
to $F_{p/d}$.

Reduced fractions with denominator divisible by four are dense: for
$k\ge2$, the fractions with odd numerator and denominator $2^k$ are
reduced and have spacing $2^{1-k}$. Choose one such base
$\rho_1=p_1/d_1\in I$ and put $F_1=F_{p_1/d_1}$. This finite set
consists of nonfocal exterior points. If it were empty, any distinct
such base would finish the proof; assume henceforth that $F_1$ is
nonempty.

It suffices to find a second reduced base of denominator $4q$, with
$q$ an arbitrarily large odd prime, whose bound misses $F_1$.
Suppose to the contrary that for every sufficiently large odd prime
$q$, every reduced $p/(4q)\in I$ satisfies
\[
 F_{p/(4q)}\cap F_1\ne\varnothing.
\]
We may take $q$ so large that these candidates differ from the first
base. We will contradict this unavoidable-intersection assumption.

\smallskip\noindent
\emph{Step 3: choose both inverse branches on a common parameter neighborhood.}
For a fixed finite point $z\in F_1$, the identities from
Lemma~\ref{arith:jacobi} give
\[
 \sn u=\frac{\mathcal P_m(u)+\mathcal P_m(u)^{-1}}2,
 \qquad
 \partial_u\mathcal P_m(u)=i\dn(u\mid m)\mathcal P_m(u).
\]
Thus a preimage of $z$ can be critical only when
\[
 1-\frac m4(z+z^{-1})^2=0,
 \qquad m=m_z:=\frac4{(z+z^{-1})^2}.
\]
The denominator is nonzero for $|z|>1$. The analytic function
$m(\rho)$ is strictly monotone on the real band, so this excludes at
most one real rotation for each $z\in F_1$. Choose a nonempty open
interval $I_0$ with compact closure in $I$, avoiding these finitely
many rotations and $1/4$.

After shrinking $I_0$ further, both preimages of every point of $F_1$
admit holomorphic lifts on a common simply connected complex
neighborhood of $\overline{I_0}$. Indeed, at a regular parameter the
degree-two map has two distinct unramified preimages of each fixed
point. Choose representatives in the lower cylinder and apply the
holomorphic implicit-function theorem to
$\mathcal P_{m(\rho)}(u)=z$ at all these finitely many points. A smaller
common neighborhood retains all branches and their distinctness.
For real $\rho\in I_0$, the representatives can be kept in
\[
 -2K'(\rho)<\operatorname{Im}u_s(\rho)<0
\]
by continuity and the strict exterior-cylinder property. Index the
finite list by $1\le s\le S$, recording both its fixed physical value
$z_s\in F_1$ and its choice of inverse branch, and define
\[
 g_s(\rho)=\frac{u_s(\rho)+iK'(\rho)}{4K(\rho)}.
\]
These functions are holomorphic on the common neighborhood. Changing
a lower-cylinder representative by $4K(\rho)k$, $k\in\Z$, changes
$g_s$ by an integer. There is no remaining imaginary-period ambiguity
between two representatives in that cylinder.

\smallskip\noindent
\emph{Step 4: express every intersection as a spacing-two lattice condition.}
Let $q$ be a sufficiently large odd prime, put $N=4q$, and let $p$ be
odd with $p/N\in I_0$. Then $0<p<2q$ and $\gcd(p,4)=1$. The only
possible multiple of $q$ in that range is $p=q$, whose rotation is
$1/4$ and was excluded. Hence $\gcd(p,N)=1$, so each such grid point
is covered by the contradiction assumption.

At an intersection with $F_1$, some point $z_s$ has a preimage
$4Kj/N-iK'$, with $0\le j<N$. The list contains both members of its
fiber. Consequently one listed representative agrees with that grid
point modulo the real period:
\[
 u_s(p/N)=\frac{4Kj}{N}-iK'+4Kk
 \quad\text{for some }k\in\Z,
\]
where the complete integrals are evaluated at $\rho=p/N$. Therefore
\[
 \boxed{N g_s(p/N)=j+Nk\in\Z.}
\]
The supplying index may depend on both $p$ and $q$. No fixed branch
has yet been asserted to supply all grid points.

\smallskip\noindent
\emph{Step 5: force one fixed normalized Abel coordinate to be affine.}
We give the spacing-two version of the argument in
Lemma~\ref{arith:grid}. Fix $\rho^\circ\in I_0$ and a compact interval
$I_1\Subset I_0$ containing it in its interior. For $N=4q$, put
\[
 p_0(N)=2\left\lfloor\frac{N\rho^\circ}{2}\right\rfloor+1-2S.
\]
Consider the $2S+1$ consecutive odd numerators
$p_0(N),p_0(N)+2,\ldots,p_0(N)+4S$. Their quotients by $N$ tend to
$\rho^\circ$ and lie in $I_1$ for large $q$. Assign to each an index
supplying the integer in Step~4. At least three have the same index;
write their numerators as $p,p+a,p+a+b$, where $a,b$ are positive even
integers and $a+b\le4S$. With $x_N=p/N$, form the integer
\[
 Q_N=N\bigl[b g_s(x_N)-(a+b)g_s(x_N+a/N)
                    +a g_s(x_N+(a+b)/N)\bigr].
\]
The weights cancel the constant and linear Taylor terms, and the
quadratic coefficient is fixed by
\[
 \begin{aligned}
 b-(a+b)+a&=0,\\
 -(a+b)a+a(a+b)&=0,\\
 -(a+b)a^2+a(a+b)^2&=ab(a+b).
 \end{aligned}
\]
Thus
\[
 Q_N=\frac{ab(a+b)}{2N}g_s''(x_N)+O(N^{-2}).
\]
The error is uniform over the finitely many possibilities for $s,a,b$:
the first three derivatives of all $g_s$ are bounded on $I_1$, the
gaps are bounded by $4S$, and the cubic Taylor error is multiplied by
$N$. In particular, $|Q_N|\le C_{I_1}/N$. The selected values are
ordinary integers, even though the $g_s$ may be complex-valued, so
integrality gives $Q_N=0$ for all sufficiently large $q$.

Pass to a subsequence of primes on which $(s,a,b)$ is fixed.
Multiplying the last expansion by $N$ and taking $x_N\to\rho^\circ$
gives $g_s''(\rho^\circ)=0$. Since $\rho^\circ$ was arbitrary, at
least one of the finitely many $g_s''$ vanishes at every point of
$I_0$. If none vanished identically, each would have only finitely
many zeros on any fixed nondegenerate compact subinterval of $I_0$,
by the identity theorem. Their finite union could not cover that
subinterval. Hence there is one fixed index $s_*$ such that
\[
 g_{s_*}(\rho)=A\rho+B\qquad(\rho\in I_0)
\]
for constants $A,B\in\C$. The physical point $z_{s_*}$ and inverse
branch are now fixed. Neither reality nor rationality of $A,B$ is
needed.

\smallskip\noindent
\emph{Step 6: contradict the marked-point Picard--Fuchs source.}
Choose a fixed complex angle $\varphi_*$ with
$z_{s_*}=e^{i\varphi_*}$. Since $J(m)=2K(m)\rho$, the affine identity
becomes
\[
 u_{s_*}+iK'=4K(A\rho+B)=2AJ+4BK.
\]
The change between $\rho$ and $m$ is locally analytic and invertible,
so this is an identity in $m$ on a nonempty parameter neighborhood.
Choose there a compatible branch of $I_{\varphi_*}$. The reference
Jacobi inverse is $K+I_{\varphi_*}$; the other inverse is obtained by
$u\mapsto2K+2iK'-u$, and representatives differ by periods $4K$ and
$4iK'$. Consequently, after restricting this neighborhood if necessary,
\[
 u_{s_*}=\epsilon I_{\varphi_*}+c_1K+c_2iK',
 \qquad \epsilon\in\{1,-1\},\quad c_1,c_2\in\Z.
\]
The sign and integer coefficients are constant: compare the branches
at one parameter, fix the corresponding deck and period choices, and
use uniqueness of the local unramified inverse to keep the equality
nearby. These are not parameter-dependent coefficients. All integral
and square-root branches in the following calculation are continued
compatibly.

Apply $4\mathcal L$, with $\eta$ and $\varphi_*$ fixed and
\[
 \mathcal L=m(1-m)\partial_m^2+(1-2m)\partial_m-\tfrac14.
\]
Lemma~\ref{arith:pf} annihilates the complete-period terms and gives
\[
 \epsilon\frac{\sin\varphi_*\cos\varphi_*}
                 {(1-m\cos^2\varphi_*)^{3/2}}
 =-\frac{2A\sqrt{\eta(1-\eta)}}{(m-\eta)^{3/2}}.
\]
The factor $2A$ comes from $2AJ$ and the identity for $4\mathcal LJ$.
Squaring removes the sign and branch powers and yields equality of
rational functions on this parameter neighborhood:
\[
 \frac{\sin^2\varphi_*\cos^2\varphi_*}
      {(1-m\cos^2\varphi_*)^3}
 =\frac{4A^2\eta(1-\eta)}{(m-\eta)^3}.
\]
For the finite exterior point $z_{s_*}$,
\[
 \cos\varphi_*=(z_{s_*}+z_{s_*}^{-1})/2\ne0,
 \qquad
 \sin\varphi_*=(z_{s_*}-z_{s_*}^{-1})/(2i)\ne0;
\]
a zero of either would put $z_{s_*}$ on the unit circle. Put
\[
 m_*:=\sec^2\varphi_*=\frac4{(z_{s_*}+z_{s_*}^{-1})^2}.
\]
It is finite and nonzero, and $m_*\ne1$ since
$\sin\varphi_*\ne0$. By Step~2 the point is nonfocal, so
$1-\eta\cos^2\varphi_*\ne0$ and therefore $m_*\ne\eta$.

Clearing denominators gives the polynomial identity
\[
 \sin^2\varphi_*\cos^2\varphi_*\,(m-\eta)^3
   =4A^2\eta(1-\eta)(1-m\cos^2\varphi_*)^3.
\]
It holds identically because it holds on a nonempty open parameter
set. At $m=m_*$ its right side is zero, whereas its left side is
\[
 \sin^2\varphi_*\cos^2\varphi_*\,(m_*-\eta)^3\ne0,
\]
a contradiction, including when $A=0$. No continuation of the inverse
branches to $m_*$ is required: only the polynomial identity is
evaluated there.

The unavoidable-intersection assumption is false. Hence a reduced
second base $\rho_2=p_2/d_2\in I$, with $d_2=4q$ for a sufficiently
large odd prime $q$, satisfies $F_1\cap F_{\rho_2}=\varnothing$.
Both denominators are divisible by four, and both finite sets are the
accumulation bounds supplied by Lemma~\ref{cont:lem:arbitrary-rational-base}.
This separates those bounds; it does not assert that the full
continuation exceptional sets are disjoint or that their monodromy is
trivial.
\end{proof}

\begin{figure}[tbp]
\centering
\begin{tikzpicture}[x=1cm,y=1cm,font=\small,
 arr/.style={-{Stealth[length=1.7mm]},line width=.7pt},
 every node/.style={inner sep=2pt}]
\path[use as bounding box] (-.1,-.65) rectangle (15.75,5.43);
\node at (3.67,5.15) {(a) Switch domains without changing the germ};
\node at (11.70,5.15) {(b) Disjoint accumulation bounds, shared point};
\fill[black!5] (2.57,3.17) circle[radius=1.40];
\fill[black!8] (4.52,3.17) circle[radius=1.40];
\draw[line width=.55pt] (2.57,3.17) circle[radius=1.40];
\draw[dashed,line width=.55pt] (4.52,3.17) circle[radius=1.40];
\node[font=\footnotesize] at (2.35,3.95) {$\mathcal U_1\subset D_1$};
\node[font=\footnotesize] at (4.77,3.95) {$\mathcal U_2\subset D_2$};
\coordinate (entry) at (1.44,2.92);
\coordinate (join) at (3.52,3.08);
\coordinate (end) at (5.54,2.95);
\draw[arr] (.25,2.45)--(entry);
\draw[line width=.85pt,postaction={decorate},decoration={markings,mark=at position .55 with {\arrow{Stealth}}}] (entry)..controls (2.0,2.08) and (2.88,2.46)..(join);

\draw[line width=.85pt,dashed,postaction={decorate},decoration={markings,mark=at position .60 with {\arrow{Stealth}}}] (join)..controls (4.04,3.40) and (4.78,2.12)..(end);

\fill (join) circle[radius=.042];
\node[above=4pt,fill=white] at (join) {$[f]_q$};
\node[below=3pt,fill=white,inner sep=1pt] at (join) {$q$};
\node[font=\footnotesize,align=center] at (.56,1.99) {previous\\history};
\node[font=\footnotesize] at (3.65,1.17)
 {$\mathscr G_1(q)=\mathscr G_{12}(q)=\mathscr G_2(q)$};
\node[font=\footnotesize,align=center] at (3.65,.49)
 {equality of sets of reachable germs,\\not a claim that the germ is unique};
\draw[black!45,line width=.35pt] (8.60,4.13)--(15.06,4.13);
\node[anchor=east] at (8.47,4.13) {$B_1$};
\draw[black!40] (9.940,4.095)--(9.940,4.165);
\node[below=5pt,font=\scriptsize] at (9.940,4.13) {$2$};
\draw[black!40] (11.280,4.095)--(11.280,4.165);
\node[below=5pt,font=\scriptsize] at (11.280,4.13) {$3$};
\draw[black!40] (13.960,4.095)--(13.960,4.165);
\node[below=5pt,font=\scriptsize] at (13.960,4.13) {$5$};
\draw[line width=.8pt] (9.885,4.075)--(9.995,4.185);
\draw[line width=.8pt] (9.885,4.185)--(9.995,4.075);
\draw[black!45,line width=.35pt] (8.60,2.77)--(15.06,2.77);
\node[anchor=east] at (8.47,2.77) {$B_2$};
\draw[black!40] (9.940,2.735)--(9.940,2.805);
\node[below=5pt,font=\scriptsize] at (9.940,2.77) {$2$};
\draw[black!40] (11.280,2.735)--(11.280,2.805);
\node[below=5pt,font=\scriptsize] at (11.280,2.77) {$3$};
\draw[black!40] (13.960,2.735)--(13.960,2.805);
\node[below=5pt,font=\scriptsize] at (13.960,2.77) {$5$};
\draw[line width=.8pt] (9.885,2.715)--(9.995,2.825);
\draw[line width=.8pt] (9.885,2.825)--(9.995,2.715);
\draw[black!45,line width=.35pt] (8.60,1.41)--(15.06,1.41);
\node[anchor=east] at (8.47,1.41) {$B_1\cap B_2$};
\draw[black!40] (9.940,1.375)--(9.940,1.445);
\node[below=5pt,font=\scriptsize] at (9.940,1.41) {$2$};
\draw[black!40] (11.280,1.375)--(11.280,1.445);
\node[below=5pt,font=\scriptsize] at (11.280,1.41) {$3$};
\draw[black!40] (13.960,1.375)--(13.960,1.445);
\node[below=5pt,font=\scriptsize] at (13.960,1.41) {$5$};
\draw[line width=.8pt] (9.885,1.355)--(9.995,1.465);
\draw[line width=.8pt] (9.885,1.465)--(9.995,1.355);
\fill (11.95000,4.13) circle[radius=.025];
\fill (11.72667,4.13) circle[radius=.025];
\fill (11.61500,4.13) circle[radius=.025];
\fill (11.54800,4.13) circle[radius=.025];
\fill (11.50333,4.13) circle[radius=.025];
\fill (11.47143,4.13) circle[radius=.025];
\fill (11.44750,4.13) circle[radius=.025];
\fill (11.42889,4.13) circle[radius=.025];
\draw[line width=.7pt,fill=white] (11.28000,4.13) circle[radius=.065];
\node[above=9pt,font=\footnotesize] at (11.530,4.13) {$F_1=\{3\}$};
\fill (14.60700,2.74700) rectangle (14.65300,2.79300);
\fill (14.38367,2.74700) rectangle (14.42967,2.79300);
\fill (14.27200,2.74700) rectangle (14.31800,2.79300);
\fill (14.20500,2.74700) rectangle (14.25100,2.79300);
\fill (14.16033,2.74700) rectangle (14.20633,2.79300);
\fill (14.12843,2.74700) rectangle (14.17443,2.79300);
\fill (14.10450,2.74700) rectangle (14.15050,2.79300);
\fill (14.08589,2.74700) rectangle (14.13189,2.79300);
\draw[line width=.7pt,fill=white] (13.96000,2.77) circle[radius=.065];
\node[above=9pt,font=\footnotesize] at (14.210,2.77) {$F_2=\{5\}$};
\node[font=\footnotesize] at (11.74,.43)
 {$F_1\cap F_2=\varnothing,\qquad B_1\cap B_2=\{2\}\ne\varnothing$};
\node[font=\footnotesize] at (7.81,-.35)
 {gluing removes points excluded by only one construction; common isolated obstructions may remain};
\end{tikzpicture}

\caption{Reached-germ gluing and the surviving intersection. Left: the
disks represent local neighborhoods $\mathcal U_i\subset D_i$, not the
entire continuation domains. Equality of the sets of reachable germs
allows the particular germ $[f]_q$ to pass into the next domain; it does
not identify all branches. Right: a set-theoretic model with
$B_1=\{2,3\}\cup\{3+1/n:n\ge2\}$ and
$B_2=\{2,5\}\cup\{5+1/n:n\ge2\}$. Only initial portions of the convergent
sequences are drawn; the open circles mark their included limit points.
The bounds $F_1=\{3\}$ and $F_2=\{5\}$ are disjoint, but the isolated
point $2$ survives in $B_1\cap B_2$. These are illustrative sets, not
computed billiard exceptions. In the theorem the intersection may still
be infinite, with points tending to infinity.}
\label{cont:fig:two-base-gluing}
\end{figure}
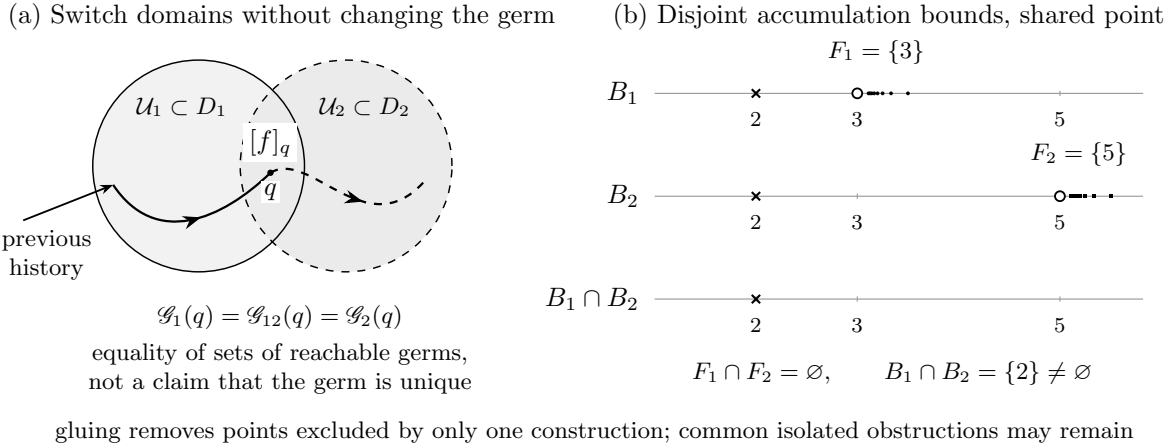

\smallskip\noindent\textbf{Exceptional sets versus accumulation bounds.}
In the common physical domain $\Omega$ used below, write $B_i$ for the
exceptional set of the $i$th base and $F_i$ for its finite accumulation
bound. A prime with subscript $\Omega$ denotes a derived set relative to
that domain. The three roles are distinct:
\begin{center}
\begin{tabular}{@{}p{.21\linewidth}p{.73\linewidth}@{}}
\toprule
Set & Role in the two-base argument\\
\midrule
$B_i$ & Relatively closed countable exceptions, with $(B_i)'_\Omega\subset F_i$.\\[3pt]
$F_i$ & A finite containing bound, not the whole exceptional set;
Lemma~\ref{cont:lem:disjoint-base-bounds} gives $F_1\cap F_2=\varnothing$.\\[3pt]
$B=B_1\cap B_2$ & Exceptions remaining after reached-germ gluing; common
isolated points may survive.\\
\bottomrule
\end{tabular}
\end{center}
The relative conclusion is
\[
 B'_\Omega\subset (B_1)'_\Omega\cap(B_2)'_\Omega
             \subset F_1\cap F_2=\varnothing.
\]
To conclude $B'_\C=\varnothing$, the proof also separates $B$ from the
inner boundary of $\Omega$ by an exception-free collar. Neither conclusion
makes $B$ finite or removes monodromy. Figure~\ref{cont:fig:two-base-gluing}
illustrates the gluing mechanism and the distinction between the sets.

\Needspace{6\baselineskip}

\begin{theorem}[Exterior continuation with no finite accumulation]
\label{cont:thm:discrete-exterior-continuation}
For every real-valued analytic displacement at a fixed noncircular
ellipse whose traces are constant at all rational rotations in a
nonempty open interior band, there is $0<r_*<1$ such that its physical
germ $V(z)$ admits
analytic continuation along every path in
\[
\{|z|>r_*\}\setminus B
\]
starting in its original collar, where $B$ is a relatively closed
discrete set disjoint from a neighborhood of the unit circle.
In particular, the possible exterior singular points have no finite
accumulation.
Arbitrary monodromy remains allowed; the statement does not yet yield
global meromorphic continuation.
\end{theorem}
\begin{proof}
Fix the reference ellipse and the open rotation band $I$, and denote
the given analytic displacement by $v$. Let $V$ be its physical germ in
$z=e^{i\varphi}$. For a subset $A$ of an open set $\mathcal O\subset\C$,
write $A'_{\mathcal O}$ for its derived set with accumulation taken
inside $\mathcal O$. All continuation paths below are continuous maps
from $[0,1]$, including their endpoints; no rectifiability is required.

\smallskip\noindent
\emph{Step 1: obtain two constructions for the same physical germ.}
Lemma~\ref{cont:lem:disjoint-base-bounds} supplies two reduced rational
bases $\rho_i=p_i/d_i\in I$, with $4\mid d_i$, whose finite physical
accumulation bounds $F_1,F_2$ satisfy
\[
 F_1\cap F_2=\varnothing.
\]
Apply Lemmas~\ref{cont:lem:arbitrary-rational-base} and
\ref{cont:lem:physical-collar-domain} at these bases. For $i=1,2$, let
\[
 \Omega_i=\{z\in\C:|z|>r_i\},\qquad 0<r_i<1,
\]
and let $E_i\subset\Omega_i$ be the resulting relatively closed
countable exceptional set. The original germ $V$ continues along every
path in $\Omega_i\setminus E_i$ starting in its original collar, and
\[
 (E_i)'_{\Omega_i}\subset F_i.
\]
Moreover, $E_i$ avoids a neighborhood of the unit circle. The collar
extension in Lemma~\ref{cont:lem:physical-collar-domain} preserves the
same containing set $F_i$; it does not add accumulation locations.

The action representatives $v\circ\Phi_{E,\rho_1}$ and
$v\circ\Phi_{E,\rho_2}$ need not be the same function of their action
variables. Nevertheless, both constructions have already been
transferred to the common physical variable $z$ and start with the
same germ $V$. No identification of the two action coordinates, and no
uniform collar for all rational bases, is needed.

\smallskip\noindent
\emph{Step 2: restrict to a common domain inside an exception-free collar.}
By compactness of the unit circle, choose $0<\delta<1$ so small that
\[
 \mathcal C_\delta=\{z:1-\delta<|z|<1+\delta\}
\]
lies in $\Omega_1\cap\Omega_2$ and in the original holomorphy annulus
of $V$, and
\[
 E_1\cap\mathcal C_\delta=E_2\cap\mathcal C_\delta=\varnothing.
\]
Choose
\[
 \max\{r_1,r_2,1-\delta\}<r_*<1,
 \qquad \Omega=\{z:|z|>r_*\},\qquad B_i=E_i\cap\Omega.
\]
Then $\Omega\subset\Omega_i$. Restriction to this open subdomain
preserves relative closedness and countability. If distinct points of
$B_i$ converge to a point of $\Omega$, the same sequence lies in $E_i$
and converges inside $\Omega_i$. Hence
\[
 (B_i)'_\Omega
 \subset (E_i)'_{\Omega_i}\cap\Omega
 \subset F_i.
\]
The all-path continuation assertion also restricts: every path in
$\Omega\setminus B_i$ lies in $\Omega_i\setminus E_i$.

The radius choice gives the additional separation
\[
 B_i\subset\{z:|z|\ge1+\delta\},\qquad i=1,2.
\]
Indeed, any point of $\Omega$ with modulus less than $1+\delta$ belongs
to $\mathcal C_\delta$ and is not in $E_i$. In particular,
\[
 \mathcal C_*:=\Omega\cap\mathcal C_\delta
     =\{z:r_*<|z|<1+\delta\}
\]
is a connected exception-free collar containing the unit circle, on
which the original $V$ is holomorphic. This separation will be needed
below to rule out accumulation at the boundary of $\Omega$, not only
inside it.

\smallskip\noindent
\emph{Step 3: glue outside the intersection while preserving reached germs.}
Take the common initial point $z_0=1\in\mathcal C_*$, with its original
germ, and put
\[
 D_i=\Omega\setminus B_i,\qquad D_{12}=D_1\cap D_2.
\]
The same germ at $z_0$ continues along every path in each $D_i$.
All hypotheses of Lemma~\ref{cont:lem:countable-gluing} are therefore
satisfied. That lemma yields continuation along every path in
\[
 D_1\cup D_2=\Omega\setminus(B_1\cap B_2).
\]
Set
\[
 B:=B_1\cap B_2.
\]

Here the gluing lemma is substantive: it proves equality of the sets
of reachable germs at every point of $D_{12}$. When a path switches
from one continuation domain to the other, the particular germ already
reached can be represented through the next domain and continued from
there. It does not replace that germ by the original branch or identify
all branches at a point. Thus the new exceptional set is the
\emph{intersection}, rather than the union, and arbitrary monodromy
remains allowed. The finite-subdivision argument in that lemma applies
to every continuous path on $[0,1]$.

\smallskip\noindent
\emph{Step 4: exclude accumulation both in the common domain and in the finite plane.}
The set $B$ is relatively closed and countable in $\Omega$. A convergent
sequence of distinct points in $B$ is simultaneously such a sequence
in both $B_i$. Consequently
\[
 B'_\Omega
   =(B_1\cap B_2)'_\Omega
   \subset (B_1)'_\Omega\cap(B_2)'_\Omega
   \subset F_1\cap F_2=\varnothing.
\]
Every point of $B$ is therefore isolated, and $B$ is relatively closed
and discrete in $\Omega$. Only the displayed inclusion for derived
sets is used; equality with the intersection of the two derived sets
is not asserted. Nor does disjointness of $F_1,F_2$ assert that the
full exceptional sets $B_1,B_2$ are disjoint: they may share isolated
points.

To obtain the stronger absence of all finite accumulation, use Step~2:
\[
 B\subset\{z:|z|\ge1+\delta\}.
\]
If pairwise distinct $b_n\in B$ converged to any finite $b_\infty\in\C$,
then $|b_\infty|\ge1+\delta>r_*$, so $b_\infty\in\Omega$. This would
contradict $B'_\Omega=\varnothing$. Thus
\[
 B'_\C=\varnothing,
 \qquad B\cap K\text{ is finite for every compact }K\subset\C.
\]
The compact-set assertion follows because an infinite subset of a
compact set has an accumulation point. In particular, no accumulation
can occur at $|z|=r_*$ or at the unit circle. Relative discreteness in
an exterior domain alone would not exclude those boundary limits;
the exception-free collar is essential. The set $B$ may still be
infinite with points tending to infinity.

\smallskip\noindent
\emph{Step 5: allow every original starting point.}
Let $z_1\in\mathcal C_*$ and let
$\gamma:[0,1]\to\Omega\setminus B$ start at $z_1$. Join $z_0$ to $z_1$
by a path $\beta$ in $\mathcal C_*$. Along $\beta$ the already
holomorphic function $V$ supplies its original germs. Applying Step~3
to $\beta*\gamma$ therefore continues along $\gamma$ from precisely
the original germ at $z_1$, by uniqueness of continuation along $\beta$.

The same conclusion holds for a starting point in any larger originally
chosen holomorphy annulus, provided it lies in $\Omega\setminus B$ as
the path requires. The intersection of that annulus with $\Omega$ is
an open annulus containing $z_0$ and $z_1$. Connect the points there by
a polygonal path. Its compact image meets $B$ in only finitely many
points by Step~4; small detours in disks contained in the annulus and
isolating those points produce a connector avoiding $B$. Along this
connector uniqueness identifies continuation with the already known
$V$. Concatenation with the requested path proves the statement for
that starting point as well. This completes the exterior assertion.

\smallskip\noindent
\emph{Step 6: record the reflected consequence and the limits of the conclusion.}
For completeness, reality also gives a punctured-interior consequence.
On a sufficiently small annulus invariant under
$z\mapsto1/\overline z$ and contained in the original holomorphy
annulus, the function
\[
 V^\sharp(z):=\overline{V(1/\overline z)}
\]
is holomorphic. On $|z|=1$ it equals $V$, since
$V(e^{i\varphi})=v(\varphi)\in\R$. The identity theorem gives
\[
 V(z)=\overline{V(1/\overline z)}
\]
on that initial annulus. Put
\[
 B^\sharp=\{1/\overline b:b\in B\},\qquad
 \Omega^\sharp=\{z:0<|z|<r_*^{-1}\}.
\]
Reflection is a homeomorphism between $\Omega$ and $\Omega^\sharp$,
so $B^\sharp$ is relatively closed and discrete in the punctured
domain and avoids a neighborhood of the unit circle. Reflecting a
path and its successive local holomorphic germs gives continuation
in $\Omega^\sharp\setminus B^\sharp$, starting with the original germ
by the displayed identity. No reality condition on an individually
continued branch is imposed. If $B$ is unbounded, $B^\sharp$ may
accumulate at zero; no extension through zero follows.

The exterior proof used only the four preceding lemmas and the common
physical germ; reality is used here for reflection. The constructed
set $B$ is a possible obstruction set, not a classification of genuine
singularities. No assertion of finiteness of $B$, regularity at
infinity, single-valuedness, finite monodromy, or meromorphicity is
made at this stage.
\end{proof}

\section{From a discrete continuation set to finitely many intrinsic obstructions}
\label{cont:sec:task6-finite-bad-set}

Write $D=\{z\in\mathbb C:|z|>1\}$. Suppose that a physical boundary germ
$V$ is holomorphic in an exterior collar of the unit circle and admits
analytic continuation along every path in $D\setminus \mathfrak E$, where $\mathfrak E$ is
relatively closed and discrete in $D$. Fix an initial point in the known
collar, outside $\mathfrak E$. An \emph{accessible germ} means any germ reached
from this initial germ by such a path. Different continuation paths
may give different accessible germs.
A point of $D\setminus \mathfrak E$ is called \emph{ordinary}.

We first remove apparent points without changing the accessible germs at
ordinary points, and then realize any prescribed source germ in a
continued rational trace. Under the standing continuation hypotheses and
an exception-free exterior collar, one odd-denominator trace makes the
intrinsic bad set finite: its companion maps have no poles at the two
source poles and avoid the provisional exceptions at sufficiently small
nonzero inverse coordinates. The final corollary preserves this same
finite set and its nonextendible witnesses on an enlarged domain containing
the unit-circle collar; regularity at infinity and elimination of
monodromy remain later tasks.

\begin{lemma}[Intrinsic obstructions and removal of apparent points]
\label{cont:lem:task6-intrinsic-bad-set}
For $b\in \mathfrak E$, choose a disk $D_b\Subset D$ centered at $b$ and
meeting $\mathfrak E$ only at $b$.
Call $b$ \emph{good} if every accessible germ at points of
$D_b\setminus\{b\}$ extends as a holomorphic function on $D_b$.
Let $B$ be the set of points of $\mathfrak E$ which are not good.
This definition is independent of the isolating disk.
The set $B$ is relatively closed and discrete in $D$, and the fixed
initial germ admits analytic continuation along every path in
$D\setminus B$ starting at its base point.
At every ordinary point, the germs reachable through $D\setminus B$
are exactly those reachable through $D\setminus \mathfrak E$.
\end{lemma}
\begin{proof}
Put $D_0=D\setminus\mathfrak E$, let $z_0\in D_0$ be the fixed
initial point, and denote its initial germ by $\mathfrak f_0$. For
$q\in D_0$, write $\mathscr G_0(q)$ for the set of germs reachable
from $\mathfrak f_0$ along paths in $D_0$, and write $[F]_q$ for the
germ of a holomorphic function $F$ at $q$. Until the enlarged-domain
continuation has been constructed, accessibility always refers to
these original sets $\mathscr G_0(q)$. All paths are continuous maps
from compact intervals; no rectifiability is required.

\smallskip\noindent
\emph{Step 1: transport accessible germs and test goodness at one point.}
If $\mathfrak f\in\mathscr G_0(q)$ is represented by a path $\alpha$
from $z_0$ to $q$, then it continues along every further path $\beta$
in $D_0$ starting at $q$: apply the original all-path hypothesis to
$\alpha*\beta$. Here $*$ means concatenation in the displayed order.
Continuation along $\beta$ followed by its reverse $\beta^{-1}$
recovers $\mathfrak f$, by uniqueness of continuation along overlapping
holomorphic elements. We also use homotopy invariance with endpoints
fixed, in the finite-chain form explained in the proof of
Lemma~\ref{cont:lem:countable-gluing}. Its hypothesis is satisfied for
any accessible germ, because a representing path can be prepended to
every path of the homotopy.

For an isolating disk $\Delta=D(b,r)\Subset D$, goodness means
\[
 \text{for every }q\in\Delta\setminus\{b\}
 \text{ and every }\mathfrak f\in\mathscr G_0(q),\quad
 \text{there is }F\in\mathcal O(\Delta)
 \text{ with }[F]_q=\mathfrak f.
\]
The extension is allowed to depend on the germ. In particular, the
condition does not assert that all accessible branches coincide.
Fix any $q_b\in\Delta\setminus\{b\}$. It suffices to test the
condition for all germs in $\mathscr G_0(q_b)$. Indeed, given
$\mathfrak f\in\mathscr G_0(q)$ at any other punctured-disk point,
transport it along a path in $\Delta\setminus\{b\}$ from $q$ to
$q_b$. The resulting germ is in $\mathscr G_0(q_b)$ and hence has a
holomorphic representative $F$ on $\Delta$. Transport back along the
reverse path. The function $F$ supplies that continuation, so
$[F]_q=\mathfrak f$. Necessity of the one-point test is immediate.
The test point can thus be fixed, but the quantifier over all its
originally accessible germs cannot be dropped.

\smallskip\noindent
\emph{Step 2: prove independence of the isolating disk.}
It suffices to compare two concentric isolating disks
$\Delta_{\rm s}\subset\Delta_{\rm l}$ about $b$. Goodness on the
larger disk implies goodness on the smaller disk by restriction.
Conversely, assume goodness on $\Delta_{\rm s}$. Take any
$q\in\Delta_{\rm l}\setminus\{b\}$ and any
$\mathfrak f\in\mathscr G_0(q)$. Choose
$a\in\Delta_{\rm s}\setminus\{b\}$ and a path $\alpha$ from $q$ to
$a$ in $\Delta_{\rm l}\setminus\{b\}$. Its transported germ
$\mathfrak h\in\mathscr G_0(a)$ has a representative
$F_{\rm s}\in\mathcal O(\Delta_{\rm s})$.

Let $\ell$ be a loop based at $a$ in the smaller punctured disk with
winding number one about $b$. Since $F_{\rm s}$ is holomorphic on the
whole smaller disk, continuation around $\ell$ fixes $\mathfrak h$.
The loop $\alpha*\ell*\alpha^{-1}$ therefore fixes $\mathfrak f$.
Its homotopy class generates the fundamental group of the larger
punctured disk. Every loop there is homotopic, with base point fixed,
to an integer power of this generator, so Step~1 shows that every
loop fixes $\mathfrak f$. Continuation of this particular germ is
therefore path independent on $\Delta_{\rm l}\setminus\{b\}$ and
defines a single-valued holomorphic function
$F_{\rm l}^{\times}$ on that punctured disk.

The functions $F_{\rm l}^{\times}$ and $F_{\rm s}$ agree near $a$,
and hence on the connected smaller punctured disk by the identity
theorem. They patch to a holomorphic function on $\Delta_{\rm l}$:
near $b$ use $F_{\rm s}$, and away from $b$ use
$F_{\rm l}^{\times}$. Its germ at $q$ is $\mathfrak f$. Since that
germ and point were arbitrary, $b$ is good on the larger disk.
This proves disk independence. The punctured disk was not assumed
simply connected: triviality of its generator's action was proved
first. The extension on the smaller disk is then needed to fill the
center; trivial local monodromy alone would not exclude a pole.

\smallskip\noindent
\emph{Step 3: verify the topology of the remaining exceptional set.}
Let $G=\mathfrak E\setminus B$ be the good points and put
\[
 D_1=D\setminus B=D_0\cup G.
\]
Each point of $B$ is isolated because $B\subset\mathfrak E$.
For relative closedness, take $z\in D\setminus B$. If
$z\notin\mathfrak E$, relative closedness of $\mathfrak E$ supplies
a neighborhood in $D$ avoiding $\mathfrak E$, and hence $B$. If
$z\in G$, an isolating neighborhood meets $\mathfrak E$ only at $z$;
since $z\notin B$, it too avoids $B$. Thus $D_1$ is open and $B$ is
relatively closed and discrete in $D$. No finiteness of $G$ or $B$
is needed.

\smallskip\noindent
\emph{Step 4: construct disk representatives while retaining original accessibility.}
Fix a path $\gamma:[0,1]\to D_1$ with $\gamma(0)=z_0$.
Cover its image by disks of two types: disks whose closures lie in
$D_0$, and isolating disks centered at good points with closures in
$D$. A good-point disk contains no member of $\mathfrak E$ other
than its center, so both types lie in $D_1$. Compactness and a
Lebesgue number for the pulled-back cover give disks $Q_1,\ldots,Q_N$
and a subdivision
\[
 0=t_0<t_1<\cdots<t_N=1,\qquad
 \gamma([t_{j-1},t_j])\subset Q_j\quad(1\le j\le N).
\]
For $1\le j<N$, put $c_j=\gamma(t_j)$ and let $O_j$ be the connected
component of $Q_j\cap Q_{j+1}$ containing $c_j$. Choose an auxiliary
point $q_j\in O_j\setminus\mathfrak E$. Such a choice is possible
because $O_j$ is open and a discrete set contains no open disk.
Keep $q_0=z_0$. If $\gamma(1)$ is ordinary, put $q_N=\gamma(1)$;
otherwise choose any $q_N\in Q_N\setminus\mathfrak E$. Thus every
$q_j$ is ordinary and $q_{j-1},q_j\in Q_j$. An original switch $c_j$
may be good; its auxiliary point is always ordinary.

We construct, simultaneously, a function $F_j\in\mathcal O(Q_j)$
and a path $\alpha_j$ in $D_0$ from $z_0$ to $q_j$, such that
continuation along $\alpha_j$ gives the germ $[F_j]_{q_j}$.
Start with the constant path $\alpha_0$ and the germ
$\mathfrak f_0$ at $q_0$. Inductively, let
$\mathfrak f_{j-1}\in\mathscr G_0(q_{j-1})$ be the germ represented
by $\alpha_{j-1}$.

If $Q_j$ avoids $\mathfrak E$, Step~1 supplies continuation of this
germ along every path in $Q_j$. Since the disk is simply connected,
homotopy invariance makes that continuation path independent and
defines $F_j\in\mathcal O(Q_j)$. If $Q_j$ isolates a good point,
its definition supplies such an $F_j$ directly, because the incoming
germ is originally accessible and based at an ordinary point of the
disk. In both cases,
\[
 [F_j]_{q_{j-1}}=\mathfrak f_{j-1}.
\]
Join $q_{j-1}$ to $q_j$ by a path $\beta_j$ in
$Q_j\setminus\mathfrak E$; a disk with at most one point removed is
path connected. Define
\[
 \alpha_j=\alpha_{j-1}*\beta_j,\qquad
 \mathfrak f_j=[F_j]_{q_j}.
\]
The path $\alpha_j$ lies entirely in $D_0$, and $F_j$ supplies its
terminal continuation along $\beta_j$. Hence
$\mathfrak f_j\in\mathscr G_0(q_j)$, closing the induction.

This original-domain representing path is the essential invariant:
we never apply the good-point condition to a new germ merely because
it has been reached through the enlarged domain. Its original
accessibility is proved before the next good-point disk is used.
Nor do these detours replace the requested path; they represent the
germs at the auxiliary ordinary points.

\smallskip\noindent
\emph{Step 5: match along the original path and include its endpoint.}
At an auxiliary switch, the construction gives equality of germs
\[
 [F_j]_{q_j}=\mathfrak f_j=[F_{j+1}]_{q_j}.
\]
The identity theorem makes $F_j$ and $F_{j+1}$ agree on $O_j$, which
also contains the original switch $c_j$. Thus they agree as germs
at $c_j$, not just in value at $q_j$. The functions $F_1,\ldots,F_N$
therefore form a chain of holomorphic elements continuing
$\mathfrak f_0$ along the original path $\gamma$.

If the endpoint is ordinary, it equals $q_N$, and the terminal germ
$[F_N]_{q_N}$ is also reached along the constructed path
$\alpha_N\subset D_0$. If the endpoint is a good point, the function
$F_N$ is holomorphic there and directly supplies its endpoint germ.
There is no claim of a path in $D_0$ ending at a point of
$\mathfrak E$.

This is a finite disk-chain construction for each continuous path.
Its validity does not depend on finite length or on a finite number
of visits to a good point: the representative on a disk covers the
entire subpath lying there. No integral along the original path and
no ordering of all good points in the domain is required. We have
proved all-path continuation on $D_1$ using only the original
all-path hypothesis on $D_0$ and the definition of goodness.

\smallskip\noindent
\emph{Step 6: identify the ordinary-point germ sets and retain all branches.}
Now let $\mathscr G_1(q)$ denote the germs reachable from
$\mathfrak f_0$ along paths in $D_1$. For $q\in D_0$, the inclusion
$\mathscr G_0(q)\subset\mathscr G_1(q)$ follows from $D_0\subset D_1$.
Conversely, represent any germ in $\mathscr G_1(q)$ by a path in
$D_1$ and apply Steps~4--5. They produce both a continuation along
that path and a path in $D_0$ with the same terminal germ. Uniqueness
of continuation along the fixed original path identifies the former
with the specified germ. Thus
\[
 \boxed{\mathscr G_1(q)=\mathscr G_0(q),\qquad q\in D\setminus\mathfrak E.}
\]
At a point $b\in B$, failure of the good-point test still supplies an
originally accessible punctured-disk germ that does not extend
holomorphically across $b$. The equality just proved preserves its
accessibility after removal of the good points. It does not assert
that every branch is singular at a bad point, or that different
holomorphic extensions at a good point coincide. The result removes
all apparent points simultaneously while preserving all ordinary-point
accessible germs. Finiteness of $B$, a classification of its
singularities, and trivial or finite global monodromy are not asserted.
\end{proof}

Figure~\ref{cont:fig:apparent-point-removal} illustrates the local extension and the original-domain accessibility retained in this construction.
\begin{figure}[htbp]
\centering
\resizebox{\linewidth}{!}{%
\begin{tikzpicture}[x=1cm,y=1cm,font=\small,
 arr/.style={-{Stealth[length=1.7mm]},line width=.75pt},
 every node/.style={inner sep=2pt}]
\path[use as bounding box] (-.10,-.10) rectangle (15.75,5.78);
\node at (3.55,5.43) {(a) A full-disk representative at a good point};
\node at (11.78,5.43) {(b) The condition concerns every accessible germ};
\fill[black!4] (3.45,3.24) circle[radius=1.48];
\draw[line width=.6pt] (3.45,3.24) circle[radius=1.48];
\node[anchor=south west] at (4.45,4.26) {$\Delta_b$};
\node at (3.45,4.16) {$F_{\mathfrak f}\in\mathcal O(\Delta_b)$};
\coordinate (qin) at (2.16,3.24);
\coordinate (qout) at (4.74,3.24);
\draw[arr] (.23,3.24)--(qin);
\node[font=\footnotesize,align=center] at (1.05,3.83) {original history\\$\alpha\subset D_0$};
\draw[line width=.85pt,postaction={decorate},decoration={markings,
 mark=at position .25 with {\arrow{Stealth}},
 mark=at position .81 with {\arrow{Stealth}}}] (qin)--(qout);
\draw[dashed,line width=.8pt,postaction={decorate},decoration={markings,
 mark=at position .55 with {\arrow{Stealth}}}]
 (qin)..controls (2.40,1.76) and (4.50,1.76)..(qout);
\fill (qin) circle[radius=.04];
\fill (qout) circle[radius=.04];
\node[above=4pt] at (qin) {$q_-$};
\node[above=4pt] at (qout) {$q_+$};
\draw[fill=white,line width=.75pt] (3.45,3.24) circle[radius=.06];
\node[above=5pt,fill=white] at (3.45,3.24) {$b$};
\node[above=5pt] at (4.12,3.24) {$\gamma$};
\node at (3.45,2.48) {$\beta$};
\node[align=center,font=\footnotesize] at (3.60,1.20)
 {$[F_{\mathfrak f}]_{q_-}\in\mathscr G_0(q_-)$\\[3pt]
  $\alpha*\beta$ reaches $[F_{\mathfrak f}]_{q_+}$ in $D_0$};
\node[font=\footnotesize,align=center] at (3.60,.38)
 {The detour represents the germ;\\the requested path may pass through $b$.};
\draw[black!30,line width=.4pt] (7.70,.13)--(7.70,5.0);
\node[font=\footnotesize] at (11.75,4.84) {$q\in\Delta_b\setminus\{b\},\qquad \mathfrak f,\mathfrak g\in\mathscr G_0(q)$};
\node at (9.10,3.90) {$\mathfrak f$};
\draw[arr] (9.53,3.90)--(10.75,3.90);
\node[draw,line width=.55pt,minimum width=3.9cm,minimum height=.67cm]
 at (12.91,3.90) {$F_{\mathfrak f}\in\mathcal O(\Delta_b)$};
\node at (9.10,2.74) {$\mathfrak g$};
\draw[arr] (9.53,2.74)--(10.75,2.74);
\node[draw,line width=.55pt,minimum width=3.9cm,minimum height=.67cm]
 at (12.91,2.74) {$F_{\mathfrak g}\in\mathcal O(\Delta_b)$};
\node[font=\footnotesize] at (11.90,1.95) {The two representatives may differ.};
\node[font=\footnotesize] at (11.80,1.08) {After removal, at every ordinary point:};
\node at (11.80,.48) {$\mathscr G_1(q)=\mathscr G_0(q)$};
\end{tikzpicture}
}
\caption{Apparent-point removal in Lemma~\ref{cont:lem:task6-intrinsic-bad-set}
(schematic). In a good-point disk, the full-disk representative continues
the incoming germ along the requested path $\gamma$, even through $b$.
The detour $\beta$ supplies an original-domain representing path for its
outgoing germ; it does not replace the requested path.
Every originally accessible germ must admit a full-disk representative.
The two representatives shown may differ, and the drawing imposes no
bound on the number of branches. Here $D_0=D\setminus\mathfrak E$,
$D_1=D\setminus B$, and $\mathscr G_i$ denotes accessibility through $D_i$.}
\label{cont:fig:apparent-point-removal}
\end{figure}

We also use an elementary fact about paths. In an open Riemann surface,
a path whose endpoints avoid a relatively closed countable set can be
perturbed, relative to its endpoints, to avoid that set, with the
perturbation homotopic to the original path in any prescribed open
neighborhood of its image. To see this, cover the compact path by a
finite chain of coordinate disks in that neighborhood. Choose switching
points outside the countable set in small connected neighborhoods of
the original switches within consecutive overlaps, and join them by
two line segments in each disk. In each coordinate disk, choose the intermediate point outside the
countable union of lines through an endpoint and an excluded point.
Move each switch to its replacement within the chosen overlap
neighborhood. With these endpoint motions included, the original
subpath and its replacement have the same endpoints in a disk and
are homotopic there. The shared endpoint motions cancel, giving a
homotopy relative to the original endpoints.

\Needspace{26\baselineskip}
\begin{lemma}[Realization of every accessible germ in a continued trace]
\label{cont:lem:task6-trace-realization}
Fix a reduced rational rotation $\rho=p/d\in(0,1/2)$ with odd denominator
$d$, and suppose its rational
trace is constant for the initial germ. Use its Jacobi parameter $m$,
put $h=4Kp/d$, and set
\[
 T_j(u)=\operatorname{sn}(u+jh\mid m)
            -i\operatorname{cn}(u+jh\mid m),\qquad 0\le j<d.
\]
Work in the lower real-period cylinder
\[
 \mathcal C_-=\{-2K'<\operatorname{Im}u<0\}/(4K\mathbb Z).
\]
Every accessible source germ at an ordinary physical point can, after
continuation to an arbitrarily close generic point if necessary, be
realized as the unshifted component of an analytically continued identity
\begin{equation}\label{cont:eq:task6-continued-trace}
 \sum_{j=0}^{d-1}V_j(T_j(u))=C_\rho.
\end{equation}
Here each $V_j$ is an accessible continuation of the original physical
germ. The unshifted lift ends at one of the two Jacobi preimages of the
possibly moved physical endpoint; neither preimage is prescribed.
The lifted trace path can avoid every shifted pullback of $\mathfrak E$ and
every shifted Jacobi pole while preserving the source germ,
transported by the stated local continuation if the endpoint moves.

Consequently, if an accessible source germ near a noncritical value
$b\in \mathfrak E$ does not extend holomorphically across $b$, then at least one
preimage $u_b$ of $b$ satisfies $T_j(u_b)\in \mathfrak E\cup\{\infty\}$ for some
$1\le j<d$. This conclusion applies separately to every
accessible source germ, including branches with nontrivial monodromy.
\end{lemma}
\begin{proof}
Put $D_0=D\setminus\mathfrak E$, where $D=\{z:|z|>1\}$, and let
$z_*\in D_0$ be the fixed initial point. Throughout this proof,
accessibility refers to continuation from the original germ along paths
in $D_0$. Every accessible germ continues along every further path in
$D_0$: concatenate that path with one representing the germ from $z_*$.
We use uniqueness along a fixed path and homotopy invariance as in
Lemma~\ref{cont:lem:countable-gluing}. All paths have compact parameter
intervals and are merely continuous. The rational caustic is fixed, as
are $m$, $K=K(m)$, and $K'=K(1-m)$; the prime in $K'$ denotes the
complementary integral, not differentiation.

\smallskip\noindent
\emph{Step 1: separate the geometric covering from the continuation domain.}
By Lemma~\ref{arith:jacobi}, the map
$T_0:\mathcal C_-\to D\cup\{\infty\}$ is proper of degree two.
Each real translation $u\mapsto u+jh$ is an automorphism of the
real-period cylinder, so the same is true of every $T_j$. In particular,
every finite value of every $T_j$ on $\mathcal C_-$ lies in $D$.
The derivative identities give
\[
 T_0'(u)=i\dn(u\mid m)T_0(u).
\]
Its poles are $-iK'$ and $2K-iK'$ modulo $4K$, and its finite critical
points in the lower cylinder are $K-iK'$ and $3K-iK'$. Their values are
\[
 C_{\mathrm{crit}}
 =\left\{\frac{1+\sqrt{1-m}}{\sqrt m},
              -\frac{1+\sqrt{1-m}}{\sqrt m}\right\}.
\]
These locations and values follow from the Jacobi identities used in
Lemmas~\ref{arith:jacobi} and \ref{arith:focal}; in particular,
$\mathcal P_m(t-iK')=(1+\dn t)/(\sqrt m\,\sn t)$, evaluated at $t=K$,
gives the positive critical value. The poles of $T_j$ are correspondingly
$-jh-iK'$ and $2K-jh-iK'$ modulo $4K$.

Define the unshifted admissible region
\[
 A_0=T_0^{-1}(D_0)
 =\mathcal C_-\setminus
   \bigl(T_0^{-1}(\mathfrak E)\cup\operatorname{Poles}(T_0)\bigr).
\]
It is open. The restriction to the full inverse image
\[
 T_0:T_0^{-1}(D_0\setminus C_{\mathrm{crit}})
           \longrightarrow D_0\setminus C_{\mathrm{crit}}
\]
is a proper unbranched two-sheeted covering. Indeed, every target point
has two regular preimages; their local inverse neighborhoods, together
with properness, give an evenly covered target disk. Deleting infinity
from the target deletes its pole preimages from the source. This
covering will be used to lift a physical path, whereas $A_0$ will be the
larger region in which to preserve its unshifted continuation germ.

\smallskip\noindent
\emph{Step 2: form the full shifted excluded set.}
Set
\[
 F=\bigcup_{j=0}^{d-1}T_j^{-1}(\mathfrak E)
       \ \cup\ \bigcup_{j=0}^{d-1}\operatorname{Poles}(T_j).
\]
The set $\mathfrak E$ is countable: a compact exhaustion of $D$ meets
this relatively closed discrete set in finite sets. Each $T_j$ has at
most two preimages of a value, and finitely many poles on the cylinder.
Consequently $F$ is countable.

It is also relatively closed in $\mathcal C_-$. To check this, take a
convergent sequence in one of the finitely many sets
$T_j^{-1}(\mathfrak E)\cup\operatorname{Poles}(T_j)$, with limit in the
cylinder. A pole limit is already included. At any other limit point,
$T_j$ is holomorphic with finite value in $D$, so relative closedness of
$\mathfrak E$ puts the limiting value back in $\mathfrak E$.
Equivalently, $\mathfrak E\cup\{\infty\}$ is relatively closed in the
spherical exterior and $F$ is the union of its inverse images.
The poles must be included: inverse images of exceptional values tending
to infinity can accumulate at them. Thus $F$ need not be discrete, but
$F\cap A_0$ is a relatively closed countable subset of the open Riemann
surface $A_0$. The path-perturbation fact preceding the lemma applies
in this precise setting.

\smallskip\noindent
\emph{Step 3: start with the initial trace and rebase the prescribed germ.}
For real $u$, the relation $u=K+2Kx/\pi$ identifies the shifts by $jh$
with action shifts by $2\pi jp/d$. Trace constancy therefore gives
\[
 \sum_{j=0}^{d-1}V(T_j(u))=C_\rho
\]
with the original branches. The finitely many maps $T_j$ send a common
complex neighborhood of the real circle into the original holomorphy
annulus of $V$. The identity theorem extends the displayed equality to
that neighborhood.

Choose $u_0$ in its lower-cylinder part, sufficiently close to the real
circle that every $T_j(u_0)$ lies in the original exterior collar, and
such that
\[
 u_0\notin F\cup T_0^{-1}(C_{\mathrm{crit}}).
\]
Such a point exists because only a countable set and a finite set have
been removed from an open collar. All starting values are ordinary and
finite, and the unshifted one is noncritical. For each $j$, choose a
path $\alpha_j$ from $z_*$ to $T_j(u_0)$ in the known physical collar
avoiding $\mathfrak E$. These paths exist by the preceding path fact,
starting with any collar paths. The original holomorphic function $V$
supplies their terminal germs, so each initial trace summand is
accessible in $D_0$.

Let $\mathfrak f$ be the prescribed accessible germ at an ordinary point
$q$, and let $\beta_*$ be a path from $z_*$ to $q$ in $D_0$ representing
it. Starting at $a=T_0(u_0)$ with the original germ, use
\[
 \beta=\alpha_0^{-1}*\beta_*.
\]
Here $*$ is concatenation in the displayed order and $\alpha_0^{-1}$
is the reversed path. Uniqueness along the collar path and its reverse
shows that $\beta$ reaches exactly $\mathfrak f$, not an independently
selected branch at $q$.

\smallskip\noindent
\emph{Step 4: avoid critical values and obtain a geometric lift.}
If $q\in C_{\mathrm{crit}}$, choose a disk $W\Subset D_0$ on which
$\mathfrak f$ has a holomorphic representative, and continue that
representative to an arbitrarily close point
$q_1\in W\setminus C_{\mathrm{crit}}$. Append a path in $W$ to $\beta$.
If $q$ is already noncritical, set $q_1=q$ and do not move it.
Denote the resulting physical path by $\beta_1$. Both of its endpoints
avoid $C_{\mathrm{crit}}$.

Apply the path-perturbation fact in $D_0$ to avoid the finite set
$C_{\mathrm{crit}}\cap D_0$, with endpoints fixed and with the homotopy
remaining in $D_0$. The resulting path $\beta_2$ is contained in
$D_0\setminus C_{\mathrm{crit}}$. Homotopy invariance, justified by the
all-path hypothesis on $D_0$, preserves the terminal source germ.
The covering in Step~1 gives a unique lift
\[
 \lambda:[0,1]\to T_0^{-1}(D_0\setminus C_{\mathrm{crit}}),
 \qquad \lambda(0)=u_0,\qquad T_0\circ\lambda=\beta_2.
\]
Its endpoint $u_1$ is one of the two preimages of $q_1$. No choice of
that terminal sheet has been imposed. At this stage the unshifted
projection avoids $\mathfrak E$, but the companion projections may
still encounter exceptional points or poles.

\smallskip\noindent
\emph{Step 5: avoid all shifted obstructions without changing the source branch.}
If $u_1\notin F$, retain the endpoint. Otherwise choose a small local
inverse neighborhood of $u_1$ contained in $A_0$, with its $T_0$-image
inside a representing neighborhood of the source germ at $q_1$ and
avoiding $C_{\mathrm{crit}}$. Choose an arbitrarily close point $u_2$
outside $F$ in that neighborhood and append a short local path from
$u_1$ to $u_2$. Countability of $F$ makes this choice possible. The
projected terminal segment continues the already specified local
representative of the source germ. When an earlier endpoint adjustment
was necessary, both adjustments can be confined to an arbitrarily small
initially chosen representing neighborhood of $\mathfrak f$ at $q$.
Write $\lambda_+$ for the resulting path, or simply $\lambda$ if no
second adjustment is needed. Its endpoints are outside $F$ and its
entire image lies in $A_0$.

Apply the same path-perturbation fact, now in the open Riemann surface
$A_0$, to the relatively closed countable set $F\cap A_0$. It gives a
path $\widehat\lambda$ avoiding $F$, homotopic to $\lambda_+$ with
endpoints fixed through a homotopy entirely in $A_0$. Composing that
homotopy with $T_0$ produces an endpoint-fixed homotopy in $D_0$,
because $T_0$ is finite and avoids $\mathfrak E$ throughout $A_0$.
Hence $T_0\circ\widehat\lambda$ reaches the prescribed source germ,
or exactly its specified local continuation when an endpoint was moved.
This homotopy need not avoid the geometric critical values: their
avoidance was needed for the first covering lift, while preservation of
the analytic germ requires only a homotopy in $D_0$.

By construction, all paths
\[
 \gamma_j:=T_j\circ\widehat\lambda,\qquad 0\le j<d,
\]
are finite-valued and lie in $D_0$. The poles in $F$ ensure finiteness,
and the lower-cylinder geometry ensures that every finite value remains
in $D$. In particular, the second perturbation has removed all shifted
obstructions without changing the prescribed unshifted branch.

\smallskip\noindent
\emph{Step 6: continue the complete trace and record exact endpoint control.}
For each $j$, continue the original germ at $T_j(u_0)$ along $\gamma_j$.
The path $\alpha_j*\gamma_j$ lies in $D_0$, so every resulting companion
germ is accessible from the fixed initial germ. There are finitely many
summands. Their continuations admit a common finite subdivision of the
parameter interval, with a physical representing neighborhood for each
summand at each step. Near the corresponding point of the lifted path,
choose a source neighborhood mapped by every $T_j$ into its physical
representing neighborhood. Each composition $V_j\circ T_j$ is
holomorphic there. No local inverse of a companion map is needed.

On consecutive source neighborhoods the summands agree as germs. The
initial constant identity from Step~3 therefore continues along the
whole chain, by the identity theorem, and gives
\eqref{cont:eq:task6-continued-trace} near the final lifted endpoint.
Its unshifted component is the prescribed source germ, with only the
permitted local continuation if the endpoint changed. Companion germs
are determined by their own projected histories; they are not separately
prescribed, and they need not coincide when their physical endpoints
coincide. Only the continued sum is asserted to remain $C_\rho$.

There is an exact-endpoint version of this construction. If $q$ is
noncritical and the endpoint of the initial lift lies outside $F$,
Steps~4--5 make no endpoint adjustment. A sufficient condition is that
both preimages of $q$ lie outside $F$. In that case any prescribed
accessible germ at $q$ is realized at $q$ itself, irrespective of which
sheet the lift reaches. This observation will be used for the final
assertion and in the subsequent applications.

\smallskip\noindent
\emph{Step 7: prove the collision criterion on both sheets.}
We prove the contrapositive. Let $b\in\mathfrak E$ be noncritical for
$T_0$, with preimages $u_b^+,u_b^-$, and suppose that
\[
 T_j(u_b^\nu)\in D_0,
 \qquad \nu\in\{+,-\},\quad 1\le j<d.
\]
Thus no companion value at either preimage is exceptional or infinite.
Choose an isolating disk $\Delta\Subset D$ about $b$ avoiding the
critical values and meeting $\mathfrak E$ only at $b$. Properness and
unramifiedness allow it to be chosen so that its full preimage is
\[
 T_0^{-1}(\Delta)=O_+\sqcup O_-,\qquad
 \psi_\nu=(T_0|_{O_\nu})^{-1}:\Delta\to O_\nu,
\]
with both restrictions biholomorphic. For each $j\ne0$ and each sign
$\nu$, choose a disk $W_{j,\nu}\Subset D_0$ about $T_j(u_b^\nu)$.
Shrinking $\Delta$ handles this finite family simultaneously and gives
\[
 T_j(O_\nu)\subset W_{j,\nu},
 \qquad 1\le j<d,\quad \nu\in\{+,-\}.
\]
In particular all these shifted maps are holomorphic and finite on the
inverse neighborhoods. The neighborhood choices depend on the fixed
maps and $\mathfrak E$, not on the branch to be realized.

Take any accessible source germ $\mathfrak f$ at any
$q\in\Delta\setminus\{b\}$. Both preimages of $q$ lie outside $F$:
their unshifted value is ordinary and all companion values lie in the
ordinary disks just chosen. Step~6 therefore realizes $\mathfrak f$ at
$q$ without moving it. Let $O_\nu$ be the inverse neighborhood reached
by the realization. Each companion germ is accessible at a point of
$W_{j,\nu}$. It has continuation along every path in that disk by the
standing all-path hypothesis and concatenation. Since the disk is
simply connected, the monodromy theorem supplies a holomorphic function
$\mathcal V_{j,\nu}\in\mathcal O(W_{j,\nu})$ representing that particular
germ. Different accessible germs may give different such functions.

The physical extension is now explicit:
\[
 \mathcal V_{\mathrm{ext}}(z)
   =C_\rho-\sum_{j=1}^{d-1}
        \mathcal V_{j,\nu}\bigl(T_j(\psi_\nu(z))\bigr),
 \qquad z\in\Delta.
\]
It is holomorphic on all of $\Delta$. The continued trace identifies its
germ at $q$ with the prescribed $\mathfrak f$, so that germ extends
holomorphically across $b$. Since $q$ and its accessible germ were
arbitrary, $b$ is good in the sense of
Lemma~\ref{cont:lem:task6-intrinsic-bad-set}. The disk-independence part
of that lemma covers a source germ initially described in a different
isolating disk. This contradicts any prescribed nonextendible germ and
proves that at least one of $u_b^+,u_b^-$ must have a companion collision.

Noncriticality is used to obtain the physical inverses $\psi_\nu$;
holomorphy in a ramified source coordinate alone would not give this
physical extension. The argument checks both sheets because the
realization preserves the source branch without prescribing the sheet.
It asserts neither a collision at every preimage nor a sufficient
condition for nonextendibility, and imposes no restriction on global
monodromy.
\end{proof}

\begin{theorem}[Finiteness from one odd-denominator trace]
\label{cont:thm:task6-finite-bad-set}
Under the standing continuation hypotheses of this section, assume additionally
that $\mathfrak E$ is disjoint from some exterior collar
$1<|z|<1+\varepsilon_0$, and that the trace of one reduced rational
rotation $p/d\in(0,1/2)$ with odd denominator $d$ is constant.
Then the intrinsic bad set $B$ from
Lemma~\ref{cont:lem:task6-intrinsic-bad-set} is finite.
No global single-valuedness, finite-monodromy, meromorphicity, or a priori
classification of the singularities is assumed.
\end{theorem}
\begin{proof}
Put $D=\{z\in\C:|z|>1\}$ and $D_0=D\setminus\mathfrak E$.
We use the standing all-path continuation hypothesis on $D_0$;
throughout this proof, an accessible germ is one reached from the fixed
initial physical germ through that original domain. By
Lemma~\ref{cont:lem:task6-intrinsic-bad-set}, a point of $\mathfrak E$ is
good precisely when every such germ in a punctured isolating disk has a
holomorphic representative on the full disk. The representative may
depend on the germ. We will choose one $R<\infty$, independently of the
accessible branch, and prove that every $b\in\mathfrak E$ with $|b|>R$
is good.

Fix the odd-denominator caustic in the statement. Use the notation of
Lemma~\ref{cont:lem:task6-trace-realization}:
\[
 K=K(m),\qquad K'=K(1-m),\qquad h=\frac{4Kp}{d},\qquad
 T_j(u)=\mathcal P_m(u+jh),\quad 0\le j<d.
\]
The parameter $m$ is fixed, and $K'$ denotes the complementary integral,
not a derivative. All source points below lie in the lower
real-period cylinder $\mathcal C_-$ from that lemma.

\smallskip\noindent
\emph{Step 1: use oddness to separate the source poles from all companion poles.}
By Lemma~\ref{arith:jacobi}, the two simple poles of $T_0$ are
\[
 a_1=-iK',\qquad a_2=2K-iK'\pmod{4K}.
\]
For $T_j$ to have a pole at either $a_\nu$, the shift $jh$ must take
one source pole to one of these two poles modulo $4K$. Their real
parts differ by $2K$, so this is equivalent to
\[
 jh\in2K\Z
 \quad\Longleftrightarrow\quad d\mid2jp.
\]
Reducedness gives $d\mid2j$, and oddness of $d$ then gives $d\mid j$.
This is impossible for $1\le j<d$. Thus
\[
 T_j(a_\nu)\ne\infty,
 \qquad 1\le j<d,\quad \nu=1,2.
\]
Real translation preserves $\mathcal C_-$, whose finite images under
$\mathcal P_m$ all lie in $D$. Consequently
$T_j(a_\nu)\in D$ for all these pairs. This is the essential use of
oddness: for even reduced $d$, the index $j=d/2$ would shift a source
pole to another source pole.

\smallskip\noindent
\emph{Step 2: construct both inverse coordinates and their companion maps.}
At each simple pole $a_\nu$, the meromorphic function $T_0$ has a local
expansion
\[
 T_0(u)=\frac{A_\nu}{u-a_\nu}+O(1),\qquad A_\nu\ne0.
\]
Hence $\xi=1/T_0(u)$ extends holomorphically across $a_\nu$, with a
simple zero there and derivative $A_\nu^{-1}\ne0$. Choose disjoint
source neighborhoods at the two poles. The holomorphic inverse-function
theorem, followed by restriction to one common small disk, gives
biholomorphic inverse maps
\[
 u_\nu:\{|\xi|<\varepsilon_{\rm inv}\}\longrightarrow\mathcal C_-,
 \qquad u_\nu(0)=a_\nu,\qquad
 T_0(u_\nu(\xi))=\xi^{-1}\quad(\xi\ne0),
\]
whose images are disjoint. Define the $2(d-1)$ companion maps
\[
 z_{j,\nu}(\xi)=T_j(u_\nu(\xi)),\qquad
 c_{j,\nu}:=z_{j,\nu}(0)=T_j(a_\nu),
 \quad 1\le j<d,\quad \nu=1,2.
\]
Step~1 shows that these maps are holomorphic and finite at zero, with
$c_{j,\nu}\in D$. Shrink the common disk so that all companion maps are
holomorphic and finite throughout it. Each is nonconstant: constancy
of $T_j\circ u_\nu$ would make the nonconstant meromorphic function
$T_j$ constant on an open source neighborhood and hence everywhere.
No nonvanishing condition on $z_{j,\nu}'(0)$ is needed.

\smallskip\noindent
\emph{Step 3: obtain one punctured avoidance radius before choosing a branch.}
Fix a pair $(j,\nu)$ and write $f=z_{j,\nu}$ and $c=f(0)\in D$.
If $c\notin\mathfrak E$, relative closedness of $\mathfrak E$ gives a
disk $W\Subset D_0$ about $c$. Continuity puts $f(\xi)$ in $W$ for
all sufficiently small $|\xi|$.

If $c\in\mathfrak E$, choose an isolating disk $W\Subset D$ with
$W\cap\mathfrak E=\{c\}$ and again shrink the parameter disk so that
$f(\xi)\in W$. Nonconstancy gives a Taylor factorization
\[
 f(\xi)-c=\xi^r a(\xi),\qquad r\ge1,\qquad a(0)\ne0.
\]
After another shrink, $a$ has no zero on the parameter disk. Therefore
$f(\xi)\ne c$ for every nonzero $\xi$ there, and these values lie in
$W\setminus\{c\}\subset D_0$.

There are only finitely many pairs $(j,\nu)$. Taking a common smaller
radius $0<\varepsilon<\min\{\varepsilon_{\rm inv},1\}$ proves
\begin{equation}\label{cont:eq:task6-endpoint-avoidance}
 z_{j,\nu}(\xi)\in D\setminus\mathfrak E
 \quad\text{for }0<|\xi|<\varepsilon,
 \qquad 1\le j<d,\quad \nu=1,2.
\end{equation}
The radius depends only on $\mathfrak E$ and the fixed caustic, not on a
branch of $V$. It does not require $c_{j,\nu}\notin\mathfrak E$; the
center $\xi=0$ is deliberately excluded. Only local discreteness near
the finitely many $c_{j,\nu}$ has been used, not finiteness of
$\mathfrak E$.

\smallskip\noindent
\emph{Step 4: control the entire large-value fiber and realize its source germs.}
Choose
\[
 R>\max\left\{\varepsilon^{-1},\ 1+\varepsilon_0,\
                   \frac{1+\sqrt{1-m}}{\sqrt m}\right\}.
\]
The last term bounds the moduli of the two finite critical values of
$T_0$ recalled in Lemma~\ref{cont:lem:task6-trace-realization}.
For every finite $z$ with $|z|>R$, both $u_1(1/z)$ and $u_2(1/z)$ are
regular preimages of $z$, in disjoint pole neighborhoods. Regularity
also follows directly by differentiating
$T_0(u_\nu(\xi))=\xi^{-1}$ at $\xi\ne0$.
The degree of $T_0$ is two, so these are its entire fiber; there are
no additional preimages outside the chosen neighborhoods. For $|z|>R$,
$1\le j<d$, and $\nu=1,2$, we therefore have
\[
 T_0^{-1}(z)=\{u_1(1/z),u_2(1/z)\},\qquad
 T_j(u_\nu(1/z))=z_{j,\nu}(1/z)\in D_0.
\]

Take any ordinary $q\in D_0$ with $|q|>R$ and any accessible source
germ there. Its two possible lifted endpoints avoid every shifted
pullback of $\mathfrak E$ and every shifted pole: the unshifted value is
$q\in D_0$, and all companion values are ordinary by
\eqref{cont:eq:task6-endpoint-avoidance}. The exact-endpoint construction
in Lemma~\ref{cont:lem:task6-trace-realization} therefore realizes the
prescribed germ at $q$ itself. On whichever pole sheet $\nu$ is reached,
the continued trace gives the identity of germs
\[
 V_0(1/\xi)=C_\rho-\sum_{j=1}^{d-1}V_j(z_{j,\nu}(\xi))
 \quad\text{near }\xi_q=1/q.
\]
Here $V_0$ is the prescribed source germ and each companion $V_j$ is
accessible through $D_0$; the companions and the sheet are not
independently prescribed.

For completeness, this identity continues the source along every path
in the full punctured disk $0<|\xi|<\varepsilon$ starting at $\xi_q$.
Indeed, each companion map sends such a path into $D_0$ by
\eqref{cont:eq:task6-endpoint-avoidance}. Prepending a representing path
for that companion germ permits its continuation along the image path.
A common finite chain for the finitely many compositions then continues
their sum and hence the source germ. This is all-path continuation of a
germ, not necessarily a single-valued holomorphic function on the
punctured disk. It may cross points where $1/\xi\in\mathfrak E$,
because the right-hand side uses only ordinary companion values.
Nothing here asserts an extension through $\xi=0$.

\smallskip\noindent
\emph{Step 5: extend every accessible branch across each large finite exception.}
Let $b\in\mathfrak E$ satisfy $|b|>R$. At both of its preimages,
\[
 T_j(u_\nu(1/b))=z_{j,\nu}(1/b)\in D_0,
 \qquad 1\le j<d,\quad \nu=1,2.
\]
Thus neither sheet has a companion collision. This is the contrapositive
case of Lemma~\ref{cont:lem:task6-trace-realization}; we spell out the
full-disk extension to verify the quantifier over all accessible germs.

For each $(j,\nu)$ choose an ordinary disk
$W_{j,\nu}\Subset D_0$ around $z_{j,\nu}(1/b)$. By continuity and the
finiteness of this family, choose one isolating disk
\[
 \Delta_b\Subset\{z:|z|>R\},\qquad
 \Delta_b\cap\mathfrak E=\{b\},
\]
centered at $b$ and small enough that
\[
 z_{j,\nu}(1/z)\in W_{j,\nu}\qquad
 (z\in\Delta_b,\ 1\le j<d,\ \nu=1,2).
\]
These disks are chosen before an accessible branch. Every point of
$\Delta_b$ is noncritical, and its two source inverses are
$z\mapsto u_\nu(1/z)$.

Now take any $q\in\Delta_b\setminus\{b\}$ and any accessible germ
$\mathfrak f$ at $q$. Step~4 realizes it at the exact endpoint $q$, on
one of the sheets $\nu$. Each resulting companion germ is based in
$W_{j,\nu}$ and is accessible through $D_0$. Concatenation with its
representing path gives continuation along every path in that disk.
Since the disk is simply connected, the monodromy theorem supplies a
holomorphic representative
$\mathcal V_{j,\nu}\in\mathcal O(W_{j,\nu})$ of that particular germ.
The representative may depend on $\mathfrak f$ and its realized trace.

Define
\[
 F_{\mathfrak f}(z)=C_\rho-
       \sum_{j=1}^{d-1}\mathcal V_{j,\nu}(z_{j,\nu}(1/z)),
       \qquad z\in\Delta_b.
\]
Every term is holomorphic on the full disk $\Delta_b$, and the continued
trace identifies $[F_{\mathfrak f}]_q$ with $\mathfrak f$. Thus
$\mathfrak f$ extends holomorphically across $b$. Since both $q$ and its
accessible germ were arbitrary, $b$ is good. The independence of the
isolating disk in Lemma~\ref{cont:lem:task6-intrinsic-bad-set} covers
germs initially described in any other isolating disk. The argument works
on either reached sheet and does not identify different accessible
branches with one another.

\smallskip\noindent
\emph{Step 6: confine the intrinsic bad set to a compact annulus.}
Step~5 gives $B\cap\{|z|>R\}=\varnothing$. Since $B\subset\mathfrak E$
and the latter set avoids the assumed exterior collar,
\[
 B\subset\mathfrak E\cap
       \{z:1+\varepsilon_0\le|z|\le R\}.
\]
The closed annulus on the right is compact and contained in $D$. An
infinite subset of $\mathfrak E$ in this annulus would have a limit point
in $D$; relative closedness would put that point in $\mathfrak E$,
contradicting its discreteness. Hence the intersection, and therefore
$B$, is finite. The collar assumption is used to keep this compact set
away from the boundary $|z|=1$.

The conclusion concerns intrinsic finite obstructions, not finiteness
of the provisional set $\mathfrak E$. At $\xi=0$ a companion may take an
exceptional value $c_{j,\nu}$, where its accessible germ need not extend.
Thus the argument proves neither holomorphic continuation through
infinity nor trivial monodromy there. It also imposes no global
single-valuedness or meromorphicity assumption.
\end{proof}

Figure~\ref{cont:fig:odd-trace-finiteness} summarizes the pole separation and punctured avoidance used above; the center $\xi=0$ is not filled by this argument.
\begin{figure}[htbp]
\centering
\resizebox{\linewidth}{!}{%
\begin{tikzpicture}[x=1cm,y=1cm,font=\small,
 arr/.style={-{Stealth[length=1.7mm]},line width=.7pt},
 every node/.style={inner sep=2pt}]
\path[use as bounding box] (-.10,-.10) rectangle (15.75,6.20);
\node at (3.58,5.93) {(a) Odd denominator: no companion pole collision};
\node at (11.83,5.93) {(b) Avoidance is only on a punctured disk};
\node[font=\footnotesize,align=center] at (3.62,5.39)
 {$\rho=2/5$: two copies of the line $\operatorname{Im}u=-K'$};
\node[anchor=east,font=\footnotesize] at (.78,4.58) {$\nu=1$};
\node[anchor=east,font=\footnotesize] at (.78,3.12) {$\nu=2$};
\draw[black!45,line width=.5pt] (1.12,4.58)--(6.74,4.58);
\draw[line width=.65pt] (1.06000,4.52000)--(1.18000,4.64000) (1.06000,4.64000)--(1.18000,4.52000);
\draw[line width=.65pt] (3.87000,4.52000)--(3.99000,4.64000) (3.87000,4.64000)--(3.99000,4.52000);
\draw[line width=.65pt] (6.68000,4.52000)--(6.80000,4.64000) (6.68000,4.64000)--(6.80000,4.52000);
\draw[fill=white,line width=.85pt] (1.04800,4.50800) rectangle (1.19200,4.65200);
\node[above=7pt,font=\scriptsize] at (1.12000,4.58000) {$j=0$};
\fill (3.36800,4.58000) circle[radius=.043];
\node[above=7pt,font=\scriptsize] at (3.36800,4.58000) {$j=1$};
\fill (5.61600,4.58000) circle[radius=.043];
\node[above=7pt,font=\scriptsize] at (5.61600,4.58000) {$j=2$};
\fill (2.24400,4.58000) circle[radius=.043];
\node[above=7pt,font=\scriptsize] at (2.24400,4.58000) {$j=3$};
\fill (4.49200,4.58000) circle[radius=.043];
\node[above=7pt,font=\scriptsize] at (4.49200,4.58000) {$j=4$};
\node[font=\footnotesize] at (3.95,4.15) {$a_1+jh$};
\draw[black!45,line width=.5pt] (1.12,3.12)--(6.74,3.12);
\draw[line width=.65pt] (1.06000,3.06000)--(1.18000,3.18000) (1.06000,3.18000)--(1.18000,3.06000);
\draw[line width=.65pt] (3.87000,3.06000)--(3.99000,3.18000) (3.87000,3.18000)--(3.99000,3.06000);
\draw[line width=.65pt] (6.68000,3.06000)--(6.80000,3.18000) (6.68000,3.18000)--(6.80000,3.06000);
\draw[fill=white,line width=.85pt] (3.85800,3.04800) rectangle (4.00200,3.19200);
\node[above=7pt,font=\scriptsize] at (3.93000,3.12000) {$j=0$};
\fill (6.17800,3.12000) circle[radius=.043];
\node[above=7pt,font=\scriptsize] at (6.17800,3.12000) {$j=1$};
\fill (2.80600,3.12000) circle[radius=.043];
\node[above=7pt,font=\scriptsize] at (2.80600,3.12000) {$j=2$};
\fill (5.05400,3.12000) circle[radius=.043];
\node[above=7pt,font=\scriptsize] at (5.05400,3.12000) {$j=3$};
\fill (1.68200,3.12000) circle[radius=.043];
\node[above=7pt,font=\scriptsize] at (1.68200,3.12000) {$j=4$};
\node[font=\footnotesize] at (3.95,2.69) {$a_2+jh$};
\node[font=\footnotesize] at (1.12,2.31) {$0$};
\node[font=\footnotesize] at (3.93,2.31) {$2K$};
\node[font=\footnotesize] at (6.74,2.31) {$4K$};

\node[font=\footnotesize] at (3.82,1.76)
 {crosses: poles of $\mathcal P_m$; squares: $j=0$};
\node at (3.80,1.15) {$d\nmid 2jp\quad(1\le j<d,\ d\text{ odd})$};
\node[font=\footnotesize,align=center] at (3.80,.43)
 {Companion centers are finite,\\but need not be ordinary.};
\draw[black!30,line width=.4pt] (7.75,.13)--(7.75,5.36);
\fill[black!4] (9.30,3.87) circle[radius=.89];
\draw[line width=.6pt] (9.30,3.87) circle[radius=.89];
\fill[black!4] (13.77,3.87) circle[radius=.95];
\draw[line width=.6pt] (13.77,3.87) circle[radius=.95];
\node[font=\footnotesize] at (9.30,5.03) {$|\xi|<\varepsilon$};
\node[font=\footnotesize] at (13.62,5.03) {$W\cap\mathfrak E=\{c\}$};
\draw[fill=white,line width=.7pt] (9.30,3.87) circle[radius=.055];
\node[below=4pt] at (9.30,3.87) {$0$};
\draw[line width=.8pt] (13.71,3.81)--(13.83,3.93) (13.71,3.93)--(13.83,3.81);
\node[below=4pt] at (13.77,3.87) {$c$};
\fill (9.59,4.25) circle[radius=.035];
\node[above right=1pt,font=\footnotesize] at (9.59,4.25) {$\xi$};
\fill (14.15,4.21) circle[radius=.035];
\node[font=\scriptsize] at (13.78,4.51) {$z_{j,\nu}(\xi)$};
\draw[arr] (10.34,3.92)--(12.64,3.92);
\node[above=5pt] at (11.47,3.92) {$z_{j,\nu}$};
\node[font=\footnotesize] at (9.31,2.55) {$0<|\xi|<\varepsilon$};
\node[font=\footnotesize] at (13.61,2.55) {$z_{j,\nu}(\xi)\in D_0$};
\node[font=\footnotesize] at (11.88,1.74)
 {$z_{j,\nu}(\xi)-c=\xi^r a(\xi),\quad a(0)\ne0$};
\node[font=\footnotesize] at (11.88,1.14) {No extension at $\xi=0$ is inferred.};
\node[font=\footnotesize,align=center] at (11.88,.43)
 {Both source sheets give holomorphic extensions\\across every sufficiently large finite exception.};
\end{tikzpicture}
}
\caption{The two local ingredients of
Theorem~\ref{cont:thm:task6-finite-bad-set}.
Left: exact positions modulo $4K$ of $a_\nu+jh$ for the sample rotation
$2/5$; the rows are two copies of the same line, not two different
imaginary heights. Crosses mark poles of $\mathcal P_m$ (the endpoints
$0$ and $4K$ are identified), and no displayed nonzero shift of either source pole
hits a pole. Right: schematic punctured avoidance in the case that the
finite center $c$ is itself exceptional; the companion is nonconstant,
so the factorization with $r\ge1$ removes this center from its sufficiently
small punctured image. The two-sheet argument makes all sufficiently
large finite exceptions good; the exception-free collar then confines
$B$ to a compact annulus. This is not an extension through infinity.}
\label{cont:fig:odd-trace-finiteness}
\end{figure}

\Needspace{12\baselineskip}
\begin{corollary}[Continuation across the collar after removal of good points]
\label{cont:cor:finite-bad-enlarged-domain}
Under the hypotheses of
Theorem~\ref{cont:thm:discrete-exterior-continuation}, the intrinsic bad
set $B\subset\{|z|>1\}$ is finite and effective: at each $b\in B$,
some accessible germ near $b$ does not extend holomorphically across $b$.
For some $0<r_*<1$, the initial physical germ admits continuation
along every path, starting at its fixed base point, in
\[
 \Omega\setminus B,\qquad \Omega=\{|z|>r_*\}.
\]
Thus the finite-set reduction preserves the all-branch collar required
by Theorem~\ref{cont:thm:task6-nonfocal-exclusion}.
\end{corollary}
\begin{proof}
Write $D=\{z\in\C:|z|>1\}$. We distinguish the geometric exceptional
set $\mathfrak E$ supplied by
Theorem~\ref{cont:thm:discrete-exterior-continuation} from the intrinsic
bad set $B$ defined by the all-accessible-germ test of
Lemma~\ref{cont:lem:task6-intrinsic-bad-set}. All paths below are
continuous maps from compact intervals; no rectifiability is required.

\smallskip\noindent
\emph{Step 1: fix an enlarged domain and the same initial physical germ.}
Let $\Omega^{(0)}=\{|z|>r_0\}$, $0<r_0<1$, and let
$\mathfrak E^{(0)}$ be the relatively closed discrete exceptional set
from Theorem~\ref{cont:thm:discrete-exterior-continuation}. The original
germ continues along every path in
$\Omega^{(0)}\setminus\mathfrak E^{(0)}$ starting in its known collar.
Choose $0<\delta<1$ so small that
\[
 \mathcal C_\delta=\{1-\delta<|z|<1+\delta\}
\]
lies in $\Omega^{(0)}$ and in the original holomorphy annulus of $V$,
and is disjoint from $\mathfrak E^{(0)}$. Choose
\[
 \max\{r_0,1-\delta\}<r_*<1,\qquad
 \Omega=\{|z|>r_*\},\qquad
 \mathfrak E=\mathfrak E^{(0)}\cap\Omega.
\]
Restriction preserves relative closedness, discreteness, and the
all-path continuation assertion. Moreover,
\[
 \mathfrak E\subset\{|z|\ge1+\delta\},\qquad
 \mathcal C_*:=\{r_*<|z|<1+\delta\}\subset\Omega\setminus\mathfrak E.
\]
Indeed, every point of $\Omega$ with modulus less than $1+\delta$ lies
in $\mathcal C_\delta$. No exterior point of the original exceptional
set is discarded, since $r_*<1$. The function $V$ is originally
holomorphic on $\mathcal C_*$.

Put $D_0=D\setminus\mathfrak E$ and
$\Omega_0=\Omega\setminus\mathfrak E$. Choose one point $z_0$ with
$1<|z_0|<1+\delta$ and use the same original germ there for both domains.
If an earlier base point was different, join it to $z_0$ within the
known physical collar while avoiding $\mathfrak E$; local detours
around the discrete set give such a connector. The original function
supplies continuation along that connector and its reverse. Thus
rebasing preserves the reachable germs, and the final assertion at the
previously fixed point follows by concatenation.

\smallskip\noindent
\emph{Step 2: compare reachable germs by an endpoint-fixed radial homotopy.}
For $q\in D_0$, let $\mathscr G_D(q)$ be the germs reachable from the
fixed germ at $z_0$ along paths in $D_0$. For $q\in\Omega_0$, define
$\mathscr G_\Omega(q)$ using paths in $\Omega_0$. Inclusion of domains
gives $\mathscr G_D(q)\subset\mathscr G_\Omega(q)$ for $q\in D_0$.

For the reverse inclusion, take a path $\gamma:[0,1]\to\Omega_0$ from
$z_0$ to an ordinary exterior point $q\in D_0$. Choose
\[
 1<r_1<\min\{|z_0|,|q|,1+\delta\}.
\]
For $0\le s\le1$ and $z\in\Omega$, define
\[
 H_s(z)=\bigl((1-s)|z|+s\max\{|z|,r_1\}\bigr)\frac{z}{|z|}.
\]
This is a continuous homotopy because $0\notin\Omega$; it need not be
holomorphic. It fixes every point of modulus at least $r_1$, in
particular both endpoints. If $|z|<r_1$, then
\[
 r_*<|z|\le |H_s(z)|\le r_1<1+\delta,
\]
so the entire trajectory lies in the exception-free collar
$\mathcal C_*$. Hence $H_s(\Omega_0)\subset\Omega_0$ for every $s$,
while $H_1\circ\gamma$ lies in $D_0$, since its modulus is at least
$r_1>1$.

The homotopy $H_s\circ\gamma$ keeps the endpoints fixed and stays in
$\Omega_0$. The all-path hypothesis on that domain therefore permits
homotopy invariance of analytic continuation, in the finite-chain form
used in Lemma~\ref{cont:lem:countable-gluing}. Its two end paths reach
the same germ at $q$. Thus every germ in $\mathscr G_\Omega(q)$ is
represented by a path in $D_0$, proving
\[
 \boxed{\mathscr G_\Omega(q)=\mathscr G_D(q),\qquad q\in D_0.}
\]
This comparison concerns exterior endpoints only. It identifies sets
of reachable germs, not all germs within each set; loops may still
have nontrivial monodromy. No removal of good points has yet been used.

\smallskip\noindent
\emph{Step 3: identify the same bad set and remove good points on $\Omega$.}
For $b\in\mathfrak E$, choose an isolating disk
$\Delta_b\Subset D$ with $\Delta_b\cap\mathfrak E=\{b\}$. In the
exterior domain the good-point test is
\[
 \begin{gathered}
 \text{for every }q\in\Delta_b\setminus\{b\}
 \text{ and every }\mathfrak f\in\mathscr G_D(q),\\
 \text{there is }F\in\mathcal O(\Delta_b)
 \text{ with }[F]_q=\mathfrak f.
 \end{gathered}
\]
Here $[F]_q$ denotes a germ, and the extension may depend on
$\mathfrak f$. Step~2 shows that replacing $\mathscr G_D(q)$ by
$\mathscr G_\Omega(q)$ gives exactly the same test, since every test
point belongs to $D_0$. Disk independence from
Lemma~\ref{cont:lem:task6-intrinsic-bad-set} applies. Consequently the
bad sets computed in the two domains coincide; denote this common set
by $B\subset\mathfrak E$. As a subset of the relatively closed discrete
set $\mathfrak E$, it is relatively closed and discrete in $\Omega$.

The disk-chain removal proof of that lemma now applies on $\Omega$.
To check the domain change explicitly, cover a path in
$\Omega\setminus B$ by a finite chain of disks avoiding
$\mathfrak E$ or isolating a good point. All good-point disks may be
chosen in $D$. Choose auxiliary switches outside $\mathfrak E$ in the
overlap components containing the original switches, and inductively
retain an $\Omega_0$-path representing the incoming germ at each
auxiliary switch. On a disk avoiding $\mathfrak E$, concatenation and
simple connectivity give a holomorphic representative on the full disk.
On a good-point disk, the incoming switch lies in $D_0$; Step~2 makes
its $\Omega_0$-accessible germ originally accessible in $D_0$, so the
same good-point test supplies a full-disk representative.

Join the auxiliary switches within the disks while avoiding
$\mathfrak E$. These detours retain an $\Omega_0$-representing path
for each outgoing germ. Equality of germs at an auxiliary switch and
the identity theorem on the corresponding overlap component give
agreement at the original switch. Thus the representatives continue
the germ along the original path, including a good endpoint. At an
endpoint in $\Omega_0$, the detours also give an $\Omega_0$-path with
the same terminal germ. This is the simultaneous construction in
Lemma~\ref{cont:lem:task6-intrinsic-bad-set}; it uses no special geometry
of the circle $|z|=1$ beyond placing good-point disks in $D$.

We have therefore proved all-path continuation on $\Omega\setminus B$.
If $\mathscr G_{\Omega,B}(q)$ denotes its reachable-germ set, the same
construction gives one inclusion below, and inclusion of domains gives
the other:
\[
 \boxed{\mathscr G_{\Omega,B}(q)=\mathscr G_\Omega(q),
        \qquad q\in\Omega_0.}
\]
In particular, the removal creates no new germs at any previously
regular point of $\Omega$. Original accessibility was checked before
each good-point extension, so this argument is not circular.

\smallskip\noindent
\emph{Step 4: apply the odd-denominator finiteness theorem.}
Choose a nondegenerate interval $(a,b)$ with $[a,b]\subset I$.
For a sufficiently large odd integer $q$, the interval $(qa,qb)$
contains an integer $p$, so $p/q\in I$. Reduction divides the odd
denominator $q$ by an odd integer; the primitive denominator remains
odd. The reduced rotation still lies in $(0,1/2)$, and its trace is
constant by the original band hypothesis.

The domain $D_0$ has the standing all-path continuation property,
$\mathfrak E$ is relatively closed and discrete in $D$, and
$\mathfrak E\cap\{1<|z|<1+\delta\}=\varnothing$. All hypotheses of
Theorem~\ref{cont:thm:task6-finite-bad-set} are therefore satisfied.
It makes the exterior intrinsic bad set finite. By Step~3 this is
exactly the same $B$ used on $\Omega$, not a newly defined bad set.

\smallskip\noindent
\emph{Step 5: preserve a nonextendible witness at every retained point.}
For $b\in B$, failure of the good-point test supplies
$q\in\Delta_b\setminus\{b\}$ and
$\mathfrak f\in\mathscr G_D(q)$ for which no holomorphic function on
$\Delta_b$ represents $\mathfrak f$ at $q$. Steps~2--3 give
\[
 \mathfrak f\in\mathscr G_D(q)
     =\mathscr G_\Omega(q)
     =\mathscr G_{\Omega,B}(q).
\]
Thus the same nonextendible germ remains accessible in the final
domain. Its nonextendibility is a property of that germ in the
physical coordinate, not of its representing path. Disk independence
in Lemma~\ref{cont:lem:task6-intrinsic-bad-set} ensures that this is a
local obstruction at $b$, rather than a choice of isolating radius.
This proves effectiveness in the stated sense: at least one accessible
germ witnesses each $b\in B$. It does not say that every branch is
singular at that point or impose any algorithmic meaning on
``effective.''

\smallskip\noindent
\emph{Step 6: retain the collar for every reachable branch.}
Since $B\subset\mathfrak E$, the collar $\mathcal C_*$ is disjoint
from $B$. Let a germ $\mathfrak h$ at $q\in\mathcal C_*$ be reached
along a path $\alpha$ in $\Omega\setminus B$ from the fixed initial
germ. For any further path $\beta$ in $\mathcal C_*$ starting at $q$,
the concatenation $\alpha*\beta$ stays in $\Omega\setminus B$.
All-path continuation along it starts on $\beta$ with precisely
$\mathfrak h$, by uniqueness along $\alpha$. The branch history is
not reset to the original branch. A connector in the original collar
also returns this assertion to any previously fixed initial point.

We have obtained the same finite, effective exterior set $B$ and
all-path continuation on a domain containing the entire unit-circle
collar. The argument allows arbitrary monodromy; it does not assert
that a reachable branch is single-valued on the annular collar, that
$B$ consists only of poles, or that any branch is regular at infinity.
Those conclusions are not used here.
\end{proof}

\section{A finite bad set is confined to the exterior focal pair}
\label{cont:sec:task6-nonfocal}

We first exclude a nonfocal bad value by turning its persistent odd-prime
trace collisions into an affine geometric gap, then showing that this gap
is a constant half-integer of a type discarded before the finite-grid
argument. After focal confinement, parity makes the limiting companion
centers ordinary, so the continued trace extends every accessible germ
through infinity. Only the bad base values form a finite family; the
number of analytic branches need not be finite.


The statements below concern all branches accessible from the physical
germ. They do not assume finite monodromy or a compact continuation
surface. Lemma~\ref{cont:lem:task6-trace-realization} includes any prescribed
accessible source germ in a continued odd-denominator rational trace,
after an arbitrarily small local continuation if necessary.
Its representing path is perturbed in the lower Jacobi cylinder to
avoid the shifted exceptional sets while preserving the source germ,
transported by that local continuation. Either Jacobi preimage of the
resulting physical endpoint may be reached.

\begin{theorem}[All-sheet exclusion of nonfocal bad values]
\label{cont:thm:task6-nonfocal-exclusion}
Fix a noncircular ellipse with semiaxes $a>b>0$ and a nonempty open
interval $I\subset(0,1/2)$. Let the physical germ
$V(e^{i\varphi})$ be holomorphic near the unit circle and satisfy every
odd-denominator rational trace identity in $I$. Suppose that, for some
$r_*<1$, it admits analytic continuation along every path in
\[
 \Omega\setminus B,\qquad \Omega=\{z\in\mathbb C:|z|>r_*\},
\]
starting in the original physical collar, where
$B\subset\{|z|>1\}$ is finite and effective. Here effectiveness means
that each $s\in B$ has an accessible germ in a punctured neighborhood of
$s$ that does not admit a holomorphic extension across $s$ in the physical
coordinate. Then, with $c=\sqrt{a^2-b^2}$,
\[
 B\subset\left\{-\frac{a+b}{c},\frac{a+b}{c}\right\}.
\]
Moreover every accessible branch near infinity extends holomorphically
in the coordinate $1/z$ at infinity. In particular every such branch
has trivial local monodromy there.
\end{theorem}

\begin{proof}[Exclusion of nonfocal finite values]
Put $D=\{z\in\C:|z|>1\}$. An accessible branch near a puncture may
initially live on the universal cover of its punctured disk. Throughout
this proof, extension at a finite puncture means extension as a
holomorphic function of the ordinary physical coordinate $z$, not merely
of a covering coordinate. We first justify using the exterior
trace-realization lemma for every branch in the hypothesis.

\smallskip\noindent
\emph{Step 1: preserve the accessible germs when restricting to the exterior.}
Choose the initial point $z_0$ in the exterior part of the original
holomorphy collar, using a connector in that collar if a different base
point was previously fixed. The original function supplies continuation
along this connector and its reverse, so rebasing preserves the
reachable germs. Since $B$ is finite and lies in $D$, choose $\delta>0$
so small that $B\cap\{|z|\le1+\delta\}=\varnothing$ and the unit-circle
collar of width $\delta$ lies in the original holomorphy annulus.

Let $q\in D\setminus B$, and let $\gamma$ be a continuous representing
path from $z_0$ to $q$ in $\Omega\setminus B$. If zero belongs to
$\Omega$, perturb the path away from zero, keeping its endpoints and its
homotopy class in $\Omega\setminus B$. The path-perturbation argument
preceding Lemma~\ref{cont:lem:task6-trace-realization} applies to this
single additional point. If zero is already outside $\Omega$, no such
perturbation is needed. Choose
\[
 1<r_1<\min\{|z_0|,|q|,1+\delta\}.
\]
On the zero-avoiding path use the continuous homotopy
\[
 H_\tau(z)=\bigl((1-\tau)|z|+\tau\max\{|z|,r_1\}\bigr)
                 \frac{z}{|z|},\qquad 0\le\tau\le1.
\]
It fixes both endpoints and every point with $|z|\ge r_1$. Every moved
point stays at radius between $|z|>r_*$ and $r_1$, in a region without
points of $B$. Hence the homotopy stays in $\Omega\setminus B$, and its
final path lies in $D\setminus B$. No holomorphy of the homotopy is
required. The all-path hypothesis and homotopy invariance, in the form
proved in Lemma~\ref{cont:lem:countable-gluing}, preserve the terminal
germ. Together with inclusion of domains, this proves
\[
 \mathscr G_{\Omega\setminus B}(q)
   =\mathscr G_{D\setminus B}(q),\qquad q\in D\setminus B,
\]
where $\mathscr G$ denotes the set of germs reachable from the fixed
initial germ. This identifies sets of branches, not the branches within
each set. All subsequent applications of
Lemma~\ref{cont:lem:task6-trace-realization} use its exceptional set
$\mathfrak E=B$. In particular every effective nonextendible witness is
eligible for that lemma. All paths have compact parameter intervals;
no rectifiability is assumed.

\smallskip\noindent
\emph{Step 2: fix a nonfocal bad value and both geometric inverse branches.}
After a uniform scaling take $a=1$, and put $\eta=1-b^2\in(0,1)$, with
$b$ now the normalized minor semiaxis. The physical parameter
$z=e^{i\varphi}$ is unchanged. Write
\[
 \begin{gathered}
 \mathcal P_m(u)=\sn(u\mid m)-i\cn(u\mid m),\qquad m\in(\eta,1),\\
 K=K(m),\qquad K'=K(1-m),\qquad J=2K\rho(m).
 \end{gathered}
\]
The prime in $K'$ denotes the complementary integral, not a derivative.
By Lemma~\ref{arith:jacobi}, $\mathcal P_m$ is a proper degree-two map
from the lower cylinder to the physical exterior with infinity adjoined.

Suppose $s\in B$ is not an exterior focal value, and choose a fixed
complex angle $\varphi_s$ with $s=e^{i\varphi_s}$. Then
\[
 \cos\varphi_s=\frac{s+s^{-1}}2\ne0,\qquad
 \sin\varphi_s=\frac{s-s^{-1}}{2i}\ne0,
 \qquad m_s:=\sec^2\varphi_s=\frac4{(s+s^{-1})^2}.
\]
A zero of either trigonometric factor would put $s$ on the unit circle.
Consequently $m_s\notin\{0,1\}$, and nonfocality gives $m_s\ne\eta$.
Fix an accessible nonextendible witness at $s$, supplied by
effectiveness. Transport to a smaller punctured isolating disk, when
needed, preserves such a witness by the disk-independence argument of
Lemma~\ref{cont:lem:task6-intrinsic-bad-set}.

For any fixed $t\in B$, the identities
\[
 \sn(u\mid m)=\frac{\mathcal P_m(u)+\mathcal P_m(u)^{-1}}2,
 \qquad \partial_u\mathcal P_m(u)=i\dn(u\mid m)\mathcal P_m(u)
\]
show that a preimage of $t$ can be critical only at
\[
 m_t=\frac4{(t+t^{-1})^2}.
\]
The denominator is nonzero for $|t|>1$. Equation~\eqref{op:eq:rhoDerivative}
and $m=\eta/(1-\lambda^2)$, with
$dm/d\lambda=2\eta\lambda/(1-\lambda^2)^2>0$, show that $m(\rho)$
is analytic with positive real derivative. Thus each $t\in B$ excludes
at most one real rotation. Choose a nonempty $I_0\Subset I$ avoiding
these finitely many parameters.

Shrink $I_0$ further so that both preimages of every $t\in B$ have
holomorphic lifts $u_{t,\beta}(\rho)$, $\beta=1,2$, on a common simply
connected complex neighborhood of $\overline{I_0}$. This follows by
applying the holomorphic implicit-function theorem at the finitely many
regular inverse points of one parameter and then taking a common
smaller neighborhood. For real $\rho\in I_0$, choose representatives
with
\[
 -2K'(\rho)<\operatorname{Im}u_{t,\beta}(\rho)<0.
\]
Include the two geometric lifts of infinity,
\[
 u_{\infty,1}(\rho)=-iK'(\rho),\qquad
 u_{\infty,2}(\rho)=2K(\rho)-iK'(\rho).
\]
These choices do not assume regularity of $V$ at infinity. The list is
finite because $B$ is finite, not because the number of analytic
branches of $V$ is finite.

\smallskip\noindent
\emph{Step 3: convert every prime-grid collision into a nonintegral lattice gap.}
Let $q$ be an odd prime and $p/q\in I_0$. Since $0<p<q/2$, the fraction
is reduced. Its trace identity holds by hypothesis. Apply the final
assertion of Lemma~\ref{cont:lem:task6-trace-realization} to the fixed
bad witness near $s$. All points of $B$ are noncritical for the current
parameter. If every companion value at both preimages of $s$ were
finite and outside $B$, the two-sheet contrapositive in that lemma
would extend every accessible source germ across $s$. Hence at least
one source preimage has a collision with a value in $B\cup\{\infty\}$
after a shift $jh$, where $h=4Kp/q$ and $1\le j<q$.

For the fixed physical source $s$, form the finite family
\[
 g_{t,\alpha,\beta}(\rho)
   =\frac{u_{t,\beta}(\rho)-u_{s,\alpha}(\rho)}{4K(\rho)},
 \qquad t\in B\cup\{\infty\},\quad \alpha,\beta\in\{1,2\}.
\]
It is holomorphic on the common neighborhood. Each index records a
fixed physical target and fixed geometric source and target lifts.
A collision at $p/q$ gives, for at least one index,
\[
 g_{t,\alpha,\beta}(p/q)=\frac{jp}{q}+n,
 \qquad 1\le j<q,\quad n\in\Z.
\]
There is only a real-period ambiguity: the shift is real and both
representatives remain in the lower cylinder, whose height is smaller
than the imaginary period $4K'$. Thus
\[
 qg_{t,\alpha,\beta}(p/q)\in\Z,
 \qquad g_{t,\alpha,\beta}(p/q)\notin\Z.
\]
The second assertion uses primitivity and $1\le j<q$. Changing either
lift by a real period changes the gap by an integer and preserves both
conditions. The supplying index, shift, and reached source sheet may
vary with the grid point; no uniform choice among them is assumed.

\smallskip\noindent
\emph{Step 4: discard the impossible constant gaps before selecting an affine one.}
Discard every gap identically constant with value in $\tfrac12\Z$.
An integer-valued constant cannot represent a genuine collision, by the
nonintegrality just proved. For a nonintegral half-integer constant,
$qg$ is not integral when $q$ is odd. Thus none of the discarded
functions can supply a genuine prime-grid collision. The remaining
fixed finite family still covers every required grid point. It is
nonempty; otherwise the collision assertion would already fail.

Lemma~\ref{arith:grid}, including its fixed-subfamily assertion, now
gives one retained member
\[
 g(\rho)=A\rho+C\qquad(\rho\in I_0),\qquad A,C\in\C.
\]
That lemma permits the supplying index to change at every grid point:
its uniform weighted second-difference argument first forces a zero of
one second derivative at each point, and its finite analytic zero-set
argument then fixes one affine member on the whole interval. Neither
reality nor rationality of $A,C$ is required, and the selected member
need not supply every collision. Fix its target and lifts, and
abbreviate them by $t$, $u_t$, and $u_s$. Since $J=2K\rho$,
\[
 u_t-u_s=2AJ+4CK.
\]
We apply the parameter equation to this unnormalized identity, not to
the quotient defining $g$.

\smallskip\noindent
\emph{Step 5: compute the marked-point sources with fixed branch choices.}
The change between $\rho$ and $m$ is locally analytically invertible.
For a fixed finite physical point $z=e^{i\varphi_z}$, choose a
compatible local branch of $I_{\varphi_z}(m)$ from
Lemma~\ref{arith:pf}. The reference inverse is $K+I_{\varphi_z}$.
The other inverse is obtained from the deck involution, and all
representatives differ by periods $4K$ and $4iK'$. Consequently,
locally in the parameter,
\[
 u_z=\varepsilon_z I_{\varphi_z}+c_{z,1}K+c_{z,2}iK',
 \qquad \varepsilon_z\in\{1,-1\},\quad c_{z,1},c_{z,2}\in\Z.
\]
These coefficients are constant: fix the deck and period choices at
one regular parameter and use uniqueness of its local unramified
inverse to keep the relation nearby. All square roots and integral
germs in the following calculation are continued compatibly.

With the ellipse and each physical angle fixed, apply
\[
 \mathcal L=m(1-m)\partial_m^2+(1-2m)\partial_m-\tfrac14.
\]
Lemma~\ref{arith:pf} gives $\mathcal LK=\mathcal LK'=0$ and
$4\mathcal LJ=-\sqrt{\eta(1-\eta)}(m-\eta)^{-3/2}$. Therefore
$4\mathcal L(u_t-u_s)=8A\mathcal LJ$. For a finite target this is
\begin{equation}\label{cont:eq:task6-nonfocal-pf}
 \varepsilon_t
 \frac{\sin\varphi_t\cos\varphi_t}
 {(1-m\cos^2\varphi_t)^{3/2}}
 -\varepsilon_s
 \frac{\sin\varphi_s\cos\varphi_s}
 {(1-m\cos^2\varphi_s)^{3/2}}
 =-\frac{2A\sqrt{\eta(1-\eta)}}{(m-\eta)^{3/2}}.
\end{equation}
When $t=\infty$, omit the target term: its two lifts are constant
linear combinations of $K$ and $iK'$ and hence have zero source.
This uses only the geometric lifts of infinity, not any extension of
$V$ there. The period terms vanish because their coefficients are
constant; the operator does not annihilate arbitrary constants or
parameter-dependent multiples of periods.

\smallskip\noindent
\emph{Step 6: force matching marked singularities and zero affine slope.}
We make the singularity comparison algebraic, so no global continuation
of $V$ or of the chosen inverse lifts is needed. For $z=s$ and, when
finite, $z=t$, set
\[
 P_z(m)=\frac{\sin^2\varphi_z\cos^2\varphi_z}
                  {(1-m\cos^2\varphi_z)^3},\qquad
 Q(m)=\frac{4A^2\eta(1-\eta)}{(m-\eta)^3},
\]
and put $P_t=0$ when $t=\infty$. Squaring
\eqref{cont:eq:task6-nonfocal-pf} twice gives the necessary identity
\[
 (P_t+P_s-Q)^2=4P_tP_s.
\]
It holds as an identity of rational functions: it first holds on the
nonempty local parameter neighborhood, and clearing denominators
makes it a polynomial identity. Only a necessary consequence of the
original source equation is used, so no converse after squaring is
being assumed.

The rational function $P_s$ has a genuine pole of order three at
$m_s=\sec^2\varphi_s$, with nonzero coefficient, whereas $Q$ is
regular there because $m_s\ne\eta$. If $P_t$ were regular there, the
left side of the last identity would have a pole of order six and the
right side a pole of order at most three. This is impossible. Hence
the target is finite and
\[
 \cos^2\varphi_t=\cos^2\varphi_s.
\]
Their squared sine factors agree as well, so $P_t=P_s$. The rational
identity reduces to $Q(Q-4P_s)=0$. The alternative $Q=4P_s$ is
impossible: $Q$ is regular at $m_s$ and $P_s$ is not. Thus $Q=0$, and
since $0<\eta<1$ we obtain $A=0$. Returning to the original compatible
source equation, its two marked terms now cancel identically.

Equivalently, on a compatible slit neighborhood of $m_s$ the source
has a nonzero $(m-m_s)^{-3/2}$ term, while the right side is regular.
The rational argument above encodes precisely why a finite target
with the same marked singularity must cancel that term. It evaluates
only rational functions outside the initial parameter neighborhood;
no inverse branch is continued to $m_s$, and no principal square root
is reset there.

\smallskip\noindent
\emph{Step 7: identify the constant gap and reach the discarded-family contradiction.}
Equality of the cosine squares gives
\[
 \frac{(t-s)(t+s)(st-1)(st+1)}{4s^2t^2}=0.
\]
Thus $t\in\{s,-s,s^{-1},-s^{-1}\}$; because both values are exterior,
$t=\pm s$. The deck involution from Lemma~\ref{arith:jacobi} and the
real half-period identity $\mathcal P_m(u+2K)=-\mathcal P_m(u)$ give
\[
 u_t=\sigma u_s+2\ell K+2n iK',\qquad
 \sigma\in\{1,-1\},\quad \ell,n\in\Z.
\]
For $t=s$, the two possibilities come from the identity and the deck
involution, followed by full periods; for $t=-s$, first add $2K$.
These formulas contain all preimages because the map has degree two.
As in Step~5, $\sigma,\ell,n$ are constant on a smaller parameter
neighborhood after they are fixed at one regular parameter.

Since $A=0$, the affine identity gives
$4\mathcal L(u_t-u_s)=4\mathcal L(4CK)=0$. Applying the same operator
to the last lift relation yields
\[
 0=(\sigma-1)4\mathcal Lu_s.
\]
But $4\mathcal Lu_s=\varepsilon_s R_{\varphi_s}$ is nonzero, since
$\sin\varphi_s\cos\varphi_s\ne0$. Hence $\sigma=1$. For real
parameters the two selected lifts lie in the lower cylinder, so
\[
 -2K'<\operatorname{Im}(u_t-u_s)<2K'.
\]
The lift relation now has imaginary part $2nK'$, forcing $n=0$.
Consequently
\[
 u_t-u_s=2\ell K,\qquad g=\frac{u_t-u_s}{4K}=\frac{\ell}{2}.
\]
The identity theorem extends this local equality to $I_0$. Thus the
retained affine member is identically a constant half-integer, contrary
to Step~4. This last step is needed: zero slope by itself would not
contradict the finite-grid conclusion.

There is therefore no nonfocal finite point of $B$. Restoring the
original semiaxes, the four zeros of $D_E(z)$ are
$\pm(a+b)/c$ and $\pm(a-b)/c$, where $c=\sqrt{a^2-b^2}$.
Since $(a-b)/c=\sqrt{(a-b)/(a+b)}<1$, only the exterior pair can belong
to $B$. This proves the finite-point conclusion. The argument has used
finiteness of bad physical values and effectiveness of the chosen
witness, not finiteness of analytic branches or regularity at infinity.
\end{proof}

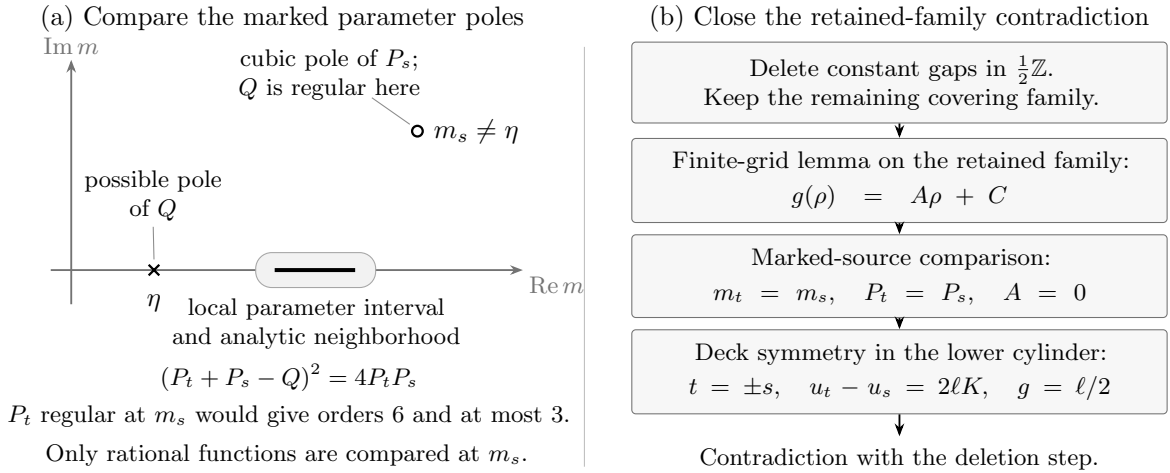
\begin{figure}[!htbp]
\centering
\resizebox{\linewidth}{!}{%
\begin{tikzpicture}[x=1cm,y=1cm,font=\small,
 arr/.style={-{Stealth[length=1.7mm]},line width=.7pt},
 box/.style={draw=black!55,rounded corners=1.5pt,fill=black!3,
   text width=6.70cm,font=\footnotesize,align=center,inner sep=5pt},
 every node/.style={inner sep=2pt}]
\path[use as bounding box] (-.06,-.10) rectangle (15.75,6.33);
\node at (3.75,6.04) {(a) Compare the marked parameter poles};
\node at (11.88,6.04) {(b) Close the retained-family contradiction};
\draw[black!25,line width=.4pt] (7.73,.14)--(7.73,5.67);
\draw[arr,black!55] (.70,2.73)--(6.92,2.73) node[below right,font=\footnotesize] {$\operatorname{Re}m$};
\draw[arr,black!55] (.97,2.28)--(.97,5.50) node[above,font=\footnotesize] {$\operatorname{Im}m$};
\draw[fill=black!5,draw=black!35,rounded corners=6pt] (3.40,2.50) rectangle (4.98,2.96);
\draw[line width=1.5pt] (3.65,2.73)--(4.71,2.73);
\node[align=center,font=\footnotesize] at (4.20,2.03) {local parameter interval\\and analytic neighborhood};
\draw[line width=.85pt] (1.99,2.66)--(2.13,2.80) (1.99,2.80)--(2.13,2.66);
\node[below=7pt] at (2.06,2.73) {$\eta$};
\node[align=center,font=\footnotesize] at (2.05,3.71) {possible pole\\of $Q$};
\draw[black!45] (2.05,3.30)--(2.06,2.87);
\draw[fill=white,line width=.85pt] (5.53,4.56) circle[radius=.07];
\node[right=4pt] at (5.53,4.56) {$m_s\ne\eta$};
\node[align=center,font=\footnotesize] at (4.37,5.31) {cubic pole of $P_s$;\\$Q$ is regular here};
\draw[black!45] (4.93,5.00)--(5.43,4.65);
\node[align=center,font=\footnotesize] at (3.84,1.08)
 {$(P_t+P_s-Q)^2=4P_tP_s$\\[3pt]
  $P_t$ regular at $m_s$ would give orders $6$ and at most $3$.};
\node[font=\footnotesize,align=center] at (3.81,.29)
 {Only rational functions are compared at $m_s$.};
\node[box] (delete) at (11.88,5.20)
 {Delete constant gaps in $\tfrac12\mathbb Z$.\\
  Keep the remaining covering family.};
\node[box] (affine) at (11.88,3.91)
 {Finite-grid lemma on the retained family:\\[2pt]
  $g(\rho)=A\rho+C$};
\draw[arr] (delete.south)--(affine.north);
\node[box] (source) at (11.88,2.65)
 {Marked-source comparison:\\[2pt]
  $m_t=m_s,\quad P_t=P_s,\quad A=0$};
\draw[arr] (affine.south)--(source.north);
\node[box] (symmetry) at (11.88,1.39)
 {Deck symmetry in the lower cylinder:\\[2pt]
  $t=\pm s,\quad u_t-u_s=2\ell K,\quad g=\ell/2$};
\draw[arr] (source.south)--(symmetry.north);
\draw[arr] (symmetry.south)--(11.88,.54);
\node[font=\footnotesize,align=center] at (11.88,.24)
 {Contradiction with the deletion step.};
\end{tikzpicture}
}%
\caption{The nonfocal contradiction in Theorem~\ref{cont:thm:task6-nonfocal-exclusion}.
(a) Schematic parameter-plane separation: the squared source $P_s$ has a cubic
pole at $m_s\ne\eta$, while $Q$ is regular there. The necessary rational
identity forces a matching target pole; no inverse branch is continued to
$m_s$. The marked points represent singularities of known parameter-dependent
sources, not physical singularities of $V$.
(b) Constant half-integer gaps are removed before the finite-grid lemma is
applied. The subsequent source comparison and deck argument force its retained
affine member into precisely that deleted class. The finite family consists
of geometric gaps, not analytic branches of the unknown.}
\label{fig:nonfocal-gap-obstruction}
\end{figure}

Figure~\ref{fig:nonfocal-gap-obstruction} summarizes the contradiction
in the first part of the proof; the zero-slope conclusion alone would
not suffice.

\Needspace{10\baselineskip}


\begin{proof}[Regularity of all branches at infinity]
The first part has proved that $B$ is contained in the exterior focal
pair. Fix one reduced $p/d\in I$ with odd $d$. Such a rotation exists:
choose a rational with a sufficiently large odd denominator in a
nondegenerate subinterval of $I$ and reduce it; its denominator remains
odd. Keep this caustic, and hence $m$, $K=K(m)$ and $K'=K(1-m)$, fixed.
Put $J=2Kp/d$, $h=2J$, and $T_j(u)=\mathcal P_m(u+jh)$.
The equality of exterior reachable-germ sets proved in the first part
allows us to use Lemma~\ref{cont:lem:task6-trace-realization} for every
accessible source germ, with $\mathfrak E=B$.

\smallskip\noindent
\emph{Step 1: exclude both poles and focal values at all companion endpoints.}
The two simple poles of $T_0$ in the lower cylinder are
\[
 a_1=-iK',\qquad a_2=2K-iK'\pmod{4K}.
\]
At either pole, a nontrivial companion $T_j$, $1\le j<d$, could have a
pole only if
\[
 jh\in2K\Z\quad\Longleftrightarrow\quad d\mid2jp.
\]
Coprimality implies $d\mid2j$, and oddness implies $d\mid j$,
impossible for $1\le j<d$. Every such companion value is therefore
finite and, by the lower-cylinder geometry, exterior.

Lemma~\ref{arith:focal} gives the full lower-cylinder fibers of the
exterior focal values as $K\pm J-iK'$ and $3K\pm J-iK'$.
A focal collision from either $a_\nu$ would require
\[
 jh\equiv K\pm J\pmod{2K}.
\]
Writing $jh=K\pm J+2Kr$ and substituting $h=4Kp/d$, $J=2Kp/d$ gives
\[
 2p(2j\mp1)=d(2r+1),\qquad r\in\Z.
\]
The left side is even and the right side is odd. Thus no such collision
occurs, on either source sheet or with either focal value. Since no
other point can remain in $B$,
\[
 c_{j,\nu}:=T_j(a_\nu)\in D\setminus B,
 \qquad 1\le j<d,\quad \nu=1,2.
\]
Unlike the limiting-value argument in
Theorem~\ref{cont:thm:task6-finite-bad-set}, these center values themselves
are now ordinary.

\smallskip\noindent
\emph{Step 2: choose both inverse coordinates and all ordinary target disks first.}
At $a_\nu$ write
\[
 T_0(u)=\frac{A_\nu}{u-a_\nu}+O(1),\qquad A_\nu\ne0.
\]
Then $1/T_0$ is holomorphic with a simple zero and nonzero derivative
at $a_\nu$. The holomorphic inverse-function theorem gives disjoint
source neighborhoods with holomorphic inverses
\[
 u_\nu:\{|\xi|<\varepsilon\}\longrightarrow\mathcal C_-,\qquad
 u_\nu(0)=a_\nu,\qquad T_0(u_\nu(\xi))=\xi^{-1}\quad(\xi\ne0),
\]
after restriction to one common small disk. For each $j\ne0$ and
$\nu$, choose an ordinary disk $W_{j,\nu}\Subset D\setminus B$
containing $c_{j,\nu}$. Shrink $\varepsilon$ so that all compositions
are holomorphic on the parameter disk and
\[
 T_j(u_\nu(\xi))\in W_{j,\nu}
 \qquad(|\xi|<\varepsilon,\ 1\le j<d,\ \nu=1,2).
\]
In particular this inclusion holds at $\xi=0$. The finite family of
target disks and the common radius are chosen before any branch of
$V$ is selected.

Choose $R>1$ so large that $R>\varepsilon^{-1}$, that $B\subset\{|z|<R\}$,
and that the two finite critical values of $T_0$ lie in $\{|z|<R\}$.
For every finite $|z|>R$, the points $u_1(1/z)$ and $u_2(1/z)$ are
distinct regular preimages. Degree two shows that they exhaust the
fiber. Thus both possible lifted endpoints lie in the controlled source
neighborhoods and avoid all shifted bad values and poles.

\smallskip\noindent
\emph{Step 3: realize an arbitrary source germ at its exact endpoint.}
Take any accessible germ $\mathfrak f$ at a finite $q$ with $|q|>R$.
It is accessible through $D\setminus B$ by the first part. The
exact-endpoint clause in
Lemma~\ref{cont:lem:task6-trace-realization} realizes it at $q$ itself,
with no endpoint movement, on whichever inverse sheet $\nu$ is reached.
The continued trace, in the coordinate $\xi=1/z$, is
\[
 V_0(1/\xi)=C_\rho-
       \sum_{j=1}^{d-1}V_j(T_j(u_\nu(\xi)))
 \quad\text{near }\xi_q=1/q.
\]
Here $V_0$ represents $\mathfrak f$, and the companion germs are
accessible at points of $W_{j,\nu}$. Concatenating a representing path
with any further path in that disk gives all-path continuation of each
companion there. Since $W_{j,\nu}$ is simply connected, the monodromy
theorem gives a holomorphic representative
$\mathcal V_{j,\nu}\in\mathcal O(W_{j,\nu})$ of the particular
companion germ. These representatives may depend on $\mathfrak f$ and
its realized trace; the neighborhoods and radius do not.

\smallskip\noindent
\emph{Step 4: construct the full-disk extension in the physical coordinate.}
Define
\[
 F_{\mathfrak f}(\xi)=C_\rho-
   \sum_{j=1}^{d-1}\mathcal V_{j,\nu}\bigl(T_j(u_\nu(\xi))\bigr),
 \qquad |\xi|<\varepsilon.
\]
Every term is holomorphic throughout this full disk, including zero.
The trace identity identifies its germ at $\xi_q$ with the prescribed
source germ expressed as $V(1/\xi)$. It is therefore a holomorphic
extension of that branch at infinity in the ordinary coordinate
$\xi=1/z$. It supplies continuation along every path in the punctured
disk from that germ, and a loop around zero returns the same germ.
Thus its local monodromy at infinity is trivial.

The argument treats every accessible source germ at every sufficiently
large finite point, and it works on either reached source sheet. Hence
every accessible branch near infinity has the asserted extension.
The conclusion is not merely holomorphy in a ramified Jacobi coordinate:
the simple zero of $1/T_0$ makes the chosen source coordinate
biholomorphic to $1/z$. Different accessible branches may give different
functions $F_{\mathfrak f}$ or values at infinity. No triviality of
monodromy at the remaining focal points, and no global meromorphicity,
is asserted here. This completes the second part of the theorem.
\end{proof}

\begin{figure}[!htbp]
\centering
\resizebox{\linewidth}{!}{%
\begin{tikzpicture}[x=1cm,y=1cm,font=\small,
 arr/.style={-{Stealth[length=1.7mm]},line width=.7pt},
 every node/.style={inner sep=2pt}]
\path[use as bounding box] (-.06,-.10) rectangle (15.75,6.12);
\node at (3.74,5.85) {(a) Ordinary centers after focal confinement};
\node at (11.85,5.85) {(b) Full-disk extension in $\xi=1/z$};
\draw[black!25,line width=.4pt] (7.73,.14)--(7.73,5.51);
\node[font=\footnotesize] at (3.72,5.22)
 {$m=16/25,\quad\rho=2/5,\quad\eta\simeq0.0936464$};
\draw[arr,black!55] (.43,3.55)--(7.17,3.55) node[below=3pt,font=\footnotesize] {$\operatorname{Re}z$};
\draw[arr,black!55] (3.79,2.84)--(3.79,4.78) node[above,font=\footnotesize] {$\operatorname{Im}z$};
\fill[black!5] (3.79,3.55) circle[radius=.46];
\draw[black!55,line width=.65pt] (3.79,3.55) circle[radius=.46];
\draw[black!45] (4.00,3.95)--(4.37,4.39);
\node[font=\footnotesize,anchor=west] at (4.40,4.43) {$|z|=1$};
\fill[black!8] (2.2411328,3.55) circle[radius=0.1840000];
\draw[black!55,line width=.55pt] (2.2411328,3.55) circle[radius=0.1840000];
\fill (2.2411328,3.55) circle[radius=.041];
\fill[black!8] (2.8241272,3.55) circle[radius=0.1840000];
\draw[black!55,line width=.55pt] (2.8241272,3.55) circle[radius=0.1840000];
\fill (2.8241272,3.55) circle[radius=.041];
\fill[black!8] (4.7558728,3.55) circle[radius=0.1840000];
\draw[black!55,line width=.55pt] (4.7558728,3.55) circle[radius=0.1840000];
\fill (4.7558728,3.55) circle[radius=.041];
\fill[black!8] (5.3388672,3.55) circle[radius=0.1840000];
\draw[black!55,line width=.55pt] (5.3388672,3.55) circle[radius=0.1840000];
\fill (5.3388672,3.55) circle[radius=.041];
\draw[line width=.9pt] (0.7907448,3.485)--(0.9207448,3.615) (0.7907448,3.615)--(0.9207448,3.485);
\node[above=7pt] at (0.8557448,3.55) {$-z_{\mathrm f}$};
\draw[line width=.9pt] (6.6592552,3.485)--(6.7892552,3.615) (6.6592552,3.615)--(6.7892552,3.485);
\node[above=7pt] at (6.7242552,3.55) {$z_{\mathrm f}$};

\node[font=\footnotesize,align=center] at (3.76,2.28)
 {four distinct centers: $\pm2.100,\ \pm3.367$\\[2pt]
  possible focal values: $\pm z_{\mathrm f}\simeq\pm6.379$};
\node[font=\footnotesize,align=center] at (3.76,1.33)
 {$c_{j,\nu}=T_j(a_\nu)\in D\setminus B$\\[2pt]
  for both $\nu=1,2$ and all $1\le j<d$};
\node[font=\footnotesize,align=center] at (3.76,.39)
 {Dots and disks avoid both possible focal values.\\
  Coincident centers share one marker.};
\node[font=\footnotesize] at (9.28,5.02) {$\Delta_\varepsilon=\{|\xi|<\varepsilon\}$};
\node[font=\footnotesize] at (13.90,5.02) {$W_{j,\nu}\Subset D\setminus B$};
\fill[black!4] (9.30,3.68) circle[radius=.85];
\draw[line width=.65pt] (9.30,3.68) circle[radius=.85];
\fill[black!4] (13.92,3.68) circle[radius=.85];
\draw[line width=.65pt] (13.92,3.68) circle[radius=.85];
\fill (9.30,3.68) circle[radius=.045];
\fill (13.92,3.68) circle[radius=.045];
\node[below=5pt] at (9.30,3.68) {$0$};
\node[below=5pt] at (13.92,3.68) {$c_{j,\nu}$};
\fill (9.62,4.07) circle[radius=.028];
\fill (14.19,4.02) circle[radius=.028];
\node[above left=2pt,font=\footnotesize] at (9.62,4.07) {$\xi$};
\draw[arr] (10.31,3.70)--(12.91,3.70);
\node[above=6pt] at (11.61,3.70) {$T_j\circ u_\nu$};
\node[font=\footnotesize] at (9.30,2.44) {the center is included};
\node[font=\footnotesize] at (13.92,2.44) {$\mathcal V_{j,\nu}\in\mathcal O(W_{j,\nu})$};
\node[font=\footnotesize,align=center] at (11.89,1.41)
 {$\displaystyle F_{\mathfrak f}(\xi)=C_\rho-
       \sum_{j=1}^{d-1}\mathcal V_{j,\nu}(T_j(u_\nu(\xi)))$};
\node[font=\footnotesize,align=center] at (11.88,.40)
 {Holomorphic at $0$ for the specified source germ;\\
  different histories need not give the same value there.};
\end{tikzpicture}
}%
\caption{Ordinary companion centers supply the full-disk extension at infinity.
(a) Computed centers for $m=16/25$ and $\rho=2/5$, with the ellipse determined
by $J=2K\rho$. There are four distinct centers among the eight companions;
the illustrated disks have radius $0.4$ in the physical $z$-plane and avoid
both the unit disk and the possible focal values. These are geometric
illustration data, not a computation of the unknown displacement.
(b) The inverse and target neighborhoods are schematic. Each center
$c_{j,\nu}$ is ordinary after focal confinement, so the solved trace is
holomorphic on a full disk in $1/z$, unlike the punctured-disk argument of
Theorem~\ref{cont:thm:task6-finite-bad-set}. The neighborhoods are chosen
before the source germ; the representative functions may depend on its
history. No uniform inverse radius is certified by the picture.}
\label{fig:ordinary-centers-infinity}
\end{figure}
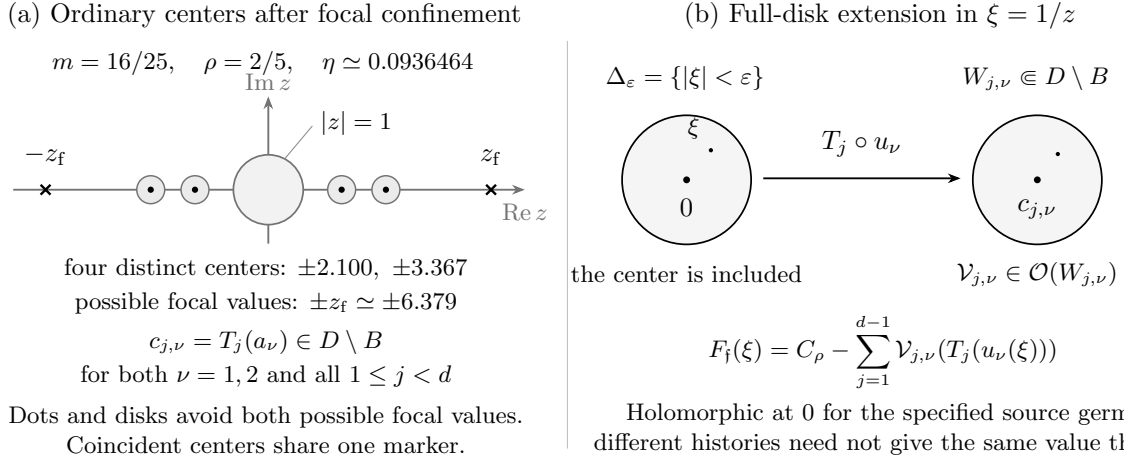

\Needspace{30\baselineskip}
Figure~\ref{fig:ordinary-centers-infinity} illustrates the additional
ordinary-center input used in the second part and the resulting extension
in the coordinate $1/z$.


\paragraph{Role of the finite obstruction set.}
Finiteness of the effective bad \emph{base values} is essential to the
finite-grid argument. No finiteness of the number of analytic branches
is used. The remaining focal germs may at this stage still have
arbitrary monodromy. They are treated by the parameter-dependent
trace argument in Section~\ref{cont:sec:logarithmic-closure}.

\section{Locked partners and elimination of monodromy}
\label{cont:sec:logarithmic-closure}

We now use the exterior continuation to eliminate monodromy. Let
$E$ be a fixed noncircular ellipse with semiaxes $a>b>0$, and let
$I\subset(0,1/2)$ be a nonempty open rotation interval. Put
$\eta=1-b^2/a^2$ and $c=\sqrt{a^2-b^2}>0$, and write
$z=e^{i\varphi}$ for the physical complex coordinate.
The two exterior zeros of $D_E$ are
\[
 A_+=\{(a+b)/c,-(a+b)/c\}.
\]
Write
\[
 \mathcal D=\{z:|z|>1\}\cup\{\infty\},\qquad
 \mathcal D^*=\mathcal D\setminus A_+.
\]
The preceding sections give the following properties of the physical
germ $V$:
\begin{enumerate}
\item $V$ is holomorphic in an annulus containing the unit circle and
satisfies $V(e^{i\varphi})\in\mathbb R$ for real $\varphi$.
\item $V$ continues along every path in $\mathcal D^*$ starting in its
original exterior collar. The continuation is not assumed single valued.
Every germ reached near infinity extends holomorphically across infinity.
\item For every reduced odd-denominator rational $\rho\in I$, the
physical rational trace is constant. Its identity continues along every
path in the corresponding lower Jacobi cylinder that starts in a
sufficiently small collar of the real boundary, where all summands use
their original annulus germs, and whose shifted physical paths avoid $A_+$.
\end{enumerate}
The extensions of every branch at infinity allow continuation along
paths through infinity, by the disk-chain and detour argument for
apparent-point removal in Lemma~\ref{cont:lem:task6-intrinsic-bad-set}.
Thus the second property applies on the domain with infinity adjoined.
The continuation assertion in the third property follows from the
second and the initial trace constancy, by the identity theorem:
begin with the constant trace in this collar and continue each summand
along its prescribed shifted physical path. These paths determine the
resulting branches, and their sum remains constant.

Starting from these inputs, we first lock a partner and its covering
determination on a parameter interval. Periodicity and fixed-contour Laurent
coefficients then give joint removability of the paired sum; parameter
differentiation descends each local derivative meromorphically. The global
topology identifies these derivatives as one differential, after which the
paired identities cancel each focal residue and the period argument yields
a single-valued meromorphic primitive. Reality is used only at the final
reflection step. Neither local paired holomorphy nor local derivative
descent alone gives the global conclusion.

\begin{theorem}[Infinite-sheet closure from the exterior continuation inputs]
\label{cont:thm:logarithmic-closure}
Under these three properties, the physical germ $V$ is single-valued and
meromorphic on $\mathcal D$, with at most simple poles at $A_+$ and no
other poles. It glues by reality to a rational function on the sphere,
with at most simple poles at the four zeros of $D_E$, regular at zero and
infinity. Consequently its displacement is an ellipse tangent.
\end{theorem}

A local paired identity alone does not exclude logarithmic branching. For example,
choose branches
\[
 f(y)=\operatorname{arcosh}(1/y)=\log(2/y)+a(y),\qquad g(y)=-f(y),
\]
where $a$ is holomorphic at zero. For every holomorphic involution with
$T(0)=0$ and $T'(0)=-1$, the function
\[
 f(y)+g(T(y))=\log(T(y)/y)+a(y)-a(T(y))
\]
is holomorphic at zero, after a consistent choice of the logarithm of
$T(y)/y$. These are opposite branches of one multivalued analytic
function, and their local monodromy is infinite. Thus the paired identities must be combined with global continuation
to eliminate the logarithmic term.

\begin{lemma}[Partner branches can be locked without a sheet bound]
\label{cont:lem:locked-partner}
Fix $a_0\in A_+$. Choose a local branch $\varphi=-i\log z$ near
$a_0$, put $\varphi_0=\varphi(a_0)$, and use the fixed physical
coordinate $y=\varphi-\varphi_0$. Let $f$ be any accessible branch of
$V$ on the universal cover of a punctured neighborhood of $a_0$.
There are a nonempty open interval $M_{\mathbb R}\subset(\eta,1)$
with $\rho(M_{\mathbb R})\subset I$, a fixed companion branch $g$
at $a_0$, and the analytic family of local involutions $T_m$ from the
focal Jacobi pair such that
\[
 f(y)+g(T_m(y))
\]
is the lift of a holomorphic germ at $y=0$ whenever
$m\in M_{\mathbb R}$ has reduced rational rotation of odd denominator.
The lifts of the arguments to the universal cover of the punctured disk
are chosen consistently throughout the interval. Neither the companion
branch nor its lift is chosen anew as the denominator varies.
\end{lemma}
\begin{proof}
We use only the three continuation inputs above. In particular, every
accessible germ continues along every further path in $\mathcal D^*$:
prepend a path representing that germ to the further path. Uniqueness
along a fixed path and homotopy invariance apply, as in
Lemma~\ref{cont:lem:countable-gluing}, without any bound on the number of
branches. We first choose one pair of marked branches $f,g$ on a real
parameter interval. Only afterward will we complete the trace separately
at each admissible rational parameter; its regular remainder and its
radius may depend on that parameter.

\smallskip\noindent
\emph{Step 1: select the source lift and one companion at a fixed parameter.}
Recall $h(m)=2J(m)=4K(m)\rho(m)$, and fix
$m_0\in\rho^{-1}(I)\subset(\eta,1)$. The notation $K'=K(1-m)$ denotes
the complementary integral, not differentiation. Write
\[
 \mathcal P_m(u)=\sn(u\mid m)-i\cn(u\mid m),\qquad
 \mathcal C_-(m)=\{-2K'(m)<\operatorname{Im}u<0\}/(4K(m)\Z).
\]
Represent the prescribed determination of $f$ by a physical path from the original exterior collar
to a punctured coordinate disk at $a_0$, followed by a radial tail
$y=re^{i\theta}$, $r\downarrow0$, with a fixed lifted argument. Every
path-lifting or homotopy argument below is applied to a compact truncation
of this tail; the focal endpoint is not an ordinary endpoint of a
continuation path.

Choose the initial point $z_{\rm in}$ sufficiently close to the unit
circle, on the inverse sheet adjacent to the real Jacobi circle, so that
all real translates of its lift have physical images in the original
holomorphy annulus. The choice can be made uniformly for real $m$ near
$m_0$: real translation is parametrized by the compact circle
$\R/(4K(m)\Z)$, and all its boundary images lie on $|z|=1$.
Compactness in this circle and in a small closed parameter interval gives
a common thin source collar. Rebasing within the original physical
annulus preserves $f$, since the original function is holomorphic and
single-valued there.

The compact physical prefix can be made piecewise smooth and perturbed
away from infinity and the two exterior critical values of
$\mathcal P_{m_0}$. These are regular points for every accessible branch:
infinity is regular by hypothesis, and no focal value is critical by
Lemma~\ref{arith:focal}. Small detours preserve the source homotopy in
$\mathcal D^*$. Lift the resulting finite-valued path through the
unbranched restriction of the proper degree-two map in
Lemma~\ref{arith:jacobi}, starting at the selected initial sheet. Its
terminal tail tends to one focal preimage $u_s(m_0)$; no terminal sheet
is prescribed in advance.

By Lemma~\ref{arith:focal}, the two members of this focal fiber differ
by $h(m_0)$. Choose $\epsilon\in\{1,-1\}$ so that the other member is
\[
 u_c(m_0)=u_s(m_0)+\epsilon h(m_0)\pmod{4K(m_0)}.
\]
Perturb the compact lifted prefix to make its translate by
$\epsilon h(m_0)$ avoid all preimages of $A_+$. The added forbidden
points are finite. Make the detours through an open neighborhood avoiding
the source preimages of $A_+$, of infinity, and of the two critical
values, keeping the prefix endpoints and a sufficiently short terminal
tail fixed. Such detours exist by the path-perturbation argument preceding
Lemma~\ref{cont:lem:task6-trace-realization}. Near the focal endpoint,
unramifiedness gives source and companion coordinate neighborhoods
whose only focal preimage is their respective center. Their punctured
tails are therefore admissible. Projection of the detour homotopy
preserves the prescribed $f$.

The companion path $\mathcal P_{m_0}(\Gamma_{m_0}+\epsilon h(m_0))$
now selects one accessible branch $g$ at $a_0$, where $\Gamma_{m_0}$
is the modified source lift. Both initial germs are the original annulus
germs, joined to a common initial germ inside that annulus. This is the
choice of $g$; no rational parameter or pigeonhole argument over branches
has been used.

\smallskip\noindent
\emph{Step 2: identify the physical involution and control its local family.}
Keep the selected focal member and sign on a small parameter neighborhood.
The explicit focal formulas $K\pm J-iK'$ and $3K\pm J-iK'$ give
analytic representatives $u_s(m)$ and $u_c(m)$, with
$u_c(m)=u_s(m)+\epsilon h(m)$ modulo $4K(m)$. Their derivatives under
$\mathcal P_m$ are nonzero, since at either focal preimage
\[
 \mathcal P_m'(u)=i\dn(u\mid m)\mathcal P_m(u),\qquad
 \dn^2(u\mid m)=1-m/\eta\ne0.
\]
In the fixed physical coordinate $z=a_0e^{iy}$, define
\[
 \psi_m(t)=-i\log\!\left(\frac{\mathcal P_m(u_s(m)+t)}{a_0}\right),
 \qquad \psi_m(0)=0,
\]
using the local logarithm equal to zero at $t=0$. Joint holomorphy and
$\psi_m'(0)\ne0$ give a jointly holomorphic local inverse after
shrinking to common parameter and coordinate disks.

The deck involution exchanges the two members of the fiber. Therefore
\[
 \mathcal P_m(u_c(m)+t)=\mathcal P_m(u_s(m)-t),\qquad
 T_m=\psi_m\circ(-\mathrm{id})\circ\psi_m^{-1}.
\]
Thus $T_m(0)=0$, $T_m'(0)=-1$, and $T_m\circ T_m=\mathrm{id}$ as
physical germs. Interchanging the two focal members replaces
$\psi_m(t)$ by $\psi_m(-t)$ and gives the same physical involution.
On a smaller common disk,
\[
 T_m(y)=-y+O(y^2),
\]
with the remainder uniform on a closed parameter neighborhood. In
particular, the companion images of a sufficiently short fixed source
ray lie in one fixed punctured sector around the opposite direction.
This gives tail control without claiming positive distance from the
limiting focal point.

\smallskip\noindent
\emph{Step 3: vary the parameter without changing either marked branch.}
Retain the modified physical source path from Step~1 and lift it for
$m$ near $m_0$, using the nearby initial inverse. Its compact prefix
avoids the moving critical values after a parameter restriction. A
finite chain of unramified inverse charts makes its lift, denoted by
$\Gamma_m$, continuous in $m$. The fixed radial tail lifts uniformly to $u_s(m)$ by the local
inverses in Step~2. The sign $\epsilon$ stays fixed.

In the common cylinder coordinate $u/(4K(m))$, the compact source and
companion prefixes at $m_0$ have positive distance from their respective
forbidden focal preimages. Continuity and a further restriction give
neighborhoods avoiding the corresponding moving forbidden points for
all parameters. The punctured tails are controlled by Step~2. Thus both
projected histories remain admissible on some nonempty open interval
$M_{\mathbb R}\Subset\rho^{-1}(I)$ about $m_0$. The source physical
path is fixed, so it still specifies $f$.

Here is the endpoint bookkeeping for the companion. Truncate its tail
at a fixed small source radius. Its starting point varies in the original
annulus, and its terminal point varies in the fixed punctured sector.
Join a fixed annulus point to the moving initial point by a continuously
chosen annulus connector. Join the moving terminal point to a fixed
sector point by a continuously chosen connector in that sector. After
shrinking the parameter interval, these choices are possible in fixed
local charts; any preceding fixed annulus connector is left unchanged.
The resulting companion histories form an endpoint-fixed homotopy in
$\mathcal D^*$. All-path continuation and homotopy invariance give the
same germ at the fixed sector point. Undoing the terminal connector
recovers its local continuations at the moving endpoints. Hence all
these histories specify the same $g$ on the universal cover.

The sector also fixes the lifted terminal arguments. There is no jump by
an integer winding as $m$ varies in this interval. We have therefore
chosen $M_{\mathbb R}$, $f$, $g$, and their covering determinations
before choosing a rational denominator. The interval may depend on the
original source history. No finiteness of the set of accessible branches
is involved.

\smallskip\noindent
\emph{Step 4: identify the two focal summands of each odd rational trace.}
Fix $m\in M_{\mathbb R}$ with $\rho(m)=p/d$ reduced and $d$ odd.
Then $d\ge3$ and $h=4Kp/d$. The full fiber of $a_0$ consists of
$u_s$ and $u_s+\epsilon h$ modulo $4K$. Since
\[
 jh\equiv kh\pmod{4K}\quad\Longleftrightarrow\quad d\mid(j-k),
\]
exactly two indices in $0\le j<d$ reach $a_0$: $j=0$ and
\[
 j_*=
 \begin{cases}
  1,&\epsilon=1,\\
  d-1,&\epsilon=-1.
 \end{cases}
\]
The opposite focal fiber is obtained by adding $2K$, because
$\mathcal P_m(u+2K)=-\mathcal P_m(u)$. A shifted endpoint in that
fiber would require
\[
 jh\equiv2K\pmod{4K}
 \quad\text{or}\quad
 (j-\epsilon)h\equiv2K\pmod{4K}.
\]
Equivalently, for some $r\in\Z$,
\[
 2jp=d(2r+1)
 \quad\text{or}\quad
 2(j-\epsilon)p=d(2r+1).
\]
The left sides are even and the right sides odd. Real shifts preserve
the lower cylinder, whose image is the spherical exterior by
Lemma~\ref{arith:jacobi}. Thus every endpoint with
$j\notin\{0,j_*\}$ lies in $\mathcal D^*$; it may be treated
as an ordinary spherical point, including infinity. Only the source and
its selected companion can contribute focal germs.

For $\epsilon=-1$, the companion path in the trace is exactly the
already selected backward path, not merely a path with the same endpoint:
\[
 (d-1)h=-h+4Kp,\qquad
 \mathcal P_m(\Gamma_m(t)+(d-1)h)
       =\mathcal P_m(\Gamma_m(t)-h).
\]
The initial branches agree as well, since both are the original annulus
branch. Thus the $j_*$ summand has the same continuation as the locked
companion even though its index depends on $d$.

\smallskip\noindent
\emph{Step 5: complete the trace by detours preserving both histories.}
The other translated paths may meet focal preimages along their compact
prefixes. For this fixed parameter, define
\[
 \mathcal A_m=\{u\in\mathcal C_-(m):
       \mathcal P_m(u)\notin A_+,\quad
       \mathcal P_m(u+\epsilon h)\notin A_+\},
\]
and the finite set
\[
 \mathcal F_{m,d}
     =\bigcup_{j=0}^{d-1}
         \{u\in\mathcal C_-(m):\mathcal P_m(u+jh)\in A_+\}.
\]
The locked source path lies in $\mathcal A_m$. Finiteness of
$\mathcal F_{m,d}$ follows from degree two and the finite number of
shifts; no bound uniform in $d$ is asserted. Step~4 permits shortening
the terminal tail so that every unpaired evaluation avoids $A_+$ there.
Leave that shortened tail and the initial point fixed.

The remaining prefix is compact in the open set $\mathcal A_m$.
Perturb it, relative to its endpoints, to avoid $\mathcal F_{m,d}$,
through a homotopy in a sufficiently small neighborhood contained in
$\mathcal A_m$. The elementary path-perturbation argument applies on
this open Riemann surface. The source projection and the projection
translated by $\epsilon h$ both send the homotopy into
$\mathcal D^*$, with endpoints fixed. Hence the detours preserve both
$f$ and $g$, including their marked lifts, not only the source germ.
All other shifted paths now avoid $A_+$ and admit continuation.
Poles of a shifted map are allowed as spherical evaluation points, since
every accessible branch is holomorphic at infinity by the second input.

This step is performed separately for each admissible rational
parameter. Its detours, regular companion branches, and terminal radius
may depend on that parameter. They do not change the pair locked in
Step~3, and no common terminal radius over unbounded denominators is
needed.

\smallskip\noindent
\emph{Step 6: express the paired sum as a holomorphic physical germ.}
For the fixed rational parameter, the original trace identity holds
near the chosen initial point of the source collar. Indeed, all real
translates there have images in the original annulus by Step~1; the
real trace identity extends into that collar by the identity theorem.
Continue its summands along the paths completed in Step~5. A common
finite subdivision of each compact truncation gives a chain of source
neighborhoods on which all compositions are defined. The identity
theorem carries their sum as the same constant $C_{\rho(m)}$.

For $j\notin\{0,j_*\}$ choose a small simply connected spherical
neighborhood $N_{j,m}$ of $\mathcal P_m(u_s+jh)$ with closure in
$\mathcal D^*$. Shrink the terminal coordinate disk so that its
$j$th image lies in $N_{j,m}$. The corresponding accessible target
germ continues along every path in $N_{j,m}$, by concatenation with
its representing history. The monodromy theorem gives a holomorphic
representative $W_{j,m}$ on that neighborhood. At infinity this uses
the ordinary coordinate $1/z$, as supplied by the continuation inputs.
The representatives may depend on the completed rational trace.

Put $t_m(y)=\psi_m^{-1}(y)$ and define
\[
 H_m(y)=C_{\rho(m)}-
  \sum_{\substack{0\le j<d\\ j\notin\{0,j_*\}}}
   W_{j,m}\!\left(\mathcal P_m(u_s(m)+t_m(y)+jh(m))\right).
\]
Each summand is holomorphic at $y=0$. If its target is infinity, write
it using the holomorphic coordinate $1/\mathcal P_m$ near the
corresponding source point. Unramifiedness at $u_s(m)$ makes $t_m$ a
holomorphic function of the fixed physical coordinate. The source term
is $f(y)$, and Steps~2 and 4 identify the paired term with
$g(T_m(y))$ on the fixed determinations. Thus the continued trace gives
\[
 f(y)+g(T_m(y))=H_m(y)
\]
on a lifted punctured sector. This proves holomorphy of the paired
\emph{germ}; its disk of extension may depend on the rational parameter.

\smallskip\noindent
\emph{Step 7: record the consistent covering determinations and their scope.}
On the universal cover write $y=e^w$, and represent the two locked
branches by holomorphic functions $F(w)$ and $G(w)$ on sufficiently
far-left half-planes. The all-path input on a punctured physical disk
justifies these representatives. The jointly holomorphic function
$T_m(y)/y$ extends to $y=0$ with value $-1$ and is nonzero on a smaller
common parameter-coordinate product. Choose its logarithm
\[
 L_m(y)=\log(T_m(y)/y)
\]
to agree at $m_0$ with the companion determination selected by the two
paths. A simply connected parameter disk and a coordinate disk give a
jointly holomorphic choice. Its value $L_m(0)=(2k+1)\pi i$ has one
fixed integer $k$, and the sector homotopies in Step~3 preserve this
choice throughout a smaller $M_{\mathbb R}$. The lifted argument is
$\widetilde T_m(w)=w+L_m(e^w)$. For each admissible rational parameter,
the paired identity reads
\[
 F(w)+G(\widetilde T_m(w))=H_m(e^w).
\]
For such a parameter the equality first holds on the selected
lifted sector. Since both sides are holomorphic on a sufficiently
far-left half-plane, the identity theorem extends it to that half-plane.
The required half-plane may depend on the extension radius of $H_m$.
The logarithm is never reset independently as $m$ or $d$ changes.

The involution assertion is in the physical coordinate. Its selected
covering lift need not be an involution: $T_m^2=\mathrm{id}$ implies
\[
 L_m(y)+L_m(T_m(y))=2(2k+1)\pi i,
 \qquad
 \widetilde T_m^2(w)=w+2(2k+1)\pi i
\]
on smaller domains where these compositions are defined. This follows
by exponentiating the first sum, which is constant with its value at
zero. Nothing in the proof requires this deck translation to act
trivially on $f$ or $g$.

We have fixed one companion and one covering determination on an open
parameter interval, and obtained the paired holomorphic germ at every
permitted odd-denominator rational parameter there. The proof does not
assert joint removability at irrational or complex parameters, a common
remainder radius over all denominators, or elimination of logarithmic or
global monodromy. These are not needed to lock the two branches.
\end{proof}

\begin{lemma}[Joint descent from the universal cover]
\label{cont:lem:universal-puncture}
In the situation of Lemma~\ref{cont:lem:locked-partner}, each source and
companion branch has a derivative which descends meromorphically in the
physical coordinate, with pole order at most two. Each branch itself
has the form
\[
 \frac{A}{y}+B\log y+h(y),
\]
where $h$ is holomorphic at zero and the logarithm is taken on its
universal covering domain. The logarithmic coefficients of the paired
branches are opposite.
\end{lemma}
\begin{proof}
Fix an exterior focal value $a_0$ and an arbitrary accessible source branch
$f$. Use the interval $M_{\mathbb R}$, the fixed companion branch $g$, and
the covering determinations supplied by
Lemma~\ref{cont:lem:locked-partner}. Throughout the proof the physical
coordinate $y=\varphi-\varphi_0$, the angle $\varphi_0$, and $\eta$ are
fixed when $m$ is differentiated. Neither branch is selected anew as $m$
varies. We first descend the paired sum, then remove its puncture jointly,
and only then differentiate in the parameter.

\smallskip\noindent
\emph{Step 1: compute the nondegenerate parameter variation of the involution.}
Let $u_s(m)$ be the selected focal preimage and let $\varphi_m$ be the local
physical-angle determination with $\varphi_m(u_s(m))=\varphi_0$. Thus
\[
 \psi_m(t)=\varphi_m(u_s(m)+t)-\varphi_0,\qquad
 T_m=\psi_m\circ(-\mathrm{id})\circ\psi_m^{-1}.
\]
Locally $\varphi_m=\operatorname{am}(\,\cdot\mid m)-\pi/2$ with a
compatible determination. The Jacobi identities give
\[
 \sn(u\mid m)=\cos\varphi_m(u),\qquad
 \cn(u\mid m)=-\sin\varphi_m(u),\qquad
 \partial_u\varphi_m=\dn(u\mid m).
\]
Consequently,
\[
 (\partial_u\varphi_m)^2=1-m\cos^2\varphi_m,\qquad
 \partial_u^2\varphi_m=m\sin\varphi_m\cos\varphi_m.
\]
Write $\psi_m(t)=a_mt+b_mt^2/2+O(t^3)$. Unramifiedness at the focal
point, proved in Lemma~\ref{arith:focal}, gives
\[
 a_m^2=1-m\cos^2\varphi_0\ne0,\qquad
 b_m=m\sin\varphi_0\cos\varphi_0.
\]
The inverse series and its conjugation by $t\mapsto-t$ are
\[
 \psi_m^{-1}(y)=\frac{y}{a_m}-\frac{b_m}{2a_m^3}y^2+O(y^3),\qquad
 T_m(y)=-y+\frac{b_m}{a_m^2}y^2+O(y^3).
\]
Thus
\[
 T_m(y)=-y+\beta(m)y^2+O(y^3),\qquad
 \beta(m)=\frac{m\sin\varphi_0\cos\varphi_0}
                   {1-m\cos^2\varphi_0},
\]
and differentiation at fixed $\varphi_0$ gives
\[
 \beta'(m)=\frac{\sin\varphi_0\cos\varphi_0}
                  {(1-m\cos^2\varphi_0)^2}\ne0.
\]
Here $\cos^2\varphi_0=1/\eta$ and
$\sin\varphi_0\cos\varphi_0\ne0$. Choose a complex parameter disk $M$
about a point of $M_{\mathbb R}$, with compact closure in a neighborhood
where the family is holomorphic, $M\cap\R\subset M_{\mathbb R}$, and
$1-m/\eta\ne0$. All coordinate disks below are centered at zero; write
$\Delta_r=\{y:|y|<r\}$. After shrinking, $T$ is jointly holomorphic on a
neighborhood of $\overline M\times\overline{\Delta_{r_T}}$ for some
$r_T>0$. Its remainder has the form $y^3E(m,y)$ with $E$ jointly
holomorphic. Therefore
\[
 \partial_mT_m(y)=\beta'(m)y^2+y^3\partial_mE(m,y),
\]
so this parameter derivative has a zero of exactly order two at $y=0$.
The motion of $u_s(m)$ is already included in $\psi_m$; its physical
value $\varphi_0$ is not varied in the formula for $\beta$.

\smallskip\noindent
\emph{Step 2: place the fixed branches on one common covering domain.}
The quotient $T_m(y)/y$ extends holomorphically across zero with value
$-1$. Shrink the product so that it is nonzero there, and choose
\[
 L_m(y)=\log\bigl(T_m(y)/y\bigr)
\]
on the simply connected parameter-coordinate product. Match its branch
to the fixed companion lift from Lemma~\ref{cont:lem:locked-partner}.
Then $L_m(0)=(2k+1)\pi i$ for one fixed $k\in\Z$. This determination
is not reset as the parameter changes. By shrinking within the preceding
product, arrange $|L_m(y)|\le C_L$ there for a finite constant $C_L$.

On the universal cover $y=e^w$, let $F(w)$ and $G(w)$ represent the
fixed branches $f$ and $g$. The all-path continuation input on a punctured
physical disk gives holomorphic representatives on half-planes
$\operatorname{Re}w<A_F$ and $\operatorname{Re}w<A_G$, respectively.
Choose $R>0$ so that
\[
 e^{-R}<r_T,\qquad -R<A_F,\qquad -R+C_L<A_G,
\]
and put $\mathbb H_R=\{w:\operatorname{Re}w<-R\}$ and $r=e^{-R}$.
The function
\[
 H(m,w)=F(w)+G\bigl(w+L_m(e^w)\bigr)
\]
is jointly holomorphic on $M\times\mathbb H_R$. Indeed, its shifted
argument lies in the half-plane of $G$. This uses no boundedness of $F$
or $G$ near the puncture or across their sheets; only the argument shift
is bounded.

\smallskip\noindent
\emph{Step 3: establish periodicity before taking Laurent coefficients.}
Let $\mathcal Q$ be the real parameters in $M$ for which $\rho(m)$ is a
reduced rational with odd denominator. This set is dense in $M\cap\R$:
odd-denominator rationals are dense, reduction preserves oddness of a
denominator, and the correspondence between $m$ and $\rho$ is a real
analytic diffeomorphism on the chosen interval, as in
\eqref{op:eq:rhoDerivative} and the relation between $m$ and $\lambda$.

Fix $m\in\mathcal Q$. The locked identity supplies a holomorphic germ
$h_m$ at zero with $H(m,w)=h_m(e^w)$ on a lifted sector. Its extension
radius may depend on $m$. On a sufficiently far-left half-plane both
sides are holomorphic, and the identity theorem in $w$ extends their
equality to that half-plane. In particular,
\[
 H(m,w+2\pi i)-H(m,w)=0
\]
there. Define
\[
 \mathscr P(m,w)=H(m,w+2\pi i)-H(m,w).
\]
For the fixed rational parameter, this difference is holomorphic on the
original common half-plane $\mathbb H_R$. It vanishes on a farther-left
open subset, and hence on all of $\mathbb H_R$. For each fixed
$w\in\mathbb H_R$, the function $m\mapsto\mathscr P(m,w)$ is holomorphic
on $M$ and vanishes on $\mathcal Q$. The identity theorem in $m$ yields
\[
 H(m,w+2\pi i)=H(m,w)
 \qquad(m\in M,\ w\in\mathbb H_R).
\]
Thus $H$ descends to a single-valued jointly holomorphic function
$\widehat H$ on $M\times\Delta_r^*$, where
$\Delta_r^*=\Delta_r\setminus\{0\}$. To see joint holomorphy, use a local
logarithm of $y$ and write $\widehat H(m,y)=H(m,\log y)$; periodicity
makes these definitions agree on overlaps. Only the paired sum has been
shown periodic. No positive lower bound on the rational extension radii,
and no periodicity of $F$ or $G$ individually, has been assumed.

\smallskip\noindent
\emph{Step 4: remove the puncture jointly by fixed-contour Laurent coefficients.}
Choose $0<r_0<r$. For each $n\in\Z$, define on the positively oriented
fixed circle
\[
 c_n(m)=\frac1{2\pi i}\int_{|\xi|=r_0}
                 \widehat H(m,\xi)\xi^{-n-1}\,d\xi.
\]
Joint holomorphy near this contour implies that every $c_n$ is holomorphic
on $M$, with parameter differentiation justified locally uniformly on
compact subsets. The coefficient is independent of the circle's radius
within the punctured disk.

For $m\in\mathcal Q$, the function $\widehat H(m,y)$ agrees near zero
with $h_m$. If $n<0$, move the coefficient contour into the extension
disk of this particular $h_m$. The integrand
$h_m(\xi)\xi^{-n-1}$ is holomorphic at zero, so its integral is zero.
Hence $c_n$ vanishes on $\mathcal Q$, and the identity theorem in $m$
gives $c_n\equiv0$ on $M$ for every $n<0$. The coefficient functions were
defined on one fixed circle; only this pointwise verification uses a
parameter-dependent smaller circle.

For a compact set $K_M\Subset M$, put
\[
 C_{K_M,r_0}=\sup_{m\in K_M,\,|\xi|=r_0}|\widehat H(m,\xi)|<\infty.
\]
Cauchy's estimate gives $|c_n(m)|\le C_{K_M,r_0}r_0^{-n}$ for $n\ge0$.
Consequently, for $0<r_1<r_0$,
\[
 \sum_{n\ge0}\sup_{m\in K_M,\,|y|\le r_1}|c_n(m)y^n|
 \le \frac{C_{K_M,r_0}}{1-r_1/r_0}.
\]
The nonnegative Laurent series therefore converges normally on smaller
parameter-coordinate products and defines a jointly holomorphic extension
across $y=0$. It agrees with the original $\widehat H$ away from zero.
Letting $r_0$ range below $r$ patches these extensions to
$M\times\Delta_r$; retain the notation $\widehat H$. In particular,
$\partial_m\widehat H$ is now holomorphic at the puncture. No bound near
zero was used to assume removability in advance.

\smallskip\noindent
\emph{Step 5: differentiate at a fixed physical coordinate on the cover.}
For $y=e^w$, set $w_c=w+L_m(y)$, so that $e^{w_c}=T_m(y)$. Before
asserting any descent for the companion derivative, denote its physical
coordinate expression on the cover by
\[
 D_G(v)=e^{-v}G'(v).
\]
The prime on $G$ differentiates its covering argument. Since $F$ and $G$
are fixed as functions of that argument, the chain rule gives
\[
 \partial_mH(m,w)=G'(w_c)\,\partial_mL_m(y),\qquad
 \partial_mT_m(y)=e^{w_c}\partial_mL_m(y).
\]
Using $H(m,w)=\widehat H(m,e^w)$, we obtain
\[
 \partial_m\widehat H(m,y)=D_G(w_c)\,\partial_mT_m(y)
\]
on the cover. In branch notation this is
$\partial_m\widehat H=g'(T_m(y))\partial_mT_m(y)$. The factor
$e^{-w_c}$ converting $G'$ to the physical derivative cancels the factor
$e^{w_c}$ from differentiating $T_m$. There is no extra factor of $y$
and no $\partial_mg$ term. Holding $w$ fixed here holds $y=e^w$ fixed;
the continuation lifts do not vary independently with $m$.

\smallskip\noindent
\emph{Step 6: descend the companion derivative through a meromorphic quotient.}
Fix $\mu\in M$ and write
\[
 \dot T_\mu(y):=\left.\partial_mT_m(y)\right|_{m=\mu}
       =y^2d_\mu(y),\qquad d_\mu(0)=\beta'(\mu)\ne0.
\]
On a smaller coordinate disk $d_\mu$ has no zero. The quotient
\[
 Q_\mu(y)=
 \frac{\left.\partial_m\widehat H(m,y)\right|_{m=\mu}}{\dot T_\mu(y)}
\]
is thus a single-valued meromorphic germ with pole order at most two.
As $T_\mu$ is a biholomorphic physical germ fixing zero, the function
\[
 R_\mu(\xi)=Q_\mu(T_\mu^{-1}(\xi))
\]
is also meromorphic with pole order at most two. Step~5 identifies it
with the selected physical derivative of $g$ on a sector. More explicitly,
$D_G(v)$ and $R_\mu(e^v)$ are holomorphic on a sufficiently far-left
connected covering half-plane and agree on a nonempty open subset,
obtained from $v=w+L_\mu(e^w)$. This change of covering coordinate is
locally biholomorphic near the puncture, since its derivative tends to
one. The identity theorem therefore gives
\[
 D_G(v)=R_\mu(e^v)
\]
on that half-plane. This proves descent of the full local derivative
branch, not merely an equality on its initial sector. The inverse used
here is the physical inverse $T_\mu^{-1}$; no assertion that its selected
covering lift squares to the identity is required. One fixed parameter
$\mu$ suffices for the conclusion.

\smallskip\noindent
\emph{Step 7: integrate and check the pole and logarithm of both branches.}
Write the descended derivative as
\[
 g'(y)=\gamma_{-2}y^{-2}+\gamma_{-1}y^{-1}+r_g(y),
\]
where $r_g$ is holomorphic at zero. Set $A=-\gamma_{-2}$ and $B=\gamma_{-1}$.
Integrating the holomorphic remainder on that disk and the two displayed
Laurent terms on the universal cover gives
\[
 g(y)=\frac Ay+B\log y+h_g(y),
\]
where $h_g$ is holomorphic at zero. Equivalently,
$G(v)=Ae^{-v}+Bv+h_g(e^v)$ on a sufficiently far-left half-plane.
The difference between this primitive and $G$ has zero derivative on a
connected covering domain, so it is one constant, absorbed into $h_g$;
there is no further sheet-dependent term.

At the fixed parameter $\mu$, the paired identity reads
$f(y)=\widehat H(\mu,y)-g(T_\mu(y))$ with the locked determination.
The identities
\[
 \log T_\mu(y)=\log y+L_\mu(y),\qquad
 \frac1{T_\mu(y)}+\frac1y\ \text{is holomorphic at }0
\]
follow from the chosen lift and $T_\mu(y)=-y+O(y^2)$. Therefore
\[
 f(y)=\frac Ay-B\log y+h_f(y),
\]
with the explicit holomorphic remainder
\begin{align*}
 h_f(y)={}&\widehat H(\mu,y)
   -A\left(\frac1{T_\mu(y)}+\frac1y\right)\\
 &-B L_\mu(y)-h_g(T_\mu(y)).
\end{align*}
In particular,
\[
 g'(y)=-\frac A{y^2}+\frac By+h_g'(y),\qquad
 f'(y)=-\frac A{y^2}-\frac By+h_f'(y).
\]
Both derivatives descend meromorphically, with pole order at most two;
the paired simple-pole coefficients agree and the logarithmic coefficients
are opposite. The coordinate change $dy/dz=1/(iz)$ is holomorphic and
nonzero near $a_0$, so the derivative with respect to $z$ has the same
pole-order bound. The original source branch and the focal value were
arbitrary, which proves the asserted local conclusion for every
accessible branch at either focal point. This argument does not eliminate
$B$ or identify derivatives obtained through different global histories;
those global conclusions are left to the following lemma.
\end{proof}

\begin{figure}[!htb]
\centering
\begin{tikzpicture}[x=1cm,y=1cm,font=\small,
 arr/.style={-{Stealth[length=1.7mm]},line width=.7pt},
 box/.style={draw=black!55,rounded corners=1.5pt,fill=black!3,
   text width=7.10cm,font=\footnotesize,align=center,inner sep=5pt},
 every node/.style={inner sep=2pt}]
\node at (3.58,6.87) {(a) Physical involution versus covering lift};
\node at (11.68,6.87) {(b) Descend, remove, then differentiate};
\draw[black!25,line width=.4pt] (7.57,.10)--(7.57,6.53);
\fill[black!3] (.45,1.59) rectangle (3.02,5.95);
\draw[dashed,black!50] (3.02,1.59)--(3.02,5.95);
\node[above,font=\footnotesize] at (1.64,6.01) {$w$-plane: $\operatorname{Re}w<\log r$};
\draw[arr] (1.16,2.03)--(1.16,3.56);
\draw[arr] (1.16,3.63)--(1.16,5.17);
\fill (1.16,2.03) circle[radius=.045];
\fill (1.16,3.60) circle[radius=.045];
\draw[fill=white,line width=.75pt] (1.16,5.17) circle[radius=.055];
\node[right=5pt,font=\footnotesize] at (1.16,2.03) {$w_0$};
\node[right=5pt,font=\footnotesize] at (1.16,3.60) {$w_0+\pi i$};
\node[right=5pt,font=\footnotesize] at (1.16,5.17) {$w_0+2\pi i$};
\node[left=3pt,font=\footnotesize] at (1.16,2.83) {$\widetilde T$};
\node[left=3pt,font=\footnotesize] at (1.16,4.43) {$\widetilde T$};
\draw[arr] (3.24,3.60)--(4.08,3.60);
\node[above,font=\footnotesize] at (3.67,3.73) {$e^w$};
\node[above,font=\footnotesize] at (5.72,6.01) {$y$-plane: $0<|y|<r$};
\draw[fill=black!3,draw=black!55] (5.66,3.60) circle[radius=1.36];
\draw[fill=white,line width=.75pt] (5.66,3.60) circle[radius=.06];
\node[below=4pt,font=\footnotesize] at (5.66,3.60) {$0$};
\fill (6.47,3.60) circle[radius=.045];
\fill (4.85,3.60) circle[radius=.045];
\node[right=3pt,font=\footnotesize] at (6.47,3.60) {$y_0$};
\node[left=3pt,font=\footnotesize] at (4.85,3.60) {$-y_0$};
\draw[arr] (6.47,3.64) arc[start angle=0,end angle=180,radius=.81];
\draw[arr] (4.85,3.56) arc[start angle=180,end angle=360,radius=.81];
\node[font=\footnotesize,fill=white] at (5.66,4.55) {$T$};
\node[font=\footnotesize,fill=white] at (5.66,2.65) {$T$};
\node[align=center,font=\footnotesize,text width=7.15cm] at (3.60,.76)
 {Exact model: $T(y)=-y$, $\widetilde T(w)=w+\pi i$.\\[3pt]
  $T^2(y)=y$, but $\widetilde T^2(w)=w+2\pi i$.\\[3pt]
  A returned physical point need not carry the same branch.};
\node[box,anchor=north] (fixed) at (11.68,6.38)
 {Fix branches $F,G$ and the lift $L_m$.\\[2pt]
  $H(m,w)=F(w)+G(w+L_m(e^w))$};
\node[box,anchor=north] (period) at ($(fixed.south)+(0,-.30)$)
 {Odd-rational paired identities imply\\[2pt]
  $H(m,w+2\pi i)=H(m,w)$;\\
  $\widehat H\in\mathcal O(M\times\Delta_r^*)$.};
\draw[arr] (fixed.south)--(period.north);
\node[box,anchor=north] (remove) at ($(period.south)+(0,-.30)$)
 {Fixed-contour Laurent coefficients:\\[2pt]
  $c_n\equiv0$ for $n<0$.\\
  Normal convergence gives $\widehat H\in\mathcal O(M\times\Delta_r)$.};
\draw[arr] (period.south)--(remove.north);
\node[box,anchor=north] (differentiate) at ($(remove.south)+(0,-.30)$)
 {Differentiate at fixed physical $y$:\\[3pt]
  $\partial_m\widehat H=g'(T_m(y))\,\partial_mT_m(y)$;\\[3pt]
  $\partial_mT_m(y)=y^2d_m(y)$, $d_m(0)\ne0$.\\
  At most a double pole for the derivative.};
\draw[arr] (remove.south)--(differentiate.north);
\end{tikzpicture}
\caption[Covering lifts and the order of joint local descent]{Local descent in Lemmas~\ref{cont:lem:locked-partner}--\ref{cont:lem:universal-puncture}.
Left: the exact linear model $T(y)=-y$, with the lift
$\widetilde T(w)=w+\pi i$ under $y=e^w$, separates the physical involution
from its non-involutive covering lift. The vertical separations are scaled;
this is not a plot of the parameter-dependent Jacobi involution. In the
actual family the fixed determination satisfies
$\widetilde T_m^2(w)=w+2(2k+1)\pi i$, and no trivial action of that deck
translation on the unknown branches is assumed. Right: only the paired sum
is first shown periodic. Fixed-contour Laurent coefficients give joint
removability before parameter differentiation yields the local derivative
bound. The separate rational extension radii need not have a common positive
lower bound; this stage does not eliminate the logarithmic coefficients.}
\label{fig:local-joint-descent}
\end{figure}
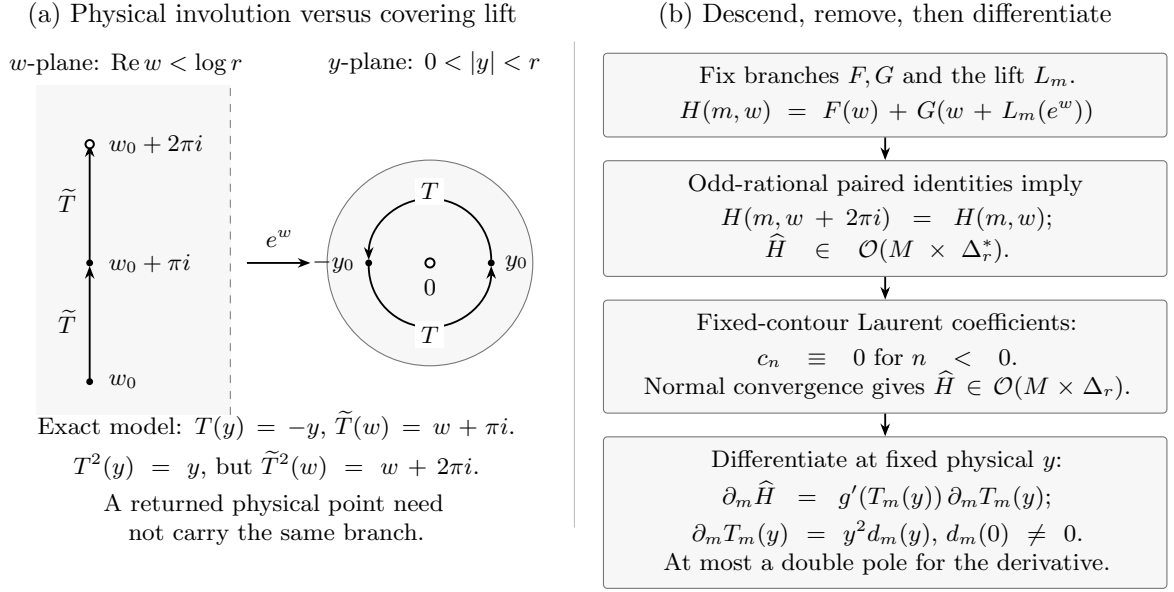

Figure~\ref{fig:local-joint-descent} separates the covering-space issue from
the order of the local parameter argument.

\begin{lemma}[Global derivative descent removes the logarithms]
\label{cont:lem:global-derivative}
Under the continuation inputs, $V$ is single-valued meromorphic on
$\mathcal D$, with at most simple focal poles.
\end{lemma}
\begin{proof}
Write $A_+=\{a_+,a_-\}$, where
$a_\pm=\pm(a+b)/\sqrt{a^2-b^2}$, and retain the domains
$\mathcal D$ and $\mathcal D^*=\mathcal D\setminus A_+$ above.
Choose a finite base point $z_0$ in the original exterior collar and a
holomorphic representative $V_0$ of the initial germ there. We first
show that the continued differential is single-valued, and only then
prove that its periods vanish and construct a global primitive.

\smallskip\noindent
\emph{Step 1: continue the differential without assuming global exactness.}
If holomorphic elements representing a continuation of $V$ agree on an
overlap, their differentials agree there as well. Differentiation
therefore gives continuation of the initial germ $dV_0$ along every
path in $\mathcal D^*$. Conversely, a continued differential germ is
the derivative of the continuation of $V_0$ along the same path, by
uniqueness of continuation. We do not yet regard $dV$ as the differential
of a globally single-valued function.

At a focal point $a_0\in A_+$ use the coordinate
$y=-i\log(z/a_0)$ with $y(a_0)=0$. The change of coordinate is
biholomorphic and $dy/dz=1/(iz)$ is holomorphic and nonzero there.
Lemma~\ref{cont:lem:universal-puncture} shows that the derivative of
every accessible branch is a single-valued meromorphic germ in this
coordinate, with pole order at most two. Thus a sufficiently small
loop around $a_0$ fixes the derivative of that branch. The conclusion
holds for every accessible branch, not merely the initial one; a common
local radius for all branches is not needed.

Infinity is an ordinary point of the continued differential. Indeed,
for any branch $V_\gamma$ reached near infinity, the continuation input
gives a holomorphic function
$\widetilde V_\gamma(\xi)=V_\gamma(1/\xi)$ near $\xi=0$. On the
punctured coordinate disk,
\[
 dV_\gamma
 =-\xi^{-2}V_\gamma'(1/\xi)\,d\xi
 =\widetilde V_\gamma'(\xi)\,d\xi,
 \qquad \xi=1/z.
\]
The right side is holomorphic at zero. Working with a one-form makes
this regularity independent of the finite coordinate $z$.

\smallskip\noindent
\emph{Step 2: eliminate the monodromy of the differential.}
Inversion $z\mapsto1/z$, with infinity mapped to zero, identifies
$\mathcal D$ with the unit disk and $\mathcal D^*$ with that disk
minus $a_+^{-1}$ and $a_-^{-1}$. Its fundamental group is generated by
two small focal loops joined to $z_0$ by connecting paths. There is no
additional puncture at infinity, by Step~1.

Consider a based generator $\alpha*\ell*\alpha^{-1}$, where $\alpha$
leads to a punctured neighborhood of a focal point and $\ell$ winds
once around that point. Continue $V_0$ and its differential along
$\alpha$. Choose the terminal part of this connector and the small
loop inside a punctured disk on which the derivative of the reached
branch has the meromorphic representative from Step~1. Such a change
of the terminal connector and loop radius preserves the generator's
homotopy class. Continuation around $\ell$ returns the same derivative
germ; reversing $\alpha$ then returns $dV_0$. This argument applies
also to a derivative branch reached after any preceding history.

Homotopy invariance applies because the initial germ has continuation
along every path in $\mathcal D^*$, as in the finite-chain argument of
Lemma~\ref{cont:lem:countable-gluing}. Each generator and its inverse
act trivially on the derivative germ, hence so does every word in the
generators and every loop. Continuation of $dV_0$ is therefore path
independent. Denote its descended holomorphic one-form on
$\mathcal D^*$ by $\omega$. Its local representatives from Step~1
extend it meromorphically across both focal points with pole order at
most two. It is holomorphic at infinity.

At this stage all accessible branches of $V$ have differential
$\omega$, but they could still differ by additive constants after
continuation around loops. Writing $\omega$ for the continued
differential does not yet assert that it has a single-valued primitive.

\smallskip\noindent
\emph{Step 3: express the paired identity as a local identity of one-forms.}
Fix $a_0\in A_+$, and write the physical chart and local differential as
\[
 \jmath(y)=a_0e^{iy},\qquad
 \omega_0=\jmath^*\omega=q(y)\,dy.
\]
Choose a pair of accessible branches $f,g$ at this same focal point
and one admissible rational parameter $m$ as in
Lemma~\ref{cont:lem:locked-partner}. The interval in that lemma contains
such parameters. By Lemma~\ref{cont:lem:universal-puncture}, their paired
identity, with its fixed covering determinations, reads
\[
 f(y)+g(T_m(y))=\widehat H(m,y),
 \qquad T_m(0)=0,\quad T_m'(0)=-1,
\]
where $\widehat H(m,\cdot)$ is holomorphic at zero. The parameter $m$
is held fixed throughout the rest of this local argument.

Global differential descent gives $df=dg=\omega_0$ on their local
covering domains, even if the functions themselves differ. Differentiating
the paired identity in the physical coordinate gives
\[
 d_y\widehat H(m,\cdot)
 =\bigl(q(y)+q(T_m(y))T_m'(y)\bigr)\,dy
 =\omega_0+T_m^*\omega_0.
\]
Initially this is an identity on a marked sector. Both differential
expressions are single-valued meromorphic germs, so the identity
theorem gives it on a sufficiently small punctured disk. The pullback
uses the physical biholomorphism $T_m$; no assertion about the square
of its chosen covering lift is needed. Both terms are germs at the
same focal point, which is essential for the residue comparison.

\smallskip\noindent
\emph{Step 4: verify the residue sign and remove the logarithmic residues.}
Write
\[
 \omega_0=
 \left(\frac{c_{-2}}{y^2}+\frac{r}{y}+h(y)\right)dy,
 \qquad h\in\mathcal O(\Delta),\qquad
 r=\operatorname{Res}_{a_0}\omega,
\]
for a sufficiently small coordinate disk $\Delta$ centered at zero.
Set $T=T_m$ and factor $T(y)=y\,v(y)$, where $v$ is holomorphic,
nonzero on a smaller disk, and $v(0)=-1$. The pulled-back simple-pole
term is
\[
 r\frac{T'(y)}{T(y)}\,dy
 =r\left(\frac1y+\frac{v'(y)}{v(y)}\right)dy.
\]
It has residue $r$, not $-r$: the factors of $-1$ in $T$ and $T'$
cancel. The double-pole term contributes no residue, since
\[
 c_{-2}\frac{T'(y)}{T(y)^2}\,dy
 =-c_{-2}\,d\!\left(\frac1{T(y)}\right),
\]
and the derivative of a meromorphic Laurent series has zero residue.
The pullback of $h(y)\,dy$ is holomorphic. Hence
\[
 \operatorname{Res}_{0}(T_m^*\omega_0)=r.
\]
Taking residues in the identity of Step~3 gives $0=r+r=2r$, because
$d_y\widehat H$ is holomorphic. Applying this at each focal point yields
\[
 \operatorname{Res}_{a_+}\omega
 =\operatorname{Res}_{a_-}\omega=0.
\]
The two residues vanish separately, not merely in their sum. Since
every accessible differential is now the same $\omega$, this removes
the logarithmic coefficient in every local branch. The double-pole
coefficient need not vanish; it gives the possible simple pole of $V$.

\smallskip\noindent
\emph{Step 5: eliminate all periods and normalize the global primitive.}
For a small positively oriented loop $\ell_\pm$ around $a_\pm$, the
residue theorem and Step~4 give
\[
 \int_{\ell_\pm}\omega
 =2\pi i\,\operatorname{Res}_{a_\pm}\omega=0.
\]
Since $\omega$ is single-valued, the two connector integrals cancel in
$\alpha_\pm*\ell_\pm*\alpha_\pm^{-1}$. Integration of a holomorphic
one-form is invariant under homotopy of piecewise smooth paths with
fixed endpoints and is additive under concatenation. The two based
loops generate the fundamental group from Step~2, so
\[
 \int_\gamma\omega=0
 \quad\text{for every piecewise smooth closed curve }
 \gamma\subset\mathcal D^*.
\]
This conclusion uses the disk-with-two-punctures topology: there are
no further independent cycles whose periods remain to be checked.

Define
\[
 P(z)=V_0(z_0)+\int_{z_0}^{z}\omega,
 \qquad z\in\mathcal D^*,
\]
using any piecewise smooth path on this Riemann surface; local charts
are used if the path passes through infinity. Vanishing of all periods
makes $P$ path independent. Local primitives of the holomorphic
one-form show that $P$ is holomorphic and $dP=\omega$. Near $z_0$,
$P-V_0$ has zero derivative and value zero at $z_0$, so it vanishes
there. Thus $P$ continues the original physical germ.

Along any continuation path the germs of this global function supply
a continuation of $V_0$. Uniqueness along that path identifies every
accessible branch with the corresponding germ of $P$. This eliminates
the remaining additive monodromy, not only derivative monodromy. The
integration construction uses piecewise smooth paths, but the resulting
global holomorphic function supplies continuation along every continuous
path; no rectifiability restriction on the original paths has been added.

\smallskip\noindent
\emph{Step 6: read off the focal pole orders and retain regularity at infinity.}
At $a_0$, the zero-residue expansion is
\[
 \omega_0=\left(\frac{c_{-2}}{y^2}+h(y)\right)dy.
\]
If $H_0$ is a holomorphic primitive of $h$ on the coordinate disk, then
$-c_{-2}/y+H_0(y)$ is a primitive of $\omega_0$ on its punctured part.
The global primitive $P\circ\jmath$ differs from it by a constant
there. Hence $P$ extends meromorphically across $a_0$ with pole order
at most one. This order is unchanged by the biholomorphic coordinate
change $z=\jmath(y)$; if $c_{-2}=0$, the focal point is removable.

At infinity write $\omega=q_\infty(\xi)\,d\xi$, where $q_\infty$ is
holomorphic at zero by Step~1. A local primitive has the form
\[
 C_\infty+\int_0^\xi q_\infty(t)\,dt,
\]
which is holomorphic at $\xi=0$. This agrees locally with $P$ after
choosing the constant, consistently with the inclusion of infinity
in $\mathcal D^*$. Thus the continuation is regular at infinity and
has no poles elsewhere in $\mathcal D\setminus A_+$. We have proved
that the original $V$ is single-valued meromorphic on $\mathcal D$,
with at most simple focal poles. No finite bound on the number of
branches was assumed. Reality and the subsequent reflection argument
are not used to obtain this conclusion.
\end{proof}

Figure~\ref{fig:global-residues-periods} distinguishes global differential
descent from the residue and period arguments needed to recover a primitive.

\begin{figure}[!htb]
\centering
\begin{tikzpicture}[x=1cm,y=1cm,font=\small,
 arr/.style={-{Stealth[length=1.7mm]},line width=.7pt},
 box/.style={draw=black!55,rounded corners=1.5pt,fill=black!3,
   text width=7.10cm,font=\footnotesize,align=center,inner sep=5pt},
 every node/.style={inner sep=2pt}]
\node at (3.55,6.94) {(a) Two generators; infinity is not a puncture};
\node at (11.68,6.94) {(b) A differential is not yet a primitive};
\draw[black!25,line width=.4pt] (7.57,.10)--(7.57,6.57);
\draw[fill=black!2,draw=black!65,line width=.7pt] (3.67,3.85) circle[radius=2.35];
\node[font=\footnotesize] at (3.67,6.43) {$\xi=1/z$, $|\xi|<1$};
\coordinate (aminus) at (2.886666667,3.85);
\coordinate (aplus) at (4.453333333,3.85);
\coordinate (seed) at (3.67,1.876);
\coordinate (lm) at (2.886666667,3.42);
\coordinate (lp) at (4.453333333,3.42);
\draw[line width=.6pt] (seed).. controls (3.09,2.25) and (2.73,2.65)..(lm);
\draw[line width=.6pt] (seed).. controls (4.23,2.25) and (4.65,2.65)..(lp);
\node[font=\footnotesize,fill=white] at (2.96,2.69) {$\alpha_-$};
\node[font=\footnotesize,fill=white] at (4.46,2.69) {$\alpha_+$};
\draw[arr] (lm) arc[start angle=-90,end angle=270,radius=.43];
\draw[arr] (lp) arc[start angle=-90,end angle=270,radius=.43];
\draw[fill=white,line width=.8pt] (aminus) circle[radius=.06];
\draw[fill=white,line width=.8pt] (aplus) circle[radius=.06];
\node[above=15pt,font=\footnotesize] at (aminus) {$a_-^{-1}$};
\node[above=15pt,font=\footnotesize] at (aplus) {$a_+^{-1}$};
\node[left=15pt,font=\footnotesize] at (aminus) {$\ell_-$};
\node[right=15pt,font=\footnotesize] at (aplus) {$\ell_+$};
\fill (3.67,3.85) circle[radius=.045];
\node[align=center,font=\footnotesize] at (3.67,5.19) {$\xi=0\ \leftrightarrow\ z=\infty$\\ordinary point};
\draw[black!40] (3.67,4.86)--(3.67,3.99);
\fill (seed) circle[radius=.045];
\node[below=4pt,font=\footnotesize] at (seed) {$\xi_0$};
\node[align=center,font=\footnotesize,text width=7.10cm] at (3.61,.75)
 {Based loops: $\alpha_\pm*\ell_\pm*\alpha_\pm^{-1}$.\\[3pt]
  Local derivative descent fixes each generator\\
  after every preceding history.};
\node[box,anchor=north] (descent) at (11.68,6.42)
 {Global differential descent:\\[2pt]
  all reached branches have the same $\omega$.\\
  Additive monodromy is still possible.};
\node[box,anchor=north] (paired) at ($(descent.south)+(0,-.28)$)
 {Now $df=dg=\omega_0$ at each focal point:\\[2pt]
  $\omega_0+T_m^*\omega_0=d_y\widehat H$,\\[2pt]
  $\operatorname{Res}(T_m^*\omega_0)=\operatorname{Res}\omega_0$.\\
  Hence $2r_+=0$ and $2r_-=0$ separately.};
\draw[arr] (descent.south)--(paired.north);
\node[box,anchor=north] (periods) at ($(paired.south)+(0,-.28)$)
 {Generating periods: $\int_{\ell_\pm}\omega=2\pi i r_\pm=0$.\\[3pt]
  Therefore every closed-curve period vanishes.};
\draw[arr] (paired.south)--(periods.north);
\node[box,anchor=north] (primitive) at ($(periods.south)+(0,-.28)$)
 {$P(z)=V_0(z_0)+\int_{z_0}^{z}\omega$ is path independent.\\[3pt]
  Every reached germ agrees with $P$.\\
  Its focal poles have order at most one.};
\draw[arr] (periods.south)--(primitive.north);
\end{tikzpicture}
\caption[Global differential descent, residues, and periods]{The distinct global steps in Lemma~\ref{cont:lem:global-derivative}.
Left: $\xi=1/z$ identifies $\mathcal D^*$ with a disk punctured at
$a_+^{-1}$ and $a_-^{-1}$, while infinity becomes the ordinary center.
The displayed punctures are at $\pm1/3$, as for $a=1$, $b=4/5$;
the connectors and loops are schematic. Each connector is traversed back
in its based generator. Right: global descent first identifies the
continued differentials, not yet their primitives. Only then does the
paired identity imply separate vanishing of the two residues. The resulting
zero generating periods give a path-independent primitive. The sign
$T_m'(0)=-1$ does not reverse a one-form residue, and the double-pole
coefficients need not vanish. The diagram imposes no bound on the number
of branches present before the descent argument.}
\label{fig:global-residues-periods}
\end{figure}

\begin{proof}[Completion of Theorem~\ref{cont:thm:logarithmic-closure}]
Lemma~\ref{cont:lem:global-derivative} gives a single-valued meromorphic
$V_{\rm ext}$ on $\mathcal D$, with possible simple poles only at
$\pm R_f$, where $R_f=(a+b)/c>1$ and $c^2=a^2-b^2$, and holomorphic
at infinity. Define
\[
 V_{\rm in}(z)=\overline{V_{\rm ext}(1/\overline z)}.
\]
Conjugating local Laurent series shows that this is meromorphic in the
interior disk, holomorphic at zero, with possible simple poles only at
$\pm R_f^{-1}=\pm(a-b)/c$.

Choose an inversion-invariant annulus about the unit circle contained
in the original holomorphy domain and avoiding the four focal values.
There the original function $V_0$ and
$V_0^\sharp(z)=\overline{V_0(1/\overline z)}$ are holomorphic and agree
on the unit circle by reality. The identity theorem gives equality on
the annulus. Its outer half identifies $V_0$ with $V_{\rm ext}$, and
its inner half identifies $V_0$ with $V_{\rm in}$. These open overlaps
glue the functions to a meromorphic $\mathcal V$ on the sphere,
regular at zero and infinity. Subtracting its finitely many simple
principal parts leaves an entire function bounded at infinity, hence
a constant; in particular $\mathcal V$ is rational.

In the original semiaxes the Laurent extension of the graph weight is
\[
 D_E(z)=\frac{a^2+b^2}{2}-\frac{c^2}{4}(z^2+z^{-2})
       =-\frac{c^2}{4z^2}(z^2-R_f^2)(z^2-R_f^{-2}).
\]
Its four zeros are simple because $R_f>1$. They cancel every possible
focal pole of $\mathcal V$, while $D_E$ has poles of order two at zero
and infinity. Thus $W=D_E\mathcal V$ has no other poles, and at these
two points its orders are at most two. Equivalently, $z^2W(z)$ is entire
and $O(|z|^4)$; Cauchy's estimates make it a polynomial of degree at
most four. Hence
\[
 D_E(z)\mathcal V(z)=\sum_{j=-2}^{2}c_jz^j,
 \qquad c_{-j}=\overline{c_j},\quad c_0\in\R,
\]
where the symmetry follows from reality of $D_Ev$ on the unit circle.
Its restriction is therefore a real trigonometric polynomial of degree
at most two divided by $D_E(\varphi)>0$. Formula~\eqref{geo:eq:tangentspace}
identifies this with $T_E$. This is an infinitesimal tangent-space
classification; no nonlinear graph assertion is used here.
If the five tangent-dual moments also vanish, expansion in the tangent
basis and its exact duality force all five coefficients to vanish.
\end{proof}

\paragraph{Completion for $L^2$ and distributional data.}
For these inputs, Corollary~\ref{band:cor:sobolev-kernels} supplies an
analytic representative from the selected equations with zero low rows
before the analytic classification is invoked.

\begin{proof}[Proof of Theorem~\ref{main:analytic}]
Put $\eta=1-b^2/a^2$, write $\Phi_\rho=\Phi_{E,\rho}$ and
$G_\rho=\Phi_\rho^{-1}$, and choose a rational anchor
$\rho_*=p_0/P\in I$. These angular coordinates are unchanged by
uniform scaling. The analytic circle diffeomorphisms and their inverses
preserve null sets and have bounded positive derivatives. Thus
$U=v\circ\Phi_{\rho_*}$ is a well-defined $L^2$ class, with
\[
 \|U\|_{L^2}^2=\frac1{2\pi}\int_0^{2\pi}
       |v(\varphi)|^2G_{\rho_*}'(\varphi)\,d\varphi
 \asymp\|v\|_{L^2}^2.
\]
For a reduced $\rho=p/d$, let $g_\rho=v\circ\Phi_\rho$.
The finite trace has Fourier multiplier
$\sum_{j=0}^{d-1}e^{2\pi injp/d}$, equal to $d$ when $d\mid n$ and
zero otherwise. Its almost-everywhere constancy therefore gives
\[
 \widehat g_\rho(n)=0\qquad(n\ne0,\ d\mid n).
\]
This includes both frequency signs and imposes no condition on the mean.

Choose the tail \eqref{band:eq:selection}. Its primitive denominators
$d_q$ divide $q$. With $H_q=G_{\rho_q}\circ\Phi_{\rho_*}$ we have
$g_{\rho_q}\circ H_q=U$ almost everywhere. For $n=\pm q$, changing
variables gives the exact selected row
\[
 0=\widehat g_{\rho_q}(n)
  =\frac1{2\pi}\int_0^{2\pi}U(x)e^{-inH_q(x)}H_q'(x)\,dx
  =(\widetilde A_\eta U)_n.
\]
The smooth-kernel integral agrees with the bounded $H^0$ realization.
Use the operator with zero low rows, not the moment-modified operator;
then $\widetilde A_\eta U=0$. Corollary~\ref{band:cor:sobolev-kernels}
provides an analytic representative $U_{\rm an}$, through its
exponential Fourier estimate. Therefore
$v_{\rm an}=U_{\rm an}\circ G_{\rho_*}$ is periodic and real analytic,
represents $v$ almost everywhere, and is real on the real circle.
For each rational rotation its continuous trace equals the original
constant almost everywhere, hence everywhere. One and the same
representative satisfies all these identities. It defines a holomorphic
physical germ $V(e^{i\varphi})=v_{\rm an}(\varphi)$ near the unit circle.

Theorem~\ref{cont:thm:discrete-exterior-continuation} applies to this
analytic germ. Corollary~\ref{cont:cor:finite-bad-enlarged-domain}
then supplies a finite effective bad set $B$ and all-path continuation
on $\Omega\setminus B$, retaining the unit-circle collar.
Theorem~\ref{cont:thm:task6-nonfocal-exclusion} confines $B$ to the
exterior focal pair and makes every accessible branch regular at
infinity. These facts verify the three inputs of
Theorem~\ref{cont:thm:logarithmic-closure}: the initial germ is real
on the circle, paths avoiding the focal pair admit continuation
including through infinity, and each initially constant odd-denominator
trace remains constant along admissible shifted paths by simultaneous
continuation of its finitely many summands. No equality of their
individual branches is required. That theorem yields a rational
continuation, regular at zero and infinity, with at most simple poles
at the four zeros of $D_E$, and $v_{\rm an}\in T_E$.
Since $v=v_{\rm an}$ almost everywhere, all the stated conclusions
follow for the original $L^2$ class, without a moment normalization
or a smallness assumption.
\end{proof}

\Needspace{14\baselineskip}

\begin{corollary}[Distributional interior-band classification]
\label{band:cor:distributional-classification}
Let $E$ be a fixed noncircular ellipse and let
$I\subset(0,1/2)$ be a nonempty open interval. Suppose a real periodic
distribution $v$ satisfies
\[
 \sum_{j=0}^{d-1}
    (v\circ\Phi_{E,\rho})(x+2\pi jp/d)
 \quad\hbox{is a constant distribution}
\]
for every reduced rational $\rho=p/d\in I$. Pullback by the analytic
circle diffeomorphism $\Phi_{E,\rho}$ defines the displayed
composition. Then $v$ is represented by an element of $T_E$.
In particular, an $L^2$ displacement whose traces are constant almost
everywhere is automatically analytic and belongs to $T_E$.
\end{corollary}
\begin{proof}
Write $\Phi_\rho=\Phi_{E,\rho}$ and $G_\rho=\Phi_\rho^{-1}$.
Use complex-linear distributional pairings, with regular functions
acting by unnormalized integration, and
$\widehat w(n)=(2\pi)^{-1}\langle w,e^{-inx}\rangle$.
For an increasing circle diffeomorphism $\Phi$ with inverse $G$,
scalar pullback is
\[
 \langle\Phi^*v,\chi\rangle
 =\langle v,(\chi\circ G)G'\rangle;
\]
see \cite[Section~10.1.3]{Dyatlov}. It agrees with composition on
regular functions and is invertible and functorial. For
$g_\rho=\Phi_\rho^*v$, the Fourier coefficient of the trace is
$\widehat g_\rho(n)\sum_{j=0}^{d-1}e^{2\pi injp/d}$.
Thus the hypothesis gives $\widehat g_\rho(n)=0$ at every nonzero
multiple of its primitive denominator, with both signs included.

Fix a rational anchor $\rho_*\in I$, choose the tail
\eqref{band:eq:selection}, and set $U=\Phi_{\rho_*}^*v$, $H_q=G_{\rho_q}\circ\Phi_{\rho_*}$.
Functoriality gives $H_q^*g_{\rho_q}=U$, hence
\[
 \widehat g_{\rho_q}(n)
 =\frac1{2\pi}\langle U,e^{-inH_q}H_q'\rangle.
\]
For $n=\pm q$ this vanishes because $d_q\mid q$. These pairings are
the selected matrix rows: polynomial growth of distributional Fourier
coefficients and exponential localization make the row sums absolutely
convergent, as proved in Corollary~\ref{band:cor:sobolev-kernels}.
With all low rows zero, $\widetilde A_\eta U=0$ as a distribution.
That corollary makes $U$ analytic. Its inverse pullback gives the
analytic representative
$v_{\rm an}=U_{\rm an}\circ G_{\rho_*}$ of the original distribution.
Reality of $v$ makes $v_{\rm an}$ real. Each analytic rational trace
represents a constant distribution and is therefore constant pointwise.
Theorem~\ref{main:analytic}, now applied to an $L^2$ function, gives
$v_{\rm an}\in T_E$. For an $L^2$ input, the conclusion concerns this
analytic representative of its equivalence class. The regularity
result used above precedes, and does not use, this classification.
\end{proof}

\section{Finite detection and nonlinear rigidity}
\label{close:sec:closure}
We now use the classification to remove the finite-dimensional kernel
of the selected operator, and then absorb the quadratic caustic error.

Fix a nonempty open interval $I\subset(0,1/2)$, a rational anchor
$\rho_*\in I$, and the selection \eqref{band:eq:selection}. In the anchor action coordinate, write
\[
 L_{\eta,\rho,n}U
 =\frac1{2\pi}\int_0^{2\pi}
 U(x)e^{-inH_\eta(\rho,x)}\partial_xH_\eta(\rho,x)\,dx.
\]
If $U=v\circ\Phi_{E,\rho_*}$, the substitution
$y=H_\eta(\rho,x)$ gives
$L_{\eta,\rho,n}U=\widehat{v\circ\Phi_{E,\rho}}(n)$.
For a reduced rational rotation $\rho$, this is a \emph{resonant row}
when $n\ne0$ is divisible by its primitive denominator.
For each fixed $\rho,n$, the input kernel is analytic and the row
is continuous in $\eta$ in the dual norm of every $H^s$.
For distributional $U$, the integral denotes the corresponding pairing.
The five moment rows below are those of \eqref{op:eq:moments},
with anchor $\rho_*$.

\begin{proposition}[Finite detection of the transverse kernel]
\label{band:prop:finite-detection}
Fix the interval and rational tail \eqref{band:eq:selection}, and let
$A^I_\eta$ be $\widetilde A_\eta$ with its five rows $|n|\le2$
replaced by the tangent-dual moments at $\rho_*$.
For each squared eccentricity $\eta_*\in(0,1)$, there are finitely
many reduced rational rotations $\sigma_j\in\mathbb Q\cap I$ and
nonzero indices $m_j$ divisible by their respective primitive
denominators such that, for each real $s$,
\begin{equation}\label{band:eq:detector-inverse}
 \|U\|_{H^s}
 \le C_s\left(\|A^I_\eta U\|_{H^s}
       +\sum_j|L_{\eta,\sigma_j,m_j}U|\right)
\end{equation}
for $\eta$ in a neighborhood of $\eta_*$ and all $U\in H^s$.
The neighborhood and detector family can be fixed independently of
$s$; the constant $C_s$ may depend on $s$.
The analogous estimate holds in every $\mathcal X_\gamma$, $\gamma>2$,
on the same neighborhood, with a constant depending on $\gamma$.
On each compact squared-eccentricity interval contained in $(0,1)$,
one finite family of detector pairs works with constants uniform in
$\eta$ for each fixed $s$ or $\gamma$.
\end{proposition}
\begin{proof}
Write $A_\eta=A^I_\eta$ and $\mathcal R_\eta=Q_\eta A_\eta-I$.
Since $N\ge3$, the five inserted rows replace zero rows and are exactly
the tangent-dual moments of \eqref{op:eq:moments}. The finite-row
versions of Proposition~\ref{band:prop:operator} and
Corollary~\ref{band:cor:sobolev-kernels} give a compact left remainder
and analytic distributional kernels. In particular,
\[
 K_\eta:=\ker(A_\eta:H^0\to H^0)
 =\ker(A_\eta:H^s\to H^s)
 =\ker(A_\eta:\mathcal X_\gamma\to\mathcal X_\gamma)
\]
as sets of distributions for every asserted exponent. These spaces
are finite dimensional, since $\mathcal R_\eta=-I$ on the kernel.
No constancy in $\eta$ is implied.

\smallskip\noindent
\emph{Finite separation and an $L^2$ estimate.}
If $U\in K_{\eta_*}$ is killed by every resonant row in $I$, all
nonzero Fourier coefficients of each rational trace of
$v=U\circ G_{E,\rho_*}$ vanish: the orbit-sum multiplier is $d$ on
$d\Z$ and zero off it. Hence these traces are constant. Apply
Theorem~\ref{main:analytic} to their real and imaginary parts to obtain
$v\in T_E\otimes_\R\C$. In the anchor coordinate,
$U=\sum_{k=-2}^2c_kt_{\eta_*,k}$, and exact moment duality gives
$c_i=\ell_{\eta_*,i}(U)=0$. Thus all resonant rows together separate
$K_{\eta_*}$. Intersecting successively with kernels of actual rows
that are nonzero on the remaining subspace selects at most
$\dim_\C K_{\eta_*}$ pairs $(\sigma_j,m_j)$. If the kernel is zero,
use the empty family. Set $\Lambda_\eta U=(L_{\eta,\sigma_j,m_j}U)_j$.
The pairs, rather than just their rotations, are retained.

If the augmented map were not bounded below on $L^2$ at $\eta_*$,
there would be norm-one $U_k$ with $A_{\eta_*}U_k\to0$ and
$\Lambda_{\eta_*}U_k\to0$. The identity
\[
 U_k=Q_{\eta_*}A_{\eta_*}U_k-\mathcal R_{\eta_*}U_k
\]
gives a strongly convergent subsequence by compactness. Its norm-one
limit contradicts the selected separation. Thus an $L^2$ lower bound
holds. The matrices of $A_\eta$ and the smooth-kernel detector rows are
norm continuous on $L^2$ and its dual. Absorbing a perturbation smaller
than half this lower bound gives one neighborhood $O$ on which the
augmented map remains uniformly bounded below. Fix
$\eta_*\in V$ with $\overline V\Subset O\cap(0,1)$.

\smallskip\noindent
\emph{A localized normal operator gives every regularity on the same set.}
On $\ell^2$ identify functions with their Fourier coefficients and put
\[
 \mathscr N_\eta=A_\eta^*A_\eta+\Lambda_\eta^*\Lambda_\eta.
\]
Finite-dimensional norm equivalence and the uniform $L^2$ estimate give
$cI\le\mathscr N_\eta\le CI$ on $\overline V$, with $c>0$.
The matrices of $A_\eta$ are uniformly exponentially localized; each
fixed detector kernel has exponentially decaying Fourier coefficients.
Hence $\mathscr N_\eta$ is uniformly exponentially localized.
Jaffard's theorem \cite[Proposition~2]{Jaffard} makes
$\mathscr N_\eta^{-1}$ exponentially localized as well.

For $b>0$ and a matrix $M=(M_{nj})_{n,j\in\mathbb Z}$, write
\[
 \|M\|_{b,\mathrm{row}}
 :=\sup_{n\in\mathbb Z}\sum_{j\in\mathbb Z}
       e^{b|n-j|}|M_{nj}|.
\]
The triangle inequality $|n-j|\le|n-k|+|k-j|$ gives
$\|M_1M_2\|_{b,\mathrm{row}}
 \le\|M_1\|_{b,\mathrm{row}}\|M_2\|_{b,\mathrm{row}}$
whenever the displayed norms are finite.

The uniform version needed here follows directly from the proof of
that theorem. Choose $a\ge\max\{C,1\}$ and write
$\mathscr N_\eta=a(I-R_\eta)$, so
$\|R_\eta\|_{\ell^2\to\ell^2}\le r_0<1$ uniformly (enlarge $a$ if
necessary). For some fixed $b>0$ and $S\ge1$, the exponentially weighted
row norm of $R_\eta$ is at most $S$. Submultiplicativity and the
$\ell^2$ bound give
\[
 |(R_\eta^k)_{nj}|\le
 \min\{r_0^k,S^ke^{-b|n-j|}\}.
\]
Split the Neumann series at $k=\lfloor\delta|n-j|\rfloor$, with
$\delta>0$ chosen so $\delta\log S<b/2$. Its first part is bounded by
$C(1+|n-j|)e^{-b|n-j|/2}$ and its second by
$Cr_0^{\delta|n-j|}$. Therefore
\[
 |(\mathscr N_\eta^{-1})_{nj}|\le C'e^{-b'|n-j|}
 \quad(\eta\in\overline V)
\]
for common $C',b'>0$. This avoids inferring uniform decay from
pointwise inverse-closedness alone.

The exact $\ell^2$ identity
\[
 U=\mathscr N_\eta^{-1}A_\eta^*(A_\eta U)
   +\mathscr N_\eta^{-1}\Lambda_\eta^*(\Lambda_\eta U)
\]
now extends to polynomially growing Fourier sequences by absolute
matrix summation. Its first coefficient matrix is exponentially
localized; the finitely many columns in the second term decay
exponentially in $|n|$. The
\hyperref[op:par:localized-calculus]{weighted-matrix estimates}
therefore prove \eqref{band:eq:detector-inverse} on every $H^s$ and its
analogue on every $\mathcal X_\gamma$, on the same $V$. Only their norm
constants depend on the regularity exponent. In particular, no
intersection of exponent-dependent neighborhoods or density in
$\mathcal X_\gamma$ is required.

Finally cover a prescribed compact eccentricity interval by finitely
many such neighborhoods and take the union of their detector pairs.
Adding detector absolute values only enlarges the right side of each
local estimate; the maximum of finitely many constants is uniform on
the interval for each fixed exponent. The augmented map is bounded
below, not asserted surjective. The finite selection is existential
and supplements the fixed infinite tail; no explicit denominator bound
or uniformity as $\eta\to0,1$ is claimed.
\end{proof}

\begin{figure}[!ht]
\centering
\resizebox{\linewidth}{!}{%
\begin{tikzpicture}[x=1cm,y=1cm,font=\footnotesize,
 arr/.style={-{Stealth[length=1.7mm]},line width=.7pt},
 box/.style={draw=black!60,rounded corners=1.5pt,fill=black!3,
 text width=6.65cm,align=center,inner sep=5pt},
 every node/.style={inner sep=2pt}]
\node[font=\small] at (3.65,7.75) {(a) Finite separation of the transverse kernel};
\node[font=\small] at (11.65,7.75) {(b) The family is fixed before the boundary};

\node[box,anchor=north] (kzero) at (3.65,7.28)
 {$K^{(0)}=K_{\eta_*}=\ker A^I_{\eta_*}$\\[2pt]
  Analytic kernel with zero tangent moments.};
\node[box,anchor=north] (kone) at ($(kzero.south)+(0,-.42)$)
 {$K^{(1)}=K^{(0)}\cap\ker L_{\eta_*,\sigma_1,m_1}$};
\draw[arr] (kzero.south)--(kone.north);
\node[box,anchor=north] (kj) at ($(kone.south)+(0,-.42)$)
 {$K^{(j)}=K^{(j-1)}\cap\ker L_{\eta_*,\sigma_j,m_j}$\\[2pt]
  Choose the row nonzero on $K^{(j-1)}$;\\
  the intersection has smaller dimension.};
\draw[arr] (kone.south)--(kj.north);
\node[box,anchor=north] (kend) at ($(kj.south)+(0,-.42)$)
 {$K^{(J_{\rm d})}=\{0\}$,\qquad
  $J_{\rm d}\le\dim_{\mathbb C}K_{\eta_*}$\\[2pt]
  Retain each rotation \emph{and} frequency.};
\draw[arr] (kj.south)--(kend.north);
\node[box,anchor=north] (normal) at ($(kend.south)+(0,-.44)$)
 {$\mathscr N_\eta=A_\eta^*A_\eta+\Lambda_\eta^*\Lambda_\eta\ge cI$\\[3pt]
  One parameter neighborhood, every regularity.};
\draw[arr] (kend.south)--(normal.north);
\node[align=center,text width=7cm,anchor=north] (limits) at ($(normal.south)+(0,-.26)$)
 {Finite rows supplement the infinite selected tail.\\
  The augmented map need not be onto.};

\node[box,anchor=north] (data) at (11.65,7.28)
 {Fix $E_0$, the band $I$, and the anchor $\rho_*$.};
\node[box,anchor=north] (family) at ($(data.south)+(0,-.32)$)
 {Choose the tail, detector pairs, and geometry.\\[2pt]
  These are fixed before $r,M$ or any boundary.};
\draw[arr] (data.south)--(family.north);
\node[box,anchor=north] (bound) at ($(family.south)+(0,-.32)$)
 {Fix $r>11/2,M$: obtain $B$ and constants.\\[2pt]
  $B$ bounds every eligible residual.};
\draw[arr] (family.south)--(bound.north);
\node[box,anchor=north] (smallness) at ($(bound.south)+(0,-.32)$)
 {Choose one $C^1$ threshold $\varepsilon>0$.\\[2pt]
  It is independent of the individual boundary.};
\draw[arr] (bound.south)--(smallness.north);
\node[box,dashed,fill=white,anchor=north] (cutoff) at ($(smallness.south)+(0,-.40)$)
 {Proof only: a hypothetical residual $0<X\le X_*$:\\[3pt]
  $\displaystyle Q=\left\lceil(K_0B/X)^{1/r}\right\rceil\ge N$.\\[3pt]
  The caustic family does not change.};
\draw[arr,dashed] (smallness.south)--(cutoff.north);
\coordinate (divtop) at (7.65,7.38);
\draw[black!25,line width=.4pt] (divtop)--(divtop |- limits.south);
\end{tikzpicture}
}
\caption{Finite detection and the order of choices. (a) The subspace chain is
schematic: $J_{\rm d}$ counts the selected rows, but no actual detector
pairs or billiard kernel dimension are computed. Here $A_\eta=A^I_\eta$
and $\Lambda_\eta$ collects the detectors. A zero kernel needs no detectors.
(b) Here $B$ is a uniform residual $H^r$ bound, $X$ is the residual $L^2$
norm, and $K_0$ is a fixed constant in the proof's cutoff choice.
The detectors are selected using full-band classification of the
\emph{reference linear equations}; the perturbed boundary is required to
have only the prescribed sparse caustic family. Neither the family nor the
geometric neighborhood depends on the hypothetical residual.}
\label{close:fig:finite-detection-choices}
\end{figure}
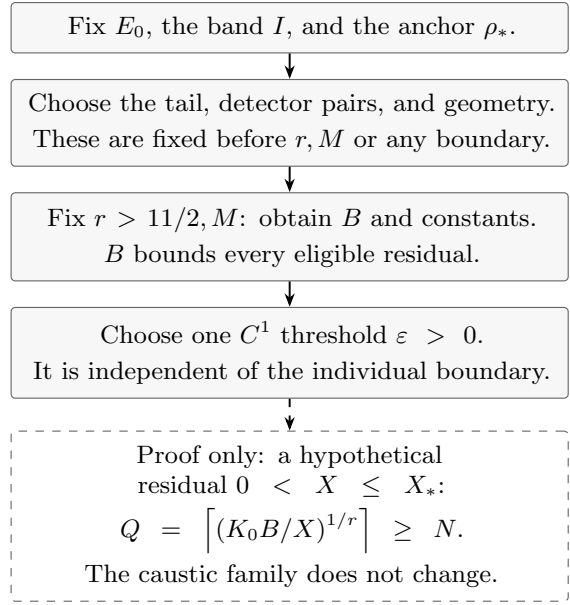

\paragraph{From reference detectors to nonlinear rigidity.}
The full-band classification selects detector rows for the reference
linear problem; it does not impose a full-band caustic hypothesis on the
perturbed table. For a fixed reference ellipse $E_0$, band $I$, and anchor
$\rho_*$, the order of choices in the next theorem is
\[
 \begin{gathered}
 (E_0,I,\rho_*)\ \longrightarrow\
   \text{tail, detector pairs, and geometric neighborhood},\\
 (r,M)\ \longrightarrow\
   \text{uniform residual bound and constants}
   \ \longrightarrow\ \varepsilon.
 \end{gathered}
\]
Only the proof cutoff may depend on a hypothetical nonzero residual;
the prescribed caustic family remains fixed. Figure~\ref{close:fig:finite-detection-choices}
summarizes this distinction and the finite separation argument.

\Needspace{20\baselineskip}

\begin{theorem}[Local rigidity from a sparse interior caustic family]
\label{band:thm:sparse-rigidity}
Fix a noncircular ellipse $E_0$, a nonempty open interval
$I\subset(0,1/2)$, and a rational $\rho_*\in I$. There are an integer
$N$ and a finite set $F\subset\mathbb Q\cap I$ with the following
property. For every real $r>11/2$ and $M>0$ there exists
$\varepsilon>0$ such that a positively curved boundary $\Gamma$ given
by an $H^r$ elliptic graph with displacement $\mu$ over $E_0$ is an
ellipse whenever
\[
 \|\mu\|_{H^r}\le M,\qquad \|\mu\|_{C^1}<\varepsilon,
\]
and its billiard has an integrable rational caustic at every rotation in
\[
 \left\{\frac{\lceil q\rho_*\rceil}{q}:q\ge N\right\}\cup F.
\]
Every rotation in this family lies in $I$, and its only accumulation
point is $\rho_*$.
\end{theorem}
\begin{proof}
Write $\rho_*=p_0/P$ in lowest terms. We first fix the prescribed caustic
family and the geometric neighborhood. Only afterward do we fix $r,M$
and choose a smallness threshold. All Fourier norms below use the
normalized measure and weights in \eqref{main:eq:sobolev}.

\smallskip\noindent
\emph{Step 1: fix the selected tail, detector pairs, and uniform geometric data.}
Choose $N\ge3$ and $J_{\rm tail}\Subset I$ as in
\eqref{band:eq:selection}, so that
\[
 a_q=\left\lceil\frac{qp_0}{P}\right\rceil,\qquad
 \rho_q=\frac{a_q}{q}\in J_{\rm tail}\quad(q\ge N),\qquad
 0\le\rho_q-\rho_*<\frac1q.
\]
Apply Proposition~\ref{band:prop:finite-detection} to this fixed tail at
the squared eccentricity $\eta_0$ of $E_0$. Retain its finite resonant
pairs $(\sigma_j,m_j)$, $1\le j\le J_{\rm d}$, and set
\[
 F=\{\sigma_j:1\le j\le J_{\rm d}\}\subset\mathbb Q\cap I.
\]
Repeated detector rotations are listed only once in $F$, but all detector
pairs are retained in the estimates. The empty detector family is
allowed. The cutoff, pairs, and detector neighborhood are independent of
$r,M$ by that proposition. They remain fixed throughout the proof.

Choose a compact interval $J\Subset I$ containing $J_{\rm tail}\cup F$.
Work in a sufficiently small compact parameter neighborhood of $E_0$,
with semiaxes positive and separated, whose squared eccentricities lie
inside the detector neighborhood. Use smoothly chosen principal-axis
frames. In these frames, rigid motions and uniform scaling do not change
the angular action formulas or rotation numbers, so the reference angular
operators are the $A^I_\eta$ of Proposition~\ref{band:prop:finite-detection}.
The normalization bounds, coordinate bounds, and constants in
Theorem~\ref{geo:thm:remainder} and Corollary~\ref{geo:cor:fourier} can be
chosen uniformly on this geometric neighborhood and on $J$.

\smallskip\noindent
\emph{Step 2: normalize the same boundary and fix the residual high-norm bound.}
Now fix $r>11/2$ and $M>0$. Since $H^r(\T)\subset C^2(\T)$, the input
is eligible for Theorem~\ref{geo:thm:gauge} at the anchor $\rho_*$.
Denote its fitting parameter vector by $\mathbf p(\mu)$, to distinguish
it from the scalar low norm used below. That theorem gives
\[
 |\mathbf p(\mu)|\le C\|\mu\|_{C^0}.
\]
Thus a preliminary bound $\|\mu\|_{C^1}<\varepsilon_{\rm pre}$ keeps the
fitted ellipse $E$ in the neighborhood of Step~1. Write the same boundary
as a residual graph $\nu$ over $E$, and put
\[
 U=\nu\circ\Phi_{E,\rho_*},\qquad
 \eta=\eta(E),\qquad X=\|U\|_{L^2},\qquad t=\|\nu\|_{C^1}.
\]
The boundary and its caustics have not changed. The increasing change of
boundary parameter preserves their rotation numbers and integrable
rationality; positive curvature is a property of the same curve, not a
consequence being inferred from $C^1$ smallness.

The low-norm bounds and norm equivalence in Theorem~\ref{geo:thm:gauge}
give fixed constants $C_{\rm g},C_{\rm eq}\ge1$ such that
\[
 t+\|U\|_{C^1}\le C_{\rm g}\|\mu\|_{C^1},\qquad
 C_{\rm eq}^{-1}t\le\|U\|_{C^1}\le C_{\rm eq}t,
 \qquad X\le C_{\rm g}\|\mu\|_{C^1}.
\]
Shrink $\varepsilon_{\rm pre}$, if needed, so that $t$ meets the uniform
caustic smallness threshold from Step~1. Its choice does not involve any
primitive denominator. The Sobolev part
(\ref{geo:gauge:sobolev}) of the normalization theorem supplies a bound
\[
 \|\nu\|_{H^r}+\|U\|_{H^r}\le M_1.
\]
Fix $B\ge\max\{1,M_1\}$ once and for all. It depends only on
$E_0,I,\rho_*,r,M$ and the data already chosen, not on the individual
boundary. No smallness of $B$, or linear dependence on $M$, is required.
In particular,
\[
 0\le X\le\|U\|_{H^r}\le B.
\]
All these bounds remain valid when the original $C^1$ threshold is later
decreased.

The normalization gives the five exact real tangent-dual moments at
$E,\rho_*$. By the equivalence recorded after \eqref{op:eq:moments},
the five complex moment rows of $A^I_\eta$ also vanish. Its remaining
rows with $2<|n|<N$ are zero by definition. Consequently every output
row with $|n|<N$ vanishes exactly; no low-row error is left over.

\smallskip\noindent
\emph{Step 3: bound all selected harmonics using primitive denominators.}
If $\rho=p/d\in J$ is reduced and the boundary has its prescribed
integrable rational caustic, Corollary~\ref{geo:cor:fourier}, applied to
$\nu$ over the fitted ellipse, gives
\[
 \big|\widehat{\nu\circ\Phi_{E,\rho}}(n)\big|
 \le C_{\rm ca}d^2t^2\qquad(n\ne0,\ d\mid n).
\]
The constant is uniform over Step~1's data and independent of the
nonzero resonant harmonic. The underlying orbit-sum estimate imposes no
condition of the form $d^\kappa t\ll1$; only the fixed bound on $t$ from
Step~2 is needed.

Write $\rho_q=p_q/d_q$ in lowest terms. Then
\[
 d_q=\frac{q}{\gcd(a_q,q)},\qquad d_q\mid q.
\]
The exact coordinate-row identity preceding
Proposition~\ref{band:prop:finite-detection} identifies the two selected
rows with harmonics of this primitive caustic. Hence
\[
 |(A^I_\eta U)_{\pm q}|
 =\big|\widehat{\nu\circ\Phi_{E,\rho_q}}(\pm q)\big|
 \le C_{\rm ca}d_q^2t^2\le C_{\rm ca}q^2t^2\qquad(q\ge N).
\]
This also covers nonreduced displayed fractions and repetitions of a
selected rotation. If $d_j$ is the primitive denominator of $\sigma_j$,
then $d_j\mid m_j$, and the same estimate gives
\[
 |L_{\eta,\sigma_j,m_j}U|\le C_{\rm ca}d_j^2t^2,\qquad
 \sum_j|L_{\eta,\sigma_j,m_j}U|\le C_{\rm det}t^2.
\]
Here $C_{\rm det}$ is fixed by the finite pairs; an empty sum is zero.
Thus one caustic for each distinct rotation in $F$ supplies all its
chosen detector harmonics. These are quadratic error estimates, not
exact vanishing of the reference high rows.

\smallskip\noindent
\emph{Step 4: split the output, using different estimates on its two parts.}
Write $W=A^I_\eta U$ and $w_n=\widehat W(n)$. Boundedness on $H^r$
from Proposition~\ref{band:prop:operator} and
Corollary~\ref{band:cor:sobolev-kernels}, including the moment rows, is
uniform on the chosen compact eccentricity set. Thus
\[
 \|W\|_{H^r}\le C_{A,r}B.
\]
For every integer $Q\ge N$, Parseval, Step~2, and the selected-row
estimate imply
\[
 \sum_{|n|\le Q}|w_n|^2
 \le 2C_{\rm ca}^2t^4\sum_{q=N}^{Q}q^4
 \le C_{\rm low}^2t^4Q^5.
\]
For the high output, use its Sobolev norm instead:
\[
 \left(\sum_{|n|>Q}|w_n|^2\right)^{1/2}
 \le (1+Q)^{-r}\|W\|_{H^r}
 \le C_{A,r}BQ^{-r}.
\]
The inverse estimate \eqref{band:eq:detector-inverse} with $s=0$,
together with the detector bound, therefore yields a fixed $C_0\ge1$
with
\begin{equation}\label{close:cutoff-estimate}
 X\le C_0\bigl(t^2Q^{5/2}+BQ^{-r}+t^2\bigr),\qquad Q\ge N.
\end{equation}
The exponent $5/2$ comes from summing the squares of the $q^2$ row
bounds over both frequency signs. That growing row estimate is not
summed over the infinite tail; the high output is controlled by the
bounded linear operator on $H^r$. The constants are independent of
$Q,X,t$ and of the individual residual.

\smallskip\noindent
\emph{Step 5: obtain the precise $C^1$ interpolation by Fourier splitting.}
Put $u_n=\widehat U(n)$ and $\langle n\rangle=1+|n|$. Since $r>3/2$,
Cauchy--Schwarz makes $\sum_n\langle n\rangle|u_n|$ finite. The Fourier
series and its first derivative converge uniformly, so
\[
 \|U\|_{C^1}\le C\sum_n\langle n\rangle|u_n|.
\]
For any real $K\ge1$, the low and high sums satisfy
\begin{align*}
 \sum_{|n|\le K}\langle n\rangle|u_n|
 &\le\left(\sum_{|n|\le K}\langle n\rangle^2\right)^{1/2}X
 \le CK^{3/2}X,\\
 \sum_{|n|>K}\langle n\rangle|u_n|
 &\le\left(\sum_{|n|>K}\langle n\rangle^{2-2r}\right)^{1/2}
                  \|U\|_{H^r}
 \le C_rK^{3/2-r}B.
\end{align*}
The second series is bounded by an integral since $2-2r<-1$.
If $X=0$, continuity gives $U=0$, hence $\nu=0$ and $\Gamma=E$.
Otherwise $0<X\le B$, so $K=(B/X)^{1/r}\ge1$ is admissible.
Both terms become $B^{3/(2r)}X^{1-3/(2r)}$. Using
$t\le C_{\rm eq}\|U\|_{C^1}$, fix $C_{\rm int}\ge1$ such that
\begin{equation}\label{close:eq:endpoint-interpolation}
 t\le C_{\rm int}B^{3/(2r)}X^{1-3/(2r)}.
\end{equation}
A real cutoff is legitimate in these index sets; rounding it would only
change constants. This is a direct two-norm estimate, not an invocation
of the endpoint embedding $H^{3/2}\hookrightarrow C^1$.

\smallskip\noindent
\emph{Step 6: choose the final threshold and absorb the nonlinear terms.}
Continue under the assumption $X>0$. Fix
\[
 K_0\ge\max\{1,4C_0\},\qquad
 Z=\left(\frac{K_0B}{X}\right)^{1/r},\qquad Q=\lceil Z\rceil.
\]
Because $X\le B$, we have $Z\ge1$ and $Z\le Q\le2Z$. If
$X\le K_0B/N^r$, then $Q\ge N$, and
\[
 C_0BQ^{-r}\le C_0BZ^{-r}=\frac{C_0}{K_0}X\le\frac X4.
\]
For the other terms, \eqref{close:eq:endpoint-interpolation} and
$Q\le2Z$ give
\begin{align*}
 C_0t^2Q^{5/2}
 &\le D_1X^{2-11/(2r)},&
 D_1&=C_0C_{\rm int}^2\,2^{5/2}K_0^{5/(2r)}B^{11/(2r)},\\
 C_0t^2
 &\le D_2X^{2-3/r},&
 D_2&=C_0C_{\rm int}^2B^{3/r}.
\end{align*}
These constants are positive and fixed before the final smallness
threshold is chosen. Put
\[
 \delta_1=1-\frac{11}{2r}>0,\qquad
 \delta_2=1-\frac3r>0,
\]
and choose, for example,
\[
 X_*=
 \min\left\{1,B,\frac{K_0B}{N^r},
                 (8D_1)^{-1/\delta_1},(8D_2)^{-1/\delta_2}\right\}>0.
\]
For $0<X\le X_*$, the two nonlinear terms have sum at most $X/4$.
Equation~\eqref{close:cutoff-estimate} would therefore imply
\[
 X\le\frac X4+D_1X^{1+\delta_1}+D_2X^{1+\delta_2}
   \le\frac X2,
\]
a contradiction. Choose the theorem's original threshold as
\[
 \varepsilon=\min\{\varepsilon_{\rm pre},X_*/(2C_{\rm g})\}.
\]
The bound $X\le C_{\rm g}\|\mu\|_{C^1}$ then rules out every nonzero
residual. Thus $X=0$, and the same reasoning as in Step~5 gives
$\nu=0$ and $\Gamma=E$ exactly.

The order of these choices is important: the caustic family and
neighborhood were fixed first; $r,M$ then fixed $B$ and the constants;
only afterward was $\varepsilon$ decreased. The cutoff depends on a
hypothetical residual, but \eqref{close:cutoff-estimate} holds for every
$Q\ge N$ and the entire selected tail was already prescribed. No
boundary-dependent caustic is added. At $r=11/2$ the first exponent
would equal one, so shrinking $X$ alone would not absorb its fixed
coefficient. No endpoint or optimality claim follows from this argument.

\smallskip\noindent
\emph{Step 7: verify the sparse family's unique accumulation point.}
All selected rotations lie in $I$ by Step~1, and $F\subset I$.
The estimate $0\le\rho_q-\rho_*<1/q$ shows that a sequence of distinct
tail values can accumulate only at $\rho_*$: its associated indices
must be unbounded. This point is actually attained as a limit of
distinct values. For sufficiently large integers $k$, take $q=kP+1$.
Since $0<p_0/P<1/2$,
\[
 \rho_{kP+1}=\frac{kp_0+1}{kP+1},\qquad
 \rho_{kP+1}-\rho_*
      =\frac{P-p_0}{P(kP+1)}>0.
\]
These differences strictly decrease to zero. Reducing the displayed
fractions does not change their values, and repetitions at other
indices are harmless. A finite set has no accumulation point, so
adjoining $F$ adds none. This completes all assertions with $N,F$
independent of $r,M$.
\end{proof}

\begin{remark}
The finite detector set supplements an infinite caustic sequence.
The theorem does not assert rigidity from finitely many caustics.
The strict condition $r>11/2$ is used in absorbing
$X^{2-11/(2r)}$; no endpoint or optimality claim is made.
\end{remark}

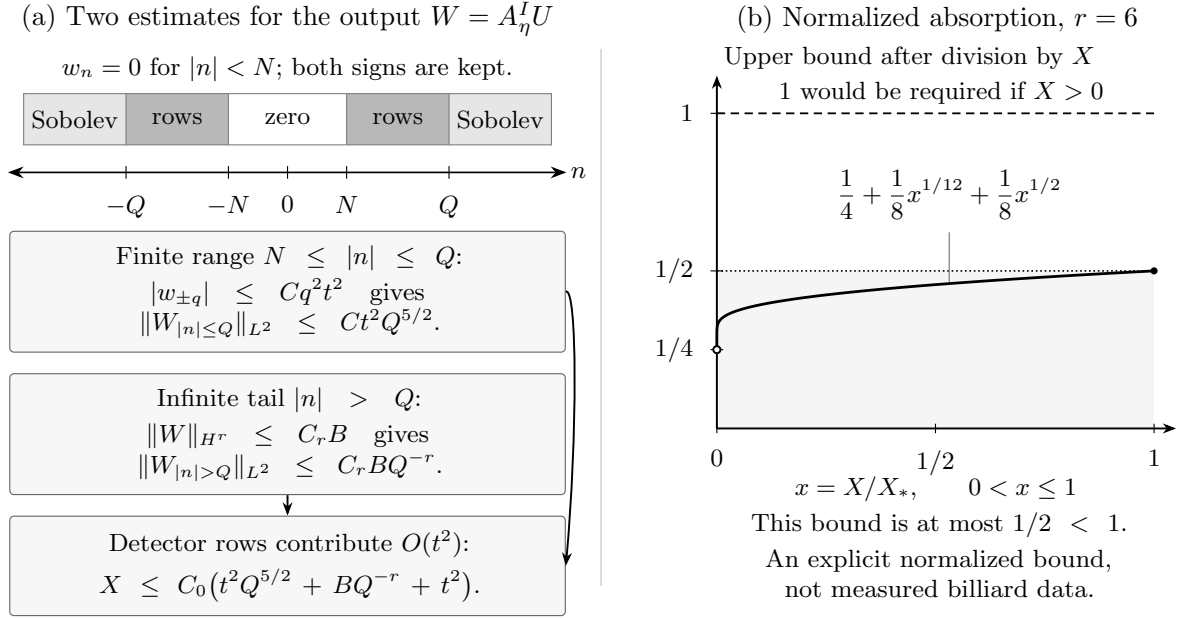
\begin{figure}[!ht]
\centering
\resizebox{\linewidth}{!}{%
\begin{tikzpicture}[x=1cm,y=1cm,font=\footnotesize,
 arr/.style={-{Stealth[length=1.7mm]},line width=.7pt},
 box/.style={draw=black!60,rounded corners=1.5pt,fill=black!3,
 text width=6.7cm,align=center,inner sep=5pt},
 every node/.style={inner sep=2pt}]
\node[font=\small] at (3.70,7.2) {(a) Two estimates for the output $W=A^I_\eta U$};
\node[font=\small] at (11.95,7.2) {(b) Normalized absorption, $r=6$};
\draw[black!25,line width=.4pt] (7.68,0)--(7.68,6.82);
\fill[black!10] (.35,5.6) rectangle (1.65,6.25);
\fill[black!28] (1.65,5.6) rectangle (2.95,6.25);
\fill[white] (2.95,5.6) rectangle (4.45,6.25);
\fill[black!28] (4.45,5.6) rectangle (5.75,6.25);
\fill[black!10] (5.75,5.6) rectangle (7.05,6.25);
\draw[black!50] (.35,5.6) rectangle (7.05,6.25);
\foreach \xx in {1.65,2.95,4.45,5.75} {\draw[black!50] (\xx,5.6)--(\xx,6.25);}
\node at (1.0,5.93) {Sobolev};
\node at (2.30,5.93) {rows};
\node at (3.70,5.93) {zero};
\node at (5.10,5.93) {rows};
\node at (6.40,5.93) {Sobolev};
\draw[<->,>=Stealth,line width=.6pt] (.15,5.25)--(7.24,5.25);
\foreach \xx/\lab in {1.65/{-Q},2.95/{-N},3.70/{0},4.45/{N},5.75/{Q}}
 {\draw (\xx,5.17)--(\xx,5.33);\node[below=4pt] at (\xx,5.17) {$\lab$};}
\node[right] at (7.22,5.25) {$n$};
\node[align=center] at (3.70,6.57)
 {$w_n=0$ for $|n|<N$; both signs are kept.};
\node[box,anchor=north] (rows) at (3.70,4.51)
 {Finite range $N\le |n|\le Q$:\\[2pt]
  $|w_{\pm q}|\le Cq^2t^2$\quad gives\quad
  $\|W_{|n|\le Q}\|_{L^2}\le Ct^2Q^{5/2}$.};
\node[box,anchor=north] (tail) at ($(rows.south)+(0,-.26)$)
 {Infinite tail $|n|>Q$:\\[2pt]
  $\|W\|_{H^r}\le C_rB$\quad gives\quad
  $\|W_{|n|>Q}\|_{L^2}\le C_rBQ^{-r}$.};
\node[box,anchor=north] (inverse) at ($(tail.south)+(0,-.26)$)
 {Detector rows contribute $O(t^2)$:\\[3pt]
  $\displaystyle X\le C_0\bigl(t^2Q^{5/2}+BQ^{-r}+t^2\bigr)$.};
\draw[arr] (rows.east).. controls (7.38,3.8) and (7.38,1.9)..(inverse.east);
\draw[arr] (tail.south)--(inverse.north);
\fill[black!4] (9.15,2.0) rectangle (14.70,4.0);
\draw[arr] (9.15,2.0)--(14.98,2.0);
\draw[arr] (9.15,2.0)--(9.15,6.48);
\node[anchor=west] at (9.16,6.70) {Upper bound after division by $X$};
\draw[densely dashed,line width=.65pt] (9.15,6.0)--(14.70,6.0);
\node[anchor=south,fill=white] at (12.00,6.04) {$1$ would be required if $X>0$};
\draw[densely dotted,line width=.6pt] (9.15,4.0)--(14.70,4.0);
\foreach \yy/\lab in {3.0/{1/4},4.0/{1/2},6.0/{1}}
 {\draw (9.07,\yy)--(9.23,\yy);\node[left=4pt] at (9.07,\yy) {$\lab$};}
\foreach \xx/\lab in {9.15/{0},11.925/{1/2},14.70/{1}}
 {\draw (\xx,1.94)--(\xx,2.06);\node[below=4pt] at (\xx,1.94) {$\lab$};}
\draw[line width=1.0pt] plot coordinates {
(9.1500000000,3.0000000000)
(9.1500000000,3.0020833333)
(9.1500000000,3.0041666667)
(9.1500000000,3.0062500000)
(9.1500000000,3.0083333333)
(9.1500000000,3.0104166667)
(9.1500000000,3.0125000001)
(9.1500000000,3.0145833336)
(9.1500000000,3.0166666674)
(9.1500000000,3.0187500014)
(9.1500000000,3.0208333359)
(9.1500000000,3.0229166713)
(9.1500000000,3.0250000078)
(9.1500000000,3.0270833460)
(9.1500000000,3.0291666864)
(9.1500000000,3.0312500298)
(9.1500000000,3.0333333772)
(9.1500000000,3.0354167298)
(9.1500000000,3.0375000890)
(9.1500000000,3.0395834564)
(9.1500000000,3.0416668341)
(9.1500000000,3.0437502244)
(9.1500000000,3.0458336300)
(9.1500000000,3.0479170540)
(9.1500000000,3.0500005000)
(9.1500000000,3.0520839721)
(9.1500000000,3.0541674749)
(9.1500000000,3.0562510136)
(9.1500000000,3.0583345941)
(9.1500000001,3.0604182230)
(9.1500000001,3.0625019073)
(9.1500000001,3.0645856554)
(9.1500000002,3.0666694760)
(9.1500000003,3.0687533790)
(9.1500000004,3.0708373751)
(9.1500000005,3.0729214763)
(9.1500000007,3.0750056953)
(9.1500000010,3.0770900463)
(9.1500000014,3.0791745445)
(9.1500000019,3.0812592064)
(9.1500000025,3.0833440501)
(9.1500000034,3.0854290948)
(9.1500000046,3.0875143615)
(9.1500000061,3.0895998725)
(9.1500000080,3.0916856520)
(9.1500000105,3.0937717259)
(9.1500000136,3.0958581218)
(9.1500000177,3.0979448693)
(9.1500000227,3.1000320000)
(9.1500000291,3.1021195475)
(9.1500000371,3.1042075478)
(9.1500000471,3.1062960388)
(9.1500000594,3.1083850610)
(9.1500000747,3.1104746573)
(9.1500000934,3.1125648732)
(9.1500001164,3.1146557567)
(9.1500001445,3.1167473587)
(9.1500001788,3.1188397329)
(9.1500002202,3.1209329359)
(9.1500002704,3.1230270274)
(9.1500003308,3.1251220703)
(9.1500004034,3.1272181307)
(9.1500004903,3.1293152783)
(9.1500005941,3.1314135859)
(9.1500007177,3.1335131303)
(9.1500008644,3.1356139919)
(9.1500010382,3.1377162550)
(9.1500012435,3.1398200078)
(9.1500014855,3.1419253426)
(9.1500017699,3.1440323561)
(9.1500021035,3.1461411491)
(9.1500024938,3.1482518270)
(9.1500029495,3.1503645000)
(9.1500034804,3.1524792828)
(9.1500040977,3.1545962950)
(9.1500048139,3.1567156613)
(9.1500056431,3.1588375116)
(9.1500066016,3.1609619810)
(9.1500077071,3.1630892101)
(9.1500089801,3.1652193450)
(9.1500104433,3.1673525377)
(9.1500121221,3.1694889459)
(9.1500140451,3.1716287335)
(9.1500162441,3.1737720703)
(9.1500187547,3.1759191328)
(9.1500216165,3.1780701037)
(9.1500248736,3.1802251725)
(9.1500285752,3.1823845354)
(9.1500327756,3.1845483957)
(9.1500375352,3.1867169638)
(9.1500429208,3.1888904572)
(9.1500490066,3.1910691010)
(9.1500558743,3.1932531281)
(9.1500636141,3.1954427789)
(9.1500723256,3.1976383020)
(9.1500821186,3.1998399541)
(9.1500931135,3.2020480000)
(9.1501054433,3.2042627133)
(9.1501192535,3.2064843762)
(9.1501347040,3.2087132797)
(9.1501519701,3.2109497238)
(9.1501712437,3.2131940178)
(9.1501927348,3.2154464806)
(9.1502166730,3.2177074404)
(9.1502433090,3.2199772353)
(9.1502729164,3.2222562134)
(9.1503057937,3.2245447332)
(9.1503422657,3.2268431632)
(9.1503826865,3.2291518828)
(9.1504274409,3.2314712820)
(9.1504769472,3.2338017619)
(9.1505316598,3.2361437348)
(9.1505920718,3.2384976241)
(9.1506587178,3.2408638653)
(9.1507321774,3.2432429054)
(9.1508130779,3.2456352034)
(9.1509020985,3.2480412307)
(9.1509999733,3.2504614711)
(9.1511074960,3.2528964210)
(9.1512255234,3.2553465898)
(9.1513549805,3.2578125000)
(9.1514968646,3.2602946873)
(9.1516522511,3.2627937012)
(9.1518222981,3.2653101048)
(9.1520082524,3.2678444753)
(9.1522114557,3.2703974041)
(9.1524333502,3.2729694972)
(9.1526754863,3.2755613751)
(9.1529395288,3.2781736735)
(9.1532272648,3.2808070432)
(9.1535406118,3.2834621503)
(9.1538816255,3.2861396767)
(9.1542525092,3.2888403203)
(9.1546556227,3.2915647950)
(9.1550934924,3.2943138312)
(9.1555688216,3.2970881760)
(9.1560845012,3.2998885933)
(9.1566436216,3.3027158642)
(9.1572494844,3.3055707872)
(9.1579056155,3.3084541787)
(9.1586157782,3.3113668727)
(9.1593839878,3.3143097216)
(9.1602145259,3.3172835963)
(9.1611119565,3.3202893862)
(9.1620811420,3.3233280000)
(9.1631272607,3.3264003654)
(9.1642558249,3.3295074298)
(9.1654726998,3.3326501602)
(9.1667841235,3.3358295440)
(9.1681967280,3.3390465886)
(9.1697175609,3.3423023224)
(9.1713541087,3.3455977943)
(9.1731143209,3.3489340748)
(9.1750066348,3.3523122555)
(9.1770400028,3.3557334500)
(9.1792239193,3.3591987940)
(9.1815684503,3.3627094453)
(9.1840842633,3.3662665845)
(9.1867826594,3.3698714152)
(9.1896756065,3.3735251639)
(9.1927757738,3.3772290809)
(9.1960965683,3.3809844403)
(9.1996521727,3.3847925402)
(9.2034575850,3.3886547030)
(9.2075286601,3.3925722762)
(9.2118821527,3.3965466321)
(9.2165357626,3.4005791681)
(9.2215081819,3.4046713077)
(9.2268191440,3.4088245000)
(9.2324894747,3.4130402206)
(9.2385411461,3.4173199715)
(9.2449973318,3.4216652818)
(9.2518824654,3.4260777076)
(9.2592223006,3.4305588327)
(9.2670439748,3.4351102687)
(9.2753760741,3.4397336554)
(9.2842487023,3.4444306612)
(9.2936935515,3.4492029832)
(9.3037439767,3.4540523477)
(9.3144350726,3.4589805107)
(9.3258037537,3.4639892578)
(9.3378888382,3.4690804049)
(9.3507311341,3.4742557984)
(9.3643735300,3.4795173156)
(9.3788610882,3.4848668649)
(9.3942411428,3.4903063863)
(9.4105634004,3.4958378517)
(9.4278800453,3.5014632652)
(9.4462458490,3.5071846636)
(9.4657182833,3.5130041165)
(9.4863576380,3.5189237268)
(9.5082271431,3.5249456311)
(9.5313930959,3.5310720000)
(9.5559249919,3.5373050385)
(9.5818956618,3.5436469862)
(9.6093814127,3.5501001179)
(9.6384621749,3.5566667438)
(9.6692216543,3.5633492100)
(9.7017474895,3.5701498986)
(9.7361314157,3.5770712283)
(9.7724694341,3.5841156550)
(9.8108619867,3.5912856716)
(9.8514141388,3.5985838087)
(9.8942357665,3.6060126350)
(9.9394417521,3.6135747578)
(9.9871521853,3.6212728230)
(10.0374925723,3.6291095157)
(10.0905940520,3.6370875608)
(10.1465936190,3.6452097229)
(10.2056343557,3.6534788071)
(10.2678656709,3.6618976593)
(10.3334435477,3.6704691664)
(10.4025307990,3.6791962570)
(10.4752973325,3.6880819014)
(10.5519204238,3.6971291125)
(10.6325849990,3.7063409458)
(10.7174839275,3.7157205000)
(10.8068183227,3.7252709172)
(10.9007978549,3.7349953836)
(10.9996410726,3.7448971296)
(11.1035757355,3.7549794306)
(11.2128391578,3.7652456069)
(11.3276785631,3.7756990245)
(11.4483514502,3.7863430956)
(11.5751259712,3.7971812785)
(11.7082813212,3.8082170784)
(11.8481081411,3.8194540481)
(11.9949089326,3.8308957876)
(12.1489984865,3.8425459453)
(12.3107043246,3.8544082181)
(12.4803671552,3.8664863518)
(12.6583413429,3.8787841417)
(12.8449953933,3.8913054328)
(13.0407124519,3.9040541204)
(13.2458908193,3.9170341505)
(13.4609444817,3.9302495203)
(13.6863036570,3.9437042786)
(13.9224153588,3.9574025259)
(14.1697439765,3.9713484157)
(14.4287718729,3.9855461539)
(14.7000000000,4.0000000000)
};
\draw[fill=white,line width=.7pt] (9.15,3.0) circle[radius=.045];
\fill (14.70,4.0) circle[radius=.045];
\node[anchor=south,align=center,fill=white] at (12.10,4.58)
 {$\displaystyle \frac14+\frac18x^{1/12}+\frac18x^{1/2}$};
\draw[black!45,line width=.5pt] (12.10,4.49)--(12.10,3.86);
\node at (11.93,1.24) {$x=X/X_*$,\qquad $0<x\le1$};
\node[align=center,text width=6.7cm] at (11.95,.40)
 {This bound is at most $1/2<1$.\\[2pt]
  An explicit normalized bound,\\not measured billiard data.};
\end{tikzpicture}
}
\caption{The output split and the nonlinear contradiction. (a) The frequency
partition is schematic. The growing $q^2$ row bound is summed only over a
finite range; the infinite tail uses the $H^r$ norm of the output.
(b) With $Q=\lceil(K_0B/X)^{1/r}\rceil\ge N$, the proof chooses $X_*$ so that
$D_iX_*^{\delta_i}\le1/8$ for $i=1,2$, where
$\delta_1=1-11/(2r)$ and $\delta_2=1-3/r$. Thus, for $x=X/X_*\in(0,1]$,
the right-hand side of \eqref{close:cutoff-estimate}, divided by $X$, is at
most $1/4+x^{\delta_1}/8+x^{\delta_2}/8\le1/2$.
The graph shows this explicit upper bound for $r=6$, not a computed
smallness threshold. It illustrates the strict range $r>11/2$ and makes no
endpoint or optimality claim.}
\label{close:fig:fourier-split-absorption}
\end{figure}

Figure~\ref{close:fig:fourier-split-absorption} separates the two output
estimates and shows how the fixed-threshold argument rules out a nonzero
residual.

\begin{proof}[Proof of Theorem~\ref{main:band}]
Fix a rational anchor $\rho_*=\rho_*(I)\in\mathbb Q\cap I$ once the
interval $I$ is fixed, independently of the boundary under consideration.
Apply Theorem~\ref{band:thm:sparse-rigidity} to $E_0,I,\rho_*$. Let
$N$ and $F\subset\mathbb Q\cap I$ be the resulting sparse-family data and,
for the prescribed $r>11/2$ and $M>0$, let
$\varepsilon_{\rm sp}=\varepsilon_{\rm sp}(E_0,I,\rho_*,r,M)>0$ be its
smallness threshold. Set $\varepsilon=\varepsilon_{\rm sp}$. The same
theorem gives
\[
 \mathcal S:=
 \left\{\frac{\lceil q\rho_*\rceil}{q}:q\ge N\right\}\cup F
 \subset\mathbb Q\cap I.
\]
These are rational rotation values: reducing a displayed fraction does
not change the required caustic, and repeated rotations impose no
additional condition. Hence the full-band hypothesis supplies every
integrable rational caustic required by
Theorem~\ref{band:thm:sparse-rigidity}. The same boundary is a positively
curved $H^r$ elliptic graph over the same $E_0$, with the same bound
$\|\mu\|_{H^r}\le M$ and the chosen threshold
$\|\mu\|_{C^1}<\varepsilon_{\rm sp}$. That theorem therefore implies
that $\Gamma$ is an ellipse. Since the anchor was fixed from $I$, the
chosen $\varepsilon$ has the stated dependence on $E_0,I,r,M$; no
uniformity over all rational anchors in $I$ is asserted or needed.
\end{proof}

\begin{corollary}[H\"older regularity]
\label{main:cor:holder-regularity}
Theorems~\ref{main:band} and \ref{band:thm:sparse-rigidity} also hold with
$\mu\in C^{5,\alpha}(\T)$ and a fixed $C^{5,\alpha}$ norm bound
replacing the Sobolev regularity assumption and bound, for every
$1/2<\alpha\le1$. Here $C^{5,1}$ means that the fifth derivative is
Lipschitz. In particular, a fixed $C^6$ bound suffices.
The smallness threshold may depend on $\alpha$.
\end{corollary}
\begin{proof}
Fix $\alpha\in(1/2,1]$ and choose
$\beta=(\alpha+1/2)/2\in(1/2,\alpha)$, so $\beta<1$.
For $g=D^5\mu$, the periodic Gagliardo characterization
\eqref{app:eq:gagliardo} and H\"older continuity give
\[
 [g]_\beta^2\le4\pi[g]_{C^\alpha}^2
       \int_0^\pi h^{2(\alpha-\beta)-1}\,dh
 =\frac{2\pi^{1+2(\alpha-\beta)}}{\alpha-\beta}[g]_{C^\alpha}^2.
\]
Together with the $L^2$ term controlling the zero frequency and
\eqref{app:eq:top-derivative-norm}, this proves the continuous embedding
$\|\mu\|_{H^{5+\beta}}\le C_{\alpha,\beta}\|\mu\|_{C^{5,\alpha}}$.
A prescribed bound $M_{\rm H}$ thus gives the fixed Sobolev bound
$M_{\rm S}=C_{\alpha,\beta}M_{\rm H}$ at
$r=5+\beta>11/2$. Apply the two rigidity theorems with these data,
retaining positive curvature, the stated caustics, and the separate
$C^1$ smallness condition. The sparse family $N,F$ is unchanged because
Theorem~\ref{band:thm:sparse-rigidity} chooses it independently of
$r,M$; only the smallness threshold depends on $\alpha,M_{\rm H}$.
At $\alpha=1$ this uses a Lipschitz fifth derivative, not a continuous
sixth derivative. Finally,
$[D^5\mu]_{\mathrm{Lip}}\le\|D^6\mu\|_\infty$ by integration along
shortest circle arcs, so a fixed $C^6$ bound gives a fixed $C^{5,1}$
bound and the claimed consequence. No endpoint assertion at
$\alpha=1/2$ is made.
\end{proof}

\appendix
\numberwithin{equation}{section}
\section{Norm estimates for changes of periodic graph coordinates}
\label{app:sec:norms}
This appendix collects the standard analytic estimates used for graph
coordinates, in the precise norms and uniformity needed here. We use
Hamilton's tame calculus \cite[Part~II, Sections~2.2--2.3]{Hamilton},
the Sobolev composition and inversion results of
Inci--Kappeler--Topalov \cite[Appendix~B]{IKT}, and the fractional
norm characterization \cite[Proposition~3.4]{DPV}. Short specializations
retain the Jacobian lower bounds and the independence from high norms.
No geometric normalization or caustic hypothesis is used; the final
proposition applies to an arbitrary small coordinate parameter.

Write $\T=\R/(2\pi\Z)$ and $D=d/dx$. For an integer $k\ge0$, use
$\|f\|_{C^k}=\max_{0\le j\le k}\|D^jf\|_{L^\infty}$. For real $s$, set
\[
 \|f\|_{H^s}^2=\sum_{n\in\Z}(1+|n|)^{2s}|\widehat f(n)|^2,
 \qquad
 \widehat f(n)=\frac1{2\pi}\int_0^{2\pi}f(x)e^{-inx}\dd x.
\]
These conventions agree with the main text, up to equivalent choices
of the $C^k$ norm. If $h(x)=x+q(x)$ is a degree-one lift, the notation
$h-\mathrm{id}\in C^k$ or $H^s$ always refers to its periodic part $q$.
All constants below are uniform for smooth coefficient families with
uniform bounds for the finitely many derivatives used on a fixed compact
parameter and amplitude neighborhood.

\Needspace{16\baselineskip}
\subsection{Integer interpolation and tame inverse composition}
\label{app:sub:integer}

\begin{lemma}[Interpolation and smooth substitution]\label{app:lem:integer}
For an integer $k\ge2$ and $g\in C^k(\T)$,
\begin{equation}\label{app:eq:integer-interpolation}
 \|D^jg\|_{L^\infty}
 \le C_k\|g\|_{C^1}^{1-(j-1)/(k-1)}
              \|g\|_{C^k}^{(j-1)/(k-1)},\qquad 1\le j\le k.
\end{equation}
Let $P(x,z)$ be smooth and periodic in $x$, defined on a neighborhood
of $\T$ times a fixed compact real interval. For inputs $u$ with range
in that interval and $\|u\|_{C^1}\le L$,
\begin{equation}\label{app:eq:integer-substitution}
 \|P(\cdot,u(\cdot))\|_{C^k}\le C_{k,L,P}(1+\|u\|_{C^k}).
\end{equation}
The constants do not depend on the high norms of the inputs.
\end{lemma}
\begin{proof}
The interpolation estimate is
\cite[Part~II, Theorem~2.2.1]{Hamilton} with lower order $1$ and upper
order $k$; finite-regularity functions follow by periodic smoothing.
Smooth substitution is the degree-zero tame estimate
\cite[Part~II, Theorems~2.2.5--2.2.6]{Hamilton}. In the present
specialization, a chain-rule product has
$\sum j m_j\le k$, so the sum of its interpolation exponents
$\sum(j-1)m_j/(k-1)$ is at most one. Its remaining factors are bounded
by the prescribed $C^1$ norm and coefficient bounds. This gives
\eqref{app:eq:integer-substitution} with no dependence of its constant
on the high norm.
\end{proof}

\begin{proposition}[Tame inverse composition]\label{app:prop:integer-inverse}
Fix an integer $k\ge2$, $c_0>0$, and $L>0$. Suppose $q,f\in C^k(\T)$,
$h(x)=x+q(x)$, $h'\ge c_0$, and $\|q\|_{C^1}\le L$. Then
\begin{equation}\label{app:eq:integer-inverse}
 \|f\circ h^{-1}\|_{C^k}
 \le C_{k,c_0,L}\bigl(\|f\|_{C^k}
                  +\|f\|_{C^1}\|q\|_{C^k}\bigr).
\end{equation}
In particular, no derivative of order $k+1$ is required.
\end{proposition}
\begin{proof}
This is the tame inverse-composition estimate of
\cite[Part~II, Section~2.3]{Hamilton}; the following specialization
records the uniformity in $c_0,L$. Composing with $h$ turns differentiation
into $(h')^{-1}D$. At order $k$ every term of
$((h')^{-1}D)^kf$ has the form
\begin{equation}\label{app:eq:inverse-monomial}
 \frac{f^{(j)}\prod_{\ell\ge2}(q^{(\ell)})^{m_\ell}}
      {(h')^{k+M}},\qquad
 j+\sum_{\ell\ge2}(\ell-1)m_\ell=k,\quad M=\sum_{\ell\ge2}m_\ell.
\end{equation}
The denominator is controlled by $c_0$. For a mixed term, interpolation
in Lemma~\ref{app:lem:integer} assigns high exponents
$\theta=(j-1)/(k-1)$ and
$\kappa=\sum(\ell-1)m_\ell/(k-1)$ with $\theta+\kappa=1$.
The $C^1$ bound on $q$ and weighted Young's inequality therefore bound
it by
$C_{k,c_0,L}(\|f\|_{C^k}+\|f\|_{C^1}\|q\|_{C^k})$.
The unmixed term and lower orders satisfy the same estimate; the
zeroth-order supremum is unchanged by composition. Only derivatives
through order $k$ occur.
\end{proof}

\subsection{Fractional product, multiplier, and composition estimates}
\label{app:sub:fractional}
For $0<\beta<1$, let $d_\T$ denote circle distance and define
\begin{equation}\label{app:eq:gagliardo}
 [g]_\beta^2=
 \int_\T\int_\T
 \frac{|g(x)-g(y)|^2}{d_\T(x,y)^{1+2\beta}}\,dx\,dy.
\end{equation}
The periodic version of the Gagliardo characterization
\cite[Proposition~3.4]{DPV} gives
$\|g\|_{H^\beta}^2\asymp_\beta\|g\|_{L^2}^2+[g]_\beta^2$.
Indeed, Parseval and Tonelli express the seminorm as
\[
 2\pi\sum_n|\widehat g(n)|^2
       \int_{-\pi}^{\pi}|e^{inh}-1|^2|h|^{-1-2\beta}\,dh;
\]
scaling makes the inner integral comparable to $|n|^{2\beta}$ for
$n\ne0$. The upper bound follows by integration on $\R$, and the lower
bound by integration on $[-\pi,\pi]$. This also checks the periodic
normalization. Unnormalized $L^2$ measure changes only fixed constants.
For an integer $m\ge1$, the Fourier weights also give
\begin{equation}\label{app:eq:top-derivative-norm}
 \|g\|_{H^{m+\beta}}
 \asymp_{m,\beta}\|g\|_{L^2}+\|D^mg\|_{H^\beta},
 \qquad 0\le\beta<1.
\end{equation}
The $L^2$ term controls the zero frequency.

\begin{lemma}[Elementary fractional estimates]\label{app:lem:fractional}
Fix $0<\beta<1$.
For $f,g\in H^\beta(\T)\cap L^\infty(\T)$,
\begin{equation}\label{app:eq:frac-product}
 \|fg\|_{H^\beta}
 \le C_\beta\bigl(\|f\|_\infty\|g\|_{H^\beta}
                 +\|g\|_\infty\|f\|_{H^\beta}\bigr).
\end{equation}
If $g\in H^\beta$ and $a$ is Lipschitz, then
\begin{equation}\label{app:eq:frac-multiplier}
 \|ag\|_{H^\beta}
 \le C_\beta\bigl(\|a\|_\infty+\operatorname{Lip}(a)\bigr)
                \|g\|_{H^\beta}.
\end{equation}
If $\kappa$ is a $C^1$ circle diffeomorphism with
$L^{-1}\le\kappa'\le L$, $L\ge1$, then
\begin{equation}\label{app:eq:frac-pullback}
 \|g\circ\kappa\|_{H^\beta}\le C_{\beta,L}\|g\|_{H^\beta}.
\end{equation}
Finally, for a smooth function $P(x,z)$, periodic in $x$, and a real
input $u\in H^\beta\cap L^\infty$ with range in a fixed compact
interval in its domain,
\begin{equation}\label{app:eq:frac-coefficient}
 \|P(\cdot,u(\cdot))\|_{H^\beta}\le C_{\beta,P}(1+\|u\|_{H^\beta}).
\end{equation}
For $\beta=0$, the corresponding $L^2$ estimates hold under the same
boundedness assumptions. In \eqref{app:eq:frac-multiplier}, an
$L^\infty$ coefficient suffices when $\beta=0$.
\end{lemma}
\begin{proof}
These are standard fractional product and composition estimates; for
the product estimate see also \cite[Lecture~II, equation~(4.20)]{Taylor}.
The difference identity
$\delta_h(fg)=f(\cdot+h)\delta_hg+g\,\delta_hf$, inserted in
\eqref{app:eq:gagliardo}, proves \eqref{app:eq:frac-product}.
For a Lipschitz multiplier the extra term is bounded by
$C\operatorname{Lip}(a)^2\|g\|_2^2\int_0^\pi h^{1-2\beta}\,dh$;
this proves \eqref{app:eq:frac-multiplier} without assuming $g$ bounded.
Changing both variables using the bi-Lipschitz map gives
$[g\circ\kappa]_\beta^2\le L^{3+2\beta}[g]_\beta^2$ and
$\|g\circ\kappa\|_2^2\le L\|g\|_2^2$.
Finally bounded first derivatives of $P$ imply
$|P(x,u(x))-P(y,u(y))|\le C(d_\T(x,y)+|u(x)-u(y)|)$, giving the last
estimate. For $\beta=0$ only the $L^2$ changes of variables and
multiplication estimates are needed.
\end{proof}

\subsection{Sobolev substitution, inversion, and composition}
\label{app:sub:sobolev}
We retain the range $s>11/2$ used in the main theorem; no sharper
threshold is claimed here. Write $s=m+\beta$, $m=\lfloor s\rfloor\ge5$,
$0\le\beta<1$. Cauchy--Schwarz on the Fourier coefficients gives
\[
 H^{m+\beta}(\T)\hookrightarrow C^{m-1}(\T),\qquad
 \|D^jg\|_\infty
 \le\left(\sum_{n\in\Z}(1+|n|)^{-2(s-j)}\right)^{1/2}\|g\|_{H^s}
 \quad (j\le m-1).
\]
In particular these norms control $C^2$. Every derivative of order at
most $m-1$ is also bounded in $H^{1+\beta}\subset H^\beta\cap L^\infty$.
Only the $m$th derivative may fail to be bounded when $\beta\le1/2$.

\begin{lemma}[Smooth substitution in $H^s$]\label{app:lem:sobolev-substitution}
Let $s>11/2$, $M>0$, and let $P(x,z)$ be smooth and periodic in $x$,
on a neighborhood of $\T$ times a fixed compact real interval.
If $u\in H^s(\T)$ is real, has range in that interval, and
$\|u\|_{H^s}\le M$, then
\[
 \|P(\cdot,u(\cdot))\|_{H^s}\le C_{s,M,P}.
\]
The bound is uniform in additional compact smooth parameters of $P$.
\end{lemma}
\begin{proof}
This is the Sobolev smooth-substitution estimate; see
\cite[Lecture~II, Section~4]{Taylor}. Its uniform version here also
follows from
$D^m(P(x,u))=P_z(x,u)D^mu+\mathcal R_m$.
The leading coefficient is Lipschitz with a bound depending only on
the fixed data and $M$. All factors in $\mathcal R_m$ involve lower
derivatives, bounded in $H^\beta\cap L^\infty$ by the embedding above.
Lemma~\ref{app:lem:fractional} and
\eqref{app:eq:top-derivative-norm} give the bound. Positive periodic
mollifiers preserve the range and bounded $H^s$ norms; weak lower
semicontinuity and convergence in $C^1$ extend it to the stated inputs.
Uniform coefficient bounds give uniformity in additional parameters.
\end{proof}

\begin{proposition}[Sobolev bounds for inverses and inverse composition]
\label{app:prop:sobolev-inverse}
Fix $s>11/2$, $M>0$, and $c_0>0$. Suppose $q,f\in H^s(\T)$ are real,
$h(x)=x+q(x)$, $h'\ge c_0$, and
$\|q\|_{H^s}+\|f\|_{H^s}\le M$. With $k=h^{-1}$, one has
\begin{equation}\label{app:eq:sobolev-inverse}
 \|k-\mathrm{id}\|_{H^s}+\|f\circ k\|_{H^s}
 \le C_{s,M,c_0}.
\end{equation}
No linear dependence on $M$ is asserted.
\end{proposition}
\begin{proof}
Sobolev inversion and composition on the circle are covered by
\cite[Appendix~B, Lemmas~B.5--B.6 and Theorem~B.2]{IKT}. We record
why the constants are uniform on the present bounded class, rather
than appeal only to continuity of these maps. The $H^s$ bound controls
$C^2$, and $h'\ge c_0$ makes $h,k=h^{-1}$ uniformly bi-Lipschitz.
For smooth inputs, differentiation gives
\begin{equation}\label{app:eq:inverse-leading}
 k^{(j)}=-\frac{h^{(j)}\circ k}{(h'\circ k)^{j+1}}+P_j,
 \qquad 2\le j\le m,
\end{equation}
where $P_j$ uses only $h^{(i)}\circ k$, $i<j$, and reciprocal powers
of $h'\circ k$. For $j=m$ all lower factors are controlled in
$H^\beta\cap L^\infty$ by Lemma~\ref{app:lem:fractional}.
The top factor $h^{(m)}\circ k$ is controlled in $H^\beta$, while
$(h'\circ k)^{-(m+1)}$ has a controlled Lipschitz norm, using only
$h''$ and $k'$. The same lemma bounds the leading product. Together
with $k-\mathrm{id}=-q\circ k$ in $L^2$, this proves the inverse bound.

In $D^m(f\circ k)$ the only top-order terms are
$(f^{(m)}\circ k)(k')^m+(f'\circ k)k^{(m)}$.
Their multipliers are Lipschitz with controlled norms; all remaining
products have lower-order factors in $H^\beta\cap L^\infty$.
The fractional estimates and \eqref{app:eq:top-derivative-norm}
therefore prove the composition bound. Their $L^2$ versions apply
when $\beta=0$. For general inputs, positive mollification preserves
$h'\ge c_0$ and the high-norm bounds; inverse maps and compositions
converge in $C^1$, and weak lower semicontinuity gives the same uniform
bounds. This is a boundedness argument, not a claim of smooth inversion
without derivative loss at a fixed Sobolev level.
\end{proof}

\begin{corollary}[Composition with a fixed smooth family]
\label{app:cor:smooth-composition}
Let $\kappa$ range over a smooth family of orientation-preserving circle
diffeomorphisms on a compact parameter set, with uniform bi-Lipschitz
bounds and uniform finite-order derivative bounds. For each integer
$k\ge1$ and each real $s>11/2$,
\[
 \|f\circ\kappa\|_{C^k}\le C_k\|f\|_{C^k},\qquad
 \|f\circ\kappa\|_{H^s}\le C_s\|f\|_{H^s}.
\]
The $C^1$ norms of $f$ and $f\circ\kappa$ are uniformly equivalent.
For the $H^s$ estimate, derivatives of $\kappa$ through order
$\lfloor s\rfloor+1$, together with the bi-Lipschitz bounds, suffice.
\end{corollary}
\begin{proof}
The integer bounds and $C^1$ equivalence follow from the chain rule for
$\kappa$ and $\kappa^{-1}$. For $s=m+\beta$, differentiate $m$ times
and apply Lemma~\ref{app:lem:fractional} to each
$f^{(j)}\circ\kappa$. Its coefficient is Lipschitz, controlled by
$\kappa$ through order $m+1$ and the bi-Lipschitz bounds. Equation
\eqref{app:eq:top-derivative-norm} gives the stated linear estimate;
this is the fixed-map specialization of \cite[Appendix~B]{IKT}.
\end{proof}

\Needspace{34\baselineskip}
\subsection{Smooth changes of graph coordinates}
\label{app:sub:graph-change}

\begin{proposition}[Norm bounds for a change of graph coordinates]
\label{app:prop:graph-change}
Let $t\in\R^d$, and suppose $\mathsf A(t,x,z)$ and $\mathsf B(t,x,z)$
are real-valued and smooth for $t,z$ near zero and $x\in\R$, with
\[
 \mathsf A(t,x+2\pi,z)=\mathsf A(t,x,z)+2\pi,\qquad
 \mathsf B(t,x+2\pi,z)=\mathsf B(t,x,z),
\]
and
\[
 \mathsf A(0,x,z)=x,\qquad \mathsf B(0,x,z)=z.
\]
There exist $\delta,\varepsilon,C>0$ such that, for $|t|<\delta$ and
real $u\in C^1(\T)$ with $\|u\|_{C^1}<\varepsilon$, the maps
\[
 h_{t,u}(x)=\mathsf A(t,x,u(x)),\qquad
 f_{t,u}(x)=\mathsf B(t,x,u(x)),\qquad
 w_{t,u}=f_{t,u}\circ h_{t,u}^{-1}
\]
are well defined, $h_{t,u}'\ge1/2$, and
\begin{equation}\label{app:eq:graph-low}
 \|h_{t,u}-\mathrm{id}\|_{C^1}\le C|t|(1+\|u\|_{C^1}),\qquad
 \|w_{t,u}\|_{C^1}\le C(|t|+\|u\|_{C^1}).
\end{equation}
For every integer $k\ge2$, if $u\in C^k$, then
\begin{equation}\label{app:eq:graph-integer}
 \|w_{t,u}\|_{C^k}\le C_k(\|u\|_{C^k}+|t|).
\end{equation}
For every real $s>11/2$ and $M>0$, if $u\in H^s$ with
$\|u\|_{H^s}\le M$, then
\begin{equation}\label{app:eq:graph-sobolev}
 \|w_{t,u}\|_{H^s}\le C_{s,M}.
\end{equation}
The choices $\delta,\varepsilon$ depend only on fixed low-order coefficient
bounds, not on $k,s,M$. For fixed indices and norm bounds, the estimates
are uniform for coefficient families with uniform finite-order bounds
on a common neighborhood. The parameter $t$ is arbitrary in this
neighborhood; no normalization equation for $t$ is assumed.
\end{proposition}
\begin{proof}
Integrating the parameter derivatives from $0$ to $t$ writes
$\mathsf A(t,x,z)-x=\sum t_\ell A_\ell(t,x,z)$ and
$\mathsf B(t,x,z)-z=\sum t_\ell B_\ell(t,x,z)$ with uniformly smooth
periodic coefficients. The $C^1$ chain rule gives the stated bound for $h-\mathrm{id}$ and
$\|f\|_{C^1}\le C(|t|+\|u\|_{C^1})$. Choosing fixed small
$\delta,\varepsilon$ ensures $h'\ge1/2$ and keeps the range in the
coordinate neighborhood. Composition preserves the zeroth-order
supremum and $(w'\circ h)=f'/h'$, proving the remaining low bound.

For $u\in C^k$, put $q=h-\mathrm{id}$. Lemma~\ref{app:lem:integer}
gives
\[
 \|q\|_{C^k}\le C_k|t|(1+\|u\|_{C^k}),\qquad
 \|f\|_{C^k}\le C_k(\|u\|_{C^k}+|t|).
\]
Proposition~\ref{app:prop:integer-inverse} then gives
\eqref{app:eq:graph-integer}, since the mixed term
$(\|u\|_{C^1}+|t|)|t|(1+\|u\|_{C^k})$ has the same bound on the fixed
low-norm neighborhood. For $H^s$ inputs, smooth substitution bounds
$q,f$ by $C_{s,M}$, and Proposition~\ref{app:prop:sobolev-inverse} with
$c_0=1/2$ proves \eqref{app:eq:graph-sobolev}. The neighborhood was
chosen solely from low-order coefficient bounds; all higher constants
are uniform for the stated families. No equation determining $t$ has
been used.
\end{proof}

\end{document}